\documentclass{memo-l}
\frontmatter
\usepackage{xurl}

\newtheorem*{introproposition}{Proposition}
\newtheorem*{introdefinition}{Definition}
\newtheorem*{introtheorem}{Theorem}
\newtheorem*{introcorollary}{Corollary}
\usepackage{tabularx,booktabs}
\usepackage[T1]{fontenc}
\usepackage[utf8]{inputenc}
\usepackage[british]{babel}
\usepackage{lmodern}
\usepackage{microtype}
\usepackage{enumitem}
\usepackage{tikz-cd}
\usepackage{placeins}
\usepackage{placeins}
\usepackage{xltabular}
\FloatBarrier
\theoremstyle{remark}
\newtheorem{discussion}{Discussion}
\usepackage{caption}
\usepackage{amsmath}
\usepackage{amsthm}
\usepackage{amssymb}
\usepackage{mathtools}
\usepackage{mathrsfs}
\usepackage{bm}
\usepackage{float}
\usepackage{enumitem}
\usepackage{graphicx}
\usepackage{booktabs}
\usepackage{tikz-cd}
\theoremstyle{remark}
\usepackage{cite}
\usepackage[
    colorlinks=true,
    linkcolor=blue,
    citecolor=blue,
    urlcolor=blue
]{hyperref}
\usepackage[
    nameinlink,
    capitalise,
    noabbrev
]{cleveref}
\usepackage{mathtools}

\allowdisplaybreaks[2]

\def \J {\mathcal{M}}
\def \pp {\mathfrak{p}}
\def \Li {\mathrm{Li}}

\def \S {\mathbb{S}}

\def \r {{g_0}}
\def \st {\mathrm{ST}}
\def \GL {\mathrm{GL}}
\def \G {\mathrm{Conj}(G)}

\def \L {\Lambda}

\def \l {\langle}

\def \rr {\rangle}

\def \Irr {\mathrm{Irr}}
\def \Gal {\mathrm{Gal}}

\def \Gm {\mathbb{G}_m}
\def \Rep {\mathrm{Rep}}
\def \Sw {\mathrm{Sw}}

\def \ff {\mathcal{F}}

\newcommand{\Q}{\mathbb{Q}}
\newcommand{\N}{\mathbb{N}}
\newcommand{\Z}{\mathbb{Z}}
\newcommand{\R}{\mathbb{R}}
\newcommand{\C}{\mathbb{C}}

\newcommand{\X}{\mathbb{M}}
\newcommand{\M}{\mathbb{M}}

\newcommand{\bb}[1]{\mathbb{#1}}

\def \e {\mathbb{R}^d}

\def \rhoo {\varrho}

\newcommand{\ind}{\mathbf{1}}

\DeclareMathOperator{\s}{\mathrm{vol}_g}

\def \r {r}
\def \u {\psi^{-1}}
\def \K {\mathcal{K}(\mathbb{M})}

\def \T {{\bb T}_{\Lambda}}

\def \f {\phi}

\DeclareMathOperator{\Hom}{Hom}

\numberwithin{equation}{chapter}

\DeclareMathOperator{\Frob}{Frob}

\def \inft {{L^\infty(\bb{M})}}
\def \s {\mathrm{vol}_g}

\allowdisplaybreaks
\theoremstyle{plain}

\newtheorem{theorem}{Theorem}[chapter]
\newtheorem{proposition}[theorem]{Proposition}
\newtheorem{lemma}[theorem]{Lemma}
\newtheorem{corollary}[theorem]{Corollary}

\newtheorem{conjecture}[theorem]{Conjecture}

\theoremstyle{definition}
\def \s {\mathrm{vol}_g}

\newtheorem{definition}[theorem]{Definition}

\newtheorem{hypothesis}[theorem]{Hypothesis}

\newtheorem{example}[theorem]{Example}
\newtheorem{remark}[theorem]{Remark}

\newtheorem{question}[theorem]{Question}

\title{
Spectral Beurling-Selberg Approximation and Applications to
Effective Motivic Chebotarev-Sato-Tate Estimates
}

\usepackage{booktabs}

\title{Spectral Beurling-Selberg approximation and application to effective motivic Chebotarev-Sato-Tate estimates}
\author{Yassine El Kaouni}
\address{Ibn Tofail University, Kenitra, Morocco}
\email{yassine.el\_kaouni@ens-paris-saclay.fr}
\date{}
\begin{document}
\pagestyle{headings}
\begin{abstract}
We establish a spectral Beurling--Selberg approximation theorem on compact connected Riemannian manifolds for sets whose boundaries have finite upper Minkowski content. Our majorants and minorants have controlled spectral coefficients and explicit stability under Hausdorff perturbations, uniformly over families satisfying a common Minkowski content bound. We extend the construction to singular spaces arising from Riemannian submetries with basic mean curvature, and apply it to equidistribution on quaternionic Shimura varieties. Assuming the conjectural formalism of the specified motivic $L$-functions, including the motivic Riemann hypothesis, we develop effective Chebotarev--Sato--Tate estimates for fixed sets and varying families with controlled boundary geometry and Hausdorff variation, and describe how additional symmetries improve the error exponents. The proofs combine classical analytic methods with new methods from differential geometry, analytic geometry, and Riemannian spectral analysis. Applications include an asymptotic formula for a motivic analogue of Atkin--Serre estimates and tubular estimates for Picard-rank jumps of $K3$ surfaces via their transcendental motives.
\end{abstract}
\maketitle
\tableofcontents
\mainmatter
\chapter{Introduction}
The problem of approximating characteristic functions by finite linear
combinations of Fourier exponentials lies at the intersection of harmonic
analysis, discrepancy theory, and analytic number theory. In its classical
periodic form, the Beurling--Selberg problem asks, for an interval
$I\subset[0,1)$, for trigonometric polynomials $P^\pm$ such that
$P^-\leq\ind_I\leq P^+$ with effective control on the Fourier coefficients and
the $L^1$-gap $\|P^+-P^-\|_{L^1(\bb{R}/\bb{Z})}$. The extremal constructions
of Beurling and Selberg, and their refinements, provide such
approximations with applications to Erd\H{o}s--Tur\'an discrepancy inequalities
\cite{Vaaler1985,ErdosTuran1948}. Tensor-product and classical multidimensional
constructions extend this principle to boxes and rectangles.

On a compact connected Riemannian manifold $(\M,g)$ of dimension $d$, the spectral analogue consists of replacing trigonometric polynomials by finite linear combinations of
eigenfunctions of the Laplace--Beltrami operator $\Delta_g$. A significant
precedent is the work of Colzani, Gigante, and Travaglini
\cite{ColzaniGiganteTravaglini2011}, which established constructions of majorants and minorants
and generalized Erd\H{o}s--Tur\'an inequalities for a wide class of domains
on compact manifolds. Our purpose is to develop a Riemannian analogue of the Beurling--Selberg construction adapted to sets with singular boundaries, with effective control
of the spectral coefficients and, crucially, of their variation under
Hausdorff perturbations of the underlying set. This control under Hausdorff perturbations makes our construction suitable for
arithmetic questions concerning equidistribution on moving sets.

The geometric quantity governing our estimates is the upper Minkowski content. More precisely, we focus on sets $A\subset\M$ for which there exist
$0<q\leq d$ and constants $C_{\partial A},r_0>0$ such that, for every
$0<r<r_0$,
$$
\operatorname{vol}_g\bigl((\partial A)^{(r)}\bigr)
\leq C_{\partial A}r^q.
$$
Unlike a smoothness assumption, this condition accommodates many cases of nonsmooth boundaries,
including suitable rectifiable, stratified, and fractal examples.
For families of sets $\mathcal{A}=(A_i)_{i\in I}$, we further require a uniform choice of the constants above, leading to the notion of \emph{uniform Minkowski content} developed below.
The underlying geometric measure
theory, including criteria for uniform bounds over families, is developed
in detail in \cref{content}, see also \cite{Mattila1995,GrayTubes}.

Our first result constructs pointwise spectral majorants and minorants for
$\ind_A$ on a connected compact Riemannian manifold $(\M,g)$, with an
explicit $L^1$-gap and spectral coefficient decay.
Write $(\varphi_j)_{j\geq0}$ for a real orthonormal basis of
Laplace-Beltrami eigenfunctions, with corresponding eigenvalues $(\lambda_j)_{j\geq 0}$.

\begin{introtheorem}[\cref{yes}]
Let $A\subset\M$ have nonempty topological boundary
and satisfy
$$
\overline{\mathcal M}_g^{d-q}(\partial A)<\infty
\qquad\text{for some }0<q\leq d.
$$
Fix $0<r_0<1$, and set
$$
C_{\partial A}
:=
\sup_{0<r<r_0}
\frac{\operatorname{vol}_g((\partial A)^{(r)})}{r^q}
<\infty.
$$
Suppose that, for some $\beta\geq0$,
$$
\|\varphi_j\|_{L^\infty(\M)}
\ll_{\M}\lambda_j^\beta,
$$
for all $j\geq 1$.
Fix an integer $s>d+4\beta$.
For every $0<r<r_0$ and $T\geq1$, there exist
spectral polynomials
$$
P_{r,T,A}^{\pm}
=
\sum_{\lambda_j<T}c_{j,A}^{\pm}\varphi_j
$$
with the following properties.
\begin{enumerate}
\item The following approximation properties hold:
\begin{enumerate}
\item
$$
P_{r,T,A}^-\leq\mathbf 1_A\leq P_{r,T,A}^+.
$$

\item For every integer $n\geq1$ and every $j\geq1$
with $\lambda_j<T$,
$$
|c_{j,A}^{\pm}|
\ll_{\M,n}
\frac{
C_{\partial A}\|\varphi_j\|_{L^\infty(\M)}
}{
r^{n-q}\lambda_j^{n/2}
}.
$$

\item There exists a constant $K_{\M,s,\beta}>0$ such that
$$
\int_{\M}(P_{r,T,A}^+-P_{r,T,A}^-)
\,d\operatorname{vol}_g
\leq
C_{\partial A}r^q
+
\frac{
K_{\M,s,\beta}C_{\partial A}
}{
r^{s-q}T^{(s-d-4\beta)/2}
}.
$$
\end{enumerate}

\item Our construction is stable under Hausdorff perturbations in the following sense:
Let $B\subset\M$ satisfy the same assumptions, with
the same Minkowski codimension $q$ and the same $r_0$, and define $C_{\partial B}$
analogously.
\begin{enumerate}
\item For every integer $n\geq1$ and every $j\geq1$
with $\lambda_j<T$,
$$
|c_{j,A}^+-c_{j,B}^+|
\ll_{\M,n}
\frac{
(C_{\partial A}+C_{\partial B})
d_H^g(\overline A,\overline B)\,
\|\varphi_j\|_{L^\infty(\M)}
}{
r^{n+1-q}\lambda_j^{n/2}
},
$$
and
$$
|c_{j,A}^--c_{j,B}^-|
\ll_{\M,n}
\frac{
(C_{\partial A}+C_{\partial B})
d_H^g\bigl(
\overline{\M\setminus A},
\overline{\M\setminus B}
\bigr)\,
\|\varphi_j\|_{L^\infty(\M)}
}{
r^{n+1-q}\lambda_j^{n/2}
}.
$$

\item For the constant coefficients corresponding
to $\lambda_0=0$,
$$
|c_{0,A}^+-c_{0,B}^+|
\ll_{\M}
\frac{d_H^g(\overline A,\overline B)}{r}
+
\frac{
K_{\M,s,\beta}|C_{\partial A}-C_{\partial B}|
}{
r^{s-q}T^{(s-d-4\beta)/2}
},
$$
and
$$
|c_{0,A}^--c_{0,B}^-|
\ll_{\M}
\frac{
d_H^g\bigl(
\overline{\M\setminus A},
\overline{\M\setminus B}
\bigr)
}{r}
+
\frac{
K_{\M,s,\beta}|C_{\partial A}-C_{\partial B}|
}{
r^{s-q}T^{(s-d-4\beta)/2}
}.
$$
\end{enumerate}
\end{enumerate}
In both cases, the dependence on $\M$, $s$, $\beta$, and $n$ is effective.
\end{introtheorem}

The Hausdorff stability of these coefficients yields the following result
for families of sets with finite uniform upper Minkowski content.

\begin{introcorollary}[\cref{family}]
Let $\mathcal{A}=(A_i)_{i\in I}$ be a family of subsets of $\M$,
each with nonempty topological boundary, and put
$
\partial\mathcal{A}:=(\partial A_i)_{i\in I}.
$
Assume that
$$
\overline{\mathcal M}_g^{d-q}(\partial\mathcal{A})<\infty
\qquad\text{for some }0<q\leq d.
$$
Fix $0<r_0<1$, and set
$$
C_{\partial\mathcal{A}}
:=
\sup_{0<r<r_0}\sup_{i\in I}
\frac{\operatorname{vol}_g((\partial A_i)^{(r)})}{r^q}.
$$
Suppose that there exists $\beta\geq0$ such that, for every $j\geq1$,
$$
\|\varphi_j\|_{L^\infty(\M)}
\ll_{\M}\lambda_j^\beta.
$$
Fix an integer $s>d+4\beta$.
For every $0<r<r_0$ and $T\geq1$, there exist,
for each $\ell\in I$, spectral polynomials
$$
P_{r,T,A_\ell}^{\pm}
=
\sum_{\lambda_j<T}c_{j,A_\ell}^{\pm}\varphi_j,
$$
constructed using the common constant $C_{\partial\mathcal{A}}$,
with the following properties.
\begin{enumerate}
\item For every $\ell\in I$, the following approximation properties hold.
\begin{enumerate}
\item
$$
P_{r,T,A_\ell}^-
\leq\mathbf 1_{A_\ell}
\leq P_{r,T,A_\ell}^+.
$$

\item For every integer $n\geq1$ and every $j\geq1$
with $\lambda_j<T$,
$$
|c_{j,A_\ell}^{\pm}|
\ll_{\M,n}
\frac{
C_{\partial\mathcal{A}}\|\varphi_j\|_{L^\infty(\M)}
}{
r^{n-q}\lambda_j^{n/2}
}.
$$

\item There exists a constant $K_{\M,s,\beta}>0$,
independent of $\ell$, such that
$$
\int_{\M}
(P_{r,T,A_\ell}^+-P_{r,T,A_\ell}^-)
\,d\operatorname{vol}_g
\leq
C_{\partial\mathcal{A}}r^q
+
\frac{
K_{\M,s,\beta}C_{\partial\mathcal{A}}
}{
r^{s-q}T^{(s-d-4\beta)/2}
}.
$$
\end{enumerate}

\item
For every $k,\ell\in I$, the following stability estimates hold.
\begin{enumerate}
\item For every integer $n\geq1$ and every $j\geq1$
with $\lambda_j<T$,
$$
|c_{j,A_k}^+-c_{j,A_\ell}^+|
\ll_{\M,n}
\frac{
C_{\partial\mathcal{A}}
d_H^g(\overline{A_k},\overline{A_\ell})\,
\|\varphi_j\|_{L^\infty(\M)}
}{
r^{n+1-q}\lambda_j^{n/2}
},
$$
and
$$
|c_{j,A_k}^--c_{j,A_\ell}^-|
\ll_{\M,n}
\frac{
C_{\partial\mathcal{A}}
d_H^g\bigl(
\overline{\M\setminus A_k},
\overline{\M\setminus A_\ell}
\bigr)\,
\|\varphi_j\|_{L^\infty(\M)}
}{
r^{n+1-q}\lambda_j^{n/2}
}.
$$

\item For the constant coefficients corresponding
to $\lambda_0=0$,
$$
|c_{0,A_k}^+-c_{0,A_\ell}^+|
\ll_{\M}
\frac{
d_H^g(\overline{A_k},\overline{A_\ell})
}{r},
$$
and
$$
|c_{0,A_k}^--c_{0,A_\ell}^-|
\ll_{\M}
\frac{
d_H^g\bigl(
\overline{\M\setminus A_k},
\overline{\M\setminus A_\ell}
\bigr)
}{r}.
$$
\end{enumerate}
\end{enumerate}
In both cases, the dependence on $\M$, $s$, $\beta$, and $n$ is effective.
\end{introcorollary}
The results in \cref{torus} extend the classical Beurling--Selberg estimates for intervals to higher dimensions, while \cref{Lie,Liefamily} give character expansions for conjugacy-invariant sets on compact connected Lie groups. We also consider
simply connected compact symmetric spaces with even restricted-root
multiplicities (see \cref{symeven}), with applications to odd-dimensional spheres $\S^{2d+1}$.

Motivated by the recent work of Lin, Mendes, and Radeschi on manifold
submetries and the Weyl law for singular Riemannian foliations
\cite{LinMendesRadeschi2025,LinMendesRadeschi2025Weyl}, we extend the
spectral construction beyond smooth quotient spaces. We introduce
\emph{SPR-spaces}, which are compact metric quotients arising from Riemannian
submetries with compact fibers and basic mean curvature. Commutation of
the basic projection with the Laplacian allows the construction of
majorants and minorants on the quotient, without requiring the quotient metric to be smooth. This includes quotients by regular and singular
Riemannian foliations. Explicit examples are recorded in \cref{what}, and related results can be found in
\cite{LytchakRadeschi2018}.

The resulting generalized Erd\H{o}s--Tur\'an inequalities apply both to fixed target sets and to moving families with effective control on the error term.

\begin{introproposition}[\cref{equidi}]
We retain the notation of \cref{basicmain,basicfamily}.
Fix integers $s>d+4\beta$ and $n\geq1$, and let
$(x_i)_{i\geq1}$ be a sequence in $X$.
\begin{enumerate}
\item
Let $A\subset X$ be a subset such that $\overline{\mathcal{M}}_{\mu_X}^{d-q}(\partial A)<\infty$ for $0<q\leq d$. Let $C_{\partial A}$ and $r_0$ be as in \cref{basicmain}.
For every $0<r<r_0$ and $T\geq1$,
\begin{align*}
&\left|
\frac1N\sum_{i=1}^N\ind_A(x_i)-\mu_X(A)
\right|
\ll
C_{\partial A}\left(
r^q+\frac{1}{r^{s-q}T^{(s-d-4\beta)/2}}
\right)\\
&\quad+
\frac{C_{\partial A}}{Nr^{n-q}}
\sum_{0<\lambda_k<T}
\frac{\|\psi_k\|_\infty}{\lambda_k^{n/2}}
\left|\sum_{i=1}^N\psi_k(x_i)\right|.
\end{align*}

\item
Let $Z$ be a $\mu_X$-continuity set.
For every $\varepsilon>0$, there exist
$C_{\partial Z,\varepsilon}>0$ and
$0<r_{0,\varepsilon}<1$ such that, for
$0<r<r_{0,\varepsilon}$ and $T\geq1$,
\begin{align*}
&\left|
\frac1N\sum_{i=1}^N\ind_Z(x_i)-\mu_X(Z)
\right|
\ll
\varepsilon+
C_{\partial Z,\varepsilon}\left(
r+\frac{1}{r^{s-1}T^{(s-d-4\beta)/2}}
\right)\\
&\quad+
\frac{C_{\partial Z,\varepsilon}}{Nr^{n-1}}
\sum_{0<\lambda_k<T}
\frac{\|\psi_k\|_\infty}{\lambda_k^{n/2}}
\left|\sum_{i=1}^N\psi_k(x_i)\right|.
\end{align*}
\item
Let $\mathcal A=(A_i)_{i\geq1}$ be a family of sets such that $\overline{\mathcal{M}}_{\mu_X}^{d-q}(\partial \mathcal{A})<\infty$ for some $0<q\leq d$. Let $r_0$ and $C_{\partial \mathcal{A}}$ be as in \cref{basicfamily}. For every $i\geq1$, define
$$
\delta_i:=
\max\left\{
d_H^X(\overline{A_i},\overline{A_{i+1}}),
d_H^X\bigl(
\overline{X\setminus A_i},
\overline{X\setminus A_{i+1}}
\bigr)
\right\}.
$$
For every $0<r<r_0$ and $T\geq1$,
\begin{align*}
&\left|
\frac1N\sum_{i=1}^N
\bigl(\ind_{A_i}(x_i)-\mu_X(A_i)\bigr)
\right|
\ll
C_{\partial\mathcal A}\left(
r^q+\frac{1}{r^{s-q}T^{(s-d-4\beta)/2}}
\right)\\
&\quad+
\frac{C_{\partial\mathcal A}}{Nr^{n-q}}
\left(1+\frac1r\sum_{i=1}^{N-1}\delta_i\right)
\sum_{0<\lambda_k<T}
\frac{\|\psi_k\|_\infty}{\lambda_k^{n/2}}
\max_{1\leq m\leq N}
\left|\sum_{i=1}^m\psi_k(x_i)\right|.
\end{align*}
\end{enumerate}
The implied constants depend effectively on
$X$, $\M$, $q$, $s$, $\beta$, and $n$.
\end{introproposition}

We apply these estimates to the quantitative equidistribution of CM
orbits on connected quaternionic Shimura orbifolds (see \cref{orbifold}),
obtaining an effective version of an equidistribution result of Zhang
\cite{Zhang2005}. Such equidistribution questions are closely related
to the Andr\'e--Oort conjecture.

The second part of the manuscript addresses effective motivic equidistribution theorems with explicit constants using the spectral and Riemannian tools developed earlier. The \emph{generalized Sato--Tate conjecture} predicts the
equidistribution of normalized Frobenius conjugacy classes attached to
a pure motive $M/K$ in the conjugacy-class space of its compact motivic
Sato--Tate group $G:=\st_K(M)$. If $y_{\pp}\in\mathrm{Conj}(G)$ denotes the
normalized Frobenius class at an unramified prime of good reduction,
then, for a Sato--Tate continuity set $A$ of positive measure, the expected asymptotic is
$$
\frac{\pi_A(x)}{\mathrm{Li}(x)}\longrightarrow \mu_{\st}(A).
$$
The effective Sato--Tate problem asks for an explicit error term, with
controlled dependence on the arithmetic and motivic invariants of $M$, the
representation theory of $G$, and the boundary geometry of $A$.

This problem is closely related to effective Chebotarev density. Let
$L/K$ be a finite Galois extension with Galois group $\Gamma$, and let
$C\subseteq\Gamma$ be a nonempty conjugacy class. The Chebotarev density
theorem asserts that
$$
\pi_C(x):=
\#\{\pp:N\pp\leq x,\ \mathrm{Frob}_{\pp}(L/K)\in C\}
\sim\frac{|C|}{|\Gamma|}\Li(x),
$$
where the finitely many ramified primes are excluded. The effective theory
was developed by Lagarias and Odlyzko \cite{LagariasOdlyzko1977}, and
refined by Serre \cite{Serre1981} and Murty, Murty, and Saradha
\cite{MurtyMurtySaradha1988}, with particularly strong estimates under the generalized Riemann hypothesis. A recent preprint develops the theory further under a
quasi-Riemann hypothesis
\cite{FreiSofosWidmer2026}.

The conditional effective Sato--Tate problem for non-CM elliptic curves was studied by V. K. Murty \cite{Murty1982}. Assuming the classical
analytic properties and the Riemann hypothesis for the symmetric-power
$L$-functions $L(s,\mathrm{sym}^mE)$, his method gives an error of order $x^{3/4}$ up to
logarithmic factors for fixed intervals. Bucur and Kedlaya
\cite{BucurKedlaya2016} extended this approach to a general motivic
setting. A substantial
generalization was later established by Bucur, Fit\'e, and Kedlaya
\cite{BucurFiteKedlaya2024} for abelian varieties with connected
Sato--Tate groups. For an abelian variety $B/K$ of dimension $g$ and
absolute conductor $N$, let $d$ be the rank of its Sato--Tate group $\mathrm{ST}(B)$, and let
$|\Phi^+|$ be the number of positive roots of its semisimple part. Put
$$
\varepsilon_{\mathfrak g}:=\frac{1}{2(d+|\Phi^+|)}.
$$
Under the Mumford--Tate conjecture and the generalized Riemann hypothesis for the motivic $L$-functions
attached to the irreducible characters of the Sato--Tate group,
\cite[Theorem 1.1]{BucurFiteKedlaya2024} gives, for an interval $I\subseteq[-2g,2g]$ of nonzero length
and an explicit range of $x$,
$$
\pi_{B,I}(x)
=
\mu_{\mathrm{ST}}(I)\Li(x)
+
O_{K,g}\left(
\frac{x^{1-\varepsilon_{\mathfrak g}}
(\log(Nx))^{2\varepsilon_{\mathfrak g}}}
{(\log x)^{1-4\varepsilon_{\mathfrak g}}}
\right).
$$

Further conditional and unconditional estimates for intervals and rectangular sets have been obtained by automorphic methods. Rouse
and Thorner \cite{RouseThorner2017} proved explicit conditional estimates under symmetric-power analytic hypotheses and the generalized Riemann hypothesis. Following the symmetric-power automorphy results of Newton and Thorne, Thorner \cite{Thorner2021} proved that, for a non-CM holomorphic newform $f$
of weight $k$, level $N$, and trivial nebentypus,
$$
\pi_{f,I}(x)
=
\mu_{\mathrm{ST}}(I)\pi(x)
+
O\left(
\pi(x)\frac{\log(kN\log x)}{\sqrt{\log x}}
\right),
$$
uniformly over intervals $I$, with an effectively computable implied
constant. Hoey, Iskander, Jin, and Trejos Su\'arez
\cite{HoeyEtAl2022} gave explicit numerical bounds in the
cases covered by their theorem. In the joint setting, Kumar, Kumari,
and Mishra \cite{KumarKumariMishra2026} recently extended unconditional
effective Sato--Tate estimates beyond rectangles to measurable planar
regions whose boundaries consist of finitely many continuous rectifiable
curves for pairs of twist-inequivalent non-CM newforms.
Their result does not impose an additional Chebotarev condition.

Our principal arithmetic question is an effective form of the Chebotarev--Sato--Tate phenomenon, the joint distribution of Artin
symbols and normalized motivic Frobenius classes. Given a motive $M/K$ and a finite
Galois extension $L/K$ with group $\Gamma$, write $c_\pp\in\Gamma$
for a representative of the Frobenius conjugacy class at $\pp$.
For a nonempty conjugacy class $C\subseteq\Gamma$ and a
conjugacy-invariant set $A\subseteq G$, define the counting function
\begin{equation}\label{introCST}
\pi_{A,C}(x):=\#\{\pp\notin S_{M,L}:N\pp\leq x,\ (y_\pp,c_\pp)\in A\times C\},
\end{equation}
where $S_{M,L}$ contains the exceptional primes specified below.
We ask whether, as $x\to\infty$,
$$\frac{\pi_{A,C}(x)}{\mathrm{Li}(x)}\longrightarrow \frac{|C|}{|\Gamma|}\mu_{\st}(A),$$
and seek an effective error term.

Effective results addressing this mixed problem remain relatively scarce.
The work of Wong \cite{Wong2019} studies Artin symbols together with
Frobenius angles of Hilbert modular forms, obtaining unconditional
results in certain cases and conditional effective statements under
Langlands functoriality, conductor bounds, and the automorphic Riemann hypothesis. The
results of \cite{Thorner2021,Thorner2025} and
\cite{KumarKumariMishra2026} provide further Sato--Tate estimates, but do
not directly address the question of a uniform motivic Chebotarev--Sato--Tate theorem for targets varying with the prime.

To treat this problem in a general motivic setting, we develop a framework inspired by the motivic Serre-group formalism of Banaszak
and Kedlaya \cite{BanaszakKedlaya2016,BanaszakKedlaya2023}. We introduce
the class of \emph{connected Sato--Tate (CST) motives}, consisting of
homogeneous polarized motives with connected Sato--Tate group subject
to the stated algebraic Sato--Tate and Serre conjectures and strong compatibility assumptions on the Weil--Deligne local parameters. The Tannakian formalism associates motivic lifts $M_\sigma$
to representations $\sigma\in\operatorname{Irr}(\st_K(M))$, with a
suitable extension of the coefficient field and up to Tate twists.
Combining these lifts with Artin motives $N_\rho$ for
$\rho\in\operatorname{Irr}(\Gamma)$ leads to mixed motivic
$L$-functions
$$
L(s,M_\sigma\otimes N_\rho),
\qquad (\sigma,\rho)\in
\operatorname{Irr}(\st_K(M))\times\operatorname{Irr}(\Gamma).
$$
We assume that these functions satisfy the analytic motivic formalism stated in Conjecture \ref{grh}, including meromorphic continuation, a functional equation, the motivic
Riemann hypothesis, and the motivic Ramanujan conjecture.
These hypotheses yield the character-sum estimates needed for an effective joint
equidistribution theorem.

Combining this motivic framework with our Riemannian approximation
method, we obtain effective Chebotarev--Sato--Tate estimates for a wide class of conjugacy-invariant sets. Let $M/K$ be a CST-motive and $G:=\st_K(M)$ the corresponding motivic Sato--Tate group of rank $d$. Let $T\subset G$ be a maximal torus, and let $m$ be the number of positive roots. Given a conjugacy-invariant set $A$, define $A_T$ to be its pullback to $T$ under the canonical map $T\to T/W\simeq\mathrm{Conj}(G)$. We equip $T$ with a fixed $W$-invariant flat Riemannian metric.
\begin{introdefinition}[\cref{min}]
Fix $0<q\leq d$. Define the class of sets
$$
\mathfrak M_g^q
:=
\left\{
A\subset \mathrm{Conj}(G):\,
0<\overline{\mathcal M}_g^{d-q}(\partial A_T)<\infty
\right\}.
$$
For $0<r_0<1$ and $A\in\mathfrak M_g^q$, define the constant
$$C_{\partial A}=C_\partial(A,r_0):=\max\left(1,\sup_{0<r<r_0}\frac{\s((\partial A_T)^{(r)})}{r^q}\right).$$
\end{introdefinition}

As discussed in \cref{equid}, this class includes intervals, boxes,
triangles, and sets with smooth boundaries, as well as rectifiable and
singular boundaries satisfying the upper Minkowski content condition on the maximal torus.
The finiteness condition is preserved under finite unions, intersections,
and complements.

Our first effective Chebotarev--Sato--Tate estimate concerns a fixed
target set, under the motivic and analytic hypotheses stated above.

\begin{introtheorem}[\cref{cst1}]
Let $A\in\mathfrak M_g^q$ and let $C\subset\Gamma$ be a nonempty union of
conjugacy classes. Fix $0<r_0<1$, and let
$C_{\partial A}$ be the constant
defined in \cref{min}. Put
$$
\delta_G:=d+m,
\qquad
S_{L/K}(M):=
3(1+C_H)A_{L/K}(M)d_L^{1/[L:K]^{1/2}},
$$
where $A_{L/K}(M)$ is defined in \cref{sum}, and set
$$
\pi_{A,C}(x):=
\#\left\{
\pp\notin S_{M,L}:
N\pp\leq x,\;
(y_\pp,c_\pp)\in A\times C
\right\}.
$$
Then there exists a threshold $x_0\geq e^4$ satisfying
$$
x_0\ll_G
n_K
\left(
\frac{|\Gamma|}{|C|^{1/2}}
\bigl(1+\log S_{L/K}(M)\bigr)
r_0^{-\delta_G}
\right)^4
$$
such that, for every $x>x_0$,
\begin{align*}
&\pi_{A,C}(x)
=
\frac{|C|}{|\Gamma|}\mu_{\mathrm{ST}}(A)\Li(x)
\\
&+
O_G\left(
n_K|C|^{1/2}\mu_{\mathrm{ST}}(A)x^{1/2}
\log\bigl(S_{L/K}(M)x\bigr)
\right)
\\
&+
O_G\left(
n_K C_{\partial A}
\frac{|C|^{1-q/(2\delta_G)}}{|\Gamma|^{1-q/\delta_G}}
x^{1-q/(2\delta_G)}
(\log x)^{2q/\delta_G-1}
\bigl[\log(S_{L/K}(M)x)\bigr]^{q/\delta_G}
\right).
\end{align*}
Moreover, the dependence on $G$ is effective.
\end{introtheorem}

The dependence of the boundary term on $q$ and $C_{\partial A}$
shows explicitly how the geometry of the target enters the arithmetic
error. We further develop a version of the Chebotarev--Sato--Tate theorem for prime-dependent targets. To control both boundary-tube volumes and Hausdorff variation, we introduce the
$(\alpha,q)$-Minkowski condition, which combines a uniform boundary-tube
estimate with a quantitative Hausdorff bound for the family.
\begin{introdefinition}[\cref{uni}]
Let
$\mathcal A=(A_{\mathfrak p})_{\mathfrak p\notin S_{M,L}}$
be a family of conjugacy-invariant Borel subsets of $G$
with nonempty boundary, and put
$$
\partial\mathcal A
:=
(\partial A_{\mathfrak p,T})_{\mathfrak p\notin S_{M,L}}.
$$
For $\mathfrak p,\mathfrak p'\notin S_{M,L}$, define
$$
d_*(\mathfrak p,\mathfrak p')
:=\max\Bigl\{
d_H^g\bigl(
\overline{A_{\mathfrak p,T}},
\overline{A_{\mathfrak p',T}}
\bigr),
d_H^g\bigl(
\overline{T\setminus A_{\mathfrak p,T}},
\overline{T\setminus A_{\mathfrak p',T}}
\bigr)
\Bigr\}.
$$
For $\alpha\geq0$ and $0<q\leq d$, we call $\mathcal A$ an
$(\alpha,q)$-Minkowski family if the following conditions hold:
\begin{enumerate}
\item
the uniform upper $(d-q)$-dimensional Minkowski content
of $\partial\mathcal A$ is finite:
$$
\overline{\mathcal M}_g^{d-q}(\partial\mathcal A)<\infty.
$$
In this case, for the fixed $0<r_0<1$, set
$$
C_{\partial\mathcal A}
:=
\max\left\{
1,\,
\sup_{0<r<r_0}
\sup_{\mathfrak p\notin S_{M,L}}
\frac{
\s\bigl((\partial A_{\mathfrak p,T})^{(r)}\bigr)
}{r^q}
\right\}.
$$

\item
there exists a constant $K_{\mathcal A}\geq0$ such that,
for every $x\geq2$, enumerating the prime ideals outside $S_{M,L}$ with
$x/2<N\mathfrak p\leq x$ as
$\mathfrak p_u,\ldots,\mathfrak p_v$ in nondecreasing
order of norm, with a fixed ordering among equal norms,
we have
$$
\sum_{k=u}^{v-1}d_*(\mathfrak p_k,\mathfrak p_{k+1})
\leq K_{\mathcal A}x^{-\alpha}.
$$
\end{enumerate}
\end{introdefinition}
\begin{introtheorem}
Let $\mathcal A=(A_\pp)_{\pp\notin S_{M,L}}$ be an
$(\alpha,q)$-Minkowski family, with $\alpha\geq0$ and
$0<q\leq d$. Fix $0<r_0<1$, and let
$C\subset\Gamma$ be a nonempty union of conjugacy classes.
Assume the hypotheses on the motivic $L$-functions used
in \cref{cst1}.
For $x\geq2$, define
$$
\pi_{\mathcal A,C}(x):=
\#\left\{
\pp\notin S_{M,L}:
N\pp\leq x,\;
(y_\pp,c_\pp)\in A_\pp\times C
\right\}.
$$
Put
$$
A_{\mathrm{ST}}
:=
\sup_{\pp\notin S_{M,L}}\mu_{\mathrm{ST}}(A_\pp),
\qquad
\beta_G:=
\frac q2
\min\left\{
\frac1{\delta_G},
\frac{1+2\alpha}{\delta_G+1}
\right\},
$$
and let $K_{\mathcal A}$ be as in \cref{uni}. Put
$$
B_0:=
\frac{|\Gamma|\bigl(1+\log S_{L/K}(M)\bigr)}
{|C|^{1/2}r_0^{\delta_G}}
\max\left\{1,\frac{K_{\mathcal A}}{r_0}\right\},
$$
and
$$
B_1:=
A_{\mathrm{ST}}|C|^{1/2}
+
C_{\partial\mathcal A}
(1+K_{\mathcal A})^{q/(\delta_G+1)}
\frac{|C|^{1-q/(2\delta_G)}}{|\Gamma|^{1-q/\delta_G}}.
$$
Then there exists a threshold $x_4\geq e^4$ satisfying
$$
x_4\ll_G
B_0^4\max\left\{1,\left(\frac{B_0^3}{B_1}\right)^4\right\}
$$
such that, for every $x>x_4$,
$$
\begin{aligned}
\pi_{\mathcal A,C}(x)
&=
\frac{|C|}{|\Gamma|}
\sum_{\substack{N\pp\leq x\\\pp\notin S_{M,L}}}
\mu_{\mathrm{ST}}(A_\pp)
\\
&\quad+n_KA_{\mathrm{ST}}|C|^{1/2}\mathcal{E}_1(x)+n_KL_{\mathcal{A}}
\frac{|C|^{1-q/(2\delta_G)}}{|\Gamma|^{1-q/\delta_G}}\mathcal{E}_2(x),
\end{aligned}
$$
with
\begin{align*}
L_{\mathcal{A}}&:=C_{\partial\mathcal A}
(1+K_{\mathcal A})^{q/(\delta_G+1)},\\
    \mathcal{E}_1(x)&:=O_G\left(
x^{1/2}
\log\bigl(S_{L/K}(M)x\bigr)\right),\\
\mathcal{E}_2(x)&:=O_G\left(
x^{1-\beta_G}
(\log x)^{2q/\delta_G-1}
\bigl[\log(S_{L/K}(M)x)\bigr]^{q/\delta_G}
\right).
\end{align*}
The dependence of the implied constants on $G$ is effective.
\end{introtheorem}

In particular, the error exponent
$$
\beta_G=\frac q2\min\left\{
\frac1{\delta_G},\frac{1+2\alpha}{\delta_G+1}
\right\}
$$
depends on the uniform Minkowski codimension $q$ and the
exponent $\alpha$ governing Hausdorff variation. Inspired  in part by Serre's treatment of
effective Chebotarev problems \cite{Serre1981}, we also discuss in \cref{tserre} of \cref{equid}
 how additional symmetries of the Sato--Tate set may
improve the resulting exponents.

We conclude with arithmetic applications illustrating the scope of the
method. We refer to \cref{sato2} for the corresponding definitions and hypotheses. For products of rank-two motives with Sato--Tate group
$\mathrm{SU}(2)^n$, the moving-family theorem gives a motivic analogue of
the Atkin--Serre problem with an additional Chebotarev restriction. It
counts primes at which at least one normalized Frobenius trace is unusually
small. In the range where the boundary error is dominated by the shrinking
main term, we obtain a conditional asymptotic equivalence, together with a bound for the least such prime; compare \cite{Rouse2007,GafniThornerWong2021}.

Let $M=\bigoplus_{i=1}^n M_i$ be a CST-motive, where the $M_i$
are rank-two polarized motives of common weight $w$, and assume we have an isomorphism of represented groups
$$
(\mathrm{ST}_K(M),H_B(M)\otimes_{\Q}\mathbb C)
\simeq
\left(\mathrm{SU}(2)^n,
\bigoplus_{i=1}^n\mathbb C^2_{\mathrm{std},i}\right).
$$
We consider the following generalized Atkin--Serre-type counting problem:
$$
\pi_{\mathrm{AS},C}(x,\varepsilon)
:=
\#\left\{
\begin{gathered}
\pp\notin S_{M,L},\quad
N\pp\leq x,\quad c_\pp\in C,
\\
|a_i(\pp)|<2(N\pp)^{w/2-\varepsilon}
\\
\text{for some }1\leq i\leq n
\end{gathered}
\right\}.
$$
Under the preceding hypotheses, applying our result for moving sets
gives the following proposition.
\begin{introproposition}[\cref{Serre}]
Put
$$
\beta_{\varepsilon}
:=
\frac12\min\left\{
\frac1{2n},\frac{1+2\varepsilon}{2n+1}
\right\}.
$$
\begin{enumerate}
\item
For every $0<\varepsilon<1$, there exists an effectively
computable $x_0>0$ such that, for $x>x_0$,
\begin{align*}
    &\pi_{AS,C}(x,\varepsilon)
=
\frac{4n}{\pi}\frac{|C|}{|\Gamma|}
\sum_{\substack{\pp\notin S_{M,L}\\N\pp\le x}}
\frac1{(N\pp)^\varepsilon}
+
O_n\left(
\frac{|C|}{|\Gamma|}
\sum_{\substack{\pp\notin S_{M,L}\\N\pp\le x}}
\frac1{(N\pp)^{2\varepsilon}}
\right)
\\
&+
O_n\left(
n_K|C|^{1/2}x^{1/2}
\log\bigl(S_{L/K}(M)x\bigr)
\right)
\\[2pt]
&+
O_n\left(
n_K
\frac{|C|^{1-1/(4n)}}{|\Gamma|^{1-1/(2n)}}
x^{1-\beta_{\varepsilon}}
(\log x)^{1/n-1}
\bigl(\log(S_{L/K}(M)x)\bigr)^{1/(2n)}
\right).
\end{align*}
In particular, for fixed
$
0<\varepsilon<\frac1{4n},
$
as $x\rightarrow\infty$ we have
$$
\pi_{AS,C}(x,\varepsilon)
\sim
\frac{|C|}{|\Gamma|}
\frac{4n}{\pi(1-\varepsilon)}
\frac{x^{1-\varepsilon}}{\log x}.
$$

\item
For every $0<\varepsilon\le1/(8n)$, there exists a prime
ideal $\mathfrak p\notin S_{M,L}$ such that
$c_{\mathfrak p}\in C$ and
$$
|a_i(\mathfrak p)|
<
2(N\mathfrak p)^{w/2-\varepsilon}
$$
for some $1\le i\le n$.
Moreover, one may choose $\mathfrak p$ with
$$
\begin{aligned}
N\mathfrak p\ll_n\max\Biggl\{
&x_5,\,
\bigl(1+n_K+|S_{M,L}|+\log d_K\bigr)^8,
\\
&
\left(\frac{n_K|\Gamma|}{|C|^{1/2}}\right)^{8n+4}
\left[
\log\left(
2S_{L/K}(M)\frac{n_K|\Gamma|}{|C|^{1/2}}
\right)
\right]^{16n+8}
\Biggr\},
\end{aligned}
$$
where $x_5=x_5(\mathcal A,C,\varepsilon)>0$ is an
effectively computable threshold.
\end{enumerate}
The implied constants depend only on $n$ and are effective.
\end{introproposition}

To illustrate the flexibility of the motivic framework beyond the
classical examples of abelian varieties, we also consider transcendental
motives of K3 surfaces satisfying our CST hypotheses. Let $X/K$ be a K3 surface. The geometry of its exceptional Picard-rank
jumps leads to a cyclotomic Frobenius locus $J$, see
\cite{ShankarShankarTangTayou2022,Tayou2025,Huang2026} for related arithmetic questions. Put $M:=t_2(X)(1)$ and $G:=\st_K(M)$. At a prime $\pp$ of good
reduction with residue field $k_\pp$, the Picard-jump condition
$\rho(X_{\overline{k_\pp}})>\rho(X_{\overline K})$ is equivalent to
$y_\pp\in J$, where $J$ is the possibly singular analytic set
described below. Assuming that $J$ is a proper subset of $G$, and
under the stated hypotheses for $M$, our Chebotarev--Sato--Tate
estimate gives the following count for the tubular neighborhoods
$J^{(\varepsilon)}:=\{x\in G:d_g(x,J)<\varepsilon\}$, where $g$ is
a fixed bi-invariant metric on $G$. In particular, it yields an
upper bound for Picard jumps.

\begin{introproposition}[\cref{jump}]
Let $S_0$ be the places of bad reduction.
Define
$$
\pi_{J,C}(x,\varepsilon)
:=
\#\{\pp\notin S_{M,L}\cup S_0:
N\pp\le x,\,
y_\pp\in J^{(\varepsilon)},\,
c_\pp\in C\},
$$
with $J^{(0)}:=J$.
There exists an effectively computable $x_0>0$,
independent of $\varepsilon$, such that, for $x>x_0$
and $0\le\varepsilon<\varepsilon_0$,
\begin{align*}
&\pi_{J,C}(x,\varepsilon)
=
2\frac{|C|}{|\Gamma|}\mathcal H_g^{\dim G-1}(J)\varepsilon\Li(x)
+O(|S_0|)+
O_{G,J}\left(
\frac{|C|}{|\Gamma|}\varepsilon^2\Li(x)
\right)
\\
&+
O_{G,J}\left(
n_K|C|^{1/2}\varepsilon x^{1/2}
\log\bigl(S_{L/K}(M)x\bigr)
\right)
\\
&+
O_{G,J}\left(
n_K
\frac{|C|^{1-1/(2\delta_G)}}
     {|\Gamma|^{1-1/\delta_G}}
x^{1-1/(2\delta_G)}
(\log x)^{2/\delta_G-1}
\bigl(\log(S_{L/K}(M)x)\bigr)^{1/\delta_G}
\right).
\end{align*}

In particular, $\pi_{J,C}(x,0)$ counts the Picard jumps
with $c_\pp\in C$.

If $\mathcal H_g^{\dim G-1}(J)>0$ and $\varepsilon=\varepsilon(x)>0$ satisfies
$$
\varepsilon(x)\longrightarrow0,
\qquad
\frac{
n_K(\log x)^{2/\delta_G}
\bigl(\log(S_{L/K}(M)x)\bigr)^{1/\delta_G}
}{
\varepsilon(x)x^{1/(2\delta_G)}
}
\longrightarrow 0,
$$
then
$$
\pi_{J,C}(x,\varepsilon(x))
\sim
2\frac{|C|}{|\Gamma|}
\mathcal H_g^{\dim G-1}(J)\varepsilon(x)\Li(x).
$$
\end{introproposition}

The proofs combine established methods from effective Chebotarev theory,
character approximation, and motivic $L$-functions, including ideas of
Serre, Zywina, Thorner, and Bucur, Fit\'e, and Kedlaya
\cite{Serre1981,ZywinaLangTrotter,Thorner2021,BucurFiteKedlaya2024}, with
regularized Riemannian distances, geometric measure theory, and spectral
analysis on compact manifolds and their singular quotients developed
throughout the manuscript. Our primary
aim is a coherent effective theory under the Riemann hypothesis for the
specified family of motivic $L$-functions. More specialized automorphic
methods can yield unconditional results, and \cref{unc} briefly discusses
this possibility.

Throughout the manuscript, we raise several questions concerning disconnected Sato-Tate groups and more general motivic settings. We leave these questions, together with problems concerning optimal spectral decay, singular quotients, and moving geometric families, as directions for future research.

\section*{Acknowledgment}
This work is the output of several months of research carried out during the author's first year of doctoral studies at Ibn Tofail University, Morocco. Although the ideas and proofs presented here are original to the author, this work would not have been possible without significant human support that helped provide a stable research environment.

The author would especially like to thank his mother for her constant attention, encouragement, and support, and Professors Youssef Ouknine and Mohamed Erraoui for their academic, institutional, and moral support.

The author is grateful to Professor Mohammed Boucetta for his useful comments on the Riemannian geometry aspects of this work, and to Gaëtan Chenevier and Bruno Kahn for their helpful remarks concerning the theory of motivic $L$-functions.

The author acknowledges with gratitude the Hassan II Academy of Sciences and Technology for its generous funding of his doctoral studies and its institutional support. Special thanks are due to the permanent secretary, Omar Fassi Fehri, as well as to Mr. Abdessattar El Amrani Jamal and Mr. Mohamed Slaoui.

Finally, the author thanks his former professor Alberto Mínguez for introducing him to number theory and the Langlands program, which helped shape his research direction.

\section*{Organization of the manuscript}
We start by fixing our notation and conventions in \cref{notation}. In \cref{spectral}, we develop the spectral Beurling--Selberg approximation, including its Hausdorff stability and a uniform version for families, and extend it to Lie groups and certain simply connected symmetric spaces. We then derive generalized Erd\H{o}s--Tur\'an inequalities and apply them to the equidistribution of CM points on certain quaternionic Shimura varieties. In \cref{class}, we introduce the class of CST-Grothendieck motives and formulate the motivic and analytic hypotheses, with detailed discussions of how to verify them and their current status in the literature. In \cref{equid}, we exploit our Riemannian tools to establish conditional effective Chebotarev-Sato-Tate estimates for fixed sets and families of sets under explicit geometric conditions, and discuss improvements to the error exponents and extensions beyond the CST hypotheses. In \cref{sato2}, we give arithmetic applications, including a random matrix application, motivic Atkin-Serre estimates, Serre-type least-prime problems, Picard-rank jumps, and exceptional Frobenius spectra of Picard curves. The appendices introduce the technical machinery and geometric tools used throughout the manuscript, together with examples of computations of certain arithmetic invariants appearing in our estimates.

\chapter{Notation and conventions}\label{notation}

\section{Euclidean spaces}

Fix a Euclidean space $\R^d$ with $d\geq1$.
We denote by $\N$ the set of nonnegative integers.
A multi-index is a tuple
$\alpha=(\alpha_1,\ldots,\alpha_d)\in\N^d$.
We write $|\alpha|:=\alpha_1+\cdots+\alpha_d$
and define the partial differential operator
$$
D^\alpha
:=
\partial_{x_1}^{\alpha_1}\cdots
\partial_{x_d}^{\alpha_d},
$$
where $\partial_{x_i}$ denotes differentiation with
respect to the $i$-th coordinate.
These notions extend naturally to the flat torus
$\bb{R}^d/\bb{Z}^d$ through its Euclidean structure.

For $s\geq0$, define
$$
\omega_s
:=
\frac{\pi^{s/2}}
     {\Gamma\left(\frac s2+1\right)}.
$$
In particular, if $s=d$ is a nonnegative integer, $\omega_d$ is the Lebesgue volume of the unit ball in $\R^d$.
\section{Metric measure spaces and Hausdorff distance}
Given a metric measure space $(X,d_X,\mu)$, where $(X,d_X)$ is a metric space and $\mu$ a Radon measure, define the distance to a set 
$$d_X(x,A):=\inf\{d_X(x,y)\,,y\in A\}.$$
We adopt the convention $d_X(x,\varnothing)=+\infty$.
Given $r>0$, define the tube
$$A^{(\r)}:=\{x\in X,d_X(x,A)<\r\}.$$

Given a subset $A\subset X$, define the topological boundary $\partial A=\overline{A}\setminus \mathrm{int}(A)$.
Let $\mathcal K(X)$ denote the space of non-empty compact subsets of
$X$. For $K_1,K_2\in\mathcal K(X)$, define the Hausdorff metric:

$$
\begin{aligned}
d_H^X(K_1,K_2)
&:=\sup_{x\in X}
\bigl|d_X(x,K_1)-d_X(x,K_2)\bigr|
\\
&=\max\left\{
\sup_{x\in K_1}d_X(x,K_2),\,
\sup_{x\in K_2}d_X(x,K_1)
\right\}.
\end{aligned}
$$

Equivalently,
$$
    d^X_H(K_1,K_2)
    =
    \inf\left\{
        {\varepsilon}>0:
        K_1\subset K_2^{({\varepsilon})}
        \ \text{and}\
        K_2\subset K_1^{({\varepsilon})}
    \right\}.
$$
The spaces $C^0(X)$ and $(L^p(X))_{1\leq p\leq \infty}$ carry their usual meanings.
If $X$ is endowed with a differentiable structure, then $C^r(X)$ and $C^\infty(X)$ denote the usual spaces of $C^r$ and smooth complex-valued functions on $X$, respectively.

\section{Riemannian geometry}
Throughout this article $(\M,g)$ denotes a smooth Riemannian manifold of dimension $d\geq 1$. Unless specified otherwise, it is assumed to be compact, connected, and without boundary. Denote by $\s$ the Riemannian volume measure. In local coordinates,
$$
    d\s
    =
    \sqrt{\det(g_{ij})}\,dx_1\cdots dx_d.
$$
Define the Riemannian distance $d_g:\M\times \M\to \R_+$ by 
$$d_g(p,q):=\inf_{\gamma}\int_{0}^1|\gamma'(t)|_g \,dt,$$
where the infimum is taken along all piecewise $C^1$-curves $\gamma:[0,1]\to \M$ such that $\gamma(0)=p$ and $\gamma(1)=q$. The existence of a minimizing geodesic in the compact case is guaranteed by Hopf--Rinow theorem \cite[Theorem 6.19 and Corollary 6.21]{Lee2018}. Without loss of generality, we normalize $g$ so that $\s(\M)=1$.

Define the diameter $$D_{\mathbb M}:=\sup_{x,y\in \M}d_g(x,y)<\infty.$$
For $x\in \M$ and $r>0$, we denote by $B_g(x,r)$ and $S_g(x,r)$,
respectively, the geodesic ball and geodesic sphere of radius $r$
centered at $x$.
The metric $g$ induces a pointwise pairing on covariant
$k$-tensors. For a $g$-orthonormal basis $(e_i)_{1\leq i\leq d}$
of $T_x\M$, we set
$$
\langle S,T\rangle_g(x)
:=
\sum_{i_1,\ldots,i_k=1}^d
S_x(e_{i_1},\ldots,e_{i_k})
T_x(e_{i_1},\ldots,e_{i_k}).
$$
This definition is independent of the orthonormal basis.
For complex-valued tensors, we extend the pairing
$\C$-bilinearly and define
$$
|S|_g:=\langle S,\overline S\rangle_g^{1/2}.
$$
In particular, $|\nabla^k f|_g$ denotes the norm
of the covariant $k$-tensor $\nabla^k f$.
Let $r_{\mathrm{inj}}$ be the global injectivity radius. Denote  by $\operatorname{Ric}_g$ the Ricci curvature $(0,2)$-tensor. Define the curvature constant
\begin{equation}
  \kappa
  :=
  \max\left\{
    0,\,
    -\frac{1}{d-1}
    \inf_{\substack{x\in \X\\|v|_g=1}}
    \operatorname{Ric}_g(v,v)
  \right\}.
\end{equation}
For $d=1$ we put $\kappa=0$.
As $\M$ is compact, we have $\kappa<\infty$ and $r_{\mathrm{inj}}>0$.
For a two-dimensional subspace
$\sigma=\operatorname{span}\{u,v\}\subseteq T_x\M$, denote $\operatorname{Sec}_g(\sigma)$ its sectional curvature. Define

\begin{equation}
K
  :=
  \sup_{\substack{x\in \X\\
                  \sigma\subseteq T_x\X,\ \dim\sigma=2}}
  \left|\operatorname{Sec}_g(\sigma)\right|<\infty.
\end{equation}
Denote by $\nabla$ the Levi-Civita connection of $g$ and by $\operatorname{div}_g$ the associated divergence operator. Let $\Delta_g:=-\operatorname{div}_g\circ\nabla:C^\infty(\X)\rightarrow C^\infty(\X)$ be the Laplace-Beltrami operator. It extends to a self-adjoint operator on the $L^2$-Sobolev space $\Delta_g:H^2(\M)\subset L^2(\X)\rightarrow L^2(\X)$ with discrete spectrum:
$$
    0=\lambda_0<\lambda_1\leq\lambda_2\leq\cdots,
$$
where the eigenvalues are repeated according to multiplicity.  We fix
an orthonormal basis $(\varphi_k)_{k\geq0}$ of $L^2(\X)$ consisting
of smooth eigenfunctions satisfying $\varphi_0=1$ and 
$
    \Delta_g\varphi_k=\lambda_k\varphi_k
$ for every $k\geq 0$.
For every $f\in L^2(\X)$, we have
$
    f
    =
    \sum_{k\geq0}\widehat f(k)\varphi_k
$
where
$$
    \widehat f(k)
    :=
    \int_{\X}f\,\overline{\varphi_k}\,d\s
$$
is the $k$-th Fourier coefficient of $f$.
\section{Submetries and foliations}
A Riemannian submetry $\sigma:(\M,g)\longrightarrow (X,d_X)$ is a surjective map of metric spaces such that 
    \begin{enumerate}
        \item for $p\in \M$ and $r>0$ we have $\sigma(B_g(p,r))= B_X(\sigma(p),r)$,
        \item for each $q\in X$ the fiber $\sigma^{-1}(q)$ is a smooth embedded manifold.
    \end{enumerate}
    We denote $\M_{\mathrm{reg}}$ the union of fibers of maximal dimension. 
    The submetry $\sigma$ is said to have basic mean curvature if for $p,q\in \M_{\mathrm{reg}}$ with $\sigma(p)=\sigma(q)=x\in X$, we have $d_p\sigma(H_p)=d_q\sigma(H_q)$, where $H$ denotes the mean curvature field of the fiber $\sigma^{-1}(x)\subset \M$.
    Given a bounded Borel function $f$ and $y\in \M$ lying in the fiber $L_p:=\sigma^{-1}(p)$ define the pointwise fiber average:
    $$P_{\sigma}f(y)=\frac{1}{\mathrm{vol}(L_p)}\int_{L^p}fd\mathrm{vol}_{g|L_p}.$$
    Here $\mathrm{vol}_{g|L_p}$ denotes the induced Riemannian measure on $L_p$ and $\mathrm{vol}(L_p)=\int_{L_p}d\mathrm{vol}_{g|L_p}$.
    A subset $A\subset \M$ is said to be $\sigma$-saturated if for any fiber $L_p$ we have $L_p\subset A$ or $L_p\cap A=\varnothing$.
    A function $f\in L^\infty(\X)$ is said to be {basic} if $f$ is constant on fibers.
For $0\leq r\leq \infty$ and $1\leq p\leq \infty$ denote by $C^r(\M,\sigma)\subset C^r(\M)$ (resp. $L^p(\M,\sigma)\subset L^p(\M)$) the subspace consisting of basic $C^r$-functions (resp. functions with a basic representative). We use the same notation for the extension
$P_\sigma:L^2(\M)\to L^2(\M,\sigma)$, which is the orthogonal
projection into the basic subspace, or equivalently, the conditional expectation with respect to the basic $\sigma$-Borel algebra generated by $\sigma$-saturated sets.
We equip $X$ with the pushforward measure $\mu_X:={\sigma_*\s}$. 
\section{Minkowski contents and the Hausdorff measure}
Let $(X,d_X,\mu_X)$ be a compact metric probability space. Fix a reference dimension $d$ and let $0\leq s\leq d$.
The lower and upper $s$-dimensional Minkowski contents
of $A\subset X$ are defined by
$$
\underline \J_{\mu_X}^{\,s}(A)
:=\liminf_{r\to0^+}
\frac{\mu_X(A^{(r)})}{\omega_{d-s}r^{d-s}},
\qquad
\overline \J_{\mu_X}^{\,s}(A)
:=\limsup_{r\to0^+}
\frac{\mu_X(A^{(r)})}{\omega_{d-s}r^{d-s}}.
$$
When these quantities agree, we call $A$
$s$-Minkowski-measurable and denote their common value
by $\J_{\mu_X}^{s}(A)$.
For a family $\mathcal A=(A_i)_{i\in I}$, define its
uniform upper and lower Minkowski contents by
$$
\overline \J_{\mu_X}^{\,s}(\mathcal A)
:=\limsup_{r\to0^+}\sup_{i\in I}
\frac{\mu_X(A_i^{(r)})}{\omega_{d-s}r^{d-s}},
$$
and
$$
\underline \J_{\mu_X}^{\,s}(\mathcal A)
:=\liminf_{r\to0^+}\inf_{i\in I}
\frac{\mu_X(A_i^{(r)})}{\omega_{d-s}r^{d-s}}.
$$
For $s>0$, the $s$-dimensional Hausdorff measure of $A$
is defined by
$$
\mathcal H_{d_X}^s(A)
:=
\lim_{\varepsilon\to0^+}
\inf\left\{
\sum_{j=1}^{\infty}
\omega_s\left(\frac{\operatorname{diam}_{d_X}U_j}{2}\right)^s:
A\subset\bigcup_{j=1}^{\infty}U_j,\;
\operatorname{diam}_{d_X}U_j\leq\varepsilon
\right\}.
$$
Define the Hausdorff dimension of $A$ as
$$
\dim_{\mathrm H}A
:=
\inf\{s>0:\mathcal H_{d_X}^s(A)=0\}.
$$
When $d_X$ is induced by a Riemannian metric $g$ and $\mu_X=\s$ is the associated volume measure, 
we write $\overline{\mathcal{M}}_g^{s}$, $  \underline{\mathcal{M}}_g^{s}$, and $\mathcal H_g^s$ respectively.
\section{Lie groups and symmetric spaces}
Let $G$ be a compact connected Lie group of
rank $d$. By \cite[Corollary 3.15]{Lee2018} $G$ admits a bi-invariant Riemannian metric such that the associated volume measure is the normalized Haar measure $\mu_G$.  Denote
$\mathrm{Conj}(G)=G/\operatorname{Ad}(G)$ the space of conjugacy classes,
and let $\pi:G\to\mathrm{Conj}(G)$ be the quotient map. Define the Sato-Tate measure $\mu_{\st}$ to be the pushforward $\pi_*\mu_{G}$ along the map $\pi$.

Fix a maximal torus $T\subset G$, with Lie algebra $\mathfrak t$. The exponential map yields a Riemannian diffeomorphism $\mathfrak t/\ker(\exp_G|_{\mathfrak t})\simeq T$. An orthonormal
identification $\mathfrak t\simeq\mathbb R^d$ induces a Riemannian
isometry
$
    \iota:T\simeq 
    {\bb T}_{\L}:=\mathbb R^d/\Lambda.
$
We use this isometry to identify the character group of $T$ with
the dual lattice
$$
    X^*(T)\simeq\Lambda^*
    :=\bigl\{\lambda\in\mathbb R^d:
              \lambda\cdot\ell\in\mathbb Z
              \text{ for every }\ell\in\Lambda\bigr\}.
$$
Thus, for $\theta=[x]\in{\bb T}_{\L}$ and
$\lambda\in X^*(T)$, the pairing
$\langle\lambda,\theta\rangle=\lambda\cdot x\pmod{\mathbb Z}$
is well defined, and the corresponding character is
$e^{2\pi i\langle\lambda,\theta\rangle}$. 
Through the above orthonormal identification, we endow $X^*(T)\subset\mathfrak t^*$ with the
Euclidean norm, denoted by $\|\cdot\|$. We denote by $X^*(T)^+$ the set of dominant weights,
equipped with the standard dominance order.
Let $\Phi:=\Phi(G,T)\subset X^*(T)$ be the root system of
$G$ relative to $T$, and let $W:=W(G,T)$ be the corresponding
Weyl group. Fix a system $\Phi^+\subset\Phi$ of positive
roots and put
$$
\rhoo:=\frac12\sum_{\alpha\in\Phi^+}\alpha.
$$
We denote by $\operatorname{Rep}(G)$ the category of
finite-dimensional continuous complex representations of $G$,
and by $\operatorname{Irr}(G)$ the set of isomorphism classes
of its irreducible representations. Given an irreducible representation $\sigma\in\operatorname{Irr}(G)$, write
$d_\sigma$ for its dimension. By the highest-weight theorem
\cite[Theorem 5.110, p. 335]{Knapp2002} continuous irreducible representations of $G$ are classified by their highest dominant weight $\lambda\in X^*(T)^+$. Given a dominant weight $\lambda\in X^*(T)^+$, let $\chi_\lambda$ be the character of the irreducible representation of highest weight $\lambda$. 

For a compact Riemannian symmetric space $X=G/K$, fix
a maximal flat torus $F$.
Our preceding root-theoretic notation extends naturally,
with $T$ replaced by $F$, $\Phi$ by the restricted root
system $\Sigma$, and $W$ by the restricted Weyl group $W_K$.
For $\alpha\in\Sigma$, we let $m_\alpha$ denote its
multiplicity. After fixing a system $\Sigma^+$ of positive
restricted roots we put
$$
\rhoo_F:=\frac12\sum_{\alpha\in\Sigma^+}m_\alpha\alpha.
$$
We write $X^*(F)^+$
for the set of dominant spherical highest weights of
$(G,K)$, equipped with the dominance order determined
by $\Sigma^+$. For each dominant spherical weight $\lambda\in X^*(F)^+$,
let $\phi_\lambda$ denote the associated spherical function,
normalized by $\phi_\lambda(eK)=1$.
\section{Galois representations}
Let $K$ be a number field. We write $n_K:=[K:\Q]$ and let $d_K$ be its absolute
discriminant. We denote $G_K:=\operatorname{Gal}(\overline K/K)$ its absolute Galois group equipped with the profinite topology. An Artin representation is a continuous finite-dimensional
complex representation
$$
\rho:G_K\longrightarrow\GL(V_\rho).
$$
Such a representation has finite image, by \cite[Section 12.3, Theorem 24 and its corollary]{Serre1977}. Every Artin representation factors through a finite Galois extension $L/K$ and admits a model
over a number field $E$. If such a model exists for an Artin representation $\rho$, we say that $\rho$ is $E$-linear. For an Artin representation $\rho$ of degree $d_\rho$, let
$\mathfrak f(\rho)$ denote its Artin conductor ideal and put
$q_\rho:=d_K^{d_\rho}N_{K/\Q}\mathfrak f(\rho)$ its Artin conductor. For a finite Galois extension $L/K$, we denote by $\mathrm{Irr}(\mathrm{Gal}(L/K))$ the finite dimensional irreducible representations of $\mathrm{Gal}(L/K)$ which we identify as irreducible Artin representations factoring through the quotient by $\mathrm{Gal}(\overline{K}/L)$.
\section{Motivic categories}
We refer the reader to \cite[Section 12]{BanaszakKedlaya2023},
\cite[pp. 447-448]{Jannsen1992}, and
\cite[Sections 2, 4]{Andre1996} for details on the
construction of the categories of motives.
Let $K,E\subset\bb C$ be number fields, and denote by
$\operatorname{SmProj}_K$ the category of smooth projective
varieties over $K$.
Fix a theory $\sim$ of correspondences with coefficients
in $E$ between smooth projective varieties over $K$.
Let $\mathcal M_\sim(K,E)$ be the associated
$E$-linear, rigid, symmetric monoidal, pseudo-abelian
category of pure motives.
Its unit object is
$
\mathbf 1:=
\bigl(\operatorname{Spec}K,\Delta_{\operatorname{Spec}K},0\bigr).
$

Let
$$
\mathbb L:=h^2(\mathbb P^1)=\mathbf 1(-1)
$$
be the Lefschetz motive.
We adopt the following convention for Tate twists:
$
M(m):=M\otimes\mathbb L^{\otimes(-m)},
$
so that if $M$ is of pure weight $w$, then $M(m)$
has weight $w-2m$.
The construction above yields a contravariant functor:
\begin{align*}
h:\operatorname{SmProj}_K^{\mathrm{op}}
&\longrightarrow\mathcal M_\sim(K,E),\\
X&\longmapsto(X,\Delta_X,0),\\
f&\longmapsto{}^t\Gamma_f.
\end{align*}

Given a motive $M\in\mathcal M_\sim(K,E)$, we denote
by $H_B(M)$ its Betti realization and by $H_\lambda(M)$
its $\lambda$-adic realization.
For a field embedding $\iota:E_\lambda\hookrightarrow\bb C$,
we identify
$$
H_B(M)\otimes_E E_\lambda
\simeq H_\lambda(M)
$$
under the canonical Betti-\'etale comparison.
For a finite extension $E'/E$ of coefficient fields, put
$M_{E'}:=M\otimes_E E'\in\mathcal M_\sim(K,E')$.
For a finite extension $L/K$, denote by
$M\times_K L\in\mathcal M_\sim(L,E)$ its base change in the ground field so that
$$
H_B(M_{E'})\simeq H_B(M)\otimes_E E',
\qquad
H_{\lambda'}(M_{E'})
\simeq H_\lambda(M)\otimes_{E_\lambda}E'_{\lambda'},
$$
for every finite place $\lambda$ of $E$ and every
place $\lambda'\mid\lambda$ of $E'$.
Following the conventions of Deligne and Milne \cite{DeligneMilne1982}
we denote by $\mathcal M_{\mathrm{Art}}(K,E)$ the
category of Artin motives over $K$ with coefficients in $E$.
Given an $E$-linear Artin representation
$\rho:G_K\rightarrow\GL_E(V)$, we denote by
$N_\rho\in\mathcal M_{\mathrm{Art}}(K,E)
\subset\mathcal M_\sim(K,E)$
the pure motive of weight $0$ associated with $(V,\rho)$.
If $\mathcal M_\sim(K,E)$ is either the category
$\mathcal M_{\mathrm{ahc}}(K,E)$ of motives of absolute
Hodge cycles or Andr\'e's category
$\mathcal M_{\mathrm{And}}(K,E)$ of motivated motives,
we let $\mathcal M_{\mathrm{Ab}}(K,E)$ be the smallest
Tannakian subcategory generated by
$\mathcal M_{\mathrm{Art}}(K,E)$ and $h(A)$
for all abelian varieties $A/K$.

\section{Landau notation}
We use the Landau notations $O$, $\ll$ and $\gg$.
The relations $f=O_{\mathcal D}(g)$ and
$f\ll_{\mathcal D}g$ mean that
$|f|\leq C_{\mathcal D}g$ for some constant
$C_{\mathcal D}>0$, where $g\geq0$.
The subscripts indicate that the implied constants
depend only on the mathematical data $\mathcal D$.
The notation $f\gg_{\mathcal D}g$ is defined analogously.
We write $f\asymp_{\mathcal D}g$ if both
$f\gg_{\mathcal D}g$ and $f\ll_{\mathcal D}g$ hold.
We say that a constant depends effectively on
$\mathcal D$ if an explicit upper bound in terms of
$\mathcal D$ can be obtained.

\chapter{Spectral Beurling-Selberg approximation}\label{spectral}
The classical Beurling-Selberg question concerns the pointwise approximation of the characteristic function of an interval by majorants and minorants with finite Fourier support as in \cite{Vaaler1985}.

On compact Riemannian manifolds, the corresponding spectral approximation problem replaces Fourier's exponentials by Laplace-Beltrami eigenfunctions as in \cite[Section 3]{ColzaniGiganteTravaglini2011}. 

\section{Construction of spectral approximations on Riemannian manifolds}
Fix a compact connected Riemannian manifold $(\M,g)$. In this section, we construct spectral pointwise majorants and minorants for the indicator function of sets whose boundaries have finite Minkowski content with explicit control on the $L^1$-gap and explicit bound on the Fourier coefficients in terms of the spectral data of the ambient manifold. We also study their stability under perturbations of
the underlying set measured by a Hausdorff distance, which allows us to formulate an analogue for a family of sets with uniform Minkowski content.

Smooth approximation of Lipschitz functions on Riemannian manifolds was studied systematically by Azagra, Ferrera, L\'opez-Mesas, and Rangel \cite{AzagraEtAl2007}. More recently, Cao and Xiao
\cite{CaoXiao2026} constructed regularized distance functions for star-shaped domains. Our goal is to construct a Riemannian analogue of Fraenkel's regularized distance \cite{Fraenkel1979} on a compact Riemannian manifold such that the assignment $A\rightarrow \rho_{A}$ is compatible with the Hausdorff metric. We refer the reader to \cref{reg} for a review of the classical Euclidean Fraenkel's regularized distance
and the technical tools used in the proof. Although our construction does not generally approximate the Riemannian distance $d_g(\cdot,F)$ with arbitrarily small error,
its stability under Hausdorff perturbations provides
a suitable framework for our arithmetic applications
to equidistribution on moving sets.
\begin{lemma}\label{regularized}
There exists a family $(\rho_F)_{F\in \K}$, with $\rho_F\in C^\infty(\M\setminus F)$, such that for every $F\in \K$, $k\geq 1$ and $x\in \M\setminus F$, we have 
\begin{enumerate}
    \item $\rho_F(x)\asymp_{\M} d_g(x,F)$,
    \item $|\nabla^k\rho_F(x)|_g
\ll_{\M,k}d_g(x,F)^{1-k}$.
\end{enumerate}
Moreover, for every $F,F'\in \K$, $k\geq 0$ and $x\in \M\setminus (F\cup F')$, we have
$$
\bigl|\nabla^k\rho_F(x)-\nabla^k\rho_{F'}(x)\bigr|_g
\ll_{\M,k}
\frac{d_H^g(F,F')}
{\min\{d_g(x,F),d_g(x,F')\}^{k}}.
$$
\end{lemma}
\begin{proof}
Choose
$
0<r<
\min\left\{
D_{\mathbb M},\frac{r_{\mathrm{inj}}}{2},
\frac{1}{\sqrt{2K}}
\right\},
$
with $1/\sqrt{2K}:=\infty$ when $K=0$.
Choose a finite family of exponential charts.
Let $(B_g(p_i,r))_{1\leq i\leq N}$ be a finite cover of $M$, and
choose a fixed smooth partition of unity $(\chi_i)_{i=1}^N$
such that
$$
\operatorname{supp}\chi_i\subset B_g(p_i,r),
\qquad
\chi_i\geq0,
\qquad
\sum_{i=1}^N\chi_i=1.
$$
For each $i$, fix a linear isometry
$B_i:\mathbb R^d\to T_{p_i}\M$, and set
$$
\psi_i:=\exp_{p_i}\circ B_i:
B_{\mathbb R^d}(0,2r)\longrightarrow B_g(p_i,2r).
$$
For a non-empty proper closed set $F\subset M$, define
$$
E_i(F):=
\psi^{-1}_i\bigl(F\cap B_g(p_i,2r)\bigr)
\cup
\bigl(\mathbb R^d\setminus B_{\mathbb R^d}(0,2r)\bigr).
$$
Let $\rho_{E_i(F)}:=\rho_{E_i(F),1/2}$ be Fraenkel's regularized distance of \cref{rgd} with the convention that $\rho_{\e}=0$. For $x\notin F$,
put
$$
\rho_{F}(x)
:=
\sum_{i=1}^N
\chi_i(x)\rho_{E_i(F)}\bigl(\u_i(x)\bigr).
$$
Each summand is extended by zero outside $B_g(p_i,r)$. The resulting function is smooth on $\M\setminus F$.

For $x\in\operatorname{supp}\chi_i\setminus F$, Lemma \ref{Stein}
provides constants $c_r,C_r>0$, independent of $i$ and $F$,
such that
$$
c_r d_g(x,F)
\leq d_{\e}\bigl(\u_i(x),E_i(F)\bigr)
\leq C_r d_g(x,F).
$$
Moreover, if $B_g(p_i,2r)\subset F$ then $\chi_i(x)=0$.
As a result
$$
\rho_{F}(x)\asymp_{\M,r}d_g(x,F).
$$

Let $k\geq1$. Fraenkel's derivative estimate of \cref{rgd} gives, for
$1\leq j\leq k$,
$$
\bigl|D^j\rho_{E_i(F)}(\u_i(x))\bigr|
\ll_{j}
d_{\e}\bigl(\u_i(x),E_i(F)\bigr)^{1-j},
$$
The covariant Leibniz rule of \cref{lbt} and covariant Faà di Bruno of \cref{fdb} give
$$
|\nabla^k\rho_{F}(x)|_g
\ll_{\M,k}
\sum_{j=0}^k d_g(x,F)^{1-j}
\ll_{\M,k}
d_g(x,F)^{1-k}.
$$
In the last step we used $d_g(x,F)\leq D_{\mathbb M}$.

Now let $F'$ be another non-empty proper closed set and let 
$x\in \M\setminus(F\cup F')$.
Define $h=d^g_H(F,F')$ and $r(x)=\min(d_g(x,F),d_g(x,F'))$.
By
Lemma \ref{Stein},
$$
d_H^{\e}\bigl(E_i(F),E_i(F')\bigr)
\leq C_rh,
$$
and $$
\min\left\{
d_{\e}(\u_i(x),E_i(F)),
d_{\e}(\u_i(x),E_i(F'))
\right\}
\geq c_rr(x),
$$
whenever $x\in\operatorname{supp}\chi_i$.
Consequently, \cref{Stein,hsd} imply, for
$0\leq j\leq k$,
$$
\left|
D^j\left(
\rho_{E_i(F)}-\rho_{E_i(F')}
\right)(\u_i(x))
\right|
\ll_{\M,j}
\frac{h}{r(x)^j}.
$$
The covariant Leibniz rule of \cref{lbt} together with the covariant Faà di Bruno of \cref{fdb} similarly give
$$
\bigl|\nabla^k\rho_{F}(x)
      -\nabla^k\rho_{F'}(x)\bigr|_g
\ll_{\M,k}
h\sum_{j=0}^kr(x)^{-j}.
$$
Since $r(x)\leq D_{\mathbb M}$, this yields for every $k\geq0$,
$$
\bigl|\nabla^k\rho_{F}(x)
      -\nabla^k\rho_{F'}(x)\bigr|_g
\ll_{\M,k}
\frac{d^g_H(F,F')}
{\min\{d_g(x,F),d_g(x,F')\}^{k}}.
$$
For $x\in F$, we have $d_g(x,F')\leq h$, and therefore
$$
|\rho_F(x)-\rho_{F'}(x)|
=\rho_{F'}(x)
\ll_{\M} d_g(x,F')
\leq h.
$$
The same argument applies to $x\in F'$.
Together with the order-zero estimate off $F\cup F'$, this gives
$$
\|\rho_F-\rho_{F'}\|_{L^\infty(\X\setminus F\cup F')}
\ll_{\M} d_H^g(F,F').
$$
\end{proof}
\begin{remark}
    The constants in the above approximation are effective in terms of
Fraenkel's parameters (which are computed explicitly in \cref{reg}),
the geometry of $\M$, the geometric data of the chosen normal atlas,
and the partition of unity. Derivatives of order $k$ involve, in addition,
bounds on derivatives through order $k$ of the normal coordinate maps
and their inverses, as discussed in
\cite[Theorem 2.5 and Lemma 3.8]{Schick2001}.
Moreover, the distance-comparison constants improve near $F$.
\end{remark}
We refer the reader to \cref{content} for a wide range of
examples, including non-smooth sets with finite, computable
upper Minkowski content. In many cases, the Minkowski content coincides with the corresponding Hausdorff measure.
We now state the main theorem of this section:
\begin{theorem}\label{yes}
Let $A\subset\M$ have nonempty topological boundary
and satisfy
$$
\overline{\mathcal M}_g^{d-q}(\partial A)<\infty
\qquad\text{for some }0<q\leq d.
$$
Fix $0<r_0<1$, and set
$$
C_{\partial A}
:=
\sup_{0<r<r_0}
\frac{\operatorname{vol}_g((\partial A)^{(r)})}{r^q}
<\infty.
$$
Suppose that, for some $\beta\geq0$,
$$
\|\varphi_j\|_{L^\infty(\M)}
\ll_{\M}\lambda_j^\beta,
$$
for all $j\geq 1$.
Fix an integer $s>d+4\beta$.
For every $0<r<r_0$ and $T\geq1$, there exist
spectral polynomials
$$
P_{r,T,A}^{\pm}
=
\sum_{\lambda_j<T}c_{j,A}^{\pm}\varphi_j
$$
with the following properties.
\begin{enumerate}
\item We have the following approximation:
\begin{enumerate}
\item
$$
P_{r,T,A}^-\leq\mathbf 1_A\leq P_{r,T,A}^+.
$$

\item For every integer $n\geq1$ and every $j\geq1$
with $\lambda_j<T$,
$$
|c_{j,A}^{\pm}|
\ll_{\M,n}
\frac{
C_{\partial A}\|\varphi_j\|_{L^\infty(\M)}
}{
r^{n-q}\lambda_j^{n/2}
}.
$$

\item There exists a constant $K_{\M,s,\beta}>0$ such that
$$
\int_{\M}(P_{r,T,A}^+-P_{r,T,A}^-)
\,d\operatorname{vol}_g
\leq
C_{\partial A}r^q
+
\frac{
K_{\M,s,\beta}C_{\partial A}
}{
r^{s-q}T^{(s-d-4\beta)/2}
}.
$$
\end{enumerate}

\item Our construction is stable under Hausdorff perturbations in the following sense:
Let $B\subset\M$ satisfy the same assumptions, with
the same Minkowski codimension $q$ and $r_0$, and define $C_{\partial B}$
analogously.
\begin{enumerate}
\item For every integer $n\geq1$ and every $j\geq1$
with $\lambda_j<T$,
$$
|c_{j,A}^+-c_{j,B}^+|
\ll_{\M,n}
\frac{
(C_{\partial A}+C_{\partial B})
d_H^g(\overline A,\overline B)\,
\|\varphi_j\|_{L^\infty(\M)}
}{
r^{n+1-q}\lambda_j^{n/2}
},
$$
and
$$
|c_{j,A}^--c_{j,B}^-|
\ll_{\M,n}
\frac{
(C_{\partial A}+C_{\partial B})
d_H^g\bigl(
\overline{\M\setminus A},
\overline{\M\setminus B}
\bigr)\,
\|\varphi_j\|_{L^\infty(\M)}
}{
r^{n+1-q}\lambda_j^{n/2}
}.
$$

\item For the constant coefficients corresponding
to $\lambda_0=0$,
$$
|c_{0,A}^+-c_{0,B}^+|
\ll_{\M}
\frac{d_H^g(\overline A,\overline B)}{r}
+
\frac{
K_{\M,s,\beta}|C_{\partial A}-C_{\partial B}|
}{
r^{s-q}T^{(s-d-4\beta)/2}
},
$$
and
$$
|c_{0,A}^--c_{0,B}^-|
\ll_{\M}
\frac{
d_H^g\bigl(
\overline{\M\setminus A},
\overline{\M\setminus B}
\bigr)
}{r}
+
\frac{
K_{\M,s,\beta}|C_{\partial A}-C_{\partial B}|
}{
r^{s-q}T^{(s-d-4\beta)/2}
}.
$$
\end{enumerate}
\end{enumerate}
In both cases, the dependence on $\M$,$s$ and $n$ is effective.
\end{theorem}
\begin{proof}
Fix $\chi\in C^\infty([0,\infty))$ such that
$0\leq\chi\leq1$, $\chi=1$ on $[0,\frac12]$, and
$\chi=0$ on $[\frac23,\infty)$.
Let $\rho_{\overline A}$ and
$\rho_{\overline{\M\setminus A}}$ be the regularized
distances of \cref{regularized}, with the convention
$\rho_\M=0$.
Let $c_\M,C_\M>0$ be constants depending only on $\M$
such that
$$
c_\M d_g(x,A)
\leq\rho_{\overline A}(x)
\leq C_\M d_g(x,A)
\qquad(x\in\M\setminus\overline A),
$$
and
$$
c_\M d_g(x,\M\setminus A)
\leq\rho_{\overline{\M\setminus A}}(x)
\leq C_\M d_g(x,\M\setminus A)
\qquad(x\in\operatorname{int}(A)).
$$

Given $0<r<r_0$, define
$$
D_{r,A}^+(x)
:=
\begin{cases}
\displaystyle
\chi\left(\frac{\rho_{\overline A}(x)}{c_\M r}\right),
& x\in\M\setminus\overline A,\\[2mm]
1,
& x\in\overline A,
\end{cases}
$$
and
$$
D_{r,A}^-(x)
:=
\begin{cases}
\displaystyle
1-\chi\left(
\frac{\rho_{\overline{\M\setminus A}}(x)}{c_\M r}
\right),
& x\in\operatorname{int}(A),\\[2mm]
0,
& x\in\M\setminus\operatorname{int}(A).
\end{cases}
$$
These functions are smooth because $\chi$ is constant
near $0$.

Define the regions
$$
\mathcal T_r^+(A)
:=
\left\{
x\in\M:
\frac{c_\M}{3C_\M}r
\leq d_g(x,A)\leq\frac23r
\right\},
$$
and
$$
\mathcal T_r^-(A)
:=
\left\{
x\in\M:
\frac{c_\M}{3C_\M}r
\leq d_g(x,\M\setminus A)\leq\frac23r
\right\}.
$$
Since $\M$ is a connected Riemannian manifold,
$$
d_g(x,\partial A)
=
\begin{cases}
d_g(x,A), & x\in\M\setminus A,\\
d_g(x,\M\setminus A), & x\in A.
\end{cases}
$$
In particular,
$\mathcal T_r^\pm(A)\subset(\partial A)^{(r)}$.
Consequently:
\begin{enumerate}
\item
$0\leq D_{r,A}^-\leq\ind_A\leq D_{r,A}^+\leq1$,

\item
$D_{r,A}^+=0$ on $\M\setminus A^{(r)}$, while
$D_{r,A}^-=1$ on
$\M\setminus(\M\setminus A)^{(r)}$;

\item
for every integer $n\geq1$,
$$
\operatorname{supp}(\nabla^nD_{r,A}^{\pm})
\subseteq\mathcal T_r^\pm(A)
\subset(\partial A)^{(r)}.
$$
\end{enumerate}

Write
$$
D_{r,A}^{\pm}
=
\sum_{j\geq0}b_{j,A}^{\pm}\varphi_j.
$$
By the covariant Faà di Bruno formula of \cref{fdb},
for every integer $n\geq1$,
\begin{align*}
|\nabla^nD_{r,A}^{\pm}(x)|_g
&\ll_{\M,n,\chi}
\sum_{m=1}^n
\sum_{\substack{j_1+\cdots+j_m=n\\j_\ell\geq1}}
r^{-m}
\prod_{\ell=1}^m
d_g(x,\partial A)^{1-j_\ell}
\mathbf1_{\mathcal T_r^\pm(A)}(x)
\\
&\ll_{\M,n,\chi}
\sum_{m=1}^n
r^{-m}r^{m-n}
\mathbf1_{\mathcal T_r^\pm(A)}(x)
\\
&\ll_{\M,n,\chi}
r^{-n}\mathbf1_{\mathcal T_r^\pm(A)}(x).
\end{align*}
Consequently,
$$
\|\nabla^nD_{r,A}^{\pm}\|_{L^1(\M)}
\ll_{\M,n,\chi}
C_{\partial A}r^{q-n}.
$$
Thus \cref{decay} implies that, for every $j\geq1$,
$$
|b_{j,A}^{\pm}|
\ll_{\M,n,\chi}
\frac{
C_{\partial A}\|\varphi_j\|_{L^\infty(\M)}
}{
r^{n-q}\lambda_j^{n/2}
}.
$$

For $\lambda\geq0$, define
$$
N(\lambda):=\#\{j\geq0:\lambda_j\leq\lambda\}.
$$
By Weyl's law \cite[Theorem 1.1]{Hormander1968},
$
N(\lambda)
\sim
\frac{\omega_d}{(2\pi)^d}
\lambda^{d/2}
$ with an effective version given in \cite[(1.5)]{HassannezhadKokarevPolterovich2016}.
Set
$$
\alpha:=\frac{s-d-4\beta}{2}>0,
\qquad
\sigma:=\frac s2-2\beta>\frac d2.
$$
Using
$\|\varphi_j\|_{L^\infty(\M)}
\ll_\M\lambda_j^\beta$,
we obtain
\begin{align*}
\left\|
\sum_{\lambda_j\geq T}b_{j,A}^{\pm}\varphi_j
\right\|_{L^\infty(\M)}
&\ll_{\M,s,\beta,\chi}
\frac{C_{\partial A}}{r^{s-q}}
\sum_{\lambda_j\geq T}\lambda_j^{-\sigma}
\\
&=
\frac{C_{\partial A}}{r^{s-q}}
\int_{[T,\infty)}
\lambda^{-\sigma}\,dN(\lambda).
\end{align*}
Integration by parts in the Stieltjes integral gives
\begin{align*}
\int_{[T,\infty)}
\lambda^{-\sigma}\,dN(\lambda)
&=
\sigma\int_T^\infty
\frac{N(\lambda)-N(T^-)}
{\lambda^{\sigma+1}}\,d\lambda
\\
&\ll_{\M,s,\beta}
T^{-(\sigma-d/2)}
\\
&=
T^{-(s-d-4\beta)/2}.
\end{align*}
For a sufficiently large constant
$C_{\M,s,\beta,\chi}>0$, set
$$
\varepsilon_{r,T,A}^{(s)}
:=
\frac{
C_{\M,s,\beta,\chi}C_{\partial A}
}{
r^{s-q}T^{(s-d-4\beta)/2}
},
$$
so that
$$
\left\|
\sum_{\lambda_j\geq T}b_{j,A}^{\pm}\varphi_j
\right\|_{L^\infty(\M)}
\leq\varepsilon_{r,T,A}^{(s)}.
$$
Define
$$
P_{r,T,A}^{\pm}(x)
:=
\sum_{\lambda_j<T}b_{j,A}^{\pm}\varphi_j(x)
\pm\varepsilon_{r,T,A}^{(s)}.
$$
Then
$$
P_{r,T,A}^-
\leq D_{r,A}^-
\leq\ind_A
\leq D_{r,A}^+
\leq P_{r,T,A}^+.
$$

Writing
$$
P_{r,T,A}^{\pm}
=
\sum_{\lambda_j<T}c_{j,A}^{\pm}\varphi_j,
$$
we have $c_{j,A}^{\pm}=b_{j,A}^{\pm}$ for $j\geq1$.
With the normalization
$\varphi_0=1$,
$$
c_{0,A}^{\pm}
=
b_{0,A}^{\pm}
\pm
\varepsilon_{r,T,A}^{(s)}.
$$
In particular, the asserted bounds for the nonconstant
coefficients follow.

Since $D_{r,A}^+-D_{r,A}^-$ is supported in
$(\partial A)^{(r)}$, and every eigenfunction associated
with a positive eigenvalue has integral zero,
\begin{align*}
\int_\M(P_{r,T,A}^+-P_{r,T,A}^-)
\,d\operatorname{vol}_g
&=
\int_\M(D_{r,A}^+-D_{r,A}^-)
\,d\operatorname{vol}_g
+
2\operatorname{vol}_g(\M)
\varepsilon_{r,T,A}^{(s)}
\\
&\leq
C_{\partial A}r^q
+
2
\varepsilon_{r,T,A}^{(s)}
\\
&\leq
C_{\partial A}r^q
+
\frac{
K_{\M,s,\beta}C_{\partial A}
}{
r^{s-q}T^{(s-d-4\beta)/2}
},
\end{align*}
for a computable constant $K_{\M,s,\beta}$.
Let
$$
h:=d_H^g(\overline A,\overline B),
\qquad
c:=\frac{c_\M}{12C_\M}.
$$
If $h\geq cr$, then, for every integer $n\geq1$,
\begin{align*}
\left|
\nabla^nD_{r,A}^+
-\nabla^nD_{r,B}^+
\right|_g
&\ll_{\M,n,\chi}r^{-n}
\\
&\ll_{\M,n,\chi}\frac{h}{r^{n+1}}.
\end{align*}

Suppose now that $h<cr$.
For $n\geq1$,
$$
\operatorname{supp}\left(
\nabla^nD_{r,A}^+
-\nabla^nD_{r,B}^+
\right)
\subseteq
\mathcal T_r^+(A)\cup\mathcal T_r^+(B).
$$
On this union,
$$
\min\{d_g(x,A),d_g(x,B)\}\geq cr,
$$
so that $x\in\M\setminus(\overline A\cup\overline B)$.
By the covariant Faà di Bruno formula of \cref{fdb},
$\nabla^nD_{r,A}^+-\nabla^nD_{r,B}^+$ is a finite
linear combination, up to permutations of the tensor
slots, of terms of the form
$$
(c_\M r)^{-m}
\chi^{(m)}
\left(\frac{\rho_{\overline A}}{c_\M r}\right)
\bigotimes_{\ell=1}^m
\nabla^{j_\ell}\rho_{\overline A}-
(c_\M r)^{-m}
\chi^{(m)}
\left(\frac{\rho_{\overline B}}{c_\M r}\right)
\bigotimes_{\ell=1}^m
\nabla^{j_\ell}\rho_{\overline B},
$$
where
$$
1\leq m\leq n,
\qquad
j_\ell\geq1,
\qquad
\sum_{\ell=1}^m j_\ell=n.
$$

By the mean value theorem and \cref{regularized},
$$
\left|
\chi^{(m)}
\left(\frac{\rho_{\overline A}}{c_\M r}\right)
-
\chi^{(m)}
\left(\frac{\rho_{\overline B}}{c_\M r}\right)
\right|
\ll_{\M,m,\chi}\frac{h}{r}.
$$
Moreover, \cref{regularized} gives
$$
\left|
\bigotimes_{\ell=1}^m
\nabla^{j_\ell}\rho_{\overline A}
\right|_g
\ll_{\M,n}r^{m-n}.
$$
Combining \cref{regularized} with the identity
\begin{align*}
&
\bigotimes_{\ell=1}^m
\nabla^{j_\ell}\rho_{\overline A}
-
\bigotimes_{\ell=1}^m
\nabla^{j_\ell}\rho_{\overline B}
\\
&\quad=
\sum_{\ell=1}^m
\left(
\bigotimes_{a<\ell}
\nabla^{j_a}\rho_{\overline A}
\right)
\otimes
\left(
\nabla^{j_\ell}\rho_{\overline A}
-
\nabla^{j_\ell}\rho_{\overline B}
\right)
\otimes
\left(
\bigotimes_{a>\ell}
\nabla^{j_a}\rho_{\overline B}
\right)
\end{align*}
yields
$$
\left|
\bigotimes_{\ell=1}^m
\nabla^{j_\ell}\rho_{\overline A}
-
\bigotimes_{\ell=1}^m
\nabla^{j_\ell}\rho_{\overline B}
\right|_g
\ll_{\M,n}
\frac{h}{r^{n-m+1}}.
$$
Consequently,
\begin{align*}
&
\left|
\nabla^nD_{r,A}^+(x)
-
\nabla^nD_{r,B}^+(x)
\right|_g
\\
&\quad\ll_{\M,n,\chi}
\sum_{m=1}^n
\sum_{\substack{j_1+\cdots+j_m=n\\j_\ell\geq1}}
r^{-m}
\left(
\frac{h}{r}r^{m-n}
+
\frac{h}{r^{n-m+1}}
\right)
\\
&\quad\ll_{\M,n,\chi}
\frac{h}{r^{n+1}}.
\end{align*}

Combining the two cases with the fact that the difference has support in $(\partial A)^{(r)}\cup (\partial B)^{(r)}\ $ we obtain, for every
integer $n\geq1$,
\begin{align*}
\|\nabla^nD_{r,A}^+
-\nabla^nD_{r,B}^+\|_{L^1(\M)}
&\ll_{\M,n,\chi}
(C_{\partial A}+C_{\partial B})
h\,r^{q-n-1}
\\
&=
\frac{
O_{\M,n,\chi}\bigl(
(C_{\partial A}+C_{\partial B})h
\bigr)
}{
r^{n+1-q}
}.
\end{align*}
Thus \cref{decay} implies, for $j\geq1$ with
$\lambda_j<T$,
$$
|c_{j,A}^+-c_{j,B}^+|
\ll_{\M,n,\chi}
\frac{
(C_{\partial A}+C_{\partial B})
d_H^g(\overline A,\overline B)
\|\varphi_j\|_{L^\infty(\M)}
}{
r^{n+1-q}\lambda_j^{n/2}
}.
$$
Applying the same argument to the complements gives
$$
|c_{j,A}^--c_{j,B}^-|
\ll_{\M,n,\chi}
\frac{
(C_{\partial A}+C_{\partial B})
d_H^g\bigl(
\overline{\M\setminus A},
\overline{\M\setminus B}
\bigr)
\|\varphi_j\|_{L^\infty(\M)}
}{
r^{n+1-q}\lambda_j^{n/2}
}.
$$

We now consider the constant coefficient corresponding
to $\lambda_0=0$.
First,
$$
\|D_{r,A}^+-D_{r,B}^+\|_{L^\infty(\M)}
\ll_{\M,\chi}
\frac{d_H^g(\overline A,\overline B)}{r}.
$$
Indeed, on $\M\setminus(\overline A\cup\overline B)$,
this follows from the Hausdorff estimate in order zero
for the regularized distances and the mean value theorem.
If $x\in\overline B\setminus\overline A$, then
$D_{r,B}^+(x)=1$ and
$$
|D_{r,A}^+(x)-1|
\leq
\|\nabla D_{r,A}^+\|_{L^\infty(\M)}
d_g(x,\overline A)
\ll_{\M,\chi}
\frac{d_H^g(\overline A,\overline B)}{r}.
$$
The case $x\in\overline A\setminus\overline B$ is
symmetric, while the difference vanishes on
$\overline A\cap\overline B$.

It follows that
\begin{align*}
|c_{0,A}^+-c_{0,B}^+|
&\leq
\int_\M|D_{r,A}^+-D_{r,B}^+|
\,d\operatorname{vol}_g
\\
&\qquad+
\left|
\varepsilon_{r,T,A}^{(s)}
-
\varepsilon_{r,T,B}^{(s)}
\right|
\\
&\leq
\frac{d_H^g(\overline A,\overline B)}{r}
+
\frac{
K_{\M,s,\beta}|C_{\partial A}-C_{\partial B}|
}{
r^{s-q}T^{(s-d-4\beta)/2}
}.
\end{align*}
The corresponding estimate for the minorant follows by applying
the same argument to the complements.
This concludes the proof.
\end{proof}
\begin{remark}\label{first}
\begin{enumerate}
\item Our method also applies to more general boundary and
eigenfunction estimates, although these lie beyond the scope of
our intended arithmetic applications. Examples include
$$
\s\bigl((\partial A)^{(r)}\bigr)
\ll C_{\partial A}r^\alpha
\bigl(\log(e/r)\bigr)^\gamma
$$
and estimates of the form
$$
\|\varphi_j\|_{L^\infty(\M)}
\ll_\M
\lambda_j^\beta
\bigl(\log(2+\lambda_j)\bigr)^\eta
$$
for the chosen eigenbasis.

\item The quality of our approximation depends on the
$L^\infty$-norms of the chosen eigenfunctions. It is therefore
closely tied to the construction of eigenbases with optimal
$L^\infty$ bounds in terms of the eigenvalues.
On the round sphere $\bb S^d$, $d\geq2$, Burq and Lebeau
\cite[Theorem 7]{BurqLebeau2013} use probabilistic Sobolev
embeddings to obtain orthonormal bases of spherical harmonics
satisfying
$$
\|\varphi_j\|_{L^\infty(\bb S^d)}
\ll \bigl(\log(2+\lambda_j)\bigr)^{1/2},
$$
with an implied constant depending on the chosen basis.
Nevertheless, the universal estimate
$$
\|\varphi_j\|_{L^\infty(\M)}
\ll_\M \lambda_j^{(d-1)/4}
$$
is sharp among compact Riemannian manifolds without boundary,
as witnessed by zonal spherical harmonics on the round sphere.
\end{enumerate}
\end{remark}

Given a family $\mathcal A=(A_i)_{i\in I}$ whose boundaries $\partial \mathcal A=(\partial A_i)_{i\in I}$ have a finite uniform Minkowski content
$\overline{\mathcal M}_g^{d-q}(\partial\mathcal A)<\infty$ for $0<q\leq d$
we may use the same constant
$C_{\partial\mathcal A}\geq\sup_{i\in I}C_{\partial A_i}$
in the preceding construction.
The additive tail correction then cancels in differences
of the constant coefficients.
The results of \cref{analytictube,tube3} provide two
classes of examples satisfying such uniform bounds that often appear in Chebotarev-Sato-Tate arithmetic applications.
We now give an analogue of \cref{yes} for families of sets.
\def \A {\mathcal A}
\begin{corollary}\label{family}
Let $\A=(A_i)_{i\in I}$ be a family of subsets of $\M$,
each with nonempty topological boundary, and put
$
\partial\A:=(\partial A_i)_{i\in I}.
$
Assume that
$$
\overline{\mathcal M}_g^{d-q}(\partial\A)<\infty
\qquad\text{for some }0<q\leq d.
$$
Fix $0<r_0<1$, and set
$$
C_{\partial\A}
:=
\sup_{0<r<r_0}\sup_{i\in I}
\frac{\operatorname{vol}_g((\partial A_i)^{(r)})}{r^q}.
$$
Suppose that there exists $\beta\geq0$ such that for $j\geq 1$,
$$
\|\varphi_j\|_{L^\infty(\M)}
\ll_{\M}\lambda_j^\beta.
$$
Fix an integer $s>d+4\beta$.
For every $0<r<r_0$ and $T\geq1$, there exist,
for each $\ell\in I$, spectral polynomials
$$
P_{r,T,A_\ell}^{\pm}
=
\sum_{\lambda_j<T}c_{j,A_\ell}^{\pm}\varphi_j,
$$
constructed using the common constant $C_{\partial\A}$,
with the following properties.
\begin{enumerate}
\item We have the following approximation:
For every $\ell\in I$:
\begin{enumerate}
\item
$$
P_{r,T,A_\ell}^-
\leq\mathbf 1_{A_\ell}
\leq P_{r,T,A_\ell}^+.
$$

\item For every integer $n\geq1$ and every $j\geq1$
with $\lambda_j<T$,
$$
|c_{j,A_\ell}^{\pm}|
\ll_{\M,n}
\frac{
C_{\partial\A}\|\varphi_j\|_{L^\infty(\M)}
}{
r^{n-q}\lambda_j^{n/2}
}.
$$

\item There exists a constant $K_{\M,s,\beta}>0$,
independent of $\ell$, such that
$$
\int_{\M}
(P_{r,T,A_\ell}^+-P_{r,T,A_\ell}^-)
\,d\operatorname{vol}_g
\leq
C_{\partial\A}r^q
+
\frac{
K_{\M,s,\beta}C_{\partial\A}
}{
r^{s-q}T^{(s-d-4\beta)/2}
}.
$$
\end{enumerate}

\item 
For every $k,\ell\in I$:
\begin{enumerate}
\item For every integer $n\geq1$ and every $j\geq1$
with $\lambda_j<T$,
$$
|c_{j,A_k}^+-c_{j,A_\ell}^+|
\ll_{\M,n}
\frac{
C_{\partial\A}
d_H^g(\overline{A_k},\overline{A_\ell})\,
\|\varphi_j\|_{L^\infty(\M)}
}{
r^{n+1-q}\lambda_j^{n/2}
},
$$
and
$$
|c_{j,A_k}^--c_{j,A_\ell}^-|
\ll_{\M,n}
\frac{
C_{\partial\A}
d_H^g\bigl(
\overline{\M\setminus A_k},
\overline{\M\setminus A_\ell}
\bigr)\,
\|\varphi_j\|_{L^\infty(\M)}
}{
r^{n+1-q}\lambda_j^{n/2}
}.
$$

\item For the constant coefficients corresponding
to $\lambda_0=0$,
$$
|c_{0,A_k}^+-c_{0,A_\ell}^+|
\ll_{\M}
\frac{
d_H^g(\overline{A_k},\overline{A_\ell})
}{r},
$$
and
$$
|c_{0,A_k}^--c_{0,A_\ell}^-|
\ll_{\M}
\frac{
d_H^g\bigl(
\overline{\M\setminus A_k},
\overline{\M\setminus A_\ell}
\bigr)
}{r}.
$$
\end{enumerate}
\end{enumerate}
In both cases, the dependence on $\M$,$s$ and $n$ is effective.

\end{corollary}
We now give a concrete application that generalizes the classical
Beurling-Selberg approximation for intervals.
Let $\Lambda\subset\mathbb R^d$ be a full lattice, and fix a
$\mathbb Z$-basis $b_1,\ldots,b_d$ of $\Lambda$. Equip
$\mathbb R^d$ with the inner product $\langle\cdot,\cdot\rangle$
for which this basis is orthonormal, and write $\|\cdot\|$
for the associated norm.
We equip $\bb T_\Lambda:=\mathbb R^d/\Lambda$ with the
associated flat metric.
Denote by $\Lambda^*\simeq\bb Z^d$ the dual lattice.
Under this identification, for $m\in\Lambda^*$ and
$\theta\in\bb T_\Lambda$, write
$\langle m,\theta\rangle=\sum_{i=1}^d m_i\theta_i
\in\mathbb R/\mathbb Z$.
With respect to normalized Haar measure, the characters
$\theta\mapsto\exp\bigl(2\pi i\langle m,\theta\rangle\bigr)$,
$m\in\Lambda^*$, form an orthonormal basis of
$L^2(\bb T_\Lambda)$ consisting of Laplace--Beltrami
eigenfunctions with corresponding eigenvalues
$4\pi^2\|m\|^2$.
\begin{corollary}\label{torus}
Let $\Lambda\subset\bb R^d$ be a lattice, let
$\T:=\bb R^d/\Lambda$ be the associated flat torus,
and denote by $\Lambda^*$ its dual lattice.
Let $A\subset\T$ have nonempty boundary and satisfy
$$
\overline{\mathcal M}_g^{d-q}(\partial A)<\infty
\qquad\text{for some }0<q\leq d.
$$
Fix $0<r_0<1$ and $s>d$, and define
$$
C_{\partial A}
:=
\sup_{0<r<r_0}
\frac{\operatorname{vol}_g((\partial A)^{(r)})}{r^q}.
$$
For every $0<r<r_0$ and $T\geq1$, there exist
spectral polynomials
$$
P_{r,T,A}^{\pm}(\theta)
=
\sum_{\substack{j\in\Lambda^*\\\|j\|<T}}
c_{j,A}^{\pm}e^{2\pi i\langle\theta,j\rangle}
$$
with the following properties.
\begin{enumerate}
\item 
\begin{enumerate}
\item
$$
P_{r,T,A}^-\leq\mathbf 1_A\leq P_{r,T,A}^+.
$$

\item For every integer $n\geq1$ and every
$j\in\Lambda^*$ with $0<\|j\|<T$,
$$
|c_{j,A}^{\pm}|
\ll_{\T,n}
\frac{C_{\partial A}}
{r^{n-q}\|j\|^n}.
$$

\item There exists a constant $K_{\T,s}>0$ such that
$$
\int_{\T}(P_{r,T,A}^+-P_{r,T,A}^-)
\,d\operatorname{vol}_g
\leq
C_{\partial A}r^q
+
\frac{K_{\T,s}C_{\partial A}}
{r^{s-q}T^{s-d}}.
$$
\end{enumerate}

\item 
Let $B\subset\T$ satisfy the same assumptions and define $C_{\partial B}$ analogously, 
\begin{enumerate}
\item For every integer $n\geq1$ and every
$j\in\Lambda^*$ with $0<\|j\|<T$,
$$
|c_{j,A}^+-c_{j,B}^+|
\ll_{\T,n}
\frac{
(C_{\partial A}+C_{\partial B})
d_H^g(\overline A,\overline B)
}{
r^{n+1-q}\|j\|^n
},
$$
and
$$
|c_{j,A}^--c_{j,B}^-|
\ll_{\T,n}
\frac{
(C_{\partial A}+C_{\partial B})
d_H^g\bigl(
\overline{\T\setminus A},
\overline{\T\setminus B}
\bigr)
}{
r^{n+1-q}\|j\|^n
}.
$$

\item For the constant coefficients,
$$
|c_{0,A}^+-c_{0,B}^+|
\ll_{\T}
\frac{d_H^g(\overline A,\overline B)}{r}
+
\frac{
K_{\T,s}|C_{\partial A}-C_{\partial B}|
}{
r^{s-q}T^{s-d}
},
$$
and
$$
|c_{0,A}^--c_{0,B}^-|
\ll_{\T}
\frac{
d_H^g\bigl(
\overline{\T\setminus A},
\overline{\T\setminus B}
\bigr)
}{r}
+
\frac{
K_{\T,s}|C_{\partial A}-C_{\partial B}|
}{
r^{s-q}T^{s-d}
}.
$$
\end{enumerate}
\end{enumerate}
\end{corollary}
We conclude the section with some remarks on how to apply our construction to a more
general class of smooth spaces.
\begin{remark}\label{remarkp}
\begin{enumerate}
    \item
    If $\M$ is disconnected, write
    $\M=\bigsqcup_{i=1}^n\M_i$ for its connected components
    and apply the construction to each $\M_i$. The construction becomes more interesting when the component space $\pi_0(\M)$ carries additional structure such as that of a finite group as in the case of disconnected compact Lie group.
    \item 
    If $\M$ has a smooth boundary $\partial\M$ and
    $A\subset\M$ satisfies
    $\overline A\cap\partial\M=\varnothing$, replace the
    Laplace-Beltrami operator by the Neumann Laplacian
    $\Delta_{\mathrm{New}}$ and use its eigenvalues
    and eigenfunctions \cite[Section 3.3]{LevitinMangoubiPolterovich2023}. Taking
    $r_0<d_g(\overline A,\partial\M)$ ensures that the
    smooth approximations vanish near $\partial\M$,
    so the integrations by parts in \cref{decay} introduce no boundary terms.
\item If $\M$ is noncompact and admits a smooth compactification
$\widetilde\M$ without boundary, we may apply the
construction on $\widetilde\M$ after choosing a smooth
Riemannian metric and verifying the required upper
Minkowski-content bound for $\partial_{\widetilde\M}A$.
Gu and Guilbault discuss in \cite{GuGuilbault2020} criteria and obstructions for
the existence of such compactifications.
\end{enumerate}
\end{remark}
We conclude this section with some questions that could serve as future research directions:
\begin{question}
\begin{enumerate}
\item
    Flat tori (and by extension flat manifolds) admit orthonormal eigenbases with uniformly bounded
$L^\infty$-norms. Whether the existence of such a basis forces
flatness remains open in general. Toth-Zelditch
\cite[Theorem 0.1(b)]{TothZelditch2002} prove it under extra hypotheses. Given the close relation between the decay of the coefficients
of our spectral majorants and minorants, the geometry of $A$,
and the $L^\infty$-norms of the eigenfunctions, it would be
interesting to determine which geometric properties of $\M$
or $A$ can be recovered if we suppose these approximations hold.
\item We use the Laplace-Beltrami operator to treat sets with irregular boundaries and varying families of sets arising in a wide range of number-theoretic problems, particularly those involving rare events. It would be interesting to identify other elliptic differential
operators with discrete spectrum that reflect the symmetries of the underlying arithmetic or geometric problem and yield sharper approximation and equidistribution estimates. A related example
involving the basic Laplacian $\Delta_{\mathrm{basic}}$ is given in
\cref{shimura}. This operator acts on functions constant along the leaves of a suitable foliation and satisfies the basic Weyl law
of \cref{basicweyl}.
\end{enumerate}
\end{question}
\section{Applications to certain homogeneous spaces}
Let $G$ be a compact connected Lie group and $T\subset G$
a maximal torus of rank $d$. For our arithmetic applications, we seek
spectral majorants and minorants expressed as finite linear
combinations of irreducible characters $\chi_\lambda$, which
fit naturally into our $L$-function formalism. The matrix
coefficients arising from the Riemannian structure on $G$
do not directly provide such expansions.

By \cite[Proposition 4.53]{Knapp2002}, the inclusion
$T\hookrightarrow G$ induces a homeomorphism
$$
T/W\simeq G/\operatorname{Ad}(G)=\mathrm{Conj}(G).
$$
Fix a Riemannian isometry $T\simeq\T$, for an appropriate $W$-invariant flat metric $g$ on $\T$. Let
$$
p:{\bb T}_{\L}\rightarrow\mathrm{Conj}(G)
$$
be the canonical projection.
By the Weyl integration formula
\cite[Theorem 8.60]{Knapp2002}, for every conjugacy-invariant
$\mu_G$-integrable function $F$ on $G$,
\begin{equation}\label{Wey}
    \int_{ G}F(g)\,d\mu_G
    =\frac1{|W|}
      \int_{{\bb T}_{\L}}
      F\bigl(p(\theta)\bigr)
      |D(\theta)|^2
      \,d\theta,
\end{equation}
where  $D(\theta)=\prod_{\alpha\in\Phi^+}
      \left(1-e^{2\pi i\langle\alpha,\theta\rangle}\right)$ is the Weyl denominator 
and $d\theta$ the normalized Haar measure on
${\bb T}_{\L}$. Before stating the main result for Lie groups we shall prove the following auxiliary lemma.

\begin{lemma}\label{Lies}
With the above notation, for $m\in\Lambda^*$ and
$\lambda\in X^*(T)^+$, define
$$
I(m,\lambda)
:=\frac{1}{|W|}
\int_{\T}e^{2\pi i\l\theta,m\rr}
\overline{\chi_\lambda(p(\theta))}
|D(\theta)|^2\,d\theta.
$$
Define the finite set
$$
W_\lambda:=\{m\in\Lambda^*:I(m,\lambda)\neq0\}.
$$
Then:
\begin{enumerate}
    \item $|I(m,\lambda)|\leq1$.
    \item $|W_\lambda|\leq|W|^2$.
    \item If $I(m,\lambda)\neq0$, then
    $
    \bigl|\|\lambda\|-\|m\|\bigr|\leq2\|\rhoo\|,
    $ and $
    \dim V_\lambda\ll_G
    (\|\rhoo\|+\|m\|)^{|\Phi^+|}.
    $
\end{enumerate}
\end{lemma}

    \begin{proof}
Given $w\in W$ let ${\varepsilon}(w):=\det(w)$. Weyl's character and denominator
formulas \cite[Corollary 5.76 and Theorem 5.77]{Knapp2002} imply that
$$
\overline{\chi_\lambda(p(\theta))}|D(\theta)|^2
=
\sum_{w,w'\in W}
{\varepsilon}(w){\varepsilon}(w')
e^{2\pi i\langle
\theta,-w(\lambda+\rhoo)+w'\rhoo
\rangle}.
$$

For $m\in\Lambda^*$ and $\lambda\in X^*(T)^+$, define
$$
W(m,\lambda)
:=
\{(w,w')\in W\times W,
m=w(\lambda+\rhoo)-w'\rhoo\}.
$$
Orthogonality of torus characters yields
$$
I(m,\lambda)
=
\frac1{|W|}
\sum_{(w,w')\in W(m,\lambda)}
{\varepsilon}(w){\varepsilon}(w').
$$

For fixed $m,\lambda,w'$, the equation
$
m+w'\rhoo=w(\lambda+\rhoo)
$
determines at most one $w$, because $\lambda+\rhoo$ is regular
and therefore has trivial Weyl stabilizer. Consequently,
$
|W(m,\lambda)|\leq |W|$ and $
|I(m,\lambda)|
\leq \frac{|W(m,\lambda)|}{|W|}
\leq1.
$ Similarly we obtain $
|W_\lambda|\leq |W|^2.
$

Suppose finally that $I(m,\lambda)\neq0$, and choose
$(w,w')\in W(m,\lambda)$. Since $W$ acts by isometries,
$$
\bigl|\|\lambda\|-\|m\|\bigr|
\leq \|w\lambda-m\|
=\|w'\rhoo-w\rhoo\|
\leq2\|\rhoo\|.
$$
Weyl's dimension formula \cite[Theorem 5.84]{Knapp2002} gives:
$$
\dim V_\lambda
=
\prod_{\alpha\in\Phi^+}
\frac{\langle\alpha,\lambda+\rhoo\rangle}
     {\langle\alpha,\rhoo\rangle}\ll_G
(\|m\|+\|\rhoo\|)^{|\Phi^+|}.
$$
This concludes the proof of the lemma.
\end{proof}
We now formulate an analogue of \cref{yes} for conjugacy-invariant sets, with majorants and minorants given by class functions.
\begin{proposition}\label{Lie}
Let $A\subset G$ be a conjugacy-invariant Borel set.
Identify $A$ with its image in $\mathrm{Conj}(G)$,
define $A_T:=p^{-1}(A)$. Assume that $\partial A_T\neq\varnothing$
and
$$
\overline{\mathcal M}_g^{d-q}(\partial A_T)<\infty,
\qquad 0<q\leq d.
$$
Fix $0<r_0<1$ and an integer $s>d$, and set
$$
C_{\partial A}
:=
\sup_{0<r<r_0}
\frac{\s((\partial A_T)^{(r)})}{r^q}.
$$
Then, for every $0<r<r_0$ and
$N\geq\max\{1+2\|\rhoo\|,3\|\rhoo\|\}$, there are class functions
$$
P_{r,N,A}^{\pm}
=
\sum_{\|\lambda\|<N}
a_{\lambda,A}^{\pm}\chi_\lambda
$$
such that:
\begin{enumerate}
    \item
    \begin{enumerate}
        \item
        $$
        P_{r,N,A}^-\leq\mathbf 1_A\leq P_{r,N,A}^+.
        $$

        \item
        For every integer $n\geq1$ and every dominant
        weight $\lambda$ with $0<\|\lambda\|<N$,
        $$
        |a_{\lambda,A}^{\pm}|
        \ll_{G,n}
        \frac{C_{\partial A}}
        {r^{n-q}\|\lambda\|^n}.
        $$

        \item
        There exists a constant $K_{G,s}>0$ such that
        $$
        \int_G(P_{r,N,A}^+-P_{r,N,A}^-)\,d\mu_G
        \ll_G
        C_{\partial A}r^q
        +
        \frac{K_{G,s}C_{\partial A}}
        {r^{s-q}N^{s-d}}.
        $$
    \end{enumerate}

    \item
If $B\subset G$ satisfies the same assumptions with the same
$q$ and $r_0$, define $C_{\partial B}$ analogously. Then:
    \begin{enumerate}
        \item
        For every integer $n\geq1$ and every dominant
        weight $\lambda$ with $0<\|\lambda\|<N$,
        $$
        |a_{\lambda,A}^+-a_{\lambda,B}^+|
        \ll_{G,n}
        \frac{
        (C_{\partial A}+C_{\partial B})
        d_H^g(\overline{A_T},\overline{B_T})
        }{
        r^{n+1-q}\|\lambda\|^n
        },
        $$
        and
        $$
        |a_{\lambda,A}^--a_{\lambda,B}^-|
        \ll_{G,n}
        \frac{
        (C_{\partial A}+C_{\partial B})
        d_H^g(
        \overline{\T\setminus A_T},
        \overline{\T\setminus B_T})
        }{
        r^{n+1-q}\|\lambda\|^n
        }.
        $$

        \item
        For $\lambda=0$, we have
        $$
        |a_{0,A}^+-a_{0,B}^+|
        \ll_G
        \frac{d_H^g(\overline{A_T},\overline{B_T})}{r}
        +
        \frac{
        K_{G,s}|C_{\partial A}-C_{\partial B}|
        }{
        r^{s-q}N^{s-d}
        },
        $$
        and
        $$
        |a_{0,A}^--a_{0,B}^-|
        \ll_{G}
        \frac{
        d_H^g(
        \overline{\T\setminus A_T},
        \overline{\T\setminus B_T})
        }{r}
        +
        \frac{
        K_{G,s}|C_{\partial A}-C_{\partial B}|
        }{
        r^{s-q}N^{s-d}
        }.
        $$
    \end{enumerate}
\end{enumerate}
In both cases, the dependence on $G$,$s$ and $n$ is effective.

\end{proposition}

\def \G {\mathcal{G}}
\begin{proof}
Define $N':=N-2\|\rhoo\|$. Since $N>0$ and
$N\geq 3\|\rhoo\|$, we have $N'\geq N/3>0$.
For $g\in G$, let $[g]\in\mathrm{Conj}(G)$ be its conjugacy
class, and choose a lift $\overline g\in\T$ under the projection
$p:\T\rightarrow\mathrm{Conj}(G)$. Let
$$
P^\pm_{r,N',A_T}(\theta)
=
\sum_{\|m\|<N'}c_m^\pm e^{2\pi i\langle m,\theta\rangle}
$$
be the spectral approximation on the flat torus $\T$
constructed in \cref{torus}. Define the class functions
$$
P^\pm_{r,N,A}(g)
=
\frac1{|W|}
\sum_{w\in W}P^\pm_{r,N',A_T}(w\overline g).
$$
Since $A_T$ is $W$-invariant, these functions are well-defined and satisfy
$$
P^-_{r,N,A}\leq\ind_A\leq P^+_{r,N,A}.
$$

The action of $W$ on $\T$ induces a dual action on $\Lambda^*$.
For $w\in W$, let $w^*:\Lambda^*\to\Lambda^*$ be characterized by
$$
\langle m,w\theta\rangle=\langle w^*m,\theta\rangle.
$$
Since $w$ is an isometry, $\|w^*m\|=\|m\|$.
By the Peter--Weyl expansion and Weyl integration,
\begin{align*}
P^\pm_{r,N,A}(g)
&=
\sum_{\lambda\in X^*(T)^+}
\left(
\int_G P^\pm_{r,N,A}(h)
\overline{\chi_\lambda(h)}\,d\mu_G(h)
\right)\chi_\lambda(g)
\\
&=
\sum_{\lambda\in X^*(T)^+}
\left(
\frac1{|W|}
\sum_{w\in W}\sum_{\|m\|<N'}
c_m^\pm I(w^*m,\lambda)
\right)\chi_\lambda(g).
\end{align*}
Here $
I(m,\lambda)
$ is as in \cref{Lies}.
Weyl invariance gives $I(w^*m,\lambda)=I(m,\lambda)$, and hence
$$
a^\pm_{\lambda,A}
:=
\frac1{|W|}
\sum_{w\in W}\sum_{\|m\|<N'}
c_m^\pm I(w^*m,\lambda)
=
\sum_{\|m\|<N'}c_m^\pm I(m,\lambda).
$$
Moreover, $I(m,\lambda)\neq0$ implies
$\|\lambda\|\leq\|m\|+2\|\rhoo\|<N$. Thus
$$
P^\pm_{r,N,A}
=
\sum_{\substack{\lambda\in X^*(T)^+\\\|\lambda\|<N}}
a^\pm_{\lambda,A}\chi_\lambda,
$$
so these class functions are smooth.

To show the estimate of Peter-Weyl coefficients, it suffices to treat the case
$\|\lambda\|\geq3\|\rhoo\|$.
Indeed, the finitely many nonzero dominant weights with
$\|\lambda\|<3\|\rhoo\|$ and their contributions can be absorbed in the ${C_{\partial A}}r^{q-n}O_G(1)$ term.
By the preceding lemma and the torus coefficient estimate, for every
integer $n\geq1$ and every nonzero $\lambda$ with
$\|\lambda\|\geq3\|\rhoo\|$,
\begin{align*}
|a^\pm_{\lambda,A}|
&\leq
|W_\lambda|
\sup_{\substack{m\in\Lambda^*,\ \|m\|<N'\\
|\|\lambda\|-\|m\||\leq2\|\rhoo\|}}
|c_m^\pm|
\\
&\ll_{G,n}
\frac{C_{\partial A}}
{r^{n-q}(\|\lambda\|-2\|\rhoo\|)^n}
\ll_{G,n}
\frac{C_{\partial A}}
{r^{n-q}\|\lambda\|^n}.
\end{align*}
Finally, Weyl integration and the $W$-invariance of
$|D(\theta)|^2\,d\theta$ give, for $s>d$:
\begin{align*}
&\int_G(P^+_{r,N,A}-P^-_{r,N,A})\,d\mu_G\\
&=
\frac1{|W|^2}
\sum_{w\in W}
\int_{\T}
\bigl(
P^+_{r,N',A_T}(w\theta)
-
P^-_{r,N',A_T}(w\theta)
\bigr)
|D(\theta)|^2\,d\theta
\\
&=
\frac1{|W|}
\int_{\T}
\bigl(
P^+_{r,N',A_T}(\theta)
-
P^-_{r,N',A_T}(\theta)
\bigr)
|D(\theta)|^2\,d\theta
\\
&\ll_{G}
C_{\partial A}r^q
+
\frac{K_{G,s}C_{\partial A}}
{r^{s-q}(N')^{s-d}}
\\
&\ll_{G}
C_{\partial A}r^q
+
\frac{K_{G,s} C_{\partial A}}
{r^{s-q}N^{s-d}}.
\end{align*}
For the second point, let $B\subset G$ satisfy the same assumptions with $B_T$, $r_0$, and $C_{\partial B}$ defined analogously.
For $\lambda\neq0$, the coefficient formula gives
$$
|a_{\lambda,A}^{\pm}-a_{\lambda,B}^{\pm}|
\ll_G
\sum_{\substack{\|m\|<N'\\m\in W_\lambda}}
|c_{m,A_T}^{\pm}-c_{m,B_T}^{\pm}|.
$$
As $W_\lambda$ is finite for $\lambda\neq 0$ and $\|\lambda\|\asymp_G \|m\|$ for $m\in W_\lambda$, it follows that
$$
|a_{\lambda,A}^{+}-a_{\lambda,B}^{+}|
\ll_{G,n}
\frac{(C_{\partial A}+C_{\partial B})
d_H^g(\overline{A_T},\overline{B_T})}
{r^{n+1-q}\|\lambda\|^n}.
$$
For $\lambda=0$,
\begin{align*}
|a_{0,A}^{+}-a_{0,B}^{+}|
&\ll_G
\|D_{r,A_T}^{+}-D_{r,B_T}^{+}\|_\infty
+
|\varepsilon_{r,N',A_T}^{(s)}
-\varepsilon_{r,N',B_T}^{(s)}|
\\
&\ll_G
\frac{d_H^g(\overline{A_T},\overline{B_T})}{r}
+
\frac{K_{G,s}|C_{\partial A}-C_{\partial B}|}
{r^{s-q}N^{s-d}}.
\end{align*}
The minus estimates follow by replacing
$\overline{A_T},\overline{B_T}$ with
$\overline{\T\setminus A_T},\overline{\T\setminus B_T}$. This concludes the proof of the statement.
\end{proof}
The same construction gives the following analogue
of Corollary 3.5 for families of conjugacy-invariant sets.
\begin{corollary}\label{Liefamily}
Let $\A=(A_i)_{i\in I}$ be a family of conjugacy-invariant
Borel subsets of $G$. Identifying each $A_i$ with its image
in $\mathrm{Conj}(G)$, define
$$
A_{i,T}:=p^{-1}(A_i),
\qquad
\partial\A:=(\partial A_{i,T})_{i\in I}.
$$
Assume that each $A_{i,T}$ has nonempty topological boundary
and that {for some } $0<q\leq d$ the $(d-q)$-th uniform Minkowski content is finite
$
\overline{\mathcal M}_g^{d-q}(\partial\A)<\infty.
$
Fix $0<r_0<1$ and an integer $s>d$, and set
$$
C_{\partial\A}
:=
\sup_{0<r<r_0}\sup_{i\in I}
\frac{\s((\partial A_{i,T})^{(r)})}{r^q}.
$$
For every $0<r<r_0$ and
$N\geq\max\{1+2\|\rhoo\|,3\|\rhoo\|\}$, there exist,
for each $\ell\in I$, class functions
$$
P_{r,N,A_\ell}^{\pm}
=
\sum_{\|\lambda\|<N}
a_{\lambda,A_\ell}^{\pm}\chi_\lambda,
$$
constructed using the common constant $C_{\partial\A}$,
with the following properties.
\begin{enumerate}
    \item
    For every $\ell\in I$:
    \begin{enumerate}
        \item
        $$
        P_{r,N,A_\ell}^-
        \leq\mathbf 1_{A_\ell}
        \leq P_{r,N,A_\ell}^+.
        $$

        \item
        For every integer $n\geq1$ and every dominant
        weight $\lambda$ with $0<\|\lambda\|<N$,
        $$
        |a_{\lambda,A_\ell}^{\pm}|
        \ll_{G,n}
        \frac{C_{\partial\A}}
        {r^{n-q}\|\lambda\|^n}.
        $$

        \item
        There exists a constant $K_{G,s}>0$, independent of $\ell$,
such that
        $$
        \int_G
        (P_{r,N,A_\ell}^+-P_{r,N,A_\ell}^-)\,d\mu_G
        \ll_G C_{\partial\A}r^q
        +
        \frac{K_{G,s}C_{\partial\A}}
        {r^{s-q}N^{s-d}}.
        $$
    \end{enumerate}

    \item
    For every $k,\ell\in I$:
    \begin{enumerate}
        \item
        For every integer $n\geq1$ and every dominant
        weight $\lambda$ with $0<\|\lambda\|<N$,
        $$
        |a_{\lambda,A_k}^+-a_{\lambda,A_\ell}^+|
        \ll_{G,n}
        \frac{
        C_{\partial\A}
        d_H^g(\overline{A_{k,T}},\overline{A_{\ell,T}})
        }{
        r^{n+1-q}\|\lambda\|^n
        },
        $$
        and
        $$
        |a_{\lambda,A_k}^--a_{\lambda,A_\ell}^-|
        \ll_{G,n}
        \frac{
        C_{\partial\A}
        d_H^g(
        \overline{\T\setminus A_{k,T}},
        \overline{\T\setminus A_{\ell,T}})
        }{
        r^{n+1-q}\|\lambda\|^n
        }.
        $$

        \item
        For $\lambda=0$, we have
        $$
        |a_{0,A_k}^+-a_{0,A_\ell}^+|
        \ll_G
        \frac{
        d_H^g(\overline{A_{k,T}},\overline{A_{\ell,T}})
        }{r},
        $$
        and
        $$
        |a_{0,A_k}^--a_{0,A_\ell}^-|
        \ll_G
        \frac{
        d_H^g(
        \overline{\T\setminus A_{k,T}},
        \overline{\T\setminus A_{\ell,T}})
        }{r}.
        $$
    \end{enumerate}
\end{enumerate}
In both cases, the dependence on $G$, $s$ and $n$ is effective.
\end{corollary}
\begin{remark}
A natural direction is to extend our construction to disconnected compact Lie groups. For $G$ compact, connected, simply connected
and simple, 
Weyl integration and character formulas 
are available in \cite[Section 2]{Wendt2001}.
Moreover, if the neutral component $G^\circ$ is compact, simply connected and simple,
and $\kappa=\operatorname{Ad}(k)|_{G^\circ}$ is a diagram
automorphism, the $\kappa$-twining characters give rise to
an orthonormal basis of
$L^2(G^\circ k)^{\operatorname{Ad}(G^\circ)}$,
see \cite[Proposition 4.5]{Zerouali2020}. 
A possible starting point for a Beurling-Selberg theory on
disconnected compact Lie groups $G$ is to pass to a finite
central isogeny
$
Z(G^\circ)^\circ\times\prod_i G_i\longrightarrow G^\circ,
$
where each $G_i$ is compact, simple and simply connected,
and then study the remaining components through twisted
conjugation.
\end{remark}
A natural question is whether our approach, which relies on
reduction to a flat manifold, extends to other classes of spaces.
We illustrate such an extension in the case of a simply connected
symmetric space whose restricted roots have even multiplicities.

We refer the reader to \cref{symmetric} for a review of simply
connected symmetric spaces and Heckman--Opdam Jacobi polynomials,
and to \cref{even} for explicit calculations of spherical functions
in the case of even restricted root multiplicities.
Let $G$ be a connected compact simply connected semisimple Lie group,
and fix a compact subgroup $K\subset G$ such that $(G,K)$ is a
symmetric pair. Denote by $X:=G/K$ the associated symmetric space,
equipped with the pushforward $\mu_X$ of the normalized Haar measure.
Let $F\subset X$ be a maximal flat torus of dimension $d$ as in
\cref{symmetric}, with restricted Weyl group $W_K$.
Let $\Sigma$ be the restricted root system and fix a positive
subsystem $\Sigma^+$. For $\alpha\in\Sigma$, denote its multiplicity
by $m_\alpha$, and put
$\rhoo_F:=\frac12\sum_{\alpha\in\Sigma^+}m_\alpha\alpha$.

For continuous $K$-invariant functions $u$ on $X$, identified with
functions on $K\setminus X$, the projection
$p_K:F\rightarrow K\setminus X$ gives the following analogue of
the Weyl integration formula \cite[p. 190]{Helgason2000}:
$$
\int_X u\,d\mu_X
=
c\int_F u(p_K(\theta))J_K(\theta)\,d\theta,
$$
where $c>0$ is a normalization constant and
$$
J_K(\theta)
=
\prod_{\alpha\in\Sigma^+}
\left|
e^{2\pi i\l\alpha,\theta\rr}
-
e^{-2\pi i\l\alpha,\theta\rr}
\right|^{m_\alpha}.
$$
If all $m_\alpha$ are even, then $J_K$ is a trigonometric polynomial.
More precisely,
$$
J_K(\theta)=
(-1)^{\frac12\sum_{\alpha>0}m_\alpha}
e^{4\pi i\langle\rhoo_F,\theta\rangle}
\prod_{\alpha>0}
(1-e^{-4\pi i\langle\alpha,\theta\rangle})^{m_\alpha},
$$
and
$$
\operatorname{supp}\widehat J_K\subseteq
\left\{2\rhoo_F-2\sum_{\alpha>0}k_\alpha\alpha:
0\leq k_\alpha\leq m_\alpha,\ k_\alpha\in\mathbb Z
\right\}.
$$
Equip $F$ with a $W_K$-invariant flat metric and fix an isometric
identification $F\simeq\T$. As in the case of Lie groups, write
$p_K:\T\rightarrow K\setminus X$ for the corresponding projection.
For spherical functions, we have the following analogue of the
Weyl character formula in terms of Heckman-Opdam Jacobi polynomials
\cite[Theorem 3.3]{RemlingRosler2015}:
$$
\f_\lambda(p_K(\theta))
=
\frac{P_\lambda(\theta)}{P_\lambda(0)}.
$$
\begin{lemma}\label{sym}
Assume that the restricted root multiplicities are even.
Under a Riemannian identification $F\simeq\T$, for $m\in \L^*$
and $\lambda\in X^*(F)^+$, define
$$
I(m,\lambda)
:=
c\int_{\T}
e^{2\pi i\l m,\theta\rr}
\overline{\phi_\lambda(p_K(\theta))}
J_K(\theta)\,d\theta.
$$
Put
$$
W_\lambda:=\{m\in \L^*:I(m,\lambda)\neq0\}.
$$
Then
\begin{enumerate}
    \item $
|I(m,\lambda)|\leq1$,
\item  $
|W_\lambda|\ll_X1.
$
\end{enumerate}

Moreover, if $I(m,\lambda)\neq0$, then
$$
\|m\|-2\|\rhoo_F\|
\leq\|\lambda\|\leq\|m\|.
$$
\end{lemma}

\begin{proof}
Since $c\int_{\T}J_K(\theta)\,d\theta=1$ and
$|\phi_\lambda|\leq1$, we have
$$
|I(m,\lambda)|
\leq c\int_{\T}J_K(\theta)\,d\theta=1.
$$

Let $m^+$ be the dominant representative of $m$.
By $W_K$-invariance,
$$
I(m,\lambda)
=
\frac{c}{|W_Km^+|P_\lambda(0)}
\int_{\T}
M_{m^+}(\theta)\overline{P_\lambda(\theta)}
J_K(\theta)\,d\theta.
$$
Triangularity and orthogonality of the Jacobi polynomials
therefore imply that $I(m,\lambda)\neq0$ requires
$\lambda\leq m^+$.
Define
$$
Q_+:=\mathbb Z_{\geq0}\Sigma^+
=
\left\{
\sum_{\alpha\in\Sigma^+}n_\alpha\alpha:
n_\alpha\in\mathbb Z_{\geq0}
\right\}.
$$
The support formula above shows that every Fourier weight $n$ of $J_K$
satisfies $2\varrho_F-n\in Q_+$. Similarly, every Fourier weight $\gamma$
of $\varphi_\lambda$ satisfies $\lambda-\gamma\in Q_+$.
Nonvanishing of $I(m,\lambda)$ therefore implies
$$
m^+=\gamma+n,
\qquad
2\varrho_F-(m^+-\lambda)
=(2\varrho_F-n)+(\lambda-\gamma)\in Q_+.
$$
Together with $\lambda\leq m^+$, this gives
$$
m^+-\lambda\in
\{\delta\in Q_+:2\varrho_F-\delta\in Q_+\}.
$$
This set is finite and depends only on $X$. Each difference determines
$m^+$, which has at most $|W_K|$ Weyl conjugates. Hence $|W_\lambda|\ll_X1$.
Finally,
$$
\|\lambda\|\leq\|m^+\|=\|m\|
\leq\|\lambda+2\varrho_F\|
\leq\|\lambda\|+2\|\varrho_F\|,
$$
which proves the asserted bounds.
\end{proof}
\begin{proposition}\label{symeven}
Assume that the symmetric space $X$ has even restricted root multiplicities.
Let $A\subset X$ be a $K$-invariant Borel set.
Identifying $A$ with its image in $K\setminus X$, define
$A_F:=p_K^{-1}(A)\subset\T$, equipped with the chosen
$W_K$-invariant flat metric. Assume that
$\partial A_F\neq\varnothing$ and
$$
\overline{\mathcal M}_g^{d-q}(\partial A_F)<\infty,
\qquad 0<q\leq d.
$$
Fix $0<r_0<1$ and an integer $s>d$, and set
$$
C_{\partial A}
:=
\sup_{0<r<r_0}
\frac{\s((\partial A_F)^{(r)})}{r^q}.
$$
Then, for every $0<r<r_0$ and
$N\geq\max\{1+2\|\rhoo_F\|,3\|\rhoo_F\|\}$, there are $K$-invariant functions
$$
P_{r,N,A}^{\pm}
=
\sum_{\|\lambda\|<N}
a_{\lambda,A}^{\pm}(\dim V_\lambda)\phi_\lambda
$$
such that:
\begin{enumerate}
\item
\begin{enumerate}
\item
$$
P_{r,N,A}^-\leq\mathbf 1_A\leq P_{r,N,A}^+.
$$

\item For every integer $n\geq 1$ and every dominant
weight $\lambda$ with $0<\|\lambda\|<N$,
$$
|a_{\lambda,A}^{\pm}|
\ll_{X,n}
\frac{C_{\partial A}}
{r^{n-q}\|\lambda\|^n}.
$$

\item There exists a constant $K_{X,s}>0$ such that
$$
\int_X(P_{r,N,A}^+-P_{r,N,A}^-)\,d\mu_X
\ll_X
C_{\partial A}
\left(
r^q+\frac {K_{X,s}}{r^{s-q}N^{s-d}}
\right).
$$
\end{enumerate}

\item
If $B\subset X$ satisfies the same assumptions with the same
$q$ and $r_0$, define $C_{\partial B}$ analogously. Then:
\begin{enumerate}
    \item
    For every integer $n\geq1$ and every dominant weight
    $\lambda$ with $0<\|\lambda\|<N$, we have
    $$
    |a_{\lambda,A}^+-a_{\lambda,B}^+|
    \ll_{X,n}
    \frac{
    (C_{\partial A}+C_{\partial B})
    d_H^g(\overline{A_F},\overline{B_F})
    }{
    r^{n+1-q}\|\lambda\|^n
    },
    $$
    and
    $$
    |a_{\lambda,A}^--a_{\lambda,B}^-|
    \ll_{X,n}
    \frac{
    (C_{\partial A}+C_{\partial B})
    d_H^g(
    \overline{\T\setminus A_F},
    \overline{\T\setminus B_F})
    }{
    r^{n+1-q}\|\lambda\|^n
    }.
    $$

    \item
    For $\lambda=0$, we have
    $$
    |a_{0,A}^+-a_{0,B}^+|
    \ll_X
    \frac{d_H^g(\overline{A_F},\overline{B_F})}{r}
    +
    \frac{
    K_{X,s}|C_{\partial A}-C_{\partial B}|
    }{
    r^{s-q}N^{s-d}
    },
    $$
    and
    $$
    |a_{0,A}^--a_{0,B}^-|
    \ll_X
    \frac{
    d_H^g(
    \overline{\T\setminus A_F},
    \overline{\T\setminus B_F})
    }{r}
    +
    \frac{
    K_{X,s}|C_{\partial A}-C_{\partial B}|
    }{
    r^{s-q}N^{s-d}
    }.
    $$
\end{enumerate}
\end{enumerate}
\end{proposition}

\begin{proof}
    The proof is identical to that for conjugacy classes, using
\cref{sym} and the spherical Peter-Weyl theorem \cite[Theorem 1.1, p. 152]{Helgason1985}.
We therefore omit it.
\end{proof}
\begin{remark}
Our even-multiplicity assumption is consistent with that of
Branson--\'Olafsson--Pasquale
\cite{BransonOlafssonPasquale2005Compact}.
In particular, this hypothesis provides, alongside classical
Fourier-analytic methods, access to a wide range of functional and differential analysis methods such as the differential shift operators
giving explicit formulas for spherical functions, including on
noncompact symmetric spaces $G/K$
\cite[Theorem 2.9]{BransonOlafssonPasquale2005Compact}.
The same paper describes the holomorphic extension of spherical
functions to the complexification
$U_{\mathbb C}/K_{\mathbb C}$ of a compact dual and identifies their restrictions
to $G/K$
\cite[Lemmas 2.4-2.5]{BransonOlafssonPasquale2005Compact}.
These connections can be a basis for future work investigating
extensions of our majorant/minorant construction to noncompact
spaces under suitable additional assumptions, as discussed
in \cref{remarkp}.
\end{remark}
\begin{example}
The uniformly bounded eigenbasis problem on the round sphere
$\bb S^d$, $d\geq2$, is a substantial, longstanding open problem
and affects the quality of our spectral approximation.
As recalled in \cref{first}(2), Burq--Lebeau construct eigenbases
probabilistically, whereas $L^2$-normalized zonal harmonics
have $L^\infty$-norms of order $\lambda_k^{(d-1)/4}$.
Han's preprint \cite[Section 1.3, Problem 1.3]{Han2022}
formulates the uniformly bounded eigenbasis question explicitly.

For odd-dimensional spheres, whose restricted root multiplicities
are even as in \cref{even}, an alternative for $K$-invariant sets
is to use spherical functions.
Let
$\bb{S}^{2d+1}=\mathrm{Spin}(2d+2)/\mathrm{Spin}(2d+1)$ be the round sphere, with $d\geq1$.
In this case, our spectral majorants and minorants take the form
$$
P_{r,N,A}^{\pm}(\theta)
=
\sum_{0\leq\ell<N}
a^\pm_{A,\ell}
\frac{\ell+d}{d}\binom{\ell+2d-1}{\ell}
{}_2F_1\left(
\begin{matrix}
-\ell,\ \ell+2d\\
d+\frac12
\end{matrix}
;\sin^2\frac{\theta}{2}
\right),
$$
where ${}_2F_1$ denotes the Gauss hypergeometric function. For radial sets, this construction reduces the expansion obtained
by directly applying \cref{yes} to the standard hyperspherical harmonic
basis to one term per degree.
\end{example}
\begin{question}
A natural question is whether an analogue of our hypergeometric approximation extends to the even-dimensional sphere
$\mathbb S^{2d}$ for general sets whose boundaries have a finite Minkowski content.
A possible approach is to use the embedding
$\mathbb S^{2d}\hookrightarrow\bb S^{2d+1}$.
However, this requires a suitable extension of $A$ and control
of the resulting boundary tube estimates, especially for
irregular boundaries satisfying
$\mathcal H_g^{2d-q}(\partial A)
\neq\overline{\mathcal M}_g^{2d-q}(\partial A)$.
\end{question}
\section{ SPR-spaces}\label{shimura}
We now consider more general spaces, which need not be smooth but can be realized as quotients of smooth Riemannian manifolds under suitable submetries or foliations. Working with basic functions allows us to construct majorants and minorants on these quotient spaces. This section is inspired by the recent work of Lin-Mendes-Radeschi on Riemannian submetries \cite{LinMendesRadeschi2025} and their previous work on the basic Weyl law for singular Foliations \cite{LinMendesRadeschi2025Weyl}.
We start with a lemma that allows us to incorporate the corresponding symmetries of the underlying set in the construction of \cref{yes}:
\begin{lemma}\label{simp}
Let $A\subset\M$ be a $\sigma$-saturated subset such that there
exist finite linear combinations of eigenfunctions
$f^\pm=\sum_{\lambda_k<T}c_k^\pm\varphi_k$
satisfying $f^-\leq\ind_A\leq f^+$.
Define
$$
f_\sigma^\pm=\sum_{\lambda_k<T}c_k^\pm P_\sigma\varphi_k.
$$
Then $f_\sigma^-\leq\ind_A\leq f_\sigma^+$ and
$$
\int_{\M}(f_\sigma^+-f_\sigma^-)\,d\s
=
\int_{\M}(f^+-f^-)\,d\s.
$$
\end{lemma}
\begin{proof}
Since $A$ is $\sigma$-saturated, $\ind_A$ is basic.
The positivity of $P_\sigma$ therefore gives
$$
f^+_\sigma:=P_\sigma f^+
\geq P_\sigma\ind_A=\ind_A.
$$
Similarly,
$$
f^-_\sigma:=P_\sigma f^-
\leq P_\sigma\ind_A=\ind_A.
$$
The second assertion follows because $P_\sigma$
preserves integrals.
\end{proof}   
\begin{definition}
    A submetric pseudo-Riemannian space (or SPR-space) $X$ is a metric space $X$ equipped with a manifold submetry $\sigma:\M\to X$ from a compact connected Riemannian manifold of dimension $d$ with basic mean curvature and compact fibers.
\end{definition}
Fix an SPR-space $X$ arising from a manifold submetry
$\sigma:\M\rightarrow X$ with compact fibers and basic
mean curvature from a connected compact Riemannian manifold of dimension $d$. Set $\mu_X:=\sigma_*\s$. The pullback induces
an isometry
$$
U:L^2(X)\longrightarrow L^2(\M),
\qquad Uu=u\circ\sigma,
$$
whose image is the closed basic subspace $L^2(\M,\sigma)$.
As per \cref{basic,no} the averaging operator $P_\sigma$ preserves pointwise smoothness, commutes with $\Delta_g$, and stabilizes the Laplace--Beltrami eigenspaces $(E_\lambda)_{\lambda}$. 
The operator $P_\sigma$ is self-adjoint and induces
an orthogonal decomposition
$$
E_\lambda=(E_\lambda\cap\operatorname{im}P_\sigma)
\mathbin{\oplus^{\perp}}
(E_\lambda\cap\ker P_\sigma).
$$
Choose orthonormal bases of these finite-dimensional
summands. Their union is an orthonormal eigenbasis
$(\varphi_j)$ of $L^2(\M)$ adapted to the basic
subspace. Hence, given $j\geq 0$ the projection $P_\sigma\varphi_j$ is either $\varphi_j$ or zero. 
The nonzero basic eigenfunctions descend to functions on $X$
which, after reindexing, form an orthonormal basis $(\psi_j)$
of $L^2(X)$. We reuse $\lambda_j$ for the corresponding
basic eigenvalues. Define the basic Laplacian and its domain
$$
\Delta_{X,\sigma}:=U^{-1}\Delta_g U,
\quad
\mathcal D(\Delta_{X,\sigma})
:=\{u\in L^2(X):u\circ\sigma\in H^2(\M)\}.
$$
It is self-adjoint and nonnegative, and
$$
(\Delta_{X,\sigma}u)\circ\sigma=\Delta_g(u\circ\sigma),
\qquad \Delta_{X,\sigma}\psi_j=\lambda_j\psi_j.
$$
The eigenvalues here are the basic eigenvalues, counted
with basic multiplicity. The operator $P_\sigma$ induces the conditional expectation with respect to the basic $\sigma$-Borel algebra and hence preserves the $L^1$-norm of positive functions.
\begin{example}\label{exampleb}
Clearly any compact Riemannian manifold is an SPR space if we take the identity map. Let $\sigma:\M\rightarrow X$ be a manifold submetry from a
compact connected Riemannian manifold. We give non-trivial examples of
SPR-spaces.
\begin{enumerate}
    \item
    Suppose that $\M$ is a normal homogeneous space.
    A theorem of Lin-Mendes-Radeschi \cite[Theorem B]{LinMendesRadeschi2025} which we recall in \cref{what}  shows $\sigma$
    has basic mean curvature.
    In particular, this applies to the quotient by any
    singular Riemannian foliation on $\M$ with closed leaves.
    Concrete examples we recorded in \cref{tableofbasic} of the Appendix include the conjugacy foliations, twisted conjugacy
    foliations, and coset decompositions associated with closed
    subgroups of a compact connected Lie group equipped with
    a bi-invariant metric.
    \item
    Suppose that the fibers of $\sigma$ are totally geodesic.
    Their second fundamental form
    $B(U,V)=(\nabla_U V)^\perp$ then vanishes for all tangent
    vector fields $U,V$. Consequently, its mean curvature vanishes and is therefore basic.
    If $\sigma$ is induced by a foliation $\ff$, then
    $\Delta_g$ preserves smooth basic functions, and its
    restriction coincides with the basic Laplacian 
    $\Delta_{\mathrm{basic}}$ and descends to $\Delta_{X,\sigma}$ on basic functions. The theory of the basic Laplacian of leaf spaces was extensively discussed in \cite{ParkRichardson1996}.
\end{enumerate}
\end{example}
We refer to \cref{foliationexamples} of \cref{what}  for more explicit examples of SPR spaces arising from Riemannian foliations on the sphere.
\begin{theorem}\label{basicmain}
Let $A\subset X$ be a Borel set with nonempty boundary.
Fix $0<q\leq d$ and $0<r_0<1$, and assume that
$$
C_{\partial A}
:=
\sup_{0<r<r_0}
\frac{\mu_X((\partial A)^{(r)})}{r^q}
<\infty.
$$
Let $\beta\geq0$, assume that $(\varphi_j)$ is an orthonormal eigenbasis of the Laplace-Beltrami operator $\Delta_g$, adapted to the
basic eigenfunctions $(\psi_j)_{j\geq 0}$, satisfying for $j>0$,
$$
\|\varphi_j\|_{L^\infty(\M)}\ll_\M \lambda_j^\beta.
$$
Fix an integer $s>d+4\beta$.
Then, for every $0<r<r_0$ and $T\geq1$, there are
spectral polynomials
$$
P_{r,T,A}^{\pm}
=
\sum_{\lambda_j<T}c_{j,A}^{\pm}\psi_j
$$
such that:
\begin{enumerate}
    \item
    \begin{enumerate}
        \item
        $$
        P_{r,T,A}^-\leq\mathbf1_A\leq P_{r,T,A}^+.
        $$

        \item
        For every integer $n\geq1$ and every $j\geq1$
        with $\lambda_j<T$,
        $$
        |c_{j,A}^{\pm}|
        \ll_{X,n}
        \frac{
        C_{\partial A}\|\psi_j\|_{L^\infty(X)}
        }{
        r^{n-q}\lambda_j^{n/2}
        }.
        $$

        \item
        There exists a constant $K_{X,s,\beta}>0$ such that
        $$
        \int_X(P_{r,T,A}^+-P_{r,T,A}^-)\,d\mu_X
        \leq
        C_{\partial A}r^q
        +
        \frac{
        K_{X,s,\beta}C_{\partial A}
        }{
        r^{s-q}T^{(s-d-4\beta)/2}
        }.
        $$
    \end{enumerate}

    \item
    If $B\subset X$ satisfies the same assumptions with
    the same $q$ and $r_0$, and $C_{\partial B}$ is defined
    analogously, then:
    \begin{enumerate}
        \item
        For every integer $n\geq1$ and every $j\geq1$
        with $\lambda_j<T$,
        $$
        |c_{j,A}^+-c_{j,B}^+|
        \ll_{X,n}
        \frac{
        (C_{\partial A}+C_{\partial B})
        d_H^X(\overline A,\overline B)\,
        \|\psi_j\|_{L^\infty(X)}
        }{
        r^{n+1-q}\lambda_j^{n/2}
        },
        $$
        and
        $$
        |c_{j,A}^--c_{j,B}^-|
        \ll_{X,n}
        \frac{
        (C_{\partial A}+C_{\partial B})
        d_H^X(
        \overline{X\setminus A},
        \overline{X\setminus B})\,
        \|\psi_j\|_{L^\infty(X)}
        }{
        r^{n+1-q}\lambda_j^{n/2}
        }.
        $$

        \item
        For $\lambda_0=0$, we have
        $$
        |c_{0,A}^+-c_{0,B}^+|
        \ll_X
        \frac{d_H^X(\overline A,\overline B)}{r}
        +
        \frac{
        K_{X,s,\beta}|C_{\partial A}-C_{\partial B}|
        }{
        r^{s-q}T^{(s-d-4\beta)/2}
        },
        $$
        and
        $$
        |c_{0,A}^--c_{0,B}^-|
        \ll_X
        \frac{
        d_H^X(
        \overline{X\setminus A},
        \overline{X\setminus B})
        }{r}
        +
        \frac{
        K_{X,s,\beta}|C_{\partial A}-C_{\partial B}|
        }{
        r^{s-q}T^{(s-d-4\beta)/2}
        }.
        $$
    \end{enumerate}
\end{enumerate}
The dependence on $n$, $s$, $X$, $\M$ and $\sigma$ are effective.
\end{theorem}

\begin{proof}
By \cref{saturated} the set $\sigma^{-1}(A)$ has a finite ($d-q$)-upper Minkowski content. Moreover we have
$$
C_{\partial A}
=
\sup_{0<r<r_0}
\frac{\mu_X((\partial A)^{(r)})}{r^q}
=
\sup_{0<r<r_0}
\frac{\s((\partial\sigma^{-1}(A))^{(r)})}{r^q}
=
C_{\partial\sigma^{-1}(A)}.
$$
Moreover, the submetry property gives
$$
d_g(z,\sigma^{-1}(\overline B))
=
d_X(\sigma(z),\overline B)
$$
Since $\sigma$ is continuous and open,
$\sigma^{-1}(\overline B)=\overline{\sigma^{-1}(B)}$.
As a result
\begin{align*}
d_H^X(\overline A,\overline B)
&=
d_H^g(
\overline{\sigma^{-1}(A)},
\overline{\sigma^{-1}(B)}),\\
d_H^X(
\overline{X\setminus A},
\overline{X\setminus B})
&=
d_H^g(
\overline{\M\setminus\sigma^{-1}(A)},
\overline{\M\setminus\sigma^{-1}(B)}).
\end{align*}

Let $Q_{r,T,\sigma^{-1}(A)}^\pm$ be the spectral
approximation constructed in \cref{yes} using the
eigenbasis $(\varphi_j)$.
The set $\sigma^{-1}(A)$ is saturated. Define
$
S_{r,T,A}^{\pm}
:=
P_\sigma Q_{r,T,\sigma^{-1}(A)}^\pm.
$
By \cref{simp} we have $S_{r,T,A}^{-} \leq \ind_{\sigma^{-1}(A)}\leq S_{r,T,A}^{+}$. As all three functions are basic they descend to the SPR-space $X$. Define $P_{r,T,A}^\pm=U^{-1}S_{r,T,A}^\pm$.
If $P_\sigma\varphi_j\neq 0$ then $P_\sigma\varphi_j=\varphi_j$.
Thus
$\|\psi_j\|_{L^\infty(X)}=\|P_\sigma\varphi_j\|_{L^{\infty}(\M)}=\|\varphi_j\|_{L^\infty(\M)}$.
The coefficient bounds and their Hausdorff estimates
follow from \cref{yes} and the identities above.
By the definition
of $\mu_X$,
\begin{align*}
\int_X(P_{r,T,A}^+-P_{r,T,A}^-)\,d\mu_X
&=
\int_\M
\left(
P_\sigma Q_{r,T,\sigma^{-1}(A)}^+
-
P_\sigma Q_{r,T,\sigma^{-1}(A)}^-
\right)\,d\s\\
&=
\int_\M
\left(
Q_{r,T,\sigma^{-1}(A)}^+
-
Q_{r,T,\sigma^{-1}(A)}^-
\right)\,d\s\\
&\leq
C_{\partial A}r^q
+
\frac{
K_{\M,s,\beta}C_{\partial A}
}{
r^{s-q}T^{(s-d-4\beta)/2}
}.
\end{align*}
This proves the statement.
\end{proof}
As in the case of Riemannian manifolds, we have an analogous estimate for families:
\begin{corollary}\label{basicfamily}
Let $X$ be an SPR space, and let
$\A=(A_i)_{i\in I}$ be a family of Borel subsets of $X$,
each with nonempty boundary and uniform $d-q$-Minkowski content. Put
$
\partial\A:=(\partial A_i)_{i\in I}.
$
Fix $0<q\leq d$ and $0<r_0<1$, and assume that
$$
C_{\partial\A}
:=
\sup_{0<r<r_0}\sup_{i\in I}
\frac{\mu_X((\partial A_i)^{(r)})}{r^q}
<\infty.
$$
Let $\beta\geq0$, and assume that $(\varphi_j)$ is an
orthonormal eigenbasis of the Laplace--Beltrami operator
$\Delta_g$, adapted to the basic eigenfunctions
$(\psi_j)_{j\geq0}$, satisfying
$$
\|\varphi_j\|_{L^\infty(\M)}
\ll_\M \lambda_j^\beta,
$$
for $j\geq 1$.
Fix an integer $s>d+4\beta$.
For every $0<r<r_0$ and $T\geq1$, there exist,
for each $\ell\in I$, spectral polynomials
$$
P_{r,T,A_\ell}^{\pm}
=
\sum_{\lambda_j<T}c_{j,A_\ell}^{\pm}\psi_j,
$$
constructed using the common constant $C_{\partial\A}$,
with the following properties.
\begin{enumerate}
\item
For every $\ell\in I$:
\begin{enumerate}
\item
$$
P_{r,T,A_\ell}^-
\leq\mathbf 1_{A_\ell}
\leq P_{r,T,A_\ell}^+.
$$

\item
For every integer $n\geq1$ and every $j\geq1$
with $\lambda_j<T$,
$$
|c_{j,A_\ell}^{\pm}|
\ll_{X,n}
\frac{
C_{\partial\A}\|\psi_j\|_{L^\infty(X)}
}{
r^{n-q}\lambda_j^{n/2}
}.
$$

\item
There exists a constant $K_{X,s,\beta}>0$,
independent of $\ell$, such that
$$
\int_X
(P_{r,T,A_\ell}^+-P_{r,T,A_\ell}^-)
\,d\mu_X
\leq
C_{\partial\A}r^q
+
\frac{
K_{X,s,\beta}C_{\partial\A}
}{
r^{s-q}T^{(s-d-4\beta)/2}
}.
$$
\end{enumerate}

\item
For every $k,\ell\in I$:
\begin{enumerate}
\item
For every integer $n\geq1$ and every $j\geq1$
with $\lambda_j<T$,
$$
|c_{j,A_k}^+-c_{j,A_\ell}^+|
\ll_{X,n}
\frac{
C_{\partial\A}
d_H^X(\overline{A_k},\overline{A_\ell})\,
\|\psi_j\|_{L^\infty(X)}
}{
r^{n+1-q}\lambda_j^{n/2}
},
$$
and
$$
|c_{j,A_k}^--c_{j,A_\ell}^-|
\ll_{X,n}
\frac{
C_{\partial\A}
d_H^X\bigl(
\overline{X\setminus A_k},
\overline{X\setminus A_\ell}
\bigr)\,
\|\psi_j\|_{L^\infty(X)}
}{
r^{n+1-q}\lambda_j^{n/2}
}.
$$

\item
For the constant coefficients corresponding to $\lambda_0=0$,
$$
|c_{0,A_k}^+-c_{0,A_\ell}^+|
\ll_X
\frac{
d_H^X(\overline{A_k},\overline{A_\ell})
}{r},
$$
and
$$
|c_{0,A_k}^--c_{0,A_\ell}^-|
\ll_X
\frac{
d_H^X\bigl(
\overline{X\setminus A_k},
\overline{X\setminus A_\ell}
\bigr)
}{r}.
$$
\end{enumerate}
\end{enumerate}
The dependence on $n$, $s$, $X$, $M$, and $\sigma$ is effective.
\end{corollary}
\begin{question}
Our construction uses basic mean curvature to obtain
a Riemannian-like framework on the SPR space $X$.
Can a Beurling-Selberg approximation theory be developed
for foliated quotients without this assumption?
More precisely, what control on the nonbasic component of the mean curvature
form suffices for such a theory, and how does the nonbasic behavior influence our effective bounds?
\end{question}
\section{Generalized Erd\H{o}s-Tur\'an inequalities}
Fix an SPR space $(X,d_X,\mu_X)$ together with a submetry $\sigma:\M\rightarrow X$ from a connected compact Riemannian manifold of dimension $d$. Let $(x_n)_{n\geq 1}$ be a sequence of elements in $X$.
It is an established fact that the linear span of the Laplace-Beltrami eigenfunctions $(\varphi_k)_{k\geq 0}$ is dense in the space $C^0(\M)$. Thus the linear span of its pushforwards under a manifold submetry $(\psi_k)_{k\geq 0}$ is dense in $C^0(X)$ as well. In particular, Portmanteau theorem implies that we have for every $n\geq 1$ that 
\begin{align*}
    \frac{1}{N}\sum_{1\leq j\leq N}\psi_n(x_j)\rightarrow 0
\end{align*}
if and only if for every bounded Borel function $f:X\rightarrow \bb{R}$ that is continuous outside a set of null measure
$$\frac{1}{N}\sum_{1\leq k\leq N}f(x_k)\underset{N\rightarrow \infty}{\longrightarrow}\int_X fd\mu_X.$$
In other words, $(x_n)_{n\geq 1}$ is equidistributed with respect to the finite measure $\mu_X$: the empirical measure $\frac{1}{N}\sum_{1\leq k\leq N}\delta_{x_k}$ converges toward $\mu_X$. 
The Erd\H{o}s--Tur\'an inequality quantifies the equidistribution
by bounding discrepancy in terms of finitely many exponential
sums \cite{ErdosTuran1948}. Colzani, Gigante, and Travaglini \cite{ColzaniGiganteTravaglini2011}
extend this approach to compact Riemannian manifolds using
Laplace-Beltrami eigenfunctions. We apply our spectral approximation results to obtain analogous estimates on SPR spaces with fully effective bounds and explicit spectral decay, including estimates for moving sets whose variation is measured by Hausdorff distance,
with an application to quaternionic Shimura varieties.
Our moving-target Erd\H{o}s--Tur\'an inequality is motivated by
Martin and Sitar \cite{MArtinSitar2011}, who developed such an inequality for varying intervals in their study of Diophantine quadruples.

Recall that a non-empty Borel set $Z\subset X$ is a $\mu_X$-continuity set if $\mu_X(\partial Z)=0$.
To extend discrepancy estimates beyond the case of a boundary with upper finite Minkowski content, we state and prove the following variant of measure regularity.
\begin{lemma}\label{approximation}
Let $Z\subset X$ be a non-empty $\mu_X$-continuity set. For every ${\varepsilon}>0$ there exist a closed subset $F_\varepsilon$ and an open subset $U_{\varepsilon}$ such that:
\begin{enumerate}
\item $
F_{\varepsilon}\subset Z\subset U_{\varepsilon}$,
    \item  $ \mu_X(U_{\varepsilon}\setminus F_{\varepsilon})<{\varepsilon},
$
 \item $0\leq \mathcal{M}_{\mu_X}^{d-1}(\partial U_\varepsilon)<\infty$ and $0\leq \mathcal{M}_{\mu_X}^{d-1}(\partial F_\varepsilon)<\infty$.
\end{enumerate}
\end{lemma}

\begin{proof}
Let $\sigma:(\M,g)\rightarrow X$ be a manifold submetry
satisfying the SPR hypothesis, and define the
$\sigma$-saturated set
$Z':=\sigma^{-1}(\overline{Z})$.
Since $\mu_X(\partial Z)=0$, choose an open set
$U'_\varepsilon\supset\overline Z$ such that
$$
\mu_X(U'_\varepsilon\setminus\overline Z)
<\varepsilon/2.
$$
By the smooth Urysohn lemma
\cite[Proposition 2.25]{lee2013smooth},
choose $0\leq f\leq1$, $f\in C^\infty(\M)$,
equal to $1$ on $Z'$ and supported in
$\sigma^{-1}(U'_\varepsilon)$.

Set $h:=P_\sigma f$. By the SPR hypothesis,
$h$ is a smooth basic function satisfying
$0\leq h\leq1$, $h=1$ on $Z'$, and
$\operatorname{supp}h\subset
\sigma^{-1}(U'_\varepsilon)$.
By Sard's theorem
\cite[Theorem 6.10]{lee2013smooth},
choose a regular value $a\in(0,1)$.
In particular, the level set $\{h=a\}$
is a smooth compact hypersurface
\cite[Corollary 5.14]{lee2013smooth},
and hence has finite $(d-1)$-dimensional
Minkowski content.

Define
$$
U_\varepsilon
:=
\sigma\left(\{x\in\M:h(x)>a\}\right).
$$
Since $\sigma$ is an open map and $h$ is basic,
$U_\varepsilon$ is open,
$\overline Z\subset U_\varepsilon$, and
$$
\sigma^{-1}(\partial U_\varepsilon)=\{h=a\}.
$$
As $\sigma$ is a Riemannian submetry,
the tubular neighborhoods of these sets
correspond under $\sigma$.
It follows that $\partial U_\varepsilon$
is $(d-1)$-Minkowski measurable with respect
to $\mu_X$. Moreover, since
$\mu_X(\partial Z)=0$, we have
$$
\begin{aligned}
\mu_X(U_\varepsilon\setminus Z)
&=\mu_X(U_\varepsilon\setminus\overline Z)\\
&\leq
\mu_X(U'_\varepsilon\setminus\overline Z)
<\varepsilon/2.
\end{aligned}
$$

We apply the same construction to
$X\setminus\operatorname{int}(Z)$, whose boundary
is contained in $\partial Z$.
Denote by $\widetilde U_\varepsilon$ the resulting
open set and put
$F_\varepsilon:=X\setminus\widetilde U_\varepsilon$.
Then $F_\varepsilon$ is closed,
$F_\varepsilon\subset\operatorname{int}(Z)$,
and $\partial F_\varepsilon
=\partial\widetilde U_\varepsilon$
has finite $(d-1)$-dimensional
Minkowski content. Moreover,
$$
\begin{aligned}
\mu_X(Z\setminus F_\varepsilon)
&=\mu_X(\operatorname{int}(Z)\setminus F_\varepsilon)\\
&=\mu_X\bigl(
\widetilde U_\varepsilon
\setminus(X\setminus\operatorname{int}(Z))
\bigr)\\
&<\varepsilon/2.
\end{aligned}
$$
\end{proof}
The preceding lemma allows us to extend the spectral discrepancy estimates to arbitrary $\mu_X$-continuity sets by approximation with sets whose boundaries have finite Minkowski content. We now formulate the resulting generalized Erdős-Turán inequalities.
\begin{proposition}\label{equidi}
We retain the notation of \cref{basicmain,basicfamily}.
Fix integers $s>d+4\beta$ and $n\geq1$, and let
$(x_i)_{i\geq1}$ be a sequence in $X$.
\begin{enumerate}
\item
Let $A\subset X$ be a subset such that $\overline{\mathcal{M}}_{\mu_X}^{d-q}(\partial A)<\infty$ for $0<q\leq d$. Let $C_{\partial A}$, $s$ and $r_0$ be as in \cref{basicmain}.
For every $0<r<r_0$ and $T\geq1$,
\begin{align*}
&\left|
\frac1N\sum_{i=1}^N\ind_A(x_i)-\mu_X(A)
\right|
\ll
C_{\partial A}\left(
r^q+\frac{1}{r^{s-q}T^{(s-d-4\beta)/2}}
\right)\\
&\quad+
\frac{C_{\partial A}}{Nr^{n-q}}
\sum_{0<\lambda_k<T}
\frac{\|\psi_k\|_\infty}{\lambda_k^{n/2}}
\left|\sum_{i=1}^N\psi_k(x_i)\right|.
\end{align*}

\item
Let $Z$ be a $\mu_X$-continuity set.
For every $\varepsilon>0$, there exist
$C_{\partial Z,\varepsilon}>0$ and
$0<r_{0,\varepsilon}<1$ such that, for
$0<r<r_{0,\varepsilon}$ and $T\geq1$,
\begin{align*}
&\left|
\frac1N\sum_{i=1}^N\ind_Z(x_i)-\mu_X(Z)
\right|
\ll
\varepsilon+
C_{\partial Z,\varepsilon}\left(
r+\frac{1}{r^{s-1}T^{(s-d-4\beta)/2}}
\right)\\
&\quad+
\frac{C_{\partial Z,\varepsilon}}{Nr^{n-1}}
\sum_{0<\lambda_k<T}
\frac{\|\psi_k\|_\infty}{\lambda_k^{n/2}}
\left|\sum_{i=1}^N\psi_k(x_i)\right|.
\end{align*}
\item
Let $\mathcal A=(A_i)_{i\geq1}$ be a family of sets such that $\mathcal{M}_{\mu_X}^{d-q}(\partial \mathcal{A})<\infty$. Let $r_0, C_{\partial \mathcal{A}}$ and $s$ be as in \cref{basicfamily}. For every $i\geq1$, define
$$
\delta_i:=
\max\left\{
d_H^X(\overline{A_i},\overline{A_{i+1}}),
d_H^X\bigl(
\overline{X\setminus A_i},
\overline{X\setminus A_{i+1}}
\bigr)
\right\}.
$$
For every $0<r<r_0$ and $T\geq1$,
\begin{align*}
&\left|
\frac1N\sum_{i=1}^N
\bigl(\ind_{A_i}(x_i)-\mu_X(A_i)\bigr)
\right|
\ll
C_{\partial\mathcal A}\left(
r^q+\frac{1}{r^{s-q}T^{(s-d-4\beta)/2}}
\right)\\
&\quad+
\frac{C_{\partial\mathcal A}}{Nr^{n-q}}
\left(1+\frac1r\sum_{i=1}^{N-1}\delta_i\right)
\sum_{0<\lambda_k<T}
\frac{\|\psi_k\|_\infty}{\lambda_k^{n/2}}
\max_{1\leq m\leq N}
\left|\sum_{i=1}^m\psi_k(x_i)\right|.
\end{align*}
\end{enumerate}
The implied constants depend effectively on
$X$, $\M$, $q$, $s$, $\beta$ and $n$.
\end{proposition}

\begin{proof}
  The first bound follows immediately from \cref{basicmain}.
For the second, let $\varepsilon>0$. By \cref{approximation},
choose $F_\varepsilon\subset Z\subset U_\varepsilon$ such that
$$
\mu_X(U_\varepsilon\setminus F_\varepsilon)<\varepsilon,
$$
and both $\partial F_\varepsilon$ and
$\partial U_\varepsilon$ have finite upper
$(d-1)$-Minkowski content.
Therefore, the inequality
$\ind_{F_\varepsilon}\leq\ind_Z\leq\ind_{U_\varepsilon}$
gives
\begin{align*}
\frac1N\sum_{k=1}^N\ind_{F_\varepsilon}(x_k)
-\mu_X(F_\varepsilon)-\varepsilon
&\leq
\frac1N\sum_{k=1}^N\ind_Z(x_k)-\mu_X(Z)\\
&\leq
\frac1N\sum_{k=1}^N\ind_{U_\varepsilon}(x_k)
-\mu_X(U_\varepsilon)+\varepsilon.
\end{align*}
Choose a radius $0<r_{0,\varepsilon}<1$ and let
$F_\varepsilon$ and $U_\varepsilon$ be as in \cref{approximation}. Set
$$
C_{\partial Z,\varepsilon}
:=\max\{C_{\partial F_\varepsilon},
C_{\partial U_\varepsilon}\}.
$$
the constants are taken to be zero if the respective boundary is empty and as in \cref{basicmain}.
Applying the first bound with $q=1$ to
$F_\varepsilon$ and $U_\varepsilon$ yields
the second point. 
\item 
To prove the third point, we start by defining
$S_k(m):=\sum_{i=1}^m\psi_k(x_i)$.
Summation by parts gives, for $0<\lambda_k<T$,
$$
\sum_{i=1}^N c_{k,A_i}^{\pm}\psi_k(x_i)
=
c_{k,A_N}^{\pm}S_k(N)
+
\sum_{i=1}^{N-1}
\bigl(c_{k,A_i}^{\pm}-c_{k,A_{i+1}}^{\pm}\bigr)S_k(i).
$$
The coefficient estimates in \cref{basicfamily} imply
$$
\left|\sum_{i=1}^N c_{k,A_i}^{\pm}\psi_k(x_i)\right|
\ll_{X,n}
\frac{C_{\partial\mathcal A}\|\psi_k\|_\infty}
{r^{n-q}\lambda_k^{n/2}}
\left(1+\frac1r\sum_{i=1}^{N-1}\delta_i\right)
\max_{1\leq m\leq N}|S_k(m)|.
$$
By \cref{basicfamily}
\begin{align*}
&\left|
\frac1N\sum_{i=1}^N
\bigl(\ind_{A_i}(x_i)-\mu_X(A_i)\bigr)
\right|
\ll 
C_{\partial\mathcal A}\left(
r^q+\frac{1}{r^{s-q}T^{(s-d-4\beta)/2}}
\right)+\\&\frac{1}{N}\sum_{0<\lambda_k<T}\left|\sum_{i=1}^N c_{k,A_i}^{\pm}\psi_k(x_i)\right| \leq C_{\partial\mathcal A}\left(
r^q+\frac{1}{r^{s-q}T^{(s-d-4\beta)/2}}
\right)\\+&\sum_{0<\lambda_k<T}\frac{C_{\partial\mathcal A}\|\psi_k\|_\infty}
{r^{n-q}\lambda_k^{n/2}}
\frac{1}{N}\left(1+\frac1r\sum_{i=1}^{N-1}\delta_i\right)
\max_{1\leq m\leq N}|S_k(m)|.
\end{align*}
This completes the proof.
\end{proof}
\begin{remark}
   Our construction is particularly useful when the increments
$\delta_i$ are summable, for instance when
$\delta_i\ll i^{-\alpha}$ for some $\alpha>1$ so that
$\sum_{N/2\leq i<N}\delta_i\ll N^{1-\alpha}$. The same argument as above gives
\begin{align*}
&\left|
\frac1N\sum_{N/2\leq i\leq N}
\bigl(\ind_{A_i}(x_i)-\mu_X(A_i)\bigr)
\right|
\ll
C_{\partial\mathcal A}\left(
r^q+\frac{1}{r^{s-q}T^{(s-d-4\beta)/2}}
\right)\\
&\quad+
\frac{C_{\partial\mathcal A}}{Nr^{n-q}}
\left(1+\frac{N^{1-\alpha}}r\right)
\sum_{0<\lambda_k<T}
\frac{\|\psi_k\|_\infty}{\lambda_k^{n/2}}
\max_{N/2\leq m\leq N}
\left|\sum_{N/2\leq i\leq m}\psi_k(x_i)\right|.
\end{align*}
This dyadic approach will be useful when studying equidistribution on moving sets as in \cref{sato} and \cref{sato2}.
\end{remark}
\begin{remark}\label{zer}
We record two ways to extend the preceding Erdős-Turán
estimate to functions.
\begin{enumerate}
\item
Let $f:X\to\mathbb R$ be bounded and Borel measurable, 
with a set of discontinuities of $\mu_X$-null measure. By \cite[Proposition 2.3.3]{MareschWinkler2009} for every $\varepsilon>0$, there exist $\mu_X$-continuity
sets $(Z_i)_{1\leq i\leq c}$ and real numbers
$(a_i)_{1\leq i\leq c}$ such that
$$
\sup_{x\in X}
\left|f(x)-\sum_{i=1}^c a_i\ind_{Z_i}(x)\right|
<\varepsilon.
$$
\item
For $X=\T$, choose lattice coordinates
$\T\simeq\mathbb R^d/\mathbb Z^d$, and let $d\theta$
denote normalized Haar measure. Given $f\in C^\infty(\T)$
and $x_1,\ldots,x_N\in[0,1)^d$, define
$$
\Delta_N(t)
:=
\frac1N\sum_{k=1}^N\ind_{[0,t)}(x_k)
-\prod_{j=1}^d t_j,
\qquad
[0,t):=\prod_{j=1}^d[0,t_j).
$$
The Hlawka-Zaremba identity
\cite[p. 91, proof of Proposition 3]{Zaremba1968} gives\begin{align*}
&\frac1N\sum_{k=1}^N f(x_k)-\int_{\T} f\,d\theta
\\&=
\sum_{\varnothing\ne u\subseteq\{1,\ldots,d\}}
(-1)^{|u|}
\int_{[0,1]^{|u|}}
\partial_u f(t_u,1)\,
\Delta_N(t_u,1)\,dt_u,
\end{align*}
where $\partial_u:=\prod_{j\in u}\partial_{t_j}$ and
$(t_u,1)$ has coordinates equal to $1$ outside $u$.
The boundaries of these rectangles have uniformly finite
upper $(d-1)$-dimensional Minkowski content. Applying the
preceding discrepancy estimate uniformly to the rectangles
therefore yields bounds for smooth averages.
\end{enumerate}
\end{remark}
As an application of our earlier estimate, we prove a quantitative
version of Zhang's result on the equidistribution of CM-points
on certain quaternionic Shimura varieties. Zhang
\cite[Corollaries 3.3 and 3.7]{Zhang2005} proves the equidistribution
of CM-points with maximal Mumford--Tate group under the
$\varepsilon$-conjecture and subconvexity bounds. As Zhang remarked, the unconditional equidistribution
of Galois orbits of strict sequences of CM-points implies the Andr\'e--Oort
conjecture.

Equidistribution of Galois orbits of strict sequences of CM-points forces any closed Shimura subvariety
containing these orbits to coincide with the ambient Shimura
variety as stated in \cite[Section 2, Remark (1)]{Zhang2005}.

Let $F$ be a totally real field, and let $B/F$ be a quaternion
algebra that splits at $d\geq1$ real places. Define
$$
G:=\mathrm{Res}_{F/\Q}\mathrm{PGL}_1(B)(\R)^0,
$$
and fix a maximal compact subgroup $K\subset G$. Let
$X=\Gamma\backslash G/K$ be a compact connected component of
Zhang's quaternionic Shimura variety $\mathcal M_U$, where $U$
is his maximal-order level \cite{Zhang2005}. Put $\M:=\Gamma\backslash G$ and $D:=\dim\M$.
Since $\Gamma$ acts freely and properly discontinuously
on $G$ by left translations, $\M$ is a compact connected
smooth manifold. The relevant metric and Laplacian are described in
\cite{RamacherWakatsuki2021}.

Although $\M$ is smooth, $X$ need not be a smooth Riemannian
manifold. For example, let $B/\Q$ be the quaternion algebra of
discriminant $6$, let $\mathcal O\subset B$ be a maximal order,
and put $\Gamma=\mathcal O^1/\{\pm1\}$. The compact quaternionic
Shimura curve $X$ has signature $(0;2,2,3,3)$
\cite[Table 4.1]{Voight2009}. At an elliptic point of order $m$,
the quotient metric has cone angle $2\pi/m$ and is therefore
not smooth, although the underlying complex curve is smooth.

Equip $G$ with the left-invariant, right $K$-invariant Cartan
metric described in
\cite[Section 2]{RamacherWakatsuki2021}.
The projection $G\to G/K$ carries a connection metric, so
its fibers are totally geodesic by
\cite[Proposition 3.1(i)]{petersenwilhelm1999}.
This property descends under the local isometry
$G\to\M$. Hence the orbit projection
$$
\sigma:\M\longrightarrow X=\Gamma\backslash G/K
$$
is a manifold submetry with compact, totally geodesic fibers.
As specified in \cref{exampleb}(2) the submetry $\sigma$ has basic mean curvature and hence $X$ is an SPR-space.

The full variety $\mathcal M_U$ has finitely many geometric
connected components, permuted by the Galois group of its
reflex field, see \cite[(3.6)]{Zhang2005}.
Choose an orthonormal eigenbasis $(\phi_\pi)$ of the basic
Laplacian on $L^2(X)$, with eigenvalues $\lambda_\pi$.

Zhang's period and Galois-orbit estimates are formulated on
the full variety $\mathcal M_U$. A Galois orbit meets each
of its visited connected components in equally many points.
Restriction to a nonempty intersection with $X$ therefore
loses at most a factor $|\pi_0(\mathcal M_U)|$ in its cardinality. We therefore work with the nonempty intersections
of these Galois orbits with $X$.
\begin{corollary}\label{orbifold}
Let $X$ be the compact connected Shimura orbifold as above, and let $(O(x_k))_{k\geq1}$ be the nonempty
intersections with $X$ of distinct Galois orbits of CM-points
whose Mumford-Tate group is the full CM torus. Put
$$
h_k:=|O(x_k)|,
\qquad
H_N:=\sum_{k=1}^N h_k.
$$
Let $0<q\leq2d$, let $A\subset X$ be a Borel set with
nonempty boundary, and fix $0<r_0<1$ such that
$$
C_{\partial A}:=
\max\left\{1,\sup_{0<r<r_0}
\frac{\mu_X((\partial A)^{(r)})}{r^q}\right\}<\infty.
$$
Assume Zhang's $\varepsilon$-conjecture
\cite[Conjecture 3.5]{Zhang2005} and the $\delta$-bounds
\cite[p. 3662-3663, (3.11)-(3.12)]{Zhang2005}
for some $0<\delta<1/2$.
There exist $1/2<\beta<1$, $t>0$ and $N_0$ such that for $N>N_0$ we have
$$
\left|
\frac1{H_N}
\sum_{k=1}^N\sum_{y\in O(x_k)}\ind_A(y)
-\mu_X(A)
\right|
\ll_{X,q,\delta}
C_{\partial A}
\left(\frac{N}{H_N}\right)^{(1-\beta)tq}.
$$
The implied constant depends on $X$, $q$, and the parameters
in Zhang’s estimates. 
\end{corollary}

\begin{proof}
Let $\phi_\pi$ be an $L^2(X)$-normalized Laplace
eigenfunction with eigenvalue $\lambda_\pi>0$.
Extend $\phi_\pi$ by zero to the whole Shimura variety
$\mathcal M_U$ and write
$$
\widetilde{\phi_\pi}
=
\sum_{i=1}^{k_\pi}a_if_{i,\pi},
$$
where $(f_{i,\pi})_{1\leq i\leq k_\pi}$ is an orthonormal
basis of simultaneous Hecke eigenforms
with Laplace eigenvalue $\lambda_\pi$.
Then
$$
\sum_{i=1}^{k_\pi}|a_i|^2
=
\|\widetilde{\phi_\pi}\|_{L^2(\mathcal M_U)}^2
\ll_X1.
$$
Since $\dim_{\mathbb R}X=2d$, the basic Weyl law of
\cref{basicweyl} implies that
$$
k_\pi\leq\sum_a
\#\{j:\lambda_j(X_a)\leq\lambda_\pi\}
\ll_X(1+\lambda_\pi)^d.$$

For a CM-point $x$, put $D_x:=\operatorname{disc}(x)$,
and let $\Gamma_x$ be the finite abelian group parametrizing
its full CM-orbit $O_{\mathrm{cm}}(x)$.
Here Galois orbits are taken over the corresponding
CM reflex fields.
After choosing a base point in $O(x)$,
we may identify $O(x)$ with a subgroup of $\Gamma_x$.
Denote its annihilator in the character group by $O(x)^\perp$.
Character orthogonality, Cauchy-Schwarz and
\cite[Proposition 5.1]{Zhang2005} imply the estimate
\begin{align*}
\left|\sum_{y\in O(x)}\phi_\pi(y)\right|
&\leq
\frac1{[\Gamma_x:O(x)]}
\sum_{i=1}^{k_\pi}
|a_i|
\left|
\sum_{\chi\in O(x)^\perp}
\sum_{y\in O_{\mathrm{cm}}(x)}
\overline{\chi(y)}f_{i,\pi}(y)
\right|\\
&\ll_X
\frac{k_\pi^{1/2}}{[\Gamma_x:O(x)]}
\sup_{1\leq i\leq k_\pi}
\left|
\sum_{\chi\in O(x)^\perp}
\sum_{y\in O_{\mathrm{cm}}(x)}
\overline{\chi(y)}f_{i,\pi}(y)
\right|\\
&\ll_{X,\eta,\delta}
(1+\lambda_\pi)^A
D_x^{1/4+\delta/2+\eta},
\end{align*}
where we used $|O(x)^\perp|=[\Gamma_x:O(x)]$ and set
$A:=A_0+d/2$, with $A_0$, $\delta$ and $\eta$
as in Zhang's estimate.

Restricting to $X$ changes the size of a nonempty Galois
orbit by at most a bounded factor depending only on $X$.
Hence \cite[Theorem 3.6]{Zhang2005} gives
$$
h_k\gg_{X,\varepsilon}D_{x_k}^{1/2-\varepsilon}.
$$
Choose
$
\frac12+\delta<\beta<1,
$
and $\eta,\varepsilon>0$ such that
$
\frac{1/4+\delta/2+\eta}{1/2-\varepsilon}<\beta.
$
Consequently, H\"older's inequality yields
\begin{align*}
\left|
\frac1{H_N}
\sum_{k=1}^N\sum_{y\in O(x_k)}\phi_\pi(y)
\right|
&\ll_{X,\delta}
(1+\lambda_\pi)^A
\frac{\sum_{k=1}^N h_k^\beta}{H_N}\\
&\leq
(1+\lambda_\pi)^A
\left(\frac{N}{H_N}\right)^{1-\beta}.
\end{align*}
Put
$$
E_N:=\left(\frac{N}{H_N}\right)^{1-\beta}.
$$
Applying \cref{equidi} to the points in these orbits,
together with the eigenfunction estimate
$$
\|\phi_\pi\|_\infty
\ll_X(1+\lambda_\pi)^{(D-1)/4},
$$
where $D=\dim\M$, we obtain, for integers
$m>D/2+2A$ and $s>2D$,
\begin{align*}
\left|
\frac1{H_N}
\sum_{k=1}^N\sum_{y\in O(x_k)}\ind_A(y)-\mu_X(A)
\right|
&\ll_{X,q,m,s,\delta}
C_{\partial A}
\left(
r^q+\frac1{r^{s-q}T^{(s-2D)/2}}
\right)\\
&\quad+
C_{\partial A}E_Nr^{q-m}
\sum_{0<\lambda_\pi<T}
\frac1{\lambda_\pi^{m/2-D/4-A}}.
\end{align*}
Fix $\varepsilon_0>0$ once and for all, and choose
$$
a:=2d+\frac D2+2A,
\qquad
m:=\lceil a\rceil-1,
\qquad
s:=2D+\left\lceil\frac{2D}{\varepsilon_0}\right\rceil.
$$
The basic Weyl law we recall in \cref{basicweyl} gives
$$
\#\{\pi:0<\lambda_\pi<T\}\ll_X T^d.
$$

Since $0<a-m\leq1$, partial summation, as in the proof
of \cref{yes}, implies, for $T\geq2$,
$$
\sum_{0<\lambda_\pi<T}
\frac1{\lambda_\pi^{m/2-D/4-A}}
\ll_{X,m,A}
T^{d+D/4+A-m/2}
=
T^{(a-m)/2}.
$$
Thus the discrepancy is bounded by
$$
\ll_{X,q,m,s,\delta}
C_{\partial A}
\left(
r^q+
\frac1{r^{s-q}T^{(s-2D)/2}}
+
E_Nr^{q-m}T^{(a-m)/2}
\right).
$$
Set
$$
B:=m+\frac{s(a-m)}{s-2D}
=
a+\frac{2D(a-m)}{s-2D}.
$$
For $E_N$ sufficiently small, take
$$
r:=E_N^{1/B},
\quad
T:=r^{-2s/(s-2D)},\quad t=\frac{1}{B}>0.
$$
Choose $r_0$ sufficiently small so that $r<r_0$ implies $T\geq2$.
Then choose $N_0$ such that $r<r_0$ for every $N>N_0$.
For such $N$, each of the three terms is $O(E_N^{q/B})$.
\end{proof}

\chapter{The class of CST-Grothendieck motives}\label{class}
The purpose of this section is to introduce the class of
CST motives, which provides the motivic framework for our
effective Chebotarev-Sato-Tate equidistribution theorems. The main idea
is to relate the irreducible representations of the
motivic Sato-Tate group to representations of the
motivic Galois group and their associated motivic
$L$-functions. This allows us to combine the spectral
approximation theory developed in the previous sections
with analytic estimates for Frobenius sums developed in \cref{formalism}.

\section{Background on motivic categories and key assumptions}
Fix a number field $K$ and a motivic category $\mathcal{M}_\sim(K,\Q)$. We recall the motivic assumptions of Banaszak and Kedlaya as in \cite[Section 13]{BanaszakKedlaya2023}:
\begin{hypothesis}\label{motiv}
\begin{enumerate}[label=\textnormal{(M\arabic*)},ref=\textnormal{M\arabic*}]

\item For every smooth projective $K$-scheme $X$, the
cohomological Künneth projectors are induced by
correspondences; that is, there exist mutually orthogonal
idempotent endomorphisms
$$
\pi_X^r\in
\operatorname{End}_{\mathcal M_{\sim}(K,\Q)}(h(X)),
\qquad 0\leq r\leq 2\dim X,
$$
such that
$$
\sum_{r=0}^{2\dim X}\pi_X^r=\operatorname{id}_{h(X)},
$$
and whose Betti, de Rham, and $\ell$-adic realizations
are the projectors onto the corresponding cohomology
groups. In particular,
$$
h(X)\simeq\bigoplus_{r=0}^{2\dim X}h^r(X),
$$
where $h^r(X):=(X,\pi_X^r,0)$.

\item The category $\mathcal M_{\sim}(K,\Q)$ is a
semisimple neutral Tannakian category over $\Q$.

\item For every smooth projective $K$-scheme $X$ of pure
dimension $d$ and every integer $r$, there exists a
natural morphism in $\mathcal M_{\sim}(K,\Q)$,
$$
\star\colon h^r(X)\longrightarrow h^{2d-r}(X)(d-r),
$$
defined by the corresponding absolute Hodge cycle,
whose Betti, $\ell$-adic, and de Rham realizations
coincide with the cohomological star operator associated
with a chosen ample class,
$$
\star\colon H^r(X)\longrightarrow H^{2d-r}(X)(d-r).
$$

\item For every pair of objects
$M,N\in\mathcal M_{\sim}(K,\Q)$, there is a canonical
monomorphism from the geometric $\operatorname{Hom}$-space
$$
\operatorname{Hom}_{\mathcal M_{\sim}(\overline K,\Q)}
(M_{\overline K},N_{\overline K}),
$$
viewed as an Artin motive via its natural $G_K$-action,
to the internal Hom object
$$
\underline{\operatorname{Hom}}_{\mathcal M_{\sim}(K,\Q)}(M,N).
$$

\end{enumerate}
\end{hypothesis}
\begin{discussion}[Current state of the literature]\label{discussion0}
We give an overview of the current state of the literature regarding the previous assumptions:
\begin{enumerate}
    \item 
Assumptions \textnormal{(M1)} through \textnormal{(M4)} are known to hold unconditionally for the categories $\mathcal M_{\mathrm{ahc}}(K,\Q)$ of motives for absolute Hodge cycles, and in the category $\mathcal M_{\mathrm{mot}}(K,\Q)$ of André motivated motives. For the case of absolute-Hodge motives, the Chow--Künneth decomposition {(M1)} is defined in \cite[p. 202]{DeligneMilne1982}. As for the semisimplicity and the Tannakian structure \textnormal{(M2)}, they are shown in \cite[Proposition 6.5]{DeligneMilne1982}. The motivic star operator {(M3)} is supplied in \cite[Proposition 6.1 and pp. 197--199]{DeligneMilne1982}, and the geometric-Hom inclusion {(M4)} can be obtained from \cite[p. 215]{DeligneMilne1982} and \cite[p. 53]{Jannsen1990}. For motivated motives, {(M1)} and {(M3)} are consequences of \cite[Proposition 2.2]{Andre1996}, and of \cite[Lemma 9.2.1.3 and Remark 9.2.2.1]{Andre2004}. The semisimplicity and Tannakian nature {(M2)} of the category is proven in \cite[Theorem 9.2.3.1]{Andre2004}, and  {(M4)} comes from \cite[Scolie 2.5 and Section 4]{Andre1996}.
\item 
In the classical categories of motives, the assumptions are mostly conditional. Assumption {(M1)} is the Chow--Künneth conjecture for rational equivalence. By Jannsen's theorem \cite[Theorem 1]{Jannsen1992}, the category of numerical motives is unconditionally semisimple abelian, and semisimplicity characterizes numerical motives among adequate equivalence relations. For homological and numerical motives, as stated in \cite[Theorem 5.4.2.1]{Andre2004}, Grothendieck’s Standard Conjecture $D$ allows one to identify homological and numerical equivalences and implies Künneth Standard Conjecture $C$. Thus, if Standard Conjecture $D$ holds, then (M2) holds for numerical motives. The Lefschetz Standard Conjecture $B$ supplies the algebraic correspondences required for {(M3)}. For homological motives, as the realization functor is faithful {(M4)} follows. For full Chow motives, (M4) can fail as geometric Hom spaces can be infinite-dimensional. For example let $E/K$ be an elliptic curve over a number field and $h(E)$ the associated abelian motive. We have the following isomorphisms
$$\mathrm{Hom}_{\mathcal M_{\mathrm{rat}(\overline{K},\Q)}}(h(E_{{\overline{K}}}),\mathbf{1})=\mathrm{CH}^1(E_{\overline{K}})_\bb{Q}=\Q\oplus (E(\overline{K})\otimes_{\bb{Z}} \bb{Q}),$$ 
which has infinite $\overline{\Q}$-dimension as $E(\overline{K})$ has infinite rank.
\end{enumerate}
\end{discussion}

\section{The main conjectures}
Fix a number field $K\subset\bb C$.
Let $\mathcal M_\sim(K,\Q)$ be a motivic category satisfying Hypothesis
\ref{motiv}.
Given a nonzero homogeneous motive
$M\in\mathcal M_\sim(K,\Q)$ of weight $n$, that is, a direct
summand of a motive of the form $h^r(X)(m)$, let
$\l M\rr^\otimes$ be the smallest strictly full Tannakian
subcategory containing $M$. Define the motivic Galois group
of $M$ by
$$
G_M:=\operatorname{Aut}^{\otimes}
\bigl(H_B|_{\l M\rr^\otimes}\bigr).
$$
By Tannakian duality, $\l M\rr^\otimes$ is equivalent to
the category of finite-dimensional algebraic representations
of $G_M$.  The semisimplicity assumption \textup{(M2)} implies that
$G_M$ is reductive.

Let
$$
\psi_\sim:M\otimes M\longrightarrow\bb L^{\otimes n}
$$
be a motivic polarization supplied by \textup{(M2)-(M3)} together with the motivic Poincar\'e duality, and let
$$
\psi:=H_B(\psi_\sim):V\otimes V\longrightarrow\Q(-n)
$$
be its Hodge realization
\cite[Remark 13.6, equations (13.4)-(13.7)]
{BanaszakKedlaya2023}.
It induces a closed immersion of algebraic groups $
G_M\hookrightarrow\operatorname{GIso}(V,\psi),
$
where $\operatorname{GIso}(V,\psi)$ is the algebraic group
of similitudes. Let
$$
\chi:\operatorname{GIso}(V,\psi)\longrightarrow\bb G_m
$$
be the corresponding similitude character. Define the motivic Serre group by
$$
\operatorname{AST}_K(M)
:=
\ker\bigl(\chi|_{G_M}\bigr)
=
G_M\cap\operatorname{Iso}(V,\psi).
$$

For every prime $\ell$, we let $G_{\ell,K}^{\mathrm{alg}}$
denote the Zariski closure of the image of the $\ell$-adic
Galois representation attached to $M$.
Under the Betti-\'etale identification
$H_\ell(M)$ with $H_B(M)\otimes_\Q\Q_\ell$, define
$$
G_{\ell,K,1}^{\mathrm{alg}}
:=
G_{\ell,K}^{\mathrm{alg}}\cap\ker(\chi_{\Q_\ell}).
$$

Let $\mathbb S=\operatorname{Res}_{\mathbb C/\mathbb R}\Gm$
be the Deligne torus, and let
$h_M:\mathbb S\longrightarrow\operatorname{GL}(H_B(M))_\R$
define the Hodge structure on $V$.
The Mumford-Tate group $\operatorname{MT}(M)$ is the
$\Q$-Zariski closure of $h_M(\mathbb S(\R))$
in $\operatorname{GL}(H_B(M))$.

Under these identifications, we recall the following
conjectures associated with $M$,
see \cite[Sections 15-16]{BanaszakKedlaya2023}.
\begin{enumerate}
\item \emph{Algebraic Sato-Tate conjecture.}
For every prime $\ell$,
$$
G_{\ell,K,1}^{\mathrm{alg}}
=
\operatorname{AST}_K(M)_{\Q_\ell}.
$$

\item \emph{Serre's conjecture.}
The natural inclusion of the Mumford--Tate group into
the motivic Galois group identifies
$$
\operatorname{MT}(M)=G_M^\circ.
$$

\item \emph{Mumford-Tate conjecture.}
For every prime $\ell$,
$$
\bigl(G_{\ell,K}^{\mathrm{alg}}\bigr)^\circ
=
\operatorname{MT}(M)_{\Q_\ell}.
$$
\end{enumerate}
\begin{discussion}[Current state of the conjectures (1)-(3)]\label{discussion1}
     The Mumford-Tate conjecture remains open in general, even for
abelian varieties over number fields. It holds for CM abelian
varieties \cite[Theorem 2.14]{Lombardo2016}, all abelian varieties
of dimension at most three \cite[Theorem 6.1]{Lombardo2016},
and absolutely simple abelian varieties of prime dimension
\cite[Theorem 2.18]{Lombardo2016}. Commelin proved that the abelian motives satisfying the
Mumford-Tate conjecture form a Tannakian subcategory
\cite[Theorem 10.3]{Commelin2019Products}.
For abelian motives, Cantoral-Farf\'an and Commelin proved
that the Mumford-Tate conjecture implies the algebraic
Sato-Tate conjecture \cite[Section 3]{CantoralFarfanCommelin2022}. Tankeev proved the Mumford-Tate conjecture for the second
cohomology motive $h^2(X)$ of a $K3$ surface over a number
field. For projective hyper-K\"ahler varieties, a recent preprint
of Ito and Zou announces a proof in every degree
over fields finitely generated over $\mathbb Q$, see \cite[Theorem A]{ItoZou2026}.
As for Serre's conjecture, it holds unconditionally for abelian motives $\mathcal{M}_{\mathrm{Ab}}(K,\Q)$, see \cite[Theorem 4.13 and Remark 4.14(1)]{Commelin2019Products}.
\end{discussion}
\section{Strong compatibility of Galois systems and CST motives}
Let $G$ be a reductive group over a number field $E$, with a fixed
embedding $G\hookrightarrow\GL_n$. Let $(\rho_\lambda)_\lambda$
be a family of continuous Galois representations such that, for
each finite place $\lambda$ of $E$, the representation
$\rho_\lambda$ has image in $G(E_\lambda)$.
Recall that $(\rho_\lambda)_\lambda$ is {strongly
$G$-compatible} if, for each finite place $v$ of $K$, there exists
a Frobenius-semisimple $G$-valued Weil--Deligne parameter
$\mathrm{WD}_v(\rho)=(r_v,N_v)$, where
$r_v\colon W_{K_v}\to G(\overline E)$ is a continuous homomorphism of the local Weil group
with $G(\overline E)$ given the discrete topology, and
$N_v\in\operatorname{Lie}(G)_{\overline E}$ is a nilpotent operator,
verifying the relation
$\operatorname{Ad}(r_v(w))N_v=|w|_vN_v$ for $w\in W_{K_v}$,
such that, for any finite place $\lambda$ of $E$ whose residue
characteristic differs from that of $v$, and any embedding
$\overline E\hookrightarrow\overline E_\lambda$ extending
$E\hookrightarrow E_\lambda$, the semi-simplified parameters satisfy
$$
\mathrm{WD}\left(\rho_\lambda|_{G_{K_v}}\right)^{\mathrm{F\text{-}ss}}
\simeq
(r_v,N_v)\otimes_{\overline E}\overline E_\lambda.
$$
Here equivalence means conjugacy by $G(\overline E_\lambda)$.
We speak of linear strong
compatibility if $G$ is the full linear group.
In the linear case, this is the local Weil--Deligne compatibility
used in Deligne's formalism of motivic $L$-functions
\cite[Section 5.2]{DeligneValeursL} and
\cite{Deligne1973Constantes}. 

Recall that we say that the associated linear system is
{strictly compatible} if there exists a finite subset $S$
of finite places of $K$ such that, for any $\lambda$, the
representation $\rho_\lambda$ is unramified outside
$S\cup S_\lambda$ and, for any $v\notin S\cup S_\lambda$,
the polynomial
$P_v(T):=\det\left(1-\rho_\lambda(\Frob_v)T\right)\in E[T]$
is independent of $\lambda$. Here $S_\lambda$ denotes the
places of $K$ having the same residue characteristic as
$\lambda$. The smallest such $S$ is the {exceptional set}
of the system \cite[Chap.I, Section 2.3]{Serre1989}.
If $(\rho_\lambda)_\lambda$ is strongly $G$-compatible and
the ramification set of its Weil--Deligne data
$$
S_M:=\{v:\ r_v|_{I_v}\neq1\ \text{or}\ N_v\neq0\}
$$
is finite, then the associated linear system is strictly
compatible and its exceptional set is contained in $S_M$.

The following lemma describes the compatibility of the
Weil-Deligne construction with pullback along
representations of $G$. It will allow us to give a global description of the local factors of the
associated motivic $L$-functions in terms of the motivic invariants of $M$.

\begin{lemma}\label{comp} Let $(\rho_\lambda)_\lambda$ be a strongly $G$-compatible system. Let $\sigma\colon G\to \operatorname{GL}(V) $ be an algebraic representation of $G$ defined over $E$. Then the system $ \bigl(\sigma_\lambda\circ\rho_\lambda\bigr)_\lambda $ is linearly strongly compatible. Its local Weil-Deligne parameter at $v$ is $\sigma_*\mathrm{WD}_v(\rho) := \bigl(\sigma\circ r_v,d\sigma(N_v)\bigr). $ \end{lemma}
\begin{proof}
For each finite place $v$ of $K$, let
$I_v\subset G_{K_v}$ denote the inertia subgroup of the local absolute Galois group. The relation for the Weil-Deligne parameter is $\sigma$-invariant since $$ \begin{aligned} \operatorname{Ad}\bigl(\sigma(r_v(w))\bigr)d\sigma(N_v) &= d\sigma\left(\operatorname{Ad}(r_v(w))N_v\right) \\ &= |w|_v\,d\sigma(N_v). \end{aligned} $$ Furthermore, Grothendieck’s monodromy theorem allows us to write $$ \rho_\lambda(\tau) = r_{v,\lambda}(\tau) \exp\bigl(t_\lambda(\tau)N_{v,\lambda}\bigr), $$
   for $\tau\in I_v$.
   Applying $\sigma_\lambda$ results in $$ \sigma_\lambda(\rho_\lambda(\tau)) = \sigma_\lambda(r_{v,\lambda}(\tau)) \exp\bigl(t_\lambda(\tau)d\sigma(N_{v,\lambda})\bigr), $$ hence $$ \mathrm{WD}\left( (\sigma_\lambda\circ\rho_\lambda)|_{G_{K_v}} \right) \simeq \sigma_* \mathrm{WD}\left( \rho_\lambda|_{G_{K_v}} \right). $$ Since homomorphisms of algebraic groups preserve Jordan
decompositions, taking the
semisimplification yields that
$$
\mathrm{WD}\left((\sigma_\lambda\circ\rho_\lambda)|_{G_{K_v}}\right)
^{\mathrm{F\text{-}ss}}
\simeq
\bigl(\sigma\circ r_v,d\sigma(N_v)\bigr)
\otimes_{\overline E}\overline E_\lambda
$$
for every embedding
$\overline E\hookrightarrow\overline E_\lambda$ over $E$, proving the assertion.
\end{proof}
    \begin{remark}
    If $(\rho_\lambda)_\lambda$ and $(\rho'_\lambda)_\lambda$ are strongly
$G$-compatible and strongly $H$-compatible, respectively, their pair
$(\rho_\lambda,\rho'_\lambda)_\lambda$ is strongly $G\times H$-compatible,
with local parameter
$
((r_v,r'_v),(N_v,N'_v)).
$
For linear representations, their tensor product is linearly strongly
compatible, with local parameter
$$
\bigl(r_v\otimes r'_v,N_v\otimes1+1\otimes N'_v\bigr).
$$
Linear strong compatibility is preserved under direct sums,
tensor products and duals.  Moreover, if $(\rho_\lambda)_\lambda$ is strongly
$G$-compatible and $W\subseteq V$ is an $E$-rational
$G$-stable subspace, then the induced systems on $W$ and
$V/W$ are linearly strongly compatible by \cref{comp}.

\end{remark}
\begin{discussion}[Known cases of strong compatibility]\label{discussion2}
   It is well known that linear strong compatibility holds
for Artin motives \cite[Section 6]{DeligneValeursL}.
Another relevant example is provided by motives attached
to Hilbert modular forms under the conditions of Carayol
\cite[Théorème A]{Carayol1986}.
In the case of the abelian motive $h^1(A)$, assuming
that the Galois representations take values in the
Mumford--Tate group of $A$, the system
$(\rho_{A,\ell})_\ell$ is strongly compatible with values
in this group if $A/K$ is semistable at every finite
place \cite[Theorem 1.2]{KisinZhou2025}.
Noot had previously established a weaker form under
additional hypotheses \cite[Introduction]{Noot2013}.
As a result for an arbitrary abelian variety $A/K$, strong compatibility holds after passing to a finite extension $L/K$ such that the base change $A_L$ has semistable reduction. Kisin and Zhou formulate the corresponding conjecture
for general smooth proper varieties
\cite[Conjecture 1.1]{KisinZhou2025}.

Now let $f=\sum_{n\ge1}a_nq^n$ be a normalized cuspidal
newform of weight $k\ge2$, level $N$, nebentypus $\chi$,
and coefficient field $E_f$. Scholl constructs a pure
Grothendieck motive $M(f)$ over $\Q$, of rank $2$ and
weight $k-1$, with coefficients in $E_f$, whose
$\lambda$-adic realization is the representation
$\rho_{f,\lambda}$ attached to $f$
\cite[Theorem 1.2.4]{Scholl1990}.
These realizations form a linearly strongly compatible
system \cite[Théorème A]{Carayol1986}; see also
\cite[Theorems 0 and 1]{Saito2009}.

More generally, Fontaine conjectures linear strong
compatibility for the local Weil-Deligne realizations
of mixed motives, including the $p$-adic realization
\cite[Section 2.4.7]{Fontaine1994PotentiallySemistable}.
A $G_M$-valued analogue for absolute Hodge motives
with semistable reduction, for $\ell\ne p$, is formulated
by Laskar \cite[Question 1.2]{Laskar2013}.

For K3 surfaces and, more generally, hyper-Kähler
varieties with second Betti number $b_2\ge4$, a recent preprint by Ito and Zou proves, after a
suitable finite extension of the base field, rational
strong compatibility of the Mumford-Tate-group-valued
Weil-Deligne representations at every finite place,
including $\ell=p$
\cite[Theorem C and Theorem 5.8]{ItoZou2026}.
In particular, for a K3 surface $X$ (with $b_2=22$) , this also gives the corresponding
Mumford-Tate-group-valued compatibility for the transcendental motive $t_2(X)(1)$
after the same base extension.
If Serre's conjecture holds, then we have a natural monomorphism $\iota:\mathrm{MT}(M)\hookrightarrow G_M$ and as a result $\mathrm{MT}(M)$-compatibility implies $G_M$-compatibility. 
\end{discussion}
We are now ready to introduce the class of CST motives.
\begin{definition}
    Let $M\in \mathcal{M}_\sim(K,\Q)$ be a homogeneous polarized motive. We call $M$ a connected Sato-Tate motive (or CST-motive) if the following hold:
    \begin{enumerate}
        \item the motivic Sato-Tate group $\st_K(M)$ is connected,
        \item the motive $M$ verifies the algebraic Sato-Tate conjecture,
        \item the motive $M$ verifies Serre's conjecture,
        \item the system of $\ell$-adic realizations
$
(\rho_{M_,\ell})_\ell$
is strongly $G_M$-compatible. 
    \end{enumerate}
\end{definition}
For abelian varieties with connected Sato-Tate group,
Bucur, Fité and Kedlaya prove, under the Mumford-Tate
conjecture, that a maximal torus of the Sato-Tate group
is generated by its Cartan Hodge circles
\cite[Theorem 3.5]{BucurFiteKedlaya2024}. We shall prove a weaker result for CST-motives,
which suffices to control the archimedean parameters of
the associated motivic $L$-functions in terms of
the highest weights of the corresponding representations:
\begin{lemma}\label{MT}
Let $M$ be a CST-motive.
    Let $\bb S=\mathrm{Res}_{\bb C/\bb R} \Gm$ be the Deligne torus and $h_M:\bb S\rightarrow \GL({H_B(M)}_\R)$ the Hodge structure of the Betti realization.
    The image of the restriction $\kappa:=h_{M| U(1)}$ lies in a conjugate of the Sato-Tate group. Moreover its algebraic Sato-Tate group $\mathrm{AST}_K(M)$ and the motivic Galois group $G_M$ are connected as well.
\end{lemma}
\begin{proof}
According to \cite[Theorem 16.8, (16.11)]{BanaszakKedlaya2023} Serre's conjecture implies that 
$$
\kappa(U(1))
 \subset H(\psi)(\mathbb R)
 \subset  \operatorname{AST}_K(M)(\mathbb C),
$$
where $H(\psi)$ denotes the Hodge group of the polarized Hodge structure.
The set $\kappa(U(1))$ is compact. According to
\cite[Proposition 1.2.3(a),(b)]{Schmid2004},
$\kappa(U(1))$ lies inside a maximal compact subgroup of
$\operatorname{AST}_K(M)(\mathbb C)$ and as a result in
a conjugate to the fixed $\operatorname{ST}_K(M)$. By \cite[Theorem 16.8, equation (16.13)]{BanaszakKedlaya2023},
the connectedness of $\st_K(M)$ implies that
$\operatorname{AST}_K(M)$ is connected.
Moreover, \cite[Corollary 16.7]{BanaszakKedlaya2023}
provides a surjection
$
\pi_0(\operatorname{AST}_K(M))
\twoheadrightarrow \pi_0(G_M),
$
so $G_M$ is connected as well.
\end{proof}
To develop an $L$-function formalism, we shall realize
the irreducible representations of the Sato-Tate group of a CST-motive $M\in \mathcal{M}_\sim(K,\Q)$ as motivic objects. More precisely, to each
$\sigma\in\Irr(\st_K(M))$, we shall associate a simple object
$
M_\sigma\in\l M\rr^\otimes\otimes_\Q E,
$
over a number field $E$, with $M_\sigma$
unique up to a Tate twist.
We begin with the following auxiliary
group-theoretic lemma:
\begin{lemma}\label{extension}
Let $k$ be a field of characteristic $0$. Let $G$ be a connected reductive split $k$-group, let $\chi:G\twoheadrightarrow\Gm$ be a surjective character, and assume that $H:=\ker\chi$ is connected. Then the restriction functor of algebraic representations
$$
\Rep_k(G)\longrightarrow\Rep_k(H)
$$
is essentially surjective and restricts to a functor of the irreducible algebraic representations
$$
\Irr(G)\longrightarrow\Irr(H).
$$
Moreover, any two extensions of an irreducible representation of $H$ differ by an integral power of $\chi$.
\end{lemma}

\begin{proof}
Let $T\subset G$ be a maximal torus. Define $T_H:=T\cap H$. We have an exact sequence
$$1\longrightarrow H\longrightarrow G\longrightarrow G/H\longrightarrow 1.$$
By \cite[Theorem 14.9]{Milne2017AlgebraicGroups} this yields an exact sequence of lattices
$$1\longrightarrow \bb Z\chi\longrightarrow X^*(T)\longrightarrow X^*(T_H)\longrightarrow 1.$$
In particular every dominant weight $\lambda\in X^*(T_H)^+$ lifts to $\tilde{\lambda}\in X^*(T)^+$, in particular the irreducible representation $\pi_\lambda$ of highest weight $\lambda$ lifts to an irreducible representation $\pi_{\tilde{\lambda}}$.
Finally, two lifts of $\lambda$ differ by $j\chi$ for some
$j\in\Z$. The corresponding irreducible representations of
$G$ therefore differ by tensoring with $\chi^j$.
\end{proof}
\begin{proposition}\label{motivic}
Let $M$ be a CST-motive, and let $E\subset\bb C$ be a
number field splitting the motivic Galois group $G_M$.
Then every irreducible continuous representation
$\sigma\in\Irr(\st(M))$ admits a motivic lift to a simple
object
$
M_\sigma\in\langle M\rangle^\otimes\otimes_\Q E
$ of the category of motives $\mathcal M_\sim(K,E)$.
Moreover, any two such lifts differ by a Tate twist,
and the $\lambda$-adic realizations
$(\rho_{M_\sigma,\lambda})_\lambda$ form a linearly
strongly compatible system with exceptional set contained
in the same finite set $S_M$.
For every finite place $\lambda$ of $E$, we have
$$
\rho_{M_\sigma,\lambda}
\simeq
\widehat{\sigma}_{E,\lambda}\circ\rho_{M_E,\lambda},
$$
where $\widehat{\sigma}_E$ is an algebraic extension
of $\sigma$ to $(G_M)_E$, unique up to tensoring with
an integral power of $\chi$.
\end{proposition}

\begin{proof}
Recall that
$
\mathrm{AST}_K(M)=\ker(\chi),
$
where $\chi:G_M\to\mathbb G_m$ is the similitude character. 
By \cref{MT}, both $\mathrm{AST}_K(M)$ and $G_M$
are connected.
Now let
$
\sigma\colon\st_K(M)\longrightarrow\GL(W_\sigma)
$
be an irreducible finite-dimensional continuous complex
representation. Since $\mathrm{AST}_K(M)$ is connected
and reductive, and $\st_K(M)$ is a compact real form of
$\mathrm{AST}_K(M)_{\bb C}$, 
Weyl's unitary trick \cite[Theorem 2.11(2)]{AdamsVanLeeuwenTrapaVogan2020} implies that $\sigma$ extends uniquely to an irreducible
algebraic representation
$
\widetilde{\sigma}_{\bb C}\colon
\mathrm{AST}_K(M)_{\bb C}\longrightarrow\GL(W_\sigma).
$
Since $E$ splits $G_M$, it also splits its connected
subgroup $\mathrm{AST}_K(M)=\ker(\chi)$.
By highest-weight theory, $\widetilde{\sigma}_{\bb C}$
therefore admits an absolutely irreducible model
$
\widetilde{\sigma}_E\colon
(\mathrm{AST}_K(M))_E\longrightarrow\GL(W_{\sigma,E})
$
such that
$
\widetilde{\sigma}_E\otimes_E\bb C
\simeq\widetilde{\sigma}_{\bb C}.
$

We now extend $\widetilde{\sigma}_E$ to a representation $\widehat{\sigma}_E$ of $G_M$. If $\chi=1$ then $G_M=\mathrm{AST}_K(M)$ and $\widetilde{\sigma}_E=\widehat{\sigma}_E$. Otherwise, since $G_M$ is connected and $\chi\ne1$, the character
$\chi:G_M\to\mathbb G_m$ is surjective and by \cref{extension}, $\widetilde{\sigma}_E$ extends
to an irreducible algebraic representation
$
\widehat{\sigma}_E\colon
(G_M)_E\longrightarrow\GL(W_{\sigma,E}).
$
Tannakian equivalence implies that $\widehat{\sigma}_E$
corresponds to a simple object
$
M_\sigma\in\langle M\rangle^\otimes\otimes_\Q E.
$

Let $M'_\sigma$ be another motivic lift in this category,
corresponding to an algebraic representation
$\widehat{\sigma}'_E$ of $(G_M)_E$.
If $\chi=1$, then $\mathrm{AST}_K(M)=G_M$, so the
two lifts are isomorphic.
Otherwise, \cref{extension} implies that, for some $a\in\Z$,
$
\widehat{\sigma}'_E
\simeq\widehat{\sigma}_E\otimes\chi^a.
$
Writing $w$ for the weight of $M$, the motivic polarization
identifies the character $\chi$ with the Tate object
$\mathbf1(-w)$. Consequently Tannakian equivalence yields
$
M'_\sigma\simeq M_\sigma(-aw).
$

Compatibility of the realization functors with Tannakian
equivalence gives, for every finite place $\lambda$ of $E$,
$$
\rho_{M_\sigma,\lambda}
\simeq
\widehat{\sigma}_{E,\lambda}\circ\rho_{M_E,\lambda}.
$$
Strong $G_M$-compatibility of $M$ and \cref{comp} therefore
give linear strong compatibility of these realizations.
For $v\notin S_M$, the common parameter
$\mathrm{WD}_v(\rho_M)=(r_v,N_v)$ satisfies
$r_v|_{I_v}=1$ and $N_v=0$, so its image under
$\widehat{\sigma}_E$ is also unramified.
Thus the exceptional set of the lifted system is
contained in $S_M$.
\end{proof}
\section{Motivic $L$-functions formalism for CST motives}
We follow the motivic $L$-function formalism of Deligne
\cite{DeligneValeursL}, together with the conjectural
analytic properties of motivic $L$-functions including the motivic Riemann hypothesis
\cite[Section 2.4]{Cohen2019LFunctions}, and refer to
\cite{Serre1994Motives} for the conjectural properties
of motivic Galois groups and their associated
$\ell$-adic representations.
Recall that the analytically normalized $L$-functions
attached to pure motives are conjectured to belong to
the Selberg class; see
\cite[Sections 4.5-5]{FarmerPitaleRyanSchmidt2019}.
For the axiomatic theory of the Selberg class,
we refer to \cite{Kaczorowski2006Selberg}.

We adapt the motivic Riemann hypothesis formulated
for abelian varieties in \cite{BucurFiteKedlaya2024}
to the class of CST motives.
We recall the definition of motivic $L$-functions
and provide proofs of the technical results used
in this section in \cref{formalism}.

\begin{conjecture}[Conjectural motivic $L$-function formalism]
\label{grh}
Let $M$ be a pure motive in $\mathcal M_\sim(K,E)$ of weight $w$,
and let $d=\dim_E\Hom(1,M(w/2))$, with $d=0$ when $w$ is odd.
Fix an embedding $\iota:E\hookrightarrow\bb{C}$.
Writing $\Lambda(\iota,M,s)=L_\infty(\iota,M,s)L(\iota,M,s)$
and letting $q_M$ denote the arithmetic conductor, define
the normalized completed $L$-function
$$
F_M(s):=q_M^{s/2}\Lambda(\iota,M,s+w/2).
$$

\begin{enumerate}
\item \emph{(Meromorphic continuation)}
The function $s^d(1-s)^dF_M(s)$ extends to an entire
function on $\mathbb C$ of order at most $1$, and is
nonzero at $s=0$ and $s=1$.

\item \emph{(Functional equation)}
There exists
$\varepsilon(\iota,M)\in\mathbb C$  determined by the local
Weil-Deligne parameters at finite places and the
archimedean Weil parameters, with
$|\varepsilon(\iota,M)|=1$, such that
$$
F_M(s)=\varepsilon(\iota,M)F_{M^\vee}(1-s)
\qquad(s\in\mathbb C).
$$

\item \emph{(Motivic Riemann hypothesis)}
Every non-trivial zero $\varrho$ of $F_M(s)$ satisfies
$\Re(\varrho)=\frac12$.

\item \emph{(Motivic Ramanujan conjecture)}
At any finite place $v$ of $K$, possibly ramified, write
$$
L_v(\iota,M,s+w/2)
=
\prod_{j=1}^{n_v}
\left(1-\alpha_{v,j}(Nv)^{-s}\right)^{-1},
$$
$$
L_v(\iota,M^\vee,s-w/2)
=
\prod_{j=1}^{n'_v}
\left(1-\alpha'_{v,j}(Nv)^{-s}\right)^{-1}.
$$
Then $|\alpha_{v,j}|\leq1$ and $|\alpha'_{v,j}|\leq1$ for every $j$.
\end{enumerate}
\end{conjecture}
\begin{discussion}[Verifying the CST conjectures]\label{discussion3}
In practice, we shall often work in the category of motives
for absolute Hodge cycles $\mathcal{M}_{\mathrm{ahc}}(K,\Q)$,
which satisfies the motivic Hypothesis \ref{motiv}.
Our main case of interest is its subcategory of abelian
motives $\mathcal{M}_{\mathrm{Ab}}$. For these motives,
Serre's conjecture, identifying the connected motivic
Galois group with the Mumford-Tate group, is known.

\begin{enumerate}
\item
As seen in Discussion \ref{discussion1}, the Mumford-Tate
conjecture implies the algebraic Sato-Tate conjecture
for abelian motives. The Mumford-Tate conjecture is known
for a wide range of abelian motives, and those satisfying
it form a Tannakian subcategory.
We shall use this fact in our arithmetic applications
to the transcendental motive $t_2(X)(1)$ of a K3 surface.  When $M$ is abelian, the motives $M_\sigma$, for
$\sigma\in\Irr(\st_K(M))$, belong by the construction in \cref{motivic} to the
Tannakian category generated by $M_E$ and Tate motives.
They are therefore abelian motives with coefficients
in $E$ and satisfy the Mumford-Tate conjecture whenever
$M$ does.
Since $N_\rho$ is an Artin motive, the motives
$$
M_\sigma\otimes N_\rho,
\qquad (\sigma,\rho)\in\Irr(\st_K(M))\times \mathrm{Irr}(\mathrm{Gal}(L/K)),
$$
are again abelian motives.

\item
As in Discussion \ref{discussion2}, strong compatibility
with values in the Mumford-Tate group is known for
abelian varieties at places of semistable reduction,
after the finite extension needed to obtain
Mumford-Tate-valued Galois representations
\cite[Theorem 1.2]{KisinZhou2025}.
As specified before, it is also known, after a suitable finite base extension,
for hyper-K\"ahler varieties with $b_2\geq4$
\cite[Theorem C]{ItoZou2026}.
Via the natural inclusion
$\operatorname{MT}(M)\hookrightarrow G_M$, these results
give the corresponding $G_M$-valued compatibility.
We shall discuss in the next section how to address
strong compatibility in our setting.

\item
Recall that a Weil-Deligne representation $(V,r,N)$
is pure of weight $w$ if the eigenvalues of $r(\Phi)$
on $\operatorname{gr}^S_jV$ are Weil $q^{w+j}$-numbers,
where $\Phi$ is geometric Frobenius at $\pp$ and $q$ is the
residue-field cardinality, and $S_\bullet V$ is the
monodromy filtration centered at $0$, characterized by
$$
N(S_jV)\subseteq S_{j-2}V,
\qquad
N^j:\operatorname{gr}^S_jV
\xrightarrow{\sim}\operatorname{gr}^S_{-j}V
\quad(j\geq0).
$$
Here a Weil $q^a$-number is a nonzero algebraic number
whose complex conjugates all have absolute value
$q^{a/2}$.
Purity implies the Ramanujan bound for the normalized
local Euler roots, while strong compatibility ensures
that the local factors are independent of $\ell$. 
Local purity for abelian varieties along with other cases is established in
\cite[Section 1.13, Proposition 4]{Blasius2006}.
It is preserved under tensor products, duals, Tate
twists and direct summands
\cite[Section 1.10, Proposition 1]{Blasius2006}.
Thus, when $M$ is abelian, every $M_\sigma$ has locally
pure Weil-Deligne parameters.
Since Artin motives have locally pure parameters of
weight zero, the same holds for $M_\sigma\otimes N_\rho$,
with the weight of $M_\sigma$.
Consequently, the normalized local Euler roots of these
twists satisfy the Ramanujan bound at every finite place.
\end{enumerate}
\end{discussion}
Fix a CST-motive $M$, together with motivic lifts
$(M_\sigma)_{\sigma\in\operatorname{Irr}(\st_K(M))}$
and Artin motives
$(N_\rho)_{\rho\in\operatorname{Irr}(\operatorname{Gal}(L/K))}$, defined over a common number field $E$ over which $G_M$ splits.
For a fixed embedding $\iota:E\hookrightarrow\mathbb C$,
\cref{motivic} shows that the normalized $L$-function
$F_{M_\sigma\otimes N_\rho}$ is independent of the
choice of motivic lift. We denote its analytic conductor
by $Q(\sigma\otimes\rho,T)$.

Fix a maximal torus $T\subset \mathrm{ST}_K(M)$ with Weyl group $W$. The next lemma illustrates how the motivic Ramanujan conjecture together with strong compatibility allows control over the conductor through an explicit constant depending on the Hodge structure. We provide explicit computations in \cref{formalism} in some special cases.
\begin{lemma}\label{conduct}
\begin{enumerate}
\item The arithmetic conductor satisfies  $q_{M_\sigma}\mid q_M^{d_\sigma}$.
\item
There exists a constant $C_H>0$, depending only on the
polarized Hodge structure $h_M$, such that
$$
H_\infty(M_\sigma)\leq1+C_H\|\lambda_\sigma\|.
$$
Moreover,
$$
\log Q(\sigma\otimes\rho,T)
\ll n_Kd_\sigma d_\rho
\log\bigl(q_Mq_\rho(3+T+C_H\|\lambda_\sigma\|)\bigr),
$$
where $\lambda_\sigma$ is the highest weight of $\sigma$
and the implied constant is absolute.
\end{enumerate}
\end{lemma}
\begin{proof}
\begin{enumerate}
    \item 
Given a place $v$ of $K$ let
$$
\mathrm{WD}_v(M)= (r_v,N_v), \qquad
\mathrm{WD}_v(M_\sigma)
=
\left(
\widehat{\sigma}_{E,\lambda}\circ r_v,
d\widehat{\sigma}_{E,\lambda}(N_v)
\right)
$$
be the Weil-Deligne parameters.
Denote by $L/K_v$ a finite Galois extension such that $I_L$ acts
trivially through $r_v$, and by $(G_{v,i})_{i\geq0}$ its
lower ramification groups. In particular, $G_{v,0}$ is the
inertia subgroup of $\operatorname{Gal}(L/K_v)$.

Then we get
$$
\begin{aligned}
\operatorname{Sw}
\left(\widehat{\sigma}_{E,\lambda}\circ r_v\right)
&=
\sum_{i\geq1}
\frac{|G_{v,i}|}{|G_{v,0}|}
\operatorname{codim}
V_\sigma^{(\widehat{\sigma}_{E,\lambda}\circ r_v)(G_{v,i})}
\\
&\leq
d_\sigma
\sum_{\substack{i\geq1 \\ r_v(G_{v,i})\neq\{1\}}}
\frac{|G_{v,i}|}{|G_{v,0}|}.
\end{aligned}
$$
If $r_v(G_{v,i})\neq\{1\}$, then
$
\operatorname{codim}
V_{M,v}^{r_v(G_{v,i})}
\geq1.
$
Thus,
$
\operatorname{Sw}
\left(\widehat{\sigma}_{E,\lambda}\circ r_v\right)
\leq
d_\sigma\operatorname{Sw}(r_v).
$

Moreover, if
$
n_M - \dim(\ker N_v)^{I_v} = 0,
$
then $N_v=0$ and $r_v|_{I_v}=1$. Hence
$
d_\sigma -
\dim(\ker d\widehat{\sigma}_{E,\lambda}(N_v))^{I_v} =0.
$
Otherwise,
$
n_M - \dim(\ker N_v)^{I_v}\geq1,
$
and the corresponding codimension for $M_\sigma$ is at most
$d_\sigma$. In any case, we have
$
\begin{aligned}
d_\sigma -
\dim(\ker d\widehat{\sigma}_{E,\lambda}(N_v))^{I_v}
\leq
d_\sigma
\left(
n_M - \dim(\ker N_v)^{I_v}
\right).
\end{aligned}
$
From these inequalities we get
$
a_v(M_\sigma)
\leq
d_\sigma a_v(M).
$
Hence
$$
\mathfrak f(M_\sigma)
\mid
\mathfrak f(M)^{d_\sigma}.
$$
Since $n_M\geq1$, we obtain
$$
\frac{q_M^{d_\sigma}}{q_{M_\sigma}}
=
|d_K|^{d_\sigma(n_M-1)}
\mathrm N_{K/\mathbb Q}
\left(
\mathfrak f(M)^{d_\sigma}
\mathfrak f(M_\sigma)^{-1}
\right)
\in\mathbb N.
$$
So $q_{M_\sigma}\mid q_M^{d_\sigma}$, which proves the first assertion.
\item The normalized archimedean parameters of $M_\sigma$ are bounded by
$H_\infty(M_\sigma)$. 
Let $g\in \mathrm{AST}_K(M)(\bb{C})$ such that the image of $g\kappa_{|U(1)} g^-1$ lies in the Sato-Tate group as in \cref{MT}. By highest weight theory
\begin{align*}
    \sum_{p,q}h^{p,q}(M_\sigma)u^{p-q}
&=\operatorname{Tr}\bigl(\sigma(g\kappa(u)g^{-1}\bigr)\
\\&=\sum_{\mu\in\mathrm{Wt}(\lambda_\sigma)}
m_\sigma(\mu)u^{\langle\mu,\kappa\rangle},
\end{align*}
where $\mathrm{Wt}(\lambda_\sigma)=\{\mu\in X^*(T),V_\lambda[\mu]\neq 0\}$ which lies in the convex hull of $W\lambda_\sigma$.
Write
$\kappa(e^{it})=\exp(tH_{\kappa,\tau})$ and set
$$
C_H:=2\pi\max_{\tau}\|H_{\kappa,\tau}\|,
$$
where $\tau$ ranges over field embeddings $K\hookrightarrow \bb{C}.$
Then $|p-q|\leq C_H\|\lambda_\sigma\|$, and hence
$H_\infty(M_\sigma)\leq1+C_H\|\lambda_\sigma\|$. From the definition of $ Q(\sigma,T)$, for every
$T\geq1$ we get
$$
\begin{aligned}
\log Q(\sigma,T)
&=
\log\mathcal Q_{\iota,M_\sigma}(T)
\\
&\leq
\log q_{M_\sigma}
+
2n_Kd_\sigma
\log\left(
3+T+H_\infty(M_\sigma)
\right)
\\
&\leq
d_\sigma\log q_M
+
2n_Kd_\sigma
\log\left(
3+T+H_\infty(M_\sigma)
\right)
\\
&\ll
n_Kd_\sigma
\log
(q_M(3+T+H_\infty(M_\sigma))) \\
&\ll n_Kd_\sigma\log(q_M(3+T+C_H\|\lambda_\sigma\|)).
\end{aligned}
$$
The same bound holds for the $L$-function $L(\iota,N_\rho,s)=L(\rho,s)$ by the classical theory of Artin-functions.
As the archimedean parameters of Artin $L$-functions are $O(1)$ bounded, \cref{tensor} implies that
$$\log Q(\sigma\otimes\rho,T)\ll d_\sigma \log Q(\rho,T)+d_\rho \log Q(\sigma,T).$$
Combining with the second point and the classical bound for the conductor of Artin $L$-functions yields the final bound in the statement.
\end{enumerate}
\end{proof}
For a finite place $\pp$ unramified for $M$, let $y_\pp$
denote the conjugacy class of the normalized semisimple Frobenius
$$
\left((N\pp)^{-w/2}
\rho_{M,\ell}(\mathrm{Frob}_\pp)\right)^{\mathrm{ss}},
\qquad \pp\nmid\ell,
$$
viewed over $\bb C$ using the common Weil--Deligne datum,
where $\mathrm{Frob}_\pp$ denotes geometric Frobenius.
Under the strong $G_M$-compatibility hypothesis, this conjugacy
class in $\mathrm{AST}_K(M)(\bb C)$ is independent of $\ell$.
It follows from \cite[Remark 7.11]{BanaszakKedlaya2023}
that the class has a representative in $\st_K(M)$; we also
write $y_\pp$ for the resulting conjugacy class in $\st_K(M)$.
We denote by $\chi_\sigma(y_\pp)$ the value of the character
of $\sigma$ on this class. For $\pp$ unramified in $L$, let
$c_\pp$ denote the geometric Frobenius conjugacy class in
$\Gal(L/K)$.

Let $\operatorname{Ram}(L/K)$ be the set of primes
$\pp\subset\mathcal O_K$ ramifying in $L$.
For the common local data
$\mathrm{WD}_v(M)=(r_v,N_v)$, recall the exceptional set of ramification of the Weil-Deligne datum
$
S_M:=\{v:r_v|_{I_v}\neq1\text{ or }N_v\neq0\}$. Define the exceptional set for the extension $L/K$ and the motive Weil--Deligne datum of the $\ell$-adic realization $M$ as 
\begin{equation}\label{set}
    S_{M,L}:=S_M\cup\operatorname{Ram}(L/K).
\end{equation}

Define
$$
A_{L/K}(M)
:=2^{|S_{M,L}|}\,d_K^{1/n_K}q_M.\\
$$

\begin{proposition}\label{sum}
Let $M\in\mathcal{M}_\sim(K,\Q)$ be a CST-motive, and let
\[
(\sigma,\rho)\in
\Irr(\st_K(M))\times\Irr(\Gal(L/K)).
\]
Assume that $(\sigma,\rho)\neq(1,1)$ and that the motivic
$L$-function formalism of Conjecture \ref{grh} holds for
$M_\sigma\otimes N_\rho$.

Then, for every $x\geq2$,
$$
\sum_{\substack{\pp\notin S_{M,L}\\N\pp\leq x}}
\chi_\sigma(y_\pp)\chi_\rho(c_\pp)
\ll
n_Kx^{1/2}d_\sigma d_\rho
\log\!\bigl(
q_\rho A_{L/K}(M)(x+C_H\|\lambda_\sigma\|)
\bigr).
$$
The implied constant in the first estimate is effective and absolute.
Moreover, for every compactly supported function $f\in C_c^\infty((0,\infty))$ and $x\geq2$,
\begin{align*}
&\sum_{\pp\notin S_{M,L}}
\chi_\sigma(y_\pp)\chi_\rho(c_\pp)
f\left(\frac{N\pp}{x}\right)
\\
&\ll_f
\frac{n_Kx^{1/2}d_\sigma d_\rho}{\log x}
\log\!\bigl(
q_\rho A_{L/K}(M)(3+C_H\|\lambda_\sigma\|)
\bigr)    
\end{align*}
The underlying constants depend effectively on the support of $f$ and
the norms $\|f^{(j)}\|_{L^\infty((0,\infty))}$ for $j\in \{0,1,2\}$.
\end{proposition}

\begin{proof}
We first show that, under our assumptions, the motivic
$L$-function $L(\iota,M_\sigma\otimes N_\rho,s)$ is entire.
Let $w_\sigma$ be the weight of $M_\sigma$.
By conjecture \ref{grh} and the simplicity of the motives involved, it suffices to prove that when $w_\sigma$ is even, there is no isomorphism
$
M_\sigma\simeq N_\rho^\vee(-w_\sigma/2).
$
Suppose that such an isomorphism exists.
If $\sigma=1$, then $M_\sigma$ is a Tate motive,
so the isomorphism forces $\rho=1$, contrary to
$(\sigma,\rho)\neq(1,1)$.
Thus $\sigma\neq1$.
The image
$
\widehat{\sigma}_E\bigl(\mathrm{AST}_K(M)_E\bigr)
$
is therefore connected and nontrivial.
On the other hand, a Tate twist of an Artin motive
has finite algebraic Sato-Tate group, which yields the desired contradiction.
The two bounds now follow by applying \cref{sumc}, together
with the conductor bounds of \cref{conduct}.
\end{proof}
An immediate consequence of Peter-Weyl density of class functions is the following equidistribution result.
\begin{corollary}
   Under Conjecture \ref{grh} for the normalized motivic $L$-functions
$F_{M_\sigma\otimes N_\rho}$ for every pair
$(\sigma,\rho)\in\Irr(\st_K(M))\times\Irr(\mathrm{Gal}(L/K))$,
the pairs $(y_\pp,c_\pp)$ are equidistributed in the conjugacy-class
space of $\st_K(M)\times\mathrm{Gal}(L/K)$ with respect to the
probability measure
\[
\sum_C \frac{|C|}{|\mathrm{Gal}(L/K)|}\,
\mu_{\st}\otimes\delta_C,
\]
where $C$ ranges over the conjugacy classes of $\mathrm{Gal}(L/K)$.
\end{corollary}
\begin{remark}
In the disconnected setting, equidistribution of the geometric Frobenius in the full Sato-Tate group might fail. Let $M=\Q(-1)$ be the Tate motive. The motivic Galois group is $\Gm$ while the algebraic Sato-Tate group is $\mu_2$. Yet every Frobenius has norm $1$ thus the Frobenius cannot be equidistributed in $\mu_2(\bb C)\simeq \{-1,1\}$. In particular, Banaszak and Kedlaya prove that if $-\mathrm{id}\notin \widetilde{\st_K(M)}$ 
then equidistribution cannot hold, where $\widetilde{\st_K(M)}$ is the extended Sato-Tate group as in \cite[Definition 17.4, Corollary 17.5 and Theorem 17.16]{BanaszakKedlaya2023}. A separate issue is that strict compatibility of the
$\ell$-adic realizations does not by itself establish
the independence of $\ell$ of the conjugacy classes
of normalized Frobenius elements in
$\operatorname{AST}_K(M)(\mathbb{C})$, see
\cite[Remarks 7.11 and 7.13]{BanaszakKedlaya2023}.

\end{remark}
We now formulate a possible direction for future research.
\begin{question}
For general motives, the Sato-Tate conjecture and its
Chebotarev analogue remain open. As discussed in the
previous section, extending our Beurling-Selberg theory
to disconnected compact groups requires an
adaptation of the Weyl integration and character formalism.

Can one develop such a theory using the representations
of a disconnected compact group $G$, and thereby obtain
effective equidistribution estimates under suitable
hypotheses on the associated $L$-functions?
Since $\mu_{\mathrm{ST}}$-continuity sets can be approximated
by sets with finite Minkowski content, as recorded in
\cref{approximation}, such a construction for
$\st_K(M)$ equipped with the singular Riemannian conjugacy foliation would
extend our approach to this setting.
The results of \cite{Wendt2001} and \cite{Zerouali2020}
on twisted conjugacy and character formulas provide
a natural starting point, particularly when the neutral component $G^\circ$
is simple and simply connected or admits a simple simply connected cover.
\end{question}
\section{Automorphic motives}\label{unc}
We briefly outline an automorphic approach which could
provide the basis for future unconditional results.
Motivic-automorphic reciprocity conjectures predict
that motivic $L$-functions arise from automorphic
$L$-functions, after the appropriate normalization
\cite{Clozel1990,Clozel2016,BuzzardGee2014},
see also \cite[Sections 3-4]{Arthur2025}.
In our setting, a suitable framework is to identify,
for each $M_\sigma$, a unitary automorphic representation
$\pi$ in an $L$-packet $\Pi$ of a connected reductive
group $H/K$, together with a finite-dimensional
representation
$
r:{}^LH\rightarrow\operatorname{GL}_n(\C),
$
algebraic on $\widehat H$, such that
$$
L^S(\iota,s+w_\sigma/2,M_\sigma)
=
L^S(s,\pi,r),
$$
where $L^S(s,\pi,r)$ is the partial Langlands
$L$-function and $S$ contains the ramified finite places.
Local-global compatibility predicts that the normalized
Frobenius classes agree with the images under $r$
of the unramified Satake parameters.
Langlands functoriality predicts the existence of an
automorphic transfer $\Pi_r$ to
$\operatorname{GL}_n(\mathbb A_K)$
whose partial standard $L$-function satisfies
$$
L^S(s,\pi,r)=L^S(s,\Pi_r);
$$
see, for example, \cite{Borel1979}.
Equality of the completed $L$-functions additionally
requires matching local factors at the remaining finite
places and at infinity.
Whenever these identifications and transfers are
established, the analytic theory of automorphic
$L$-functions becomes available for our
equidistribution problems.
This approach applies, in particular, to modular motives,
for which the relevant $L$-functions include the symmetric-power
$L$-functions $L(s,\operatorname{Sym}^m f)$ and the
Rankin-Selberg $L$-functions
$L(s,\operatorname{Sym}^m f\times\operatorname{Sym}^n f')$,
all in unitary normalization.
Their analytic properties and associated prime character
sums have been extensively studied, see, for example,
\cite[Propositions 2.1-2.2]{Thorner2021} and
\cite[Corollary 1.2 and Theorem 2.1]{Thorner2025}.

Motivated by the extensive literature on symmetric-power
$L$-functions, we ask whether the recent work of OpenAI
on the quasi-Riemann hypothesis, and its extension to
Dedekind zeta functions by Cureton and Thorner,
can be extended to the automorphic setting
\cite{CuretonThorner2026}.
Frei, Sofos and Widmer have already established an
effective Chebotarev density theorem under a
quasi-Riemann hypothesis
\cite{FreiSofosWidmer2026}.
\begin{question}
Let $f$ and $f'$ be non-CM holomorphic cuspidal newforms.
Can these methods establish a quasi-Riemann hypothesis
for the symmetric-power and Rankin-Selberg $L$-functions
$$
L(s,\operatorname{Sym}^m f),
\qquad
L(s,\operatorname{Sym}^m f
\times\operatorname{Sym}^n f'),
\qquad m,n\geq1.
$$
More precisely, can one obtain zero-free half-planes
$\Re(s)>1-\delta$, for an effectively computable $\delta>0$ independent of $m,n$?
\end{question}

\chapter{Motivic equidistribution theorems}\label{equid}
In this section, we develop conditional effective
Chebotarev-Sato-Tate estimates for CST-motives by applying the Beurling-Selberg approximation on Lie
groups from \cref{Lie} to conjugacy-invariant sets.
We express the error terms in terms of their boundary
geometry and the relevant arithmetic invariants.
The assumed motivic $L$-function conjectural package \ref{grh} supplies
the analytic estimates through control of the zeros
of the associated $L$-functions. The uniformity of the construction also allows the target set to vary with the prime, with its variation
controlled by Hausdorff distance and uniform Minkowski
content. We derive several corollaries and formulate questions for future research projects. Since the error exponent deteriorates as the relevant Lie-theoretic parameters increase, we explain at the
end of the section how additional symmetries can be
exploited to improve it. We also recover several results from the literature,
including conditional effective Sato-Tate estimates
for abelian varieties and trace intervals
\cite[Theorem 3.8]{BucurFiteKedlaya2024}.
\section{Effective Chebotarev-Sato-Tate for a fixed set}
Throughout this section, fix a number field $K$ of degree
$n_K$ and absolute discriminant $d_K$, and let
$M\in\mathcal M_\sim(K,\Q)$ be a CST-motive. Let $G:=\st_K(M)$ be the Sato-Tate group.
Equip $\mathrm{Conj}(G)$ with the Sato-Tate
measure $\mu_{\st}$. Fix a maximal torus $T\subset G$ of dimension
$d=\operatorname{rank}(G)$ equipped with the underlying $W$-invariant flat Metric $g$. Fix a system of positive
roots $\Phi^+$ with cardinality $m:=|\Phi^+|$.
For a conjugacy-invariant subset $A\subset G$
with nonempty boundary, denote by $A_T\subset\T$
the inverse image of its image in
$\mathrm{Conj}(G)$ under the map
$$
T\longrightarrow T/W
\simeq\mathrm{Conj}(G).
$$
Let $L/K$ be a finite Galois extension with Galois group
$\Gamma:=\operatorname{Gal}(L/K)$.
Fix a number field $E$ over which $G_M$ splits and
all irreducible complex representations of $\Gamma$
admit models. Fix an embedding
$\iota:E\hookrightarrow\C$.
For each pair
$(\sigma,\rho)\in
\operatorname{Irr}(G)
\times\operatorname{Irr}(\Gamma)$,
let
$
M_\sigma\in\l M\rr^\otimes\otimes_\Q E
$
be a motivic lift supplied by \cref{motivic},
defined up to a Tate twist, and let
$N_\rho\in\mathcal M_\sim(K,E)$ be the corresponding
Artin motive. The exceptional set $S_{M,L}$ is defined in \cref{set}.

\begin{hypothesis}\label{trh}
We assume that all the normalized motivic $L$-functions $F_{M_\sigma\otimes N_\rho}(s)$ satisfy Conjecture \ref{grh}.
\end{hypothesis}
Our assumptions imply that the Dedekind zeta function $\zeta_K$ identified with the normalized $L$-function of the weight $0$ unit motive $h(\operatorname{Spec} K)$ satisfies the Generalized Riemann hypothesis.
Define the counting functions:
\begin{align*}
\pi_{K}(x)&:=\#\{\pp\subset \mathcal{O}_K,\,N\pp\leq x\}\\
    \pi_{K,M,\Gamma}(x)&:=\#\{\pp\notin S_{M,L},\,N\pp\leq x\}.
\end{align*}

The Rosser--Schoenfeld bound \cite{RosserSchoenfeld1962} implies that 
\begin{equation}\label{primebound}
    \pi_{K,M,\Gamma}(x)\leq n_K\pi_{\Q}(x)
 <1.25506\,n_K\frac{x}{\log x}.
\end{equation}
We shall make use of the following formulas
\cite[Chapter VII, Section 11, (11.9)]{Neukirch1999},
\begin{equation}\label{cond}
    \sum_{\rho\in\Irr(\Gamma)}
d_\rho \log q_\rho
=
\log d_L, \qquad \sum_{\rho\in \Irr(\Gamma)}d_\rho^2=|\Gamma|.
\end{equation}
We start by defining the class of sets on which our effective estimate holds. 
\begin{definition}\label{min}
\begin{enumerate}
    \item  Fix $0<q\leq d$, and define
$$
\mathfrak M_g^q
:=\{A\subset\operatorname{Conj}(G):
0<\overline \J_g^{d-q}(\partial A_T)<\infty\}.
$$
\item
For $A\in\mathfrak M_g^q$ and $0<r_0<1$, put
$$
C_{\partial A}=C_\partial(A,r_0)
:=\max\left\{
1,\sup_{0<r<r_0}
\frac{\operatorname{vol}_g((\partial A_T)^{(r)})}{r^q}
\right\}.
$$
   \end{enumerate}
\end{definition} 
As discussed in \cref{content}, our class of sets
extends beyond the classical cases considered in the literature,
such as intervals, boxes, and sets with smooth or rectifiable
boundaries, to sets whose boundaries have finite Minkowski content. Moreover the condition, $\overline{\mathcal M}_g^{d-q}(\partial A_T)<\infty$ is stable under finite unions, intersections, and complements and independent of the choice of the flat metric $g$ on the torus $T$. 
We refer to \cref{content} for a detailed discussion of
Minkowski content.

We start with an effective Chebotarev-Sato-Tate for CST motives.
\begin{theorem}\label{cst1}
Let $A\in\mathfrak M_g^q$ and let $C\subset\Gamma$ be a non-empty union of
conjugacy classes. Fix $0<r_0<1$, and let
$C_{\partial A}$ be the constant
defined in \cref{min}. Put
$$
\delta_G:=d+m,
\qquad
S_{L/K}(M):=
3(1+C_H)A_{L/K}(M)d_L^{1/[L:K]^{1/2}},
$$
where $A_{L/K}(M)$ is defined in \cref{sum}, and set
$$
\pi_{A,C}(x):=
\#\left\{
\pp\notin S_{M,L}:
N\pp\leq x,\;
(y_\pp,c_\pp)\in A\times C
\right\}.
$$
Then there exists a threshold $x_0\geq e^4$ satisfying
$$
x_0\ll_G
n_K
\left(
\frac{|\Gamma|}{|C|^{1/2}}
\bigl(1+\log S_{L/K}(M)\bigr)
r_0^{-\delta_G}
\right)^4
$$
such that, for every $x>x_0$,
\begin{align*}
&\pi_{A,C}(x)
=
\frac{|C|}{|\Gamma|}\mu_{\mathrm{ST}}(A)\Li(x)
\\
&+
O_G\left(
n_K|C|^{1/2}\mu_{\mathrm{ST}}(A)x^{1/2}
\log\bigl(S_{L/K}(M)x\bigr)
\right)
\\
&+
O_G\left(
n_K C_{\partial A}
\frac{|C|^{1-q/(2\delta_G)}}{|\Gamma|^{1-q/\delta_G}}
x^{1-q/(2\delta_G)}
(\log x)^{2q/\delta_G-1}
\bigl[\log(S_{L/K}(M)x)\bigr]^{q/\delta_G}
\right).
\end{align*}
Moreover, the dependence on $G$ is effective.
\end{theorem}
\begin{proof}
Let
$
P_{r,N}^{\pm}
=
\sum_{\|\lambda\|<N}
a_\lambda^\pm\chi_\lambda
$
be the spectral Beurling-Selberg polynomials associated with
$A$ constructed in \cref{Lie}, so that
$
P_{r,N}^{-}\leq\mathbf 1_A\leq P_{r,N}^{+}.
$

By a standard orthogonality argument
(cf. \cite[Proof of Lemma 7.2]{ZywinaLangTrotter}), we have
\begin{align}
&\sum_{D'\subset\Gamma}
\frac{1}{|D'|}
\left|
\sum_{\substack{N\mathfrak p\leq x\\
\mathfrak p\notin S_{M,L}}}
\mathbf 1_{D'}(c_{\mathfrak p})
P_{r,N}^{\pm}(y_{\mathfrak p})
-
\frac{|D'|}{|\Gamma|}
\sum_{\substack{N\mathfrak p\leq x\\
\mathfrak p\notin S_{M,L}}}
P_{r,N}^{\pm}(y_{\mathfrak p})
\right|^2
\nonumber\\
&\qquad=
\frac1{|\Gamma|}
\sum_{\substack{\rho\in\Irr(\Gamma)\\
\rho\neq 1}}
\left|
\sum_{\substack{N\mathfrak p\leq x\\
\mathfrak p\notin S_{M,L}}}
P_{r,N}^{\pm}(y_{\mathfrak p})
\chi_\rho(c_{\mathfrak p})
\right|^2,
\label{1}
\end{align}
where the sum on the left is over the conjugacy classes
$D'$ of $\Gamma$. Hence
\begin{align}
&
\left|
\sum_{\substack{N\mathfrak p\leq x\\
\mathfrak p\notin S_{M,L}}}
\mathbf 1_C(c_{\mathfrak p})
P_{r,N}^{\pm}(y_{\mathfrak p})
-
\frac{|C|}{|\Gamma|}
\sum_{\substack{N\mathfrak p\leq x\\
\mathfrak p\notin S_{M,L}}}
P_{r,N}^{\pm}(y_{\mathfrak p})
\right|
\nonumber\\
&\qquad\leq
\left(\frac{|C|}{|\Gamma|}\right)^{1/2}
\left(
\sum_{\rho\neq 1}
\left|
\sum_{\substack{N\mathfrak p\leq x\\
\mathfrak p\notin S_{M,L}}}
P_{r,N}^{\pm}(y_{\mathfrak p})
\chi_\rho(c_{\mathfrak p})
\right|^2
\right)^{1/2}.
\label{2}
\end{align}
Since
$
P_{r,N}^{-}\leq\mathbf 1_A\leq P_{r,N}^{+},
$
it follows that
\begin{equation}
\begin{aligned}\label{3}
&\left|
\pi_{A,C}(x)
-
\frac{|C|}{|\Gamma|}
\mu_{\mathrm{ST}}(A)\,
\pi_{K,M,\Gamma}(x)
\right|
\leq
\frac{|C|}{|\Gamma|}
\sum_{\substack{\lambda\neq0\\ \|\lambda\|<N}}
\max\{|a_\lambda^+|,|a_\lambda^-|\}
\left|
\sum_{\substack{N\mathfrak p\leq x\\
\mathfrak p\notin S_{M,L}}}
\chi_\lambda(y_{\mathfrak p})
\right|
\\&+
\frac{|C|}{|\Gamma|}
\pi_{K,M,\Gamma}(x)
\int_G
\left(P_{r,N}^{+}-P_{r,N}^{-}\right)
\,d\mu_{\mathrm{ST}}
\\&
+
\left(\frac{|C|}{|\Gamma|}\right)^{1/2}
\max_{\pm}
\left(
\sum_{\rho\neq 1}
\left|
\sum_{\substack{N\mathfrak p\leq x\\
\mathfrak p\notin S_{M,L}}}
P_{r,N}^{\pm}(y_{\mathfrak p})
\chi_\rho(c_{\mathfrak p})
\right|^2
\right)^{1/2}.
\end{aligned}
\end{equation}
The coefficient bound in \cref{Lie}, together with the
character estimate in \cref{sum}, gives
$$
\begin{aligned}
&
\sum_{0<\|\lambda\|<N}
\max\{|a_\lambda^+|,|a_\lambda^-|\}
\left|
\sum_{\substack{N\mathfrak p\leq x\\
\mathfrak p\notin S_{M,L}}}
\chi_\lambda(y_{\mathfrak p})
\right|
\\
&\qquad\ll_G
\frac{C_{\partial A}n_Kx^{1/2}}
{r^{\delta_G-q}}
\sum_{0<\|\lambda\|<N}
\frac{
d_\lambda
\log\bigl(A_{L/K}(M)(x+C_H\|\lambda\|)\bigr)
}{
\|\lambda\|^{\delta_G}
}.
\end{aligned}
$$

There are $O_G(k^{\dim T-1})$ dominant weights whose
norm lies in $[k,k+1)$. Moreover, by the Weyl dimension
formula (\cref{Lies}),
$
d_\lambda\ll_G
(1+\|\lambda\|)^{\#\Phi^+}.
$
Since $\delta_G=\dim T+\#\Phi^+$, absorbing the finitely
many nonzero weights of norm below $1$ into the implied
constant, we obtain, for $N\geq3$,
\begin{align}
&
\sum_{0<\|\lambda\|<N}
\frac{
d_\lambda
\log\bigl(A_{L/K}(M)(x+C_H\|\lambda\|)\bigr)
}{
\|\lambda\|^{\delta_G}
}
\nonumber\\
&\qquad\ll_G
\log\bigl(A_{L/K}(M)(x+C_HN)\bigr)
\left(1+\sum_{1\leq k<N}\frac1k\right)
\nonumber\\
&\qquad\ll_G
\log N\,
\log\bigl(A_{L/K}(M)(x+C_HN)\bigr).
\label{sCST-lambda-sum}
\end{align}

On the other hand, the prime bound \cref{primebound},
together with the $L^1$-gap estimate in \cref{Lie}, implies
\begin{equation}
\pi_{K,M,\Gamma}(x)
\int_G
(P_{r,N}^{+}-P_{r,N}^{-})
\,d\mu_{\mathrm{ST}}
\ll_G
n_K\frac{x}{\log x}
\left(
\frac{C_{\partial A}}
{r^{s-q}N^{s-d}}
+
C_{\partial A}r^q
\right).
\label{4}
\end{equation}

For $\rho\neq 1$, \cref{sum} yields
\begin{equation}\label{5}
\begin{aligned}
&\left|
\sum_{\substack{N\mathfrak p\leq x\\
\mathfrak p\notin S_{M,L}}}
P_{r,N}^{\pm}(y_{\mathfrak p})
\chi_\rho(c_{\mathfrak p})
\right|
\ll_G
n_Kx^{1/2}d_\rho
\Biggl[
\left(
\mu_{\mathrm{ST}}(A)
+
C_{\partial A}r^q
+
\frac{C_{\partial A}}
{r^{s-q}N^{s-d}}
\right)\\&\times
\log(q_\rho A_{L/K}(M)x)
+
\frac{C_{\partial A}\log N}
{r^{\delta_G-q}}
\log\bigl(
A_{L/K}(M)q_\rho(x+C_HN)
\bigr)
\Biggr].
\end{aligned}
\end{equation}
By the conductor-discriminant formulas \eqref{cond} it follows that:
$$
\left(
\sum_\rho d_\rho^2(\log q_\rho)^2
\right)^{1/2}
\leq
\sum_\rho d_\rho\log q_\rho
=
\log d_L.
$$
Minkowski's inequality therefore gives
\begin{equation}\label{6}
\begin{aligned}
&
\left(
\sum_\rho d_\rho^2
\Biggl[
\left(
\mu_{\mathrm{ST}}(A)
+
C_{\partial A}r^q
+
\frac{C_{\partial A}}
{r^{s-q}N^{s-d}}
\right)
\log(A_{L/K}(M)q_\rho x)
\right.\\&
+
\frac{C_{\partial A}\log N}
{r^{\delta_G-q}}
\log\bigl(
A_{L/K}(M)q_\rho(x+C_HN)
\bigr)
\Biggr]^2
)^{1/2}
\\
&\leq
|\Gamma|^{1/2}
\left[
\mu_{\mathrm{ST}}(A)
+
C_{\partial A}
\left(
r^q
+
\frac1{r^{s-q}N^{s-d}}
+
\frac{\log N}{r^{\delta_G-q}}
\right)
\right]\\&\times
\Biggl[
\log\bigl(A_{L/K}(M)(x+C_HN)\bigr)
+
\frac{\log d_L}{|\Gamma|^{1/2}}
\Biggr].
\end{aligned}
\end{equation}
Combining
\eqref{3},
\eqref{sCST-lambda-sum},
\eqref{4},
\eqref{5},
and \eqref{6}, we obtain
\begin{equation}
    \begin{aligned}
&
\left|
\pi_{A,C}(x)
-
\frac{|C|}{|\Gamma|}
\mu_{\mathrm{ST}}(A)
\pi_{K,M,\Gamma}(x)
\right|
\nonumber\ll_G
\frac{|C|}{|\Gamma|}
n_K\frac{x}{\log x}
\left(
\frac{C_{\partial A}}
{r^{s-q}N^{s-d}}
+
C_{\partial A}r^q
\right)
\nonumber\\
&+
n_K|C|^{1/2}x^{1/2}
\left[
\mu_{\mathrm{ST}}(A)
+
C_{\partial A}
\left(
r^q
+
\frac1{r^{s-q}N^{s-d}}
+
\frac{\log N}{r^{\delta_G-q}}
\right)
\right]
\\
&\times \left[
\log\bigl(A_{L/K}(M)(x+C_HN)\bigr)
+
\frac{\log d_L}{|\Gamma|^{1/2}}
\right]+\\&
n_KC_{\partial A}\frac{|C|}{|\Gamma|}
\frac{x^{1/2}\log N}{r^{\delta_G-q}}
\log\bigl(A_{L/K}(M)(x+C_HN)\bigr).
\label{8}
\end{aligned}
\end{equation}

Choose $s=d+1$,
$
N=
\lceil\max\{1+2\|\rhoo\|,
3,\,
3\|\rhoo\|,\,
 r^{-(d+1)}
\}\rceil,
$
and
\begin{equation}
r=
\left(
\frac{|\Gamma|}{|C|^{1/2}}\,
\frac{(\log x)^2}{x^{1/2}}\,
\log\bigl(S_{L/K}(M)x\bigr)
\right)^{1/\delta_G}.
\label{7}
\end{equation}
For $x\geq e^4$, the elementary estimate
$$
(\log x)^2
\bigl(\log x+\log S_{L/K}(M)\bigr)
\leq
128\bigl(1+\log S_{L/K}(M)\bigr)x^{1/4}
$$
shows that $r<r_0$ whenever $x>x_0$, where one may take
$$
x_0:=
\max\left\{
e^4,\,
\left(
128\,\frac{|\Gamma|}{|C|^{1/2}}\,
\bigl(1+\log S_{L/K}(M)\bigr)
r_0^{-\delta_G}
\right)^4
\right\}.
$$
This choice satisfies the asserted bound for $x_0$.

Moreover, \eqref{7} gives
$$
r^{-\delta_G}\leq x^{1/2},
\qquad
N\ll_G x^{(d+1)/(2\delta_G)}\ll_G x,
\qquad
\log N\ll_G\log x.
$$
Consequently, the definition of $S_{L/K}(M)$ implies
$$
\log\bigl(A_{L/K}(M)(x+C_HN)\bigr)
+
\frac{\log d_L}{|\Gamma|^{1/2}}
\ll_G
\log\bigl(S_{L/K}(M)x\bigr).
$$
Our choice of $N$ also gives
$$
\frac1{r^{s-q}N^{s-d}}\leq r^q.
$$
Using
$
|C|/|\Gamma|\leq |C|^{1/2}
$
and
$
r^q\ll r^{q-\delta_G}\log N,
$
we deduce from \eqref{8} that
\begin{align*}
&
\left|
\pi_{A,C}(x)
-
\frac{|C|}{|\Gamma|}
\mu_{\mathrm{ST}}(A)
\pi_{K,M,\Gamma}(x)
\right|\ll_G
n_K|C|^{1/2}\mu_{\mathrm{ST}}(A)x^{1/2}
\log\bigl(S_{L/K}(M)x\bigr)
\\
&\qquad+
n_KC_{\partial A}
\left[
\frac{|C|}{|\Gamma|}\frac{x}{\log x}r^q
+
|C|^{1/2}x^{1/2}
\frac{\log N}{r^{\delta_G-q}}
\log\bigl(S_{L/K}(M)x\bigr)
\right].
\end{align*}
By \eqref{7},
$$
|C|^{1/2}x^{1/2}
\frac{\log N}{r^{\delta_G-q}}
\log\bigl(S_{L/K}(M)x\bigr)
=
\frac{|C|}{|\Gamma|}
\frac{x}{\log x}r^q
\frac{\log N}{\log x}.
$$
Thus the second term in brackets is
$\ll_G$ the first. Furthermore,
$$
\frac{|C|}{|\Gamma|}
\frac{x}{\log x}r^q
={}
\frac{|C|^{\,1-q/(2\delta_G)}}
{|\Gamma|^{\,1-q/\delta_G}\log x}
x^{\,1-q/(2\delta_G)}
(\log x)^{2q/\delta_G}
\bigl(\log(S_{L/K}(M)x)\bigr)^{q/\delta_G}.
$$

Finally, under the Riemann hypothesis for $\zeta_K$, which is included in
Hypothesis \ref{trh},
$$
\pi_{K,M,\Gamma}(x)
=
\Li(x)
+
O\left(
x^{1/2}\bigl(\log d_K+n_K\log x\bigr)
+
|S_{M,L}|
\right).
$$
Using the bound
$
|S_{M,L}|\ll\log A_{L/K}(M)
$
incorporated in the definition of $A_{L/K}(M)$,
and the definition of $S_{L/K}(M)$, we obtain
$$
\pi_{K,M,\Gamma}(x)
=
\Li(x)
+
O\left(
n_Kx^{1/2}\log\bigl(S_{L/K}(M)x\bigr)
\right).
$$
As a result,
$$
\begin{aligned}
&\left|
\pi_{A,C}(x)
-
\frac{|C|}{|\Gamma|}
\mu_{\mathrm{ST}}(A)\Li(x)
\right|
\ll_G
n_K|C|^{1/2}\,
\mu_{\mathrm{ST}}(A)x^{1/2}\,
\log\bigl(S_{L/K}(M)x\bigr)
+\\&
\frac{
C_{\partial A}n_K
|C|^{\,1-q/(2\delta_G)}
}{
|\Gamma|^{\,1-q/\delta_G}\log x
}
x^{\,1-q/(2\delta_G)}
(\log x)^{2q/\delta_G}
\bigl(\log(S_{L/K}(M)x)\bigr)^{q/\delta_G}.
\end{aligned}
$$
This completes the proof.
\end{proof}
\begin{remark}
\begin{enumerate}
    \item For $C=\Gamma$, $M=h^1(B)$ with $B/K$ an abelian variety,
and $A$ the trace interval,
our estimate recovers the exponent of $x$ obtained in
\cite[Theorem 1.1]{BucurFiteKedlaya2024}.
\item Our estimate improves as the boundary becomes thinner,
as measured by the tubular-neighborhood bound.
In particular, a smooth boundary $\partial A_T\subset T$
of codimension $q$ satisfies the required bound with
exponent $q$.
\item  By choosing $r_0>0$ sufficiently small,
the tube constant $C_{\partial A}$ can be taken
arbitrarily close to $\max(1,\omega_q\overline{\mathcal{M}}_g^{d-q}(\partial A))$.
For a smooth boundary of codimension $q$, Gray's tube
formula gives
$$
\lim_{r\to0^+}
\frac{\operatorname{vol}_g((\partial A_T)^{(r)})}{r^q}
=
\omega_q\,\mathcal H_g^{d-q}(\partial A_T).
$$

\item 
The construction is uniform over families $\mathcal{A}$ whose boundary has a uniform $d-q$-Minkowski content $\mathcal{M}_g^{d-q}(\partial \mathcal{A})$.
Its explicit dependence on these geometric quantities
provides a useful framework for studying Frobenius
classes near rare events, as illustrated in the next
section.
\end{enumerate}
\end{remark}
Following the smoothed Chebotarev estimate of
\cite[Theorem 2.2]{ZywinaLangTrotter}, we state a smoothed
version of our estimate, which will be useful for obtaining
logarithmic savings and bounding the least prime satisfying
prescribed Chebotarev-Sato-Tate conditions.
\begin{theorem}\label{cst2}
Retain the notation of
\cref{cst1}. Let $0\neq f\in C_c^\infty((0,\infty))$ be nonnegative
with support in an interval $(a,b)$.
Let
$
S_{L/K}(M)
$
be as in \cref{cst1} and define
$$
\pi_{A,C,f}(x):=
\sum_{\pp\notin S_{M,L}}
f\left(\frac{N\pp}{x}\right)
\ind_A(y_\pp)\ind_C(c_\pp).
$$
There exists a threshold $x_1\geq \max(e^4,a^{-2})$ satisfying
$$
x_1\ll_{G,f} \left(n_K\frac{|\Gamma|}{|C|^{1/2}}
\bigl(1+\log S_{L/K}(M)\bigr)
r_0^{-\delta_G}
\right)^4
$$
such that for every $x>x_1$,
\begin{align*}
&\pi_{A,C,f}(x)
=
\frac{|C|}{|\Gamma|}\mu_{\mathrm{ST}}(A)
\int_2^\infty\frac{f(u/x)}{\log u}\,du
\\
&+
O_{G,f}\left(
n_K|C|^{1/2}\mu_{\mathrm{ST}}(A)x^{1/2}
\frac{\log\bigl(S_{L/K}(M)x\bigr)}{\log x}
\right)
\\
&+
O_{G,f}\left(
n_KC_{\partial A}
\frac{|C|^{1-q/(2\delta_G)}}{|\Gamma|^{1-q/\delta_G}}
x^{1-q/(2\delta_G)}
(\log x)^{q/\delta_G-1}
\bigl[\log(S_{L/K}(M)x)\bigr]^{q/\delta_G}
\right).
\end{align*}
\end{theorem}

\begin{proof}
Define the smoothed prime sum
$$
\pi_{K,f,M,\Gamma}(x):=
\sum_{\pp\notin S_{M,L}}
f\left(\frac{N\pp}{x}\right).
$$
Put $a:=\inf\operatorname{supp}f$.
For $x\geq\max\{e^4,a^{-2}\}$, the functions
$$
g_x(t):=\frac{f(t)}{\log(xt)},
$$
defined near $\operatorname{supp}f$ and extended by zero,
satisfy
$$
\|g_x^{(j)}\|_\infty
\ll_{f,j}\frac1{\log x}
\qquad(j\geq0).
$$
Applying the argument of
\cite[Theorem 2.2]{ZywinaLangTrotter} to the unit motive
with the test function $g_x$, and treating the excluded
primes as in \cref{sum}, gives
$$
\pi_{K,f,M,\Gamma}(x)
=
\int_2^\infty\frac{f(u/x)}{\log u}\,du
+
O_f\left(
n_Kx^{1/2}
\frac{\log\bigl(S_{L/K}(M)x\bigr)}{\log x}
\right).
$$
The argument of \cref{cst1}, using character orthogonality,
the estimates in \cref{Lie,sum}, and the
conductor-discriminant formulas \eqref{cond}, yields
\begin{equation}\label{9}
\begin{aligned}
&
\left|
\pi_{A,C,f}(x)
-
\frac{|C|}{|\Gamma|}
\mu_{\mathrm{ST}}(A)
\pi_{K,f,M,\Gamma}(x)
\right|\ll_{G,f}
\frac{|C|}{|\Gamma|}
n_K\frac{x}{\log x}
\left(
\frac{C_{\partial A}}
{r^{s-q}N^{s-d}}
+
C_{\partial A}r^q
\right)
\nonumber\\
&+
n_K|C|^{1/2}\frac{x^{1/2}}{\log x}
\left[
\mu_{\mathrm{ST}}(A)
+
C_{\partial A}
\left(
r^q
+
\frac1{r^{s-q}N^{s-d}}
+
\frac{\log N}{r^{\delta_G-q}}
\right)
\right]
\nonumber\\
&\times
\left[
\log\bigl(A_{L/K}(M)(3+C_HN)\bigr)
+
\frac{\log d_L}{|\Gamma|^{1/2}}
\right].
\end{aligned}
\end{equation}
Here the untwisted nonconstant contribution is absorbed
using $|C|/|\Gamma|\leq|C|^{1/2}$.

Choose $s=d+1$, $N=\left\lceil
\max\{1+2\|\varrho\|,\,3,\,3\|\varrho\|,\,r^{-(d+1)}\}
\right\rceil$
and 
$$r
=
\left(
\frac{|\Gamma|}{|C|^{1/2}}\,
\frac{\log x\,
\log\bigl(S_{L/K}(M)x\bigr)}
{x^{1/2}}
\right)^{1/\delta_G}.$$
For $x\geq e^4$, we have
$$
\log x\,
\bigl(\log x+\log S_{L/K}(M)\bigr)
\leq
256\bigl(1+\log S_{L/K}(M)\bigr)x^{1/4}.
$$
Hence $x>x_1$ ensures that $r<r_0$.
Moreover,
$$
\frac1{r^{s-q}N^{s-d}}\leq r^q,
\qquad
N\ll_G x,
\qquad
\log N\ll_G\log x,
$$
and therefore
$$
\log\bigl(A_{L/K}(M)(3+C_HN)\bigr)
+
\frac{\log d_L}{|\Gamma|^{1/2}}
\ll_G
\log\bigl(S_{L/K}(M)x\bigr).
$$
As in the proof of \cref{cst1}, the choice of $r$ balances
the geometric and spectral errors, since
$$
|C|^{1/2}\frac{x^{1/2}}{\log x}
\frac{\log N}{r^{\delta_G-q}}
\log\bigl(S_{L/K}(M)x\bigr)
=
\frac{|C|}{|\Gamma|}
\frac{x}{\log x}r^q
\frac{\log N}{\log x}.
$$
Substituting these estimates into the preceding bound, together with the estimate for $\pi_{K,f,M,\Gamma}(x)$,
proves the assertion.
\end{proof}
We recover a useful upper bound for $\pi_{A,C}$  with logarithmic saving in comparison with \cref{cst1}.
\begin{corollary}\label{cst3}
Let $x_1$ be the threshold in \cref{cst2} associated with
a fixed nonnegative function
$f\in C_c^\infty((1/4,2))$ satisfying
$f\geq1$ on $[1/2,1]$. Put
$$
x_2:=
\max\left\{
2x_1,\,
\left(
\frac{x_1}{
|C|^{1/2}\mu_{\mathrm{ST}}(A)
+
C_{\partial A}
|C|^{1-q/(2\delta_G)}
|\Gamma|^{q/\delta_G-1}
}
\right)^2
\right\},
$$
with the convention $x_2:=2x_1$ if the denominator vanishes.
Then, for every $x>x_2$,
\begin{align*}
&\pi_{A,C}(x)
\ll_{G,q}
\frac{|C|}{|\Gamma|}
\mu_{\mathrm{ST}}(A)\frac{x}{\log x}
\\
&+
n_K|C|^{1/2}\mu_{\mathrm{ST}}(A)x^{1/2}
\frac{\log\bigl(S_{L/K}(M)x\bigr)}{\log x}
\\
&+
n_KC_{\partial A}
\frac{|C|^{1-q/(2\delta_G)}}{|\Gamma|^{1-q/\delta_G}}
x^{1-q/(2\delta_G)}
(\log x)^{q/\delta_G-1}
\bigl[\log(S_{L/K}(M)x)\bigr]^{q/\delta_G}.
\end{align*}
\end{corollary}
\begin{proof}
Fix a nonnegative compactly supported function
$f\in C_c^\infty((1/4,2))$ satisfying
$f\geq1$ on $[1/2,1]$.
Write $g(t)=t/\log t$, and let $h(t)$ denote the sum
of the two error bounds in \cref{cst2}. For $t>x_1$,
$$
\pi_{A,C}(t)-\pi_{A,C}(t/2)
\leq \pi_{A,C,f}(t)
\ll_G
\mu_{\mathrm{ST}}(A)\frac{|C|}{|\Gamma|}g(t)+h(t).
$$

Since $x_1\geq e^4>16$ and $0<q\leq\delta_G$, we have,
for $t>2x_1$,
$$
g(t/2)\leq\frac23g(t),
\qquad
h(t/2)\leq\frac{4}{3\sqrt2}\,h(t).
$$
Consequently,
$$
\sum_{\substack{j\geq0\\x/2^j>x_1}}
g(x/2^j)\ll g(x),
\qquad
\sum_{\substack{j\geq0\\x/2^j>x_1}}
h(x/2^j)\ll h(x).
$$
Telescoping and using $\pi_{A,C}(x_1)\leq n_Kx_1$,
we obtain
$$
\pi_{A,C}(x)
\ll_G
n_Kx_1+
\mu_{\mathrm{ST}}(A)\frac{|C|}{|\Gamma|}g(x)+h(x).
$$

Finally, since $x^{1/2}\geq\log x$ and
$\log(S_{L/K}(M)x)\geq\log x\geq1$ for $x\geq x_2$,
and $0<q\leq\delta_G$, we have
$$
h(x)\geq
n_Kx^{1/2}
\left(
\mu_{\mathrm{ST}}(A)|C|^{1/2}
+
C_{\partial A}
\frac{|C|^{1-q/(2\delta_G)}}
     {|\Gamma|^{1-q/\delta_G}}
\right).
$$
By the definition of $x_2$, it follows that
$h(x)\geq n_Kx_1$ for $x\geq x_2$.
Absorbing the initial term into $h(x)$ proves the statement.
\end{proof}
\begin{corollary}
Retain the notation and assumptions of \cref{cst1}.
Assume that $\mu_{\mathrm{ST}}(A)>0$ and
$C\neq\varnothing$. Let $x_1$ be the threshold in
\cref{cst2} corresponding to a fixed nonzero nonnegative
function $f\in C_c^\infty((1/2,1))$.
The least prime ideal $\mathfrak p\notin S_{M,L}$
satisfying $y_{\mathfrak p}\in A$ and
$c_{\mathfrak p}\in C$ has norm
$$
\begin{aligned}
&N\mathfrak p
\leq \max\Biggl\{x_1,\,
\frac{|\Gamma|^2}{|C|}
\left[
n_K+
\left(
\frac{n_K C_{\partial A}}{\mu_{\mathrm{ST}}(A)}
\right)^{\delta_G/q}
\right]^2
\\
&\qquad{}\times
O_{G,q}\!\left(
\log^4\!\left(
S_{L/K}(M)
\left[
2+\frac{|\Gamma|}{|C|^{1/2}}
\left(
n_K+
\left(
\frac{n_K C_{\partial A}}{\mu_{\mathrm{ST}}(A)}
\right)^{\delta_G/q}
\right)
\right]
\right)
\right)
\Biggr\}.
\end{aligned}
$$
The dependence on $G$ and $q$ is effective.
\end{corollary}

\begin{proof}
With the fixed choice of $f$ above, we have, for $x>x_1$,
$$
\int_2^\infty\frac{f(u/x)}{\log u}\,du
\geq
\frac{x}{\log x}\int_0^\infty f(u)\,du.
$$
Hence \cref{cst2} implies $\pi_{A,C,f}(x)>0$ provided
$$
x^{1/2}\gg_G
\frac{n_K|\Gamma|}{|C|^{1/2}}
\log\bigl(S_{L/K}(M)x\bigr)
$$
and, since $q>0$,
$$
x^{1/2}\gg_{G,q}
\frac{|\Gamma|}{|C|^{1/2}}
\left(
\frac{n_KC_{\partial A}}{\mu_{\mathrm{ST}}(A)}
\right)^{\delta_G/q}
\log x\,\log\bigl(S_{L/K}(M)x\bigr).
$$
Put
$$
B:=
2+\frac{n_K|\Gamma|}{|C|^{1/2}}
+\frac{|\Gamma|}{|C|^{1/2}}
\left(
\frac{n_KC_{\partial A}}{\mu_{\mathrm{ST}}(A)}
\right)^{\delta_G/q}.
$$
Both inequalities hold for
$$
x>
\max\left\{
x_1,\,
c_{G,q}B^2
\bigl(\log(S_{L/K}(M)B)\bigr)^4
\right\},
$$
with a sufficiently large constant $c_{G,q}$.
Since $f\geq0$ and
$\operatorname{supp}f\subset(1/2,1)$, positivity of
$\pi_{A,C,f}(x)$ gives a prime
$\mathfrak p\notin S_{M,L}$ with $N\mathfrak p<x$,
$y_{\mathfrak p}\in A$, and $c_{\mathfrak p}\in C$.
The stated bound follows.
\end{proof}
\begin{remark}
\begin{enumerate}
    \item Our estimates apply more generally to equidistribution
problems admitting a compact Sato-Tate group and the
character-sum and conductor bounds used in our argument.
For instance, for automorphic representations, the relevant $L$-function formalism is related to the Langlands $L$-functions. Local-global
compatibility identifies the corresponding normalized
Frobenius classes with the unramified Satake parameters
\cite{Borel1979,BuzzardGee2014}.
\item
Let $H/K$ be a connected reductive group, and choose
a finite Galois extension $L/K$ over which $H$ splits.
Put $\Gamma:=\operatorname{Gal}(L/K)$.
Using the finite Galois form of the $L$-group, write
$
{}^LH=\widehat H\rtimes\Gamma,
$
where $\Gamma$ acts through pinned automorphisms.
Outside the ramified places, the Satake parameter
$c(\pi_v)$ is represented by an element of the coset
$\widehat H\gamma_v$, where $\gamma_v$ represents
the Frobenius class of $v$ in $\Gamma$.
\item
If the Galois action is trivial and the joint
Sato-Tate group is a direct product $G\times\Gamma$,
the resulting distribution problem takes precisely
the Chebotarev-Sato-Tate form considered here.
\item
For a nontrivial action, conjugation by $U$ on the
coset $U\gamma$ corresponds to twisted conjugation
on $U$:
$
h\cdot_\gamma u:=hu\gamma(h)^{-1}.
$
The associated basic projection is
$
(P_\gamma f)(u)
:=
\int_U f\bigl(hu\gamma(h)^{-1}\bigr)\,d\mu_U(h).
$
For a $\Gamma$-invariant bi-invariant metric,
twisted Weyl integration formulas and expansions
in twining characters offer a stArting point as recorded in \cref{tableofbasic}. Such methods are used in
\cite[Section 5.1, Lemma 5.2]{ShinTemplier2016}
to describe the Sato-Tate measures associated with
non-split reductive groups via twisted Weyl
integration.
\end{enumerate}

\end{remark}
\section{Effective Chebotarev-Sato-Tate for a family of sets}\label{sato}
We keep the notations of the previous section.
The main novelty of this article is an effective
Chebotarev-Sato-Tate theorem for prime-dependent
families of sets. Its formulation requires a geometric
condition combining uniform control of the boundaries
with control of the variation of the sets in Hausdorff
distance, and compatible with the stability estimates for
our Beurling-Selberg spectral approximations.
This motivates the notion of an $(\alpha,q)$-Minkowski family.
\begin{definition}\label{uni}
Let
$\mathcal A=(A_{\mathfrak p})_{\mathfrak p\notin S_{M,L}}$
be a family of conjugacy-invariant Borel subsets of $G$
with nonempty boundary, and put
$$
\partial\mathcal A
:=
(\partial A_{\mathfrak p,T})_{\mathfrak p\notin S_{M,L}}.
$$
For $\mathfrak p,\mathfrak p'\notin S_{M,L}$, define
$$
d_*(\mathfrak p,\mathfrak p')
:=\max\Bigl\{
d_H^g\bigl(
\overline{A_{\mathfrak p,T}},
\overline{A_{\mathfrak p',T}}
\bigr),
d_H^g\bigl(
\overline{ T_{}\setminus A_{\mathfrak p,T}},
\overline{T_{}\setminus A_{\mathfrak p',T}}
\bigr)
\Bigr\}.
$$
We call $\mathcal A$ an $(\alpha,q)$-Minkowski family if the following conditions hold:
\begin{enumerate}
\item
the uniform upper $(d-q)$-dimensional Minkowski content
of $\partial\mathcal A$ is finite:
$$
\overline{\mathcal M}_g^{d-q}(\partial\mathcal A)<\infty.
$$
In this case, for the fixed $0<r_0<1$, set
$$
C_{\partial\mathcal A}
:=
\max\left\{
1,\,
\sup_{0<r<r_0}
\sup_{\mathfrak p\notin S_{M,L}}
\frac{
\s\bigl((\partial A_{\mathfrak p,T})^{(r)}\bigr)
}{r^q}
\right\}.
$$

\item
there exists a constant $K_{\mathcal A}\geq0$ such that,
for every $x\geq2$, enumerating the prime ideals outside $S_{M,L}$ with
$x/2<N\mathfrak p\leq x$ as
$\mathfrak p_u,\ldots,\mathfrak p_v$ in nondecreasing
order of norm, with a fixed ordering among equal norms,
we have
$$
\sum_{k=u}^{v-1}d_*(\mathfrak p_k,\mathfrak p_{k+1})
\leq K_{\mathcal A}x^{-\alpha}.
$$
\end{enumerate}
\end{definition}
\begin{remark}
To verify the uniform Minkowski content condition,
\cref{tube3,analytictube} provide criteria covering a wide
range of smooth and subanalytic families appearing in the Sato-Tate literature.
A standard way to control the variation of a family
$(A_t)_{t\in I}$ of subsets of $T$ is to construct an
 isotopy
$$
\Phi_{t,t'}:T\rightarrow T,
\qquad
\Phi_{t,t'}(A_t)=A_{t'}.
$$
Its restrictions give homeomorphisms
$$
\varphi_{t,t'}:A_t\longrightarrow A_{t'},
\qquad
\psi_{t,t'}:T\setminus A_t\longrightarrow T\setminus A_{t'}.
$$
Since these maps are surjective, we have
\begin{align*}
&\max\Bigl\{
d_H^{g_T}(\overline{A_t},\overline{A_{t'}}),
d_H^{g_T}\bigl(
\overline{T\setminus A_t},
\overline{T\setminus A_{t'}}
\bigr)
\Bigr\}\\
&\qquad\leq
\sup_{\theta\in T}
d_{g_T}\bigl(\theta,\Phi_{t,t'}(\theta)\bigr).
\end{align*}
Ehresmann's fibration theorem provides local trivializations for smooth proper submersions. For singular analytic families, Thom's first isotopy lemma applies to submersions compatible with a Whitney stratification \cite[Proposition 11.1 and Corollary 10.3]{Mather2012}.
We discuss these tools in greater detail when addressing
Serre's least-prime problem in \cref{minimalserre}.
\end{remark}
We now establish an effective Chebotarev-Sato-Tate
estimate for $(\alpha,q)$-Minkowski families.
\begin{theorem}\label{cst4}
Let $\mathcal A=(A_\pp)_{\pp\notin S_{M,L}}$ be an
$(\alpha,q)$-Minkowski family, with $\alpha\geq0$ and
$0<q\leq d$. Fix $0<r_0<1$, and let
$C\subset\Gamma$ be a non-empty union of conjugacy classes.
Assume the hypotheses on the motivic $L$-functions used
in \cref{cst1}.
For $x\geq2$, define
$$
\pi_{\mathcal A,C}(x):=
\#\left\{
\pp\notin S_{M,L}:
N\pp\leq x,\;
(y_\pp,c_\pp)\in A_\pp\times C
\right\}.
$$
Put
$$
A_{\mathrm{ST}}
:=
\sup_{\pp\notin S_{M,L}}\mu_{\mathrm{ST}}(A_\pp),
\qquad
\beta_G:=
\frac q2
\min\left\{
\frac1{\delta_G},
\frac{1+2\alpha}{\delta_G+1}
\right\},
$$
and let $K_{\mathcal A}$ be as in \cref{uni}. Put
$$
B_0:=
\frac{|\Gamma|\bigl(1+\log S_{L/K}(M)\bigr)}
{|C|^{1/2}r_0^{\delta_G}}
\max\left\{1,\frac{K_{\mathcal A}}{r_0}\right\},
$$
and
$$
B_1:=
A_{\mathrm{ST}}|C|^{1/2}
+
C_{\partial\mathcal A}
(1+K_{\mathcal A})^{q/(\delta_G+1)}
\frac{|C|^{1-q/(2\delta_G)}}{|\Gamma|^{1-q/\delta_G}}.
$$
Then there exists a threshold $x_4\geq e^4$ satisfying
$$
x_4\ll_G
B_0^4\max\left\{1,\left(\frac{B_0^3}{B_1}\right)^4\right\}
$$
such that, for every $x>x_4$,
$$\pi_{\mathcal A,C}(x)
=
\frac{|C|}{|\Gamma|}
\sum_{\substack{N\pp\leq x\\\pp\notin S_{M,L}}}
\mu_{\mathrm{ST}}(A_\pp)+n_KA_{\mathrm{ST}}|C|^{1/2}\mathcal{E}_1(x)+n_KL_{\mathcal{A}}
\frac{|C|^{1-q/(2\delta_G)}}{|\Gamma|^{1-q/\delta_G}}\mathcal{E}_2(x),
$$
with 
\begin{align*}
L_{\mathcal{A}}&:=C_{\partial\mathcal A}
(1+K_{\mathcal A})^{q/(\delta_G+1)},\\
    \mathcal{E}_1(x)&:=O_G\left(
x^{1/2}
\log\bigl(S_{L/K}(M)x\bigr)\right),\\
\mathcal{E}_2(x)&:=O_G\left(
x^{1-\beta_G}
(\log x)^{2q/\delta_G-1}
\bigl[\log(S_{L/K}(M)x)\bigr]^{q/\delta_G}
\right).
\end{align*}
The dependence of the implied constants on $G$ is effective.
\end{theorem}

\begin{proof}
Choose a sufficiently large constant $c_G\geq256$, depending
only on $G$, and set
$$
x'_4:=(c_GB_0)^4,
\qquad
x_4:=
\max\left\{
2x'_4,\,
\left(\frac{4x'_4}{B_1}\right)^4
\right\}.
$$
Since $B_0\geq1$, these choices satisfy $x'_4\geq e^4$
and the stated upper bound for $x_4$.

For each $\pp\notin S_{M,L}$, let
$$
P_{r,N,\pp}^{\pm}
=
\sum_{\|\lambda\|<N}
c_{\lambda,\pp}^{\pm}\chi_\lambda
$$
be the spectral majorants and minorants associated with
$A_\pp$ in \cref{Lie}. We use the common boundary constant
$C_{\partial\mathcal A}$ and the same additive tail corrections
for every $\pp$.
Throughout the dyadic estimate below, sums over
$\mathfrak p$ are restricted to
$\mathfrak p\notin S_{M,L}$ with
$x/2<N\mathfrak p\leq x$.

As in the proof of \cref{cst1}, we have
\begin{align}
&\left|
\pi_{\mathcal A,C}(x)-\pi_{\mathcal A,C}(x/2)
-\frac{|C|}{|\Gamma|}
\sum_{\mathfrak p}\mu_{\mathrm{ST}}(A_{\mathfrak p})
\right|
\nonumber \\&\ll_{G,s}
n_KC_{\partial\mathcal A}\frac{|C|}{|\Gamma|}
\frac{x}{\log x}
\left(
r^q+\frac1{r^{s-q}N^{s-d}}
\right)
\nonumber\\
&\qquad+
\frac{|C|}{|\Gamma|}
\max_\pm
\sum_{0<\|\lambda\|<N}
\left|
\sum_{\mathfrak p}
c_{\lambda,\mathfrak p}^{\pm}
\chi_\lambda(y_{\mathfrak p})
\right|\nonumber\\&+
\left(\frac{|C|}{|\Gamma|}\right)^{1/2}
\max_\pm
\left(
\sum_{\rho\neq 1}
\left|
\sum_{\mathfrak p}
P_{r,N,\mathfrak p}^{\pm}(y_{\mathfrak p})
\chi_\rho(c_{\mathfrak p})
\right|^2
\right)^{1/2}.
\label{eqvarying-sandwich}
\end{align}

Enumerate the prime ideals in the interval as
$\mathfrak p_u,\ldots,\mathfrak p_v$ in the fixed norm
order. For $\rho\in\Irr(\Gamma)$, Abel summation gives
\begin{align*}
&
\sum_{0<\|\lambda\|<N}
\left|
\sum_{k=u}^v
c_{\lambda,\mathfrak p_k}^{\pm}
\chi_\lambda(y_{\mathfrak p_k})
\chi_\rho(c_{\mathfrak p_k})
\right|\\&\leq
\sum_{0<\|\lambda\|<N}
\left(
|c_{\lambda,\mathfrak p_v}^{\pm}|
+
\sum_{k=u}^{v-1}
|c_{\lambda,\mathfrak p_k}^{\pm}
-c_{\lambda,\mathfrak p_{k+1}}^{\pm}|
\right)
\max_{u\leq j\leq v}
\left|
\sum_{k=u}^j
\chi_\lambda(y_{\mathfrak p_k})
\chi_\rho(c_{\mathfrak p_k})
\right|
\\
&\quad\ll_G
n_Kd_\rho x^{1/2}
\frac{C_{\partial\mathcal A}}{r^{\delta_G-q}}
\left(
1+\frac{K_{\mathcal A}}{r x^\alpha}
\right)
\log N\,
\log\bigl(q_\rho A_{L/K}(M)(x+C_HN)\bigr).
\end{align*}
Here we used \cref{sum}, the coefficient and variation
bounds from the construction in \cref{Lie}, and
\cref{uni}, followed by the same weight summation as in
the proof of \cref{cst1}.

For $\lambda=0$ and $\rho\neq 1$, similarly
$$
\begin{aligned}
&
\left|
\sum_{k=u}^{v}
c_{0,\mathfrak p_k}^{\pm}
\chi_\rho(c_{\mathfrak p_k})
\right|
\\
&\quad\leq
\left(
|c_{0,\mathfrak p_v}^{\pm}|
+
\sum_{k=u}^{v-1}
|c_{0,\mathfrak p_{k+1}}^{\pm}
-c_{0,\mathfrak p_k}^{\pm}|
\right)
\max_{u\leq j\leq v}
\left|
\sum_{k=u}^{j}\chi_\rho(c_{\mathfrak p_k})
\right|
\\
&\quad\ll_G
n_Kd_\rho x^{1/2}
\log\bigl(q_\rho A_{L/K}(M)x\bigr)
\left[
A_{\mathrm{ST}}
+
C_{\partial\mathcal A}
\left(
r^q+
\frac1{r^{s-q}N^{s-d}}
+
\frac{K_{\mathcal A}}{r x^\alpha}
\right)
\right].
\end{aligned}
$$

Finally, the conductor-discriminant formula
\eqref{cond} and Minkowski's inequality give
$$\left(
\sum_\rho d_\rho^2
\log^2\bigl(A_{L/K}(M)q_\rho(x+C_HN)\bigr)
\right)^{1/2}\leq
|\Gamma|^{1/2}
\log\bigl(S_{L/K}(M)(x+C_HN)\bigr).
$$
Inserting these bounds into
\eqref{eqvarying-sandwich} yields
\begin{equation}\label{eqvarying-unbalanced}
\begin{aligned}
&\left|
\pi_{\mathcal A,C}(x)-\pi_{\mathcal A,C}(x/2)
-\frac{|C|}{|\Gamma|}
\sum_{\mathfrak p}\mu_{\mathrm{ST}}(A_{\mathfrak p})
\right|\\
&\quad\ll_{G,s}
n_KC_{\partial\mathcal A}\frac{|C|}{|\Gamma|}
\frac{x}{\log x}
\left(
r^q+\frac1{r^{s-q}N^{s-d}}
\right)\\&+
n_K|C|^{1/2}x^{1/2}
\log\bigl(S_{L/K}(M)(x+C_HN)\bigr)
\nonumber\\
&\times
\Biggl[
A_{\mathrm{ST}}
+
C_{\partial\mathcal A}
\left(
r^q+\frac1{r^{s-q}N^{s-d}}
+\frac{K_{\mathcal A}}{r x^\alpha}
\right)
+
\frac{C_{\partial\mathcal A}}{r^{\delta_G-q}}
\left(
1+\frac{K_{\mathcal A}}{r x^\alpha}
\right)
\log N
\Biggr].
\end{aligned}
\end{equation}
Choose $s=d+1$, $N=\left\lceil
\max\{1+2\|\varrho\|,\,3,\,3\|\varrho\|,\,r^{-(d+1)}\}
\right\rceil$, and 
$$
\begin{aligned}
r
&:=
\max\Biggl\{
\left(
\frac{|\Gamma|}{|C|^{1/2}}
\frac{(\log x)^2\log(S_{L/K}(M)x)}{x^{1/2}}
\right)^{1/\delta_G},
\\
&\hspace{24mm}
\left(
\frac{K_{\mathcal A}|\Gamma|}{|C|^{1/2}}
\frac{(\log x)^2\log(S_{L/K}(M)x)}
{x^{1/2+\alpha}}
\right)^{1/(\delta_G+1)}
\Biggr\},
\\
\end{aligned}
$$
For $x>x'_4$, these choices satisfy $r<r_0$,
$N\ll_G x$, and $\log N\ll_G\log x$.
Moreover,
$$
\frac1{r^{s-q}N^{s-d}}\leq r^q,
\qquad
\log\bigl(S_{L/K}(M)(x+C_HN)\bigr)
\ll_G\log\bigl(S_{L/K}(M)x\bigr).
$$
Since $q\leq\delta_G$ and $N\geq3$, the term
$K_{\mathcal A}/(r x^\alpha)$ in
the preceding dyadic bound is absorbed by
$$
\frac1{r^{\delta_G-q}}
\left(
1+\frac{K_{\mathcal A}}{r x^\alpha}
\right)\log N.
$$
The two branches in the definition of $r$ balance this
spectral contribution against the geometric error.
They also give
$$
r^q\leq{}
(1+K_{\mathcal A})^{q/(\delta_G+1)}
\left(\frac{|\Gamma|}{|C|^{1/2}}\right)^{q/\delta_G}
x^{-\beta_G}
(\log x)^{2q/\delta_G}
\bigl[\log(S_{L/K}(M)x)\bigr]^{q/\delta_G}.
$$
Consequently, the dyadic discrepancy is $\ll_G E(x)$, where
\begin{align*}
&E(x):=
n_KA_{\mathrm{ST}}|C|^{1/2}x^{1/2}
\log\bigl(S_{L/K}(M)x\bigr)
\\&+
n_K L_{\mathcal{A}}
\frac{|C|^{1-q/(2\delta_G)}}{|\Gamma|^{1-q/\delta_G}}
x^{1-\beta_G}
(\log x)^{2q/\delta_G-1}
\bigl[\log(S_{L/K}(M)x)\bigr]^{q/\delta_G}.
\end{align*}

Since $\beta_G\leq q/(2\delta_G)\leq1/2$, comparison
of the powers and logarithms gives, for
$x\geq t\geq x'_4$,
$$
E(t)\leq
\left(\frac{t}{x}\right)^{1/4}E(x).
$$
Apply the dyadic estimate at $x/2^j$ for every $j$ such
that $x/2^j>x'_4$.
The main terms telescope, while
$$
\sum_{\substack{j\geq0\\x/2^j>x'_4}}
E(x/2^j)
\leq
E(x)\sum_{j\geq0}2^{-j/4}
\ll E(x).
$$
The remaining primes contribute at most $n_Kx'_4$ to
the discrepancy. Finally, for $x\geq x'_4$, 
$$
\begin{aligned}
E(x)\geq n_Kx^{1/4}\Biggl[
&A_{\mathrm{ST}}|C|^{1/2}+
C_{\partial\mathcal A}
(1+K_{\mathcal A})^{q/(\delta_G+1)}
\frac{|C|^{1-q/(2\delta_G)}}{|\Gamma|^{1-q/\delta_G}}
\Biggr].
\end{aligned}
$$
Thus $n_Kx'_4=O(E(x))$ for $x>x_4$, by the definition
of $x_4$. This proves the asserted estimate.
\end{proof}
For Chebotarev-Sato-Tate problems with varying sets,
the least-prime bound depends strongly on the global asymptotic
behavior of the Sato-Tate measure of the family of sets.
In \cref{ST} , we record the small-radius asymptotics
of the Sato-Tate measure of geodesic balls for several groups.
\begin{corollary}\label{smallest}
Retain the notation and assumptions of \cref{cst4}.
Suppose that there exists $c_{\mathcal A}>0$ such that
$$
\mu_{\mathrm{ST}}(A_{\mathfrak p})
\geq
\frac{c_{\mathcal A}}{(N\mathfrak p)^\beta}$$
for
$
\mathfrak p\notin S_{M,L},
$
where $0\leq\beta<\beta_G$.
Put
$$
\eta:=\beta_G-\beta,
\qquad
\kappa:=\max\left\{2,\frac{3q}{\delta_G}\right\},
$$
$$
B:=
2+\frac{n_K}{c_{\mathcal A}}
\left[
\frac{A_{\mathrm{ST}}|\Gamma|}{|C|^{1/2}}
+
C_{\partial\mathcal A}
(1+K_{\mathcal A})^{q/(\delta_G+1)}
\left(
\frac{|\Gamma|}{|C|^{1/2}}
\right)^{q/\delta_G}
\right],
$$
and
$$
x_5:=
\max\left\{
x_4,\,
c_0\bigl(1+n_K+|S_{M,L}|+\log d_K\bigr)^8
\right\},
$$
where $c_0>0$ is a sufficiently large effective absolute
constant.
Then there exists an effective constant $C_G>0$,
depending only on $G$ and independent of $\beta$, such that
$$
\min_{\substack{\mathfrak p\notin S_{M,L}\\
y_{\mathfrak p}\in A_{\mathfrak p},\
c_{\mathfrak p}\in C}}
N\mathfrak p
\leq
\max\left\{
x_5,\,
\left[
C_G B
\left(
\frac{1+\log(S_{L/K}(M)B/\eta)}{\eta}
\right)^\kappa
\right]^{1/\eta}
\right\}.
$$
\end{corollary}

\begin{proof}
Under the generalized Riemann hypothesis for $\zeta_K$, a sufficiently large choice of
$c_0$ ensures that
$$
\#\{\mathfrak p\notin S_{M,L}:N\mathfrak p\leq x\}
\geq\frac{x}{2\log x}
\qquad (x>x_5).
$$
Since $\beta\geq0$, it follows that
$$
\frac{|C|}{|\Gamma|}
\sum_{\substack{N\mathfrak p\leq x\\
\mathfrak p\notin S_{M,L}}}
\mu_{\mathrm{ST}}(A_{\mathfrak p})
\geq
\frac{c_{\mathcal A}|C|}{2|\Gamma|}
\frac{x^{1-\beta}}{\log x}.
$$

Put $L(x):=\log(S_{L/K}(M)x)$.
Using $\beta_G\leq1/2$ and $0<q/\delta_G\leq1$,
the sum of the absolute values of the two error terms
in \cref{cst4}, divided by this lower bound, is
$$
\begin{aligned}
&\ll_G
Bx^{-\eta}
\left(
\log x\,L(x)
+
(\log x)^{2q/\delta_G}L(x)^{q/\delta_G}
\right)\\
&\ll_G
Bx^{-\eta}L(x)^\kappa.
\end{aligned}
$$
Consequently, $\pi_{\mathcal A,C}(x)>0$ whenever
$$
x^\eta
>
C'_G B\bigl[\log(S_{L/K}(M)x)\bigr]^\kappa
$$
for a sufficiently large effective constant $C'_G$
depending only on $G$, which holds for
$$
x>
\left[
C_G B
\left(
\frac{1+\log(S_{L/K}(M)B/\eta)}{\eta}
\right)^\kappa
\right]^{1/\eta},
$$
with a sufficiently large effective constant $C_G$
depending only on $G$.
Thus $\pi_{\mathcal A,C}(x)>0$ for every $x$ exceeding
the asserted bound, which proves the result.
\end{proof}

\begin{remark}
Our estimates also yield bounds for sums of the form
$$
\sum_{\substack{N\pp\leq x\\ \pp\notin S_{M,L}}}
f(y_\pp,c_\pp),
$$
where $f$ is a class function on $G\times\Gamma$
whose pullback to a maximal torus is smooth in the
first variable.
One may apply the generalized Koksma--Hlawka inequality
of Aistleitner and Dick
\cite[Theorem 1]{AistleitnerDick2015}
to the probability measure
$$
\frac{|D(\theta)|^2}{|W|}\,d\theta
\otimes
\sum_{C\in\operatorname{Conj}(\Gamma)}
\frac{|C|}{|\Gamma|}\,\delta_{C},
$$
where $d\theta$ is normalized Haar measure on $T$.
This gives the corresponding Chebotarev-Sato-Tate
formalism considered in \cite[Hypothesis 4.4]{Yang2026}. \end{remark}

\section{Improvements in the error exponents}\label{tserre}
As seen in the previous estimates, the exponent in the error term
deteriorates significantly as the dimension of the motivic
Sato-Tate group increases. We illustrate an approach to improving the error term by exploiting the available symmetries of the Chebotarev-Sato-Tate class. Our idea is similar in principle to Serre's treatment of
effective Chebotarev estimates
\cite[Section 5, Theorem 12, and Section 7.3]
{Serre1981}.

Let $N\subset G$ be a closed normal subgroup and $T_N$ the maximal torus of $G/N$ of dimension $d_N$ and denote by $m_N$ the number of its positive roots. 
\begin{definition}
Let $g_{T,N}$ be the induced quotient metric on $T_N$.
Let $p_N:G\to G/N$ be the quotient map and equip $G/N$
with the quotient metric. Define $d_N:=\mathrm{rank}(G/N)$ and for $0<q\leq d_N$ define
$$
\mathfrak M^q_{g,N}:=
\{p_N^{-1}(B):B\in\mathfrak M^q_{g_{G/N}}\}=\{A\in \mathfrak M^q_{g}\,, NA\subset A\}.
$$
Let $A_{T,N}$ be the inverse image in the quotient torus $T/T\cap N$.
Define $$C_{\partial A,N}=C_\partial(A,r_0,N):=\max\left(1,\sup_{0<r<r_0}\frac{\mathrm{vol}_{g_{T,N}}((\partial A_{T,N})^{(r)})}{r^q}\right).$$
\end{definition}
\begin{proposition}\label{imp}
For $A\in\mathfrak M_{g,N}^q$,
the exponent $\delta_G$ in the preceding estimates
can be replaced by
$$
\delta_N:=d_N+m_N.
$$ The implied constants may depend on $G$ and $N$.

Moreover, if $H\triangleleft\Gamma$ satisfies $HC=C$, the Galois parameters become
$$
\left(
|\Gamma/H|,\ |C/H|,\ d_{L^H}
\right).
$$
\end{proposition}

\begin{proof}
Identify $A$ with a subset of $G/N$. We apply the Lie-group approximation \cref{Lie} to the compact connected group $G/N$
with its quotient metric. This yields majorant/minorant $P_{r,T,A}^{N\pm}$ such that 
$P_{r,T,A}^{N-}\circ \pi\leq \ind_A\leq P_{r,T,A}^{N+}\circ \pi $ where $\pi:G\rightarrow G/N$ is the projection.
There are $O_{G,N}(k^{d_N-1})$ surviving dominant weights
associated to characters $\chi_\lambda$ trivial on $N$, and the Weyl dimension formula gives
$$
d_\lambda\ll_{G,N}(1+\|\lambda\|)^{m_N}.
$$
Using the coefficient estimate in \cref{Lie} with
order $\delta_N=d_N+m_N$, the corresponding spectral
sum is therefore bounded by
$$
\sum_{\substack{0<\|\lambda\|<R\\
N\subseteq\ker\pi_\lambda}}
\frac{d_\lambda}{\|\lambda\|^{\delta_N}}
\ll_{G,N}
1+\sum_{1\leq k<R}\frac1k
\ll_{G,N}\log(2R).
$$
The remaining arguments are unchanged, with
$\delta_G$ replaced by $\delta_N$.
If we assume moreover that $HC=C$ for a normal subgroup, the indicator $\ind_C$ factors through
$\Gamma/H$. In particular, the characters can be identified with characters of the extension $L^H/K$. We may therefore apply the orthogonality
argument of \cite[Proof of Lemma 7.2]{ZywinaLangTrotter}
directly to $\Gamma/H$, whose irreducible characters
identify with those of $\Gamma$ trivial on $H$.
\end{proof}
\section{Beyond the CST hypotheses}
We conclude this section by discussing a possible reduction
when the CST hypotheses hold after a finite base extension.
Fix a motivic category $\mathcal M_\sim(K,\Q)$ satisfying the Hypothesis
\ref{motiv}. 

Let $M\in\mathcal M_\sim(K,\Q)$ be a homogeneous motive.
Assume that there exists a finite extension $K_0/K$,
independent of $\ell$, such that
$$
G_{\ell,K_0}^{\mathrm{alg}}
=
(G_{\ell,K}^{\mathrm{alg}})^\circ,
\qquad
\st_{K_0}(M\times_K {K_0})=G^\circ.
$$
For the corresponding connectedness statements in odd weight,
see \cite[Theorem 5.5 and Remark 5.9]{BanaszakKedlaya2023}.
For example, if $M=h^1(C)$ for an elliptic curve $C/K$
with complex multiplication by an imaginary quadratic
field $F$, one may take $K_0=KF$.

Assume that there exists a finite Galois extension
$E/K$, containing $K_0$, such that the base change $M\times _K E\in \mathcal{M}_\sim(E,\Q)$ is
a CST-motive and
$$
G^\circ:=\st_E(M\times _K E)=G^\circ.
$$

Let $T\subset G^\circ$ be a maximal torus. For every $n\geq1$,
the power map
$$
f_n:T\longrightarrow T,\qquad g\longmapsto g^n,
$$
is a smooth covering of degree $n^{\dim T}$, with differential
$n\,\mathrm{id}$.

Let $A\subset G$ be conjugacy-invariant, and put
$$
\widetilde A_s
:=
\{g^s:g\in\overline A\}\cap G^\circ,
\qquad
\widetilde A
:=
\bigcup_{s\mid[E:K]}\widetilde A_s.
$$
Let $L/K$ be the finite Galois extension associated with the
Chebotarev condition $C\subset\Gal(L/K)$, and put $E':=EL$.
Define
$$
\widetilde C_s
:=
\left\{
\gamma\in\Gal(E'/E):
\gamma|_L\in\{c^s:c\in C\}
\right\},
\qquad
\widetilde C
:=
\bigcup_{s\mid[E:K]}\widetilde C_s.
$$

Outside a fixed finite set of bad or ramified primes,
let $\mathfrak P$ be a prime of $E$ above $\pp$, with residue
degree $f=f(\mathfrak P/\pp)$. The normalized Frobenius classes
$y'_{\mathfrak P}$ of $M\times_KE$ is conjugate to $y_\pp^f$, while
$$
N\mathfrak P=(N\pp)^f,
\qquad
c_{\mathfrak P}|_L\in c_\pp^f.
$$
Since $f\mid[E:K]$, we obtain
$$
\pi_{A,C}^{K}(x)
\leq
\sum_{s\mid[E:K]}
\frac{s}{[E:K]}\,
\pi_{\widetilde A_s,\widetilde C_s}^{E}(x^s)
+
O_{M,E,L}(1),
$$
where the counting functions on the right are taken over
$E$, for $M_E$ and the extension $E'/E$. Indeed, every
unramified prime of residue degree $s$ has exactly
$[E:K]/s$ primes above it.

Our estimates over $E$ therefore provide upper bounds
whenever the new sets satisfy the required
geometric hypotheses. For example, this holds if the power maps $f_s$ are local
diffeomorphisms on a neighborhood of $\overline A$.

\chapter{Arithmetic applications}\label{sato2}
In this section, we exhibit several arithmetic applications of the
theorems of Section\ref{equid}, together with open questions for
future research. Throughout this section, $K/\Q$ denotes a number
field. The corresponding motivic categories are assumed to satisfy
the Banaszak--Kedlaya motivic hypotheses\ref{motiv}.
Given the connected Sato--Tate group $G:=\st_K(M)$ of a CST-motive
$M$, we equip it with a suitable bi-invariant Riemannian metric $g$
whose volume measure is the normalized Haar measure $\mu_G$.

The maximal torus $T$ carries the induced $W$-invariant flat metric
$g_T:=g|_T$. Recall that $T$ is geodesically connected, as stated
in \cite[Lemma 8.2]{SecoSanMArtin2018}. In particular, the restriction
of the ambient distance $d_g$ to $T$ agrees with the intrinsic
distance $d_{g_T}$. We fix an isometric identification
$T\simeq\bb{T}^{\operatorname{rank}G}
:=(\bb{R}/2\pi\bb{Z})^{\operatorname{rank}G}$
with the corresponding flat $W$-invariant metric.
Given a conjugacy-invariant subset $A$, we denote by $A_T$ its
pullback to $\bb{T}^{\operatorname{rank}G}$ under the Weyl isomorphism
\begin{equation}\label{we}
    \bb{T}^{\operatorname{rank}G}/W\simeq\mathrm{Conj}(G).
\end{equation}
The Weyl isomorphism \eqref{we} is an isometry for the quotient
metrics \cite[Remarks 4.43 and 4.10]{AlexandrinoBettiol2015}
and restricts to an isometry between $A^{(\varepsilon)}$ and
$A_T^{(\varepsilon)}/W$. Their volumes are related by the Weyl
integration formula
\begin{equation}\label{int}
    \mu_{\st}(A^{(\varepsilon)})
    =\frac{1}{|W|}\int_{A_T^{(\varepsilon)}}|D(\theta)|^2d\theta,
\end{equation}
where $d\theta=\frac{d\s}{\s(T)}$.

We equip $\mathrm{Conj}(G)$ with the Sato--Tate measure
$\mu_{\mathrm{ST}}$, defined as the pushforward of the normalized
Haar measure $\mu_G$. We fix a finite Galois extension $L/K$ and
write $\Gamma=\mathrm{Gal}(L/K)$. Here $C\subset\Gamma$ denotes a
nonempty union of conjugacy classes. For a prime $\pp$ unramified
in $L$, let $c_\pp$ denote its Frobenius conjugacy class.

Throughout this section, we fix $0<r_0<1$. For a set $A\subset G$
with nonempty boundary such that $\partial A_T$ has finite
Minkowski content, we denote by $C_{\partial A}$ the constant
in \cref{min}. Similarly, for a family
$\mathcal{A}=(A_i)_{i\in I}$ such that
$\partial\mathcal{A}_T=(\partial A_{i,T})_{i\in I}$
has finite uniform Minkowski content, we denote by
$C_{\partial\mathcal{A}}$ the constant in \cref{uni}.
Throughout this section, we assume the following for each CST-motive $M$:

\begin{hypothesis}
There exist a number field $E$ splitting the motivic Galois
group $G_M$ and an embedding $\iota:E\hookrightarrow\bb{C}$
such that the Artin motives $N_\rho$ are defined with
coefficients in $E$ and the $L$-functions
$L(\iota,M_\sigma\otimes N_\rho,s)$ satisfy the motivic $L$-functions formalism of Conjecture \ref{grh} for every
$
(\sigma,\rho)\in
\operatorname{Irr}(G)\times
\operatorname{Irr}(\Gamma).
$
\end{hypothesis}
\section{Random matrices and eigenvalue repulsion for generic motives}
Fix an integer $n\geq 2$. Call a motive $M$ $n$-generic if there is an isomorphism of represented groups $$
\bigl(\st_K(M),H_B(M)\otimes_\Q\bb{C}\bigr)
\simeq
\bigl(\mathrm{USp}(2n),\bb{C}_{\mathrm{std}}^{2n}\bigr),
$$
where $\mathrm{USp}(2n)$ is the compact unitary symplectic group
and $\bb{C}_{\mathrm{std}}^{2n}$ denotes its standard representation.
Its Sato-Tate measure is given by
    $$d\mu_{\st}(\mathbf{\theta})=\frac{2^{n(n-1)}}{n!\pi^n}\prod_{1\leq i<j\leq n}(\cos\theta_i-\cos\theta_j)^2\prod_{j=1}^n\sin^2\theta_jd\mathbf{\theta}.$$

Examples of generic motives include
$$
h^1(A)\in\mathcal M_{\mathrm{ahc}}(K,\mathbb Q),
$$
where $A/K$ is an abelian variety with potentially
semistable reduction of toric rank one at some finite place
and $\operatorname{End}_{\overline K}(A)=\mathbb Z$ (cf. \cite[Theorem 1]{Hall2011}).
Let $M$ be a CST-$n$-generic motive in a suitable motivic category. Recall that a complex random matrix is a measurable map $X:\Omega\to\operatorname{Mat}_n(\mathbb C)$ from a probability space ($\Omega,\mathcal F,\mathbb P$) with values in an $n\times n$ complex matrix space. For a compact group $G$, a Haar-distributed random matrix takes values in $G$ and has normalized Haar measure as its law.
In \cite[Theorem 3]{Mezzadri2007}, a Haar-distributed random matrix $U$ taking values in
$\mathrm{Sp}(n)\simeq\mathrm{USp}(2n)$ is constructed by applying the
quaternionic QR decomposition to an $n\times n$ Gaussian
quaternionic matrix.
The motivic character estimates of \cref{sum} together with the discussion of Section \ref{shimura} applied to the conjugacy Foliation imply that 
$$\frac{1}{\pi_K(x)}\sum_{\substack{N\pp}\leq x,\,\pp\notin S_M}f(y_\pp)\rightarrow \bb{E}_{\mu_{\st}}[f(U)],$$
for any bounded Borel function $f:\mathrm{Conj}(\mathrm{USp}(2n))\rightarrow{ \bb{R}}$ whose set of discontinuities has zero Sato-Tate measure. Informally, the normalized Frobenius conjugacy classes $(y_\pp)_{\pp\notin S_{M,L}}$ is asymptotically distributed in $\mathrm{USp}(2n)$ according to the law of the random Haar matrix $U$. Eigenvalue repulsion for Haar-distributed matrices in
compact classical groups is manifested in the collision
factors of the Weyl eigenangle densities
see \cite[Section 3.1]{Meckes2019}. In the classical Gaussian situation, the
joint eigenvalue density vanishes at eigenvalue collisions
\cite{Trefethen2026}.
\def \t {\mathrm{\theta}}
Define the function
\begin{align*}
    S&:\mathbb{T}^n\rightarrow \bb{R}\\
    &\mathbf{\theta}\rightarrow \min_{\pm}\min_{1\leq a<b\leq n}|e^{i\theta_a}-e^{\pm i\theta_b}|.
\end{align*}
Fix $0<{\varepsilon}<\varepsilon_0:=\min\left(1,2\sin\left(\frac{\pi}{2(n-1)}\right)\right)$. Define the set 
$$A({\varepsilon}):=\{\theta\in \mathbb{T}^n, |S(\mathbf{\theta})|<{\varepsilon}\}.$$
\def  \A {A({\varepsilon})}
The family of boundaries
$$
\partial A(\varepsilon)
=
\{\boldsymbol{\theta}\in\mathbb T^n:
S(\boldsymbol{\theta})=\varepsilon\}\neq \varnothing,
\qquad 0<\varepsilon<\varepsilon_0,
$$
is a subfamily of the subanalytic family
$$
S_{\theta'}
:=
\{\mathbf t\in\mathbb T^n:
S(\mathbf t)=|\cos\theta'|\},
\qquad \theta'\in\mathbb T.
$$
Moreover,
$$
\sup_{\theta'\in\mathbb T}\dim_{\mathrm{H}} S_{\theta'}=n-1.
$$
Thus \cref{analytictube} implies that
$(\partial A(\varepsilon))_{0<\varepsilon<\varepsilon_0}$
has uniformly bounded $(n-1)$-dimensional Minkowski content.
The Weyl integration formula \cref{Wey} gives
$$
\mu_{\mathrm{ST}}(A(\varepsilon))
=
c_n\varepsilon^3+O_n(\varepsilon^5),
\qquad
c_n:=\frac{n(n-1)(n+1)(2n+1)}{18\pi}.
$$
Identify $S$ with a function of $\mathrm{Conj}(\mathrm{USp(2n))}$ according to the convention of the section.
$$\pi_{A(\varepsilon)}(x)
:=\#\{\mathfrak p\notin S_M:N\mathfrak p\leq x,
S(y_{\mathfrak p})<\varepsilon\}.$$
Applying \cref{cst1} yields the following estimate:

\begin{proposition}
Let $n\ge2$.
There exists an effectively computable $x_0>0$,
independent of $\varepsilon$, such that, for every
$0<\varepsilon<\varepsilon_0$ and $x>x_0$,
$$
\begin{aligned}
\pi_{A(\varepsilon)}(x)
={}&c_n\varepsilon^3\Li(x)
+O_n\bigl(\varepsilon^5\Li(x)\bigr)
\\
&+
O_n\left(
n_Kx^{1-\frac{1}{2n(n+1)}}
(\log x)^{\frac{2}{n(n+1)}-1}
\bigl(\log(S_{K/K}(M)x)\bigr)^{\frac{1}{n(n+1)}}
\right).
\end{aligned}
$$
The implied constants are effective in terms of $n$ and independent of
$\varepsilon$.

In particular, if $\varepsilon=\varepsilon(x)>0$ satisfies
$$
\varepsilon(x)\longrightarrow0,
\qquad
\frac{
n_K(\log x)^{\frac{2}{n(n+1)}}
\bigl(\log(S_{K/K}(M)x)\bigr)^{\frac{1}{n(n+1)}}
}{
\varepsilon(x)^3x^{\frac{1}{2n(n+1)}}
}
\longrightarrow0,
$$
then as $x\rightarrow +\infty$ we have
$$
\frac{\pi_{A(\varepsilon(x))}(x)}
{\varepsilon(x)^3\Li(x)}
\longrightarrow c_n.
$$
\end{proposition}

\begin{question}
The eigenvalue repulsion obtained above invites comparison
with the GUE models for the local statistics of zeros of
$L$-functions. These models originate in Montgomery's pair
correlation conjecture \cite{Montgomery1973}. Related
correlation results for primitive automorphic $L$-functions
were obtained, under restrictions on the test functions, by
Rudnick and Sarnak \cite{RudnickSarnak1996}.
See also \cite{KatzSarnak1999Symmetry}.

The explicit formula in the proof of \cref{sumc}, obtained through Mellin inversion
and contour shifting, relates sums over zeros of the
corresponding $L$-functions to prime-power sums involving
Artin-Sato-Tate characters evaluated at powers of
Frobenius conjugacy classes. Although the two statistics
concern different limiting regimes, we conjecture that
the random matrix behavior of normalized Frobenius
eigenvalues and the conjectural GUE behavior of high-lying
zeros reflect a common arithmetic mechanism. A precise
relationship between these statistics remains to be
formulated. A quantitative form of such a relationship
might lead to refined conditional estimates in effective
Chebotarev-Sato-Tate theory. Improvement of the classical effective Chebotarev density theorem was obtained by \cite{MurtyMurtyWong2018} under a joint pair correlation conjecture for the family of the Artin $L$-functions $(L(s,\rho))_{\rho\in \mathrm{Irr}{(\Gamma)}}$.
\end{question}

\section{Generalized Atkin-Serre for families of rank $2$ motives}
The Atkin-Serre conjecture \cite[p. 244]{Serre1976} predicts that the coefficients of a non-CM newform $f$ of even weight $k\geq 4$ satisfy $
|a_f(p)|\gg_{f,{\varepsilon}}p^{(k-3)/2-{\varepsilon}}
$ for every sufficiently large prime $p$. Rouse \cite{Rouse2007} obtained qualitative bounds for primes conditioned on the automorphic Riemann hypothesis of the $L$-functions $(L(s,\mathrm{sym}^mf))_{m\geq 0}$. Subsequently, Gafni, Thorner, and Wong proved unconditionally that the conjectural lower bound holds for a set of primes of density $1$, with an effective estimate for the
exceptional set \cite{GafniThornerWong2021}. In this subsection, we shall formulate a motivic variant of the Atkin-Serre problem with an additional Chebotarev condition by applying our effective Chebotarev-Sato-Tate for a family of sets. In particular, these estimates apply to modular forms arising from modular motives satisfying
the CST hypotheses.

Let $M=\bigoplus_{i=1}^n M_i$ be a CST-motive, where the $M_i$'s
are rank-two polarized motives of common weight $w$, and assume we have an isomorphism of represented groups
$$
(\mathrm{ST}_K(M),H_B(M)\otimes_{\Q}\mathbb C)
\simeq
\left(\mathrm{SU}(2)^n,
\bigoplus_{i=1}^n\mathbb C^2_{\mathrm{std},i}\right).
$$
For every prime ideal $\mathfrak p\notin S_{M,L}$ define $$
a_i(\mathfrak p)
:=\operatorname{tr}\left(
\operatorname{Frob}_{\mathfrak p}\mid H_{\ell}(M_{i})
\right)
=2(N\mathfrak p)^{w/2}\cos\theta_{i,\mathfrak p},
$$
where $\theta_{i,\mathfrak p}\in \bb{T}^1$.
Let $0<{\varepsilon}<1$. We consider the following generalized Atkin-Serre-type counting problem:
$$
\pi_{\mathrm{AS},C}(x,\varepsilon)
:=
\#\left\{
\begin{gathered}
\pp\notin S_{M,L},\
N\pp\leq x,\quad c_\pp\in C,
\\
|a_i(\pp)|<2(N\pp)^{w/2-\varepsilon}
\\
\text{for some }1\leq i\leq n
\end{gathered}
\right\}.
$$
Define the family of sets $\mathcal{A}=(A_\pp)_{\pp\notin S_{M,L}}$, given by
$$A_\pp:=\bigcup_{1\leq i\leq n} \left\{\t\in \mathbb{T}^n,|\cos\theta_i|< \frac{1}{(N\pp)^{\varepsilon}}\right\}.$$ 
\begin{proposition}\label{Serre}
Put
$$
\beta_{\varepsilon}
:=
\frac12\min\left\{
\frac1{2n},\frac{1+2\varepsilon}{2n+1}
\right\}.
$$
\begin{enumerate}
\item
For every $0<\varepsilon<1$, there exists an effectively
computable $x_0>0$ such that, for $x>x_0$,
\begin{align*}
&\pi_{AS,C}(x,\varepsilon)
=
\frac{4n}{\pi}\frac{|C|}{|\Gamma|}
\sum_{\substack{\pp\notin S_{M,L}\\N\pp\le x}}
\frac1{(N\pp)^\varepsilon}
+
O_n\left(
\frac{|C|}{|\Gamma|}
\sum_{\substack{\pp\notin S_{M,L}\\N\pp\le x}}
\frac1{(N\pp)^{2\varepsilon}}
\right)
\\
&+
O_n\left(
n_K|C|^{1/2}x^{1/2}
\log\bigl(S_{L/K}(M)x\bigr)
\right)
\\[2pt]
&+
O_n\left(
n_K
\frac{|C|^{1-1/(4n)}}{|\Gamma|^{1-1/(2n)}}
x^{1-\beta_{\varepsilon}}
(\log x)^{1/n-1}
\bigl(\log(S_{L/K}(M)x)\bigr)^{1/(2n)}
\right).
\end{align*}
In particular, for fixed
$
0<\varepsilon<\frac1{4n},
$
as $x\rightarrow\infty$ we have
$$
\pi_{AS,C}(x,\varepsilon)
\sim
\frac{|C|}{|\Gamma|}
\frac{4n}{\pi(1-\varepsilon)}
\frac{x^{1-\varepsilon}}{\log x}.
$$

\item
For every $0<\varepsilon\le1/(8n)$, there exists a prime
ideal $\mathfrak p\notin S_{M,L}$ such that
$c_{\mathfrak p}\in C$ and
$$
|a_i(\mathfrak p)|
<
2(N\mathfrak p)^{w/2-\varepsilon}
$$
for some $1\le i\le n$.
Moreover, one may choose $\mathfrak p$ with
$$
\begin{aligned}
N\mathfrak p\ll_n\max\Biggl\{
&x_5,\,
\bigl(1+n_K+|S_{M,L}|+\log d_K\bigr)^8,
\\
&
\left(\frac{n_K|\Gamma|}{|C|^{1/2}}\right)^{8n+4}
\left[
\log\left(
2S_{L/K}(M)\frac{n_K|\Gamma|}{|C|^{1/2}}
\right)
\right]^{16n+8}
\Biggr\},
\end{aligned}
$$
where $x_5=x_5(\mathcal A,C,\varepsilon)>0$ is an
effectively computable threshold.
\end{enumerate}
The implied constants depend only on $n$ and are effective.
\end{proposition}
\begin{proof}
For $\pp\notin S_{M,L}$, the topological boundary satisfies
$$
\partial A_\pp
\subset
\bigcup_{1\le i\le n}
\left\{
\boldsymbol{\theta}\in\mathbb T^n:
|\cos\theta_i|=(N\pp)^{-\varepsilon}
\right\}.
$$
Thus the family
$\partial\mathcal A=(\partial A_\pp)_{\pp\notin S_{M,L}}$
has uniformly bounded $(n-1)$-dimensional upper Minkowski
content:
$$
\overline{\mathcal M}_{g_T}^{n-1}
(\partial\mathcal A)<\infty.
$$
Direct computations give
\begin{align*}
d_H^{g_T}
\left(\overline{A_{\mathfrak p}},
      \overline{A_{\mathfrak p'}}\right)
&\asymp
\left|
\arcsin\bigl((N\mathfrak p)^{-\varepsilon}\bigr)
-
\arcsin\bigl((N\mathfrak p')^{-\varepsilon}\bigr)
\right|,
\\
d_H^{g_T}
\left(
\overline{\mathbb T^n\setminus A_{\mathfrak p}},
\overline{\mathbb T^n\setminus A_{\mathfrak p'}}
\right)
&\asymp
\sqrt n\,
\left|
\arcsin\bigl((N\mathfrak p)^{-\varepsilon}\bigr)
-
\arcsin\bigl((N\mathfrak p')^{-\varepsilon}\bigr)
\right|.
\end{align*}
Enumerating the prime ideals with
$x/2<N\mathfrak p\le x$ in nondecreasing order of norm together with the inequality
$\arcsin t\le(\pi/2)t$ give
$$
\begin{aligned}
\sum_{k=u}^{v-1}
d_*(\mathfrak p_k,\mathfrak p_{k+1})
&\ll
\sqrt n\,
\arcsin\bigl((x/2)^{-\varepsilon}\bigr)
\\
&\ll
\frac{\pi\sqrt n}{2}\,2^\varepsilon x^{-\varepsilon}
\ll_n x^{-\varepsilon}.
\end{aligned}
$$
Hence $\mathcal A$ is an $(\varepsilon,1)$-Minkowski
family, with uniformly bounded constants for
$0<\varepsilon<1$.
The first estimate follows from \cref{cst4}, since
$\delta_G=2n$, together with
$$
\mu_{\mathrm{ST}}(A_{\mathfrak p})
=
\frac{4n}{\pi}(N\mathfrak p)^{-\varepsilon}
+
O_n\bigl((N\mathfrak p)^{-2\varepsilon}\bigr).
$$

For fixed $0<\varepsilon<1$, partial summation and the
prime ideal theorem give
$$
\sum_{\substack{\pp\notin S_{M,L}\\N\pp\le x}}
\frac1{(N\pp)^\varepsilon}
\sim
\frac{x^{1-\varepsilon}}{(1-\varepsilon)\log x}.
$$
If $0<\varepsilon<1/(4n)$, then
$\beta_\varepsilon>\varepsilon$, so all the error terms
are of lower order than the main term. This proves the
asymptotic statement.

For the second assertion, the elementary inequality
$$
\mu_{\mathrm{ST}}(A_\pp)
\ge (N\pp)^{-\varepsilon}
$$
allows us to apply \cref{smallest} with
$c_{\mathcal A}=1$ and $\beta=\varepsilon$.
For $0<\varepsilon\le1/(8n)$, we have
$$
\eta
=
\beta_\varepsilon-\varepsilon
=
\frac{1-4n\varepsilon}{4n+2}
\ge\frac1{8n+4}.
$$
Moreover, $A_{\mathrm{ST}}\le1$ and the uniform bounds
above imply
$$
B\ll_n\frac{n_K|\Gamma|}{|C|^{1/2}}.
$$
Substituting these bounds into \cref{smallest} gives the
asserted estimate.
\end{proof}

\section{A Serre-type least-prime problem}\label{minimalserre}
We remain within the setting of the previous subsection.
Serre studied the effective nonvanishing of polynomial expressions
$F(p,a_p,a'_p)$ in the Frobenius traces of two elliptic curves and obtained a conditional bound on the generalized Riemann hypothesis for the Dedekind zeta function of the division fields \cite[8.3, Remark 2]{Serre1981}. He further indicated the possibility of analogous conditional extensions to other compatible systems
of Galois representations of multiple elliptic curves. We confirm Serre's intuition.

We consider a motivic version of this
nonvanishing problem with an additional Chebotarev condition.
Fix $\alpha>0$ and put $a:=2^{-\alpha}$.
Let $F\in C^\infty([0,a]\times\mathbb T^n)^W$
be a real-valued smooth function, invariant under
the Weyl group in the second variable.
For $t\in[0,a]$, define
$$
Z_t:=\{\theta\in\mathbb T^n:F(t,\theta)=0\}.
$$
Assume that the fibers $(Z_t)_{t\in [0,a]}$ are non-empty
and that
\begin{equation}\label{Fregular}
\inf_{\substack{
t\in[0,a],\,\theta\in Z_t
}}
|D_\theta F(t,\theta)|_{g_T}>0.
\end{equation}

By Ehresmann's fibration theorem, there is a smooth family of
diffeomorphisms $\Phi_{t,t'}:Z_{t'}\to Z_t$. We give an explicit
description. In a neighborhood of the zero set, define the vector field
$$
V_t(\theta)
=
-\frac{\partial_tF(t,\theta)}
{|\nabla_\theta F(t,\theta)|_{g_T}^2}
\nabla_\theta F(t,\theta).
$$
The diffeomorphisms $(\Phi_{t,t'})$ are obtained by solving the time-dependent ODE below (see
\cite[Theorem 9.48 and Problem 9-20]{lee2013smooth}) .
For $p\in Z_{t'}$, they satisfy
\begin{equation}\label{ODE}
    \frac{\partial}{\partial t}\Phi_{t,t'}(p)
=
V_t\bigl(\Phi_{t,t'}(p)\bigr),
\qquad
\Phi_{t',t'}(p)=p.
\end{equation}

This regularity hypothesis gives the following conditional
bound for the least prime ideal satisfying
$F((N\pp)^{-\alpha},\theta_\pp)\ne0$
and $c_\pp\in C$.

\begin{proposition}
Put
$$
\eta_\alpha
:=
\frac12\min\left\{
\frac1{2n},\frac{1+2\alpha}{2n+1}
\right\}>0.
$$
There exists an effectively computable constant
$B=B(C,F,\alpha,L/K)>0$ and a prime ideal
$\mathfrak p\notin S_{M,L}$ such that
$$
F\bigl((N\mathfrak p)^{-\alpha},
\theta_{\mathfrak p}\bigr)\ne0,
\qquad c_{\mathfrak p}\in C,
$$
with
$$
N\mathfrak p
\ll_n
\max\left\{
x_5,\,
B^{1/\eta_\alpha}
\bigl(1+\log(S_{L/K}(M)B)\bigr)^{2/\eta_\alpha}
\right\}.
$$
Here $x_5$ is the effectively computable threshold defined in \cref{smallest}.  
\end{proposition}

\begin{proof}
The map
$$
[0,a]\longrightarrow\mathcal K(\mathbb T^n),
\qquad t\longmapsto Z_t,
$$
is Lipschitz. Indeed, solving \cref{ODE} gives diffeomorphisms
$\Phi_{t,t'}:Z_{t'}\to Z_t$. By compactness and
\eqref{Fregular}, the speed of this flow is uniformly bounded.
Consequently, we have
$$
\begin{aligned}
d_H^{g_T}(Z_t,Z_{t'})
&\leq
\sup_{\theta\in Z_{t'}}
d_{g_T}\bigl(\theta,\Phi_{t,t'}(\theta)\bigr)
\\
&\leq K_\mathcal{A}|t-t'|,
\end{aligned}
$$
with 
$$
K_\mathcal{A}:=
\sup_{\substack{t,t'\in[0,a]\\ \theta\in Z_{t'}}}
\left|
\left.\frac{d}{ds}\Phi_{s,t'}(\theta)\right|_{s=t}
\right|_{g_T}
<\infty.
$$
For $\pp\notin S_{M,L}$, define
$$
A_{\mathfrak p}
:=
\{\theta\in\mathbb T^n:
F((N\mathfrak p)^{-\alpha},\theta)\ne0\}.
$$
The regularity hypothesis gives
$$
\partial A_{\mathfrak p}
=Z_{(N\mathfrak p)^{-\alpha}},
\qquad
\overline{A_{\mathfrak p}}=\mathbb T^n.
$$
Moreover, compactness and \cref{tube3} imply that
$\partial\mathcal A=(\partial A_\pp)_{\pp\notin S_{M,L}}$
has uniformly bounded $(n-1)$-dimensional upper Minkowski
content.

Since each $A_\pp$ is dense in $\mathbb T^n$ and
$\mathbb T^n\setminus A_\pp=Z_{(N\pp)^{-\alpha}}$,
we obtain
$$
d_*(\mathfrak p,\mathfrak p')
=
d_H^{g_T}\left(
Z_{(N\mathfrak p)^{-\alpha}},
Z_{(N\mathfrak p')^{-\alpha}}
\right).
$$

Enumerating the prime ideals with
$x/2<N\mathfrak p\le x$ in nondecreasing order of norm,
we obtain
$$
\begin{aligned}
\sum_{k=u}^{v-1}d_*(\mathfrak p_k,\mathfrak p_{k+1})
&\le K_\mathcal{A}\sum_{k=u}^{v-1}
\left|
(N\mathfrak p_k)^{-\alpha}
-(N\mathfrak p_{k+1})^{-\alpha}
\right|
\\
&\le K_\mathcal{A}(2^\alpha-1)x^{-\alpha}.
\end{aligned}
$$
Thus $\mathcal A$ is an $(\alpha,1)$-Minkowski family.

Since each $Z_t$ is a smooth hypersurface, the Weyl
integration formula gives
$$
\mu_{\mathrm{ST}}(A_{\mathfrak p})=1.
$$
We may therefore apply \cref{smallest} with
$c_{\mathcal A}=1$, $\beta=0$, and $q=1$.
As $\delta_G=2n$, we have
$$
\eta=\beta_G
=
\frac12\min\left\{
\frac1{2n},\frac{1+2\alpha}{2n+1}
\right\}
=\eta_\alpha.
$$
Moreover, $\kappa=2$. Since
$$
\frac1{4n+2}\le\eta_\alpha\le\frac1{4n},
$$
the factors depending only on $\eta_\alpha$ in
\cref{smallest} can be absorbed into the implied
constant depending on $n$.
The asserted bound follows.
\end{proof}

We now record some examples of functions satisfying
the regularity hypothesis \eqref{Fregular}.

\begin{example}
We record several examples of functions satisfying our smoothness assumption:

\begin{enumerate}
\item 
Consider
$$
F(t,\theta)
:=
\cos\theta_1-\frac{1+t}{4}.
$$
The zero fibers are given by
$$
Z_t
=
\left\{
\theta\in\mathbb T^n:
\cos\theta_1=\frac{1+t}{4}
\right\}.
$$
For a rank-two motive with normalized Frobenius trace
$2\cos\theta_{1,\pp}$, the corresponding nonvanishing
condition becomes
$$
2\cos\theta_{1,\pp}
\ne
\frac{1+(N\pp)^{-\alpha}}{2}.
$$

\item 
For $n\geq2$, consider
$$
F(t,\theta)
:=
\cos\theta_1-\frac12
+\frac t8\cos\theta_2.
$$
This gives a family of moving hypersurfaces
involving two Frobenius angles.

\item 
Consider
$$
F(t,\theta)
:=
\prod_{j=1}^n\cos\theta_j-\frac{1+t}{4}.
$$
This gives a polynomial nonvanishing condition
involving the normalized Frobenius traces of
a family of rank-two motives.

\item
Consider
$$
F(t,\theta)
:=
\sum_{j=1}^n\cos\theta_j-c(t),
$$
where
$$
c(t):=
\begin{cases}
\displaystyle\frac12+\frac t4,
&\text{if }n\text{ is even},\\[2mm]
\displaystyle\frac t4,
&\text{if }n\text{ is odd}.
\end{cases}
$$
The values of $c(t)$ avoid the critical values
$$
\{-n,-n+2,\ldots,n-2,n\}
$$
of the function
$\theta\mapsto\sum_{j=1}^n\cos\theta_j$.

\item
More generally, let
$H\in C^\infty(\mathbb T^n)^W$ be real-valued,
and let $J\subset\mathbb R$ be a compact interval
consisting entirely of regular values of $H$,
with nonempty level sets.
For any smooth function
$h:[0,1]\rightarrow J$, define
$$
F(t,\theta):=H(\theta)-h(t).
$$
Then
$$
Z_t=H^{-1}(h(t))
$$
is a nonempty smooth hypersurface for every
$t\in[0,1]$, and \eqref{Fregular} holds.
\end{enumerate}
\end{example}
We no longer assume \eqref{Fregular}. 
For background on subanalytic and semianalytic sets and on the associated monomialization algorithm,
we refer the reader to \cite{bierstone1988semianalytic} and \cite{bierstone1997canonical}.

View $I=[0,a]$ as a compact
subanalytic arc in $\bb{T}^1=\bb{R}/\bb{Z}$. Assume instead that $F$ is a subanalytic function,
that every $Z_t$ is nonempty, and that $t\mapsto Z_t$ is continuous
in the Hausdorff topology.

The function
$(t,t')\mapsto d_H^{g_T}(Z_t,Z_{t'})$
is subanalytic with compact graph, and its zero set contains
the zero set of $(t,t')\mapsto d_{\bb{T}^1}(t,t')$.
The subanalytic L{}ojasiewicz inequality
\cite[Theorem 6.4]{bierstone1988semianalytic} gives an exponent $0<\beta\leq1$ and a constant $K_\mathcal{A}>0$ such that
$$
d_H^{g_T}(Z_t,Z_{t'})
\leq K_\mathcal{A} d_{\bb{T}^1}(t,t')^\beta
\leq  K_\mathcal{A} |t-t'|^\beta.
$$
The exponent $\beta$ and the constant $K_{\mathcal A}$
are effective and arise from the monomialization algorithm:
after simultaneous monomialization, on finitely many compact charts
with $|x_j|\le1$, we may write
$$
d_H(Z_t,Z_{t'})\circ\pi_i
=
|u_i(x)|\prod_j |x_j|^{a_{ij}},
\quad
d_{\bb{T}^1}(t,t')\circ\pi_i
=
|v_i(x)|\prod_j |x_j|^{b_{ij}},
$$
where $u_i,v_i$ are analytic units. In that case one have
$$
\beta
:=
\min\left(1,\min_{i,j:b_{ij}>0}\frac{a_{ij}}{b_{ij}}\right)>0,
\qquad
K_{\mathcal A}
:=
\max_i
\frac{\sup |u_i|}{(\inf |v_i|)^\beta}.
$$
See
\cite[Theorem 4.4, Lemma 5.3, and Theorem 6.4]
{bierstone1988semianalytic}, also \cite[Theorems 2.4 and 2.5]{bierstone1997canonical}.
\begin{proposition}
Put
$$
q:=n-\sup_{t\in I}\dim_{\mathrm{H}} Z_t,
$$
and assume that $q>0$ and $\beta(1+\alpha)>1$. Put
$$
\eta_{\beta,\alpha}
:=
\frac q2
\min\left\{
\frac1{2n},
\frac{2\beta(1+\alpha)-1}{2n+1}
\right\}>0.
$$
There exists a computable constant $B=B(C,F,\alpha,\beta,L/K)>0$
and a prime ideal $\pp\notin S_{M,L}$ such that
$$
F((N\pp)^{-\alpha},\theta_\pp)\neq0,
\qquad c_\pp\in C,
$$
with
$$
N\pp
\ll_n
\max\left\{
x_5,\,
B^{1/\eta_{\beta,\alpha}}
\bigl(1+\log(S_{L/K}(M)B)\bigr)^{2/\eta_{\beta,\alpha}}
\right\}.
$$
Here $x_5$ effectively computable as in \cref{smallest}.
\end{proposition}

\begin{proof}
For $\pp\notin S_{M,L}$, put
$$
A_\pp
:=
\{\theta\in \bb{T}^n:F((N\pp)^{-\alpha},\theta)\neq0\}.
$$
Since $q>0$, every $Z_t$ is closed and nowhere dense. Thus
$$
\partial A_\pp=T^n\setminus A_\pp=Z_{(N\pp)^{-\alpha}},
\qquad
\mu_{\mathrm{ST}}(A_\pp)=1.
$$
By \cref{analytictube}, these boundaries have uniformly bounded
$(n-q)$-dimensional upper Minkowski content.

Enumerate the prime ideals with $x/2<N\pp\leq x$
in nondecreasing order of norm. For $u<v$, H\"older's inequality gives
$$
\begin{aligned}
\sum_{k=u}^{v-1}d_*(\pp_k,\pp_{k+1})
&\ll_F
(v-u)^{1-\beta}
\left(
\sum_{k=u}^{v-1}
\left|(N\pp_k)^{-\alpha}
-(N\pp_{k+1})^{-\alpha}\right|
\right)^\beta
\\[3pt]
&\ll_{F,K,\alpha}
x^{1-\beta(1+\alpha)}.
\end{aligned}
$$
Hence $\mathcal A=(A_\pp)_{\pp\notin S_{M,L}}$ is a
$(\beta(1+\alpha)-1,q)$-Minkowski family.
Since $\mu_{\mathrm{ST}}(A_\pp)=1$ and $1\leq q\leq n$,
the conclusion follows from \cref{smallest}, as in the smooth case.
\end{proof}
    \begin{example}
The following real-analytic, Weyl-invariant functions satisfy
the preceding hypotheses.

\begin{enumerate}
\item For $n=2$, consider
$$
F(t,\theta)
:=
\prod_{\varepsilon=\pm1}
\left(2-2\cos(\theta_1+\varepsilon\theta_2)-t^2\right).
$$
Here
$
q=1$ and  $\beta=1$.
The zero set describes the chordal distance $t$ between eigenvalues
of two normalized Frobenius matrices.
\item For $n=3$, consider
$$
F(t,\theta)
:=
4\cos^2\theta_1+4\cos^2\theta_2+4\cos^2\theta_3-t.
$$
Here
$
q=1, \beta=\frac12.
$
This prescribes the sum of the squares of three normalized
Frobenius traces.
\end{enumerate}
\end{example}
\begin{remark}
\begin{enumerate}
\item
The map $t\mapsto Z_t$ need not be Hausdorff continuous on $I$.
Define
$
\mathcal Z:=\{(t,\theta)\in I\times\bb T^n:F(t,\theta)=0\},
$
and let
$
\pi:\mathcal Z\rightarrow I
$
be the first projection.

As $\pi$ is a proper subanalytic map, the stratification theorem
for subanalytic maps
\cite[Part I, Section 1.7, p.43]{GoreskyMacPherson1988}
gives Whitney stratifications
\[
\mathcal Z=\bigsqcup_\nu S_\nu,
\qquad
I=\bigsqcup_\eta T_\eta,
\]
such that, for every $\nu$, the restriction
$\pi|_{S_\nu}:S_\nu\rightarrow T_\eta$ is a submersion for some
$\eta$. Thom's first isotopy lemma
\cite[Proposition11.1]{Mather2012} then gives a local
stratified trivialization over each $T_\eta$. In particular,
$t\mapsto Z_t$ is Hausdorff continuous within each target stratum,
so any discontinuities lie in lower-dimensional strata of $I$.

\item
The case $0<\beta\leq 1/(1+\alpha)$, or the case in which global
Hausdorff continuity fails, remains unresolved in our approach.
A possible direction is to analyze the family stratum by stratum.
\end{enumerate}
\end{remark}

\section{Exceptional Picard-rank jump}
Let $X/K$ be a $K3$ surface over a number field $K$. We work within the category of absolute Hodge cycles $\mathcal{M}_{\mathrm{ahc}}(K,\Q)$ which satisfies the motivic hypotheses \ref{motiv} as stated in Discussion \ref{discussion0}. As in \cite[Proposition 7.2.3]{KMP2007} the refined Chow--Künneth decomposition is given by
$$h^2(X)(1)=h^2_{\mathrm{alg}}(X)(1)\oplus M,$$
where $M:=t_2(X)(1)$ is the transcendental motive of weight $0$. 

Following our detailed explanations of how to verify the CST hypothesis in
Discussions \ref{discussion2} and \ref{discussion3}, we now examine the CST
assumptions for the motive $M$.
The motive $h^2(X)$ belongs to the Tannakian subcategory
$\mathcal{M}_{\mathrm{Ab}}(K,\Q)$ of
$\mathcal{M}_{\mathrm{ahc}}(K,\Q)$.
Consequently, $M$ is an abelian motive, Serre's conjecture
holds for $M$, and the motivic Ramanujan conjecture holds
for the corresponding Artin-Sato-Tate $L$-functions.
Moreover, Tankeev's theorem together
with the Tannakian stability of the Mumford-Tate conjecture
implies that $M$ satisfies the Mumford-Tate conjecture,
and hence the algebraic Sato-Tate conjecture.

After replacing $K$ by a finite extension,
we assume that the Sato-Tate group of $M$ is
connected and its realizations are strongly
$G_M$-compatible.

If we retain the original field $K$, we instead
assume that $\st_K(M)$ is connected and that the
realizations of $M$ form a strongly $G_M$-compatible
system, as predicted by the strong compatibility
conjecture formulated by Kisin and Zhou. Let $G:=\st_K(M)$ be the motivic Sato-Tate group viewed as a smooth compact connected Riemannian manifold together with a fixed choice of bi-invariant Riemannian metric $g$.

Define $\rho(X_{\overline{K}}):=\mathrm{rank}(\mathrm{NS}(X_{\overline{K}}))$ to be the rank of the N\'eron-Severi group over $\overline{K}$. Let $\pp$ be a prime of good reduction with residue field $k$ and define $\rho(X_{\overline{k}}):=\mathrm{rank}(\mathrm{NS}(X_{\overline{k}}))$. The specialization map $\mathrm{NS}(X_{\overline{K}})\hookrightarrow \mathrm{NS}(X_{{\overline{k}}})$, together with $K3$-Tate conjecture (which is known to hold for finite fields) implies that 
$$\rho(X_{\overline{k}})=\rho(X_{\overline{K}})+\sum_{\zeta\in \mu_{\infty}}m_\zeta.$$ 
Here $\mu_{\infty}$ denotes the roots of unity and $m_\zeta=m{(\zeta|H_\ell(M))}$ the multiplicity of $\zeta\in \mu_{\infty}$ in the characteristic polynomial of $y_\pp$. Put $D=22-\rho(X_{\overline{K}})$.

Define
$$
Y_D:=\{\zeta\in\mu_\infty:
[\Q(\zeta):\Q]\le D\}.
$$
For a prime of good reduction $\pp$, recall that a Picard
jump occurs at $\pp$ if
$\rho(X_{\overline{k(\pp)}})>\rho(X_{\overline{K}})$.
Define the cyclotomic locus
$$
J:=\{g\in G:
\det(g-\zeta\mathrm{Id})=0
\text{ for some }\zeta\in Y_D\}.
$$
By the Tate conjecture for K3 surfaces over finite fields,
a Picard jump occurs at $\pp$ if and only if $y_\pp\in J$.
Motivated by Serre's conjecture on Picard jumps, we study
the distribution of Frobenius classes in tubular neighborhoods of the jump locus $J^{(\varepsilon)}$. Our estimates yield, in particular, a conditional upper bound for the number of primes at which Picard jumps occur.
Assume that $J\ne G$, and put
$n:=\operatorname{rank}(G)$.
Since $J$ is compact and analytic it has a finite
$(\dim G-1)$-dimensional upper Minkowski content,
although it need not be smooth.
Moreover, $J_T:=J\cap T$ is a finite union of translated
subtori of codimension one. In particular,
$\partial J_T=J_T$, and by Gray's tube formula recorded in \cref{gray} we have
$$
\operatorname{vol}_{g_T}\bigl(J_T^{(\varepsilon)}\bigr)
=
2\varepsilon\mathcal H_{g_T}^{n-1}(J_T)
+O_J(\varepsilon^2).
$$
More generally, \cref{analytictube} implies that there
exists $\varepsilon_0>0$ such that
$$
\bigl(\partial(J_T^{(\varepsilon)})\bigr)_
{0<\varepsilon<\varepsilon_0}
$$
has uniformly bounded $(n-1)$-dimensional upper Minkowski
content: these boundaries are contained in fibers of
the subanalytic family
$$
\{(\theta,y)\in\mathbb T^{n+1},
d_{g_T}(\theta,J_T)=|\cos y|\}.
$$

Write $J_T=\bigcup_i H_i$, with each
$H_i$ a translated subtorus. Distinct components are either disjoint or intersect transversely, so
the overlaps $H_i^{(\varepsilon)}\cap H_j^{(\varepsilon)}$ for distinct components have volumes
$O_{G,J}(\varepsilon^2)$. On each component $H_i$ the coarea formula implies
$$
\int_{H_i^{(\varepsilon)}}|D|^2\,d\operatorname{vol}_{g_T}
=
2\varepsilon\int_{H_i}|D|^2\,d\mathcal H_{g_T}^{n-1}
+O_{G,J}(\varepsilon^2).
$$
Weyl integration therefore gives
$$
\mu_{\st}\bigl(J^{(\varepsilon)}\bigr)
=
\frac{2\varepsilon}{|W|\s(T)}
\int_{J_T}|D(\theta)|^2\,
d\mathcal H_{g_T}^{n-1}(\theta)
+O_{G,J}(\varepsilon^2).
$$
Since $J$ is compact and subanalytic and as a result its
$(\dim G-1)$-dimensional Minkowski content equals
$\mathcal H_g^{\dim G-1}(J)$ (cf. Discussion in the beginning of \cref{content}), hence the equality
$$
\begin{aligned}
\frac{2}{|W|\s(T)}
\int_{J_T}|D(\theta)|^2\,
d\mathcal H_{g_T}^{n-1}(\theta)=
2\mathcal H_g^{\dim G-1}(J).
\end{aligned}
$$

Applying \cref{cst1} with $q=1$ gives the following:

\begin{proposition}\label{jump}
Let $S_0$ be the places of bad reduction.
Define
$$
\pi_{J,C}(x,\varepsilon)
:=
\#\{\pp\notin S_{M,L}\cup S_0:
N\pp\le x,\,
y_\pp\in J^{(\varepsilon)},\,
c_\pp\in C\},
$$
with $J^{(0)}:=J$.
There exists an effectively computable $x_0>0$,
independent of $\varepsilon$, such that, for $x>x_0$
and $0\le\varepsilon<\varepsilon_0$,
\begin{align*}
&\pi_{J,C}(x,\varepsilon)
=
2\frac{|C|}{|\Gamma|}\mathcal H_g^{\dim G-1}(J)\varepsilon\Li(x)
+O(|S_0|)+
O_{G,J}\left(
\frac{|C|}{|\Gamma|}\varepsilon^2\Li(x)
\right)
\\
&+
O_{G,J}\left(
n_K|C|^{1/2}\varepsilon x^{1/2}
\log\bigl(S_{L/K}(M)x\bigr)
\right)
\\
&+
O_{G,J}\left(
n_K
\frac{|C|^{1-1/(2\delta_G)}}
     {|\Gamma|^{1-1/\delta_G}}
x^{1-1/(2\delta_G)}
(\log x)^{2/\delta_G-1}
\bigl(\log(S_{L/K}(M)x)\bigr)^{1/\delta_G}
\right).
\end{align*}

In particular, $\pi_{J,C}(x,0)$ counts the Picard jumps
with $c_\pp\in C$.

If $\mathcal H_g^{\dim G-1}(J)>0$ and $\varepsilon=\varepsilon(x)>0$ satisfies
$$
\varepsilon(x)\longrightarrow0,
\qquad
\frac{
n_K(\log x)^{2/\delta_G}
\bigl(\log(S_{L/K}(M)x)\bigr)^{1/\delta_G}
}{
\varepsilon(x)x^{1/(2\delta_G)}
}
\longrightarrow 0,
$$
then
$$
\pi_{J,C}(x,\varepsilon(x))
\sim
2\frac{|C|}{|\Gamma|}
\mathcal H_g^{\dim G-1}(J)\varepsilon(x)\Li(x).
$$
\end{proposition}
\begin{remark}
\begin{enumerate}
Define $E=\mathrm{End}_{\mathcal{M}_{\mathrm{ahc}}(\bb C,\bb Q)}(M \times_K \C)$ to be the Hodge endomorphism field. Equivalently, $E$ is the endomorphism field of Hodge structures on the rational transcendental cohomology $T(X)_\Q$. By Zarhin's theorem, $E$ is either a totally real field or a CM field. 
\item
Our connectedness hypothesis, together with $J\neq G$,
excludes forced jumps. A nontrivial action of
$\operatorname{Gal}(K_0/K)$ can lead to Picard jumps of positive density. This was shown by \cite{CostaTschinkel2014} in the case of Kummer surface $X=\mathrm{Kum}(C\times C)$ for an elliptic curve $C$ with complex multiplication that jumps occur with density $1/2$. That density is computed explicitly in \cite[Theorem 1]{Huang2026}. If $E$ is totally real and $\dim_E T(X)_\Q$ is odd, every element of the connected Sato--Tate group has eigenvalue $1$, thus $J=G$ see \cite[Proposition 15]{Charles2014}.
An effective Chebotarev-Sato-Tate estimate might be a suitable application for future work on disconnected Sato-Tate groups.
\item By \cite[Definition 1.1]{Valloni2021} and the discussion immediately following it, $X$ has complex multiplication whenever $E$ has complex multiplication and $m=1$. Moreover, $$m:=\dim_ET(X)_\Q=\frac{22-\rho(X_{\overline{K}})}{[E:\Q]}.$$ Similar Chebotarev-Sato-Tate estimates hold for the locus
$$
J_E=\left\{g\in G:
\sum_{\zeta\in Y_D}m_\zeta(g)>[E:\Q]\right\},
$$
where the sum ranges over roots of unity and $m_\zeta(g)$ denotes
the dimension of the $\zeta$-eigenspace of $g$. This locus corresponds
to prime ideals $\pp$ such that
$$
\rho(X_{\overline{k(\pp)}})
>
\rho(X_{\overline{K}})+[E:\Q].
$$
Our work complements the result of \cite{Huang2026} in the number
field case. In the forced-jump case discussed above, it would be more interesting to study the locus $J_E\neq G$. If $J_E$ is nonempty, its torus inverse image is a compact
subanalytic set of positive codimension. Our Chebotarev-Sato-Tate estimates of \cref{equid} apply.
\end{enumerate}
\end{remark}
\begin{question}
     Tayou proved in \cite{Tayou2025} that there are infinitely many Picard jumps. It removes the assumption of everywhere potentially good reduction from the earlier work with Shankar--Shankar--Tang \cite{ShankarShankarTangTayou2022}.
\cref{cst3} implies that they are conditionally bounded: $$
\pi_{J,C}(x,0)
\ll_{X,K,J}x^{1-\frac1{2\delta_G}}
(\log x)^{\frac2{\delta_G}-1}.
$$ It would be interesting to obtain a bound for the first Picard jump or even an asymptotic for the number of Picard jumps following Lang-Trotter type methods in the connected Sato-Tate group case.
\end{question}
\section{Exceptional Frobenius spectra of Picard curves}
The distribution of exceptional split reductions is a natural
extension of the study of Frobenius conditions by effective
Chebotarev methods, see, in particular,
\cite[Theorem 12]{Serre1981}.
For Picard curves, a specific splitting condition admits a
description in terms of the three-dimensional
Frobenius factor. Let $\mathcal{J}$ be the Jacobian of the Picard curve $\mathscr C:\ y^3=x^4+x+1$ over the field $\Q(\zeta_3)$. It is an abelian threefold, a range in which the Mumford--Tate conjecture for abelian varieties is known. It has Sato-Tate group $\mathrm U(3)$ as in
\cite[Example 8.3]{FiteKedlayaSutherland2025}.
We now give a justification of the strong Mumford-Tate compatibility of the $\ell$-adic realizations.
We have an inclusion of the motivic Galois group $G_M\subseteq\mathrm{GU}(V,h)$. Since $\mathrm{ST}_K(M)=\mathrm U(3)$, $G_M$ is of Zariski dimension $10$. In particular over the splitting field $E$ we have $(G_M)_E=\GL_3\times \bb{G}_m.$
By \cite[Corollary 2.7]{Noot2013}, the associated rank-three
system with $\Q(\zeta_3)$-coefficients is strongly compatible. Under this identification, the first factor is the strongly
compatible rank-three system, while the second is the
Tate similitude character. Hence the full system is
strongly $G_M$-compatible.

Let $Q_\pp(T)=T^3-a_\pp T^2+b_\pp T-d_\pp \in \bb{Z}[\zeta_3][T]$ be the Frobenius factor.
Denote $S:=S_M\cup S_0$ where $S_0$ are the places of bad reduction of $\mathcal{J}$. The associated motive is the abelian motive $h^1(\mathcal{J})\simeq h^1(\mathcal C)$ which is a rank six polarized abelian motive of weight $1$. In particular, in the category $\mathcal{M}_{\mathrm{ahc}}(\Q(\zeta_3),\Q)$ it is a CST-motive of Sato-Tate group $U(3)$.
Define the counting function 
$$\pi_{\mathcal{J}}(x)=\#\{\pp \notin S, N\pp\leq x, \pp\,\text{of good reduction},\, N_{K/\Q}(a_\pp)=N\pp \}.$$
Given a prime $\pp$ of good reduction satisfying the conditions of the counting function, $\mathcal{J}_{\pp}$ is nonordinary and geometrically nonsimple. More precisely by Tate's isogeny theorem, $\mathcal{J}_{\pp}$ is isogeneous to a product $E\times B$ of an ordinary elliptic curve with CM by $\Q(\zeta_3)$ and a supersingular abelian surface.
Under the generalized Riemann hypothesis for the Dedekind zeta functions of the relevant
division fields, Serre's Theorem \cite[Théorème 12(ii), p. 157]{Serre1981}  gives the upper bound
$\pi_J(x)\ll x^{11/12}(\log x)^{-2/3}$.
\begin{proposition}
We have the following estimate:
$$
\pi_{\mathcal{J}}(x)\ll_{\mathcal{J},K}x^{9/10}(\log x)^{-3/5}.
$$
\end{proposition}

\begin{proof}
Let $y_{\mathfrak p}\in\mathrm U(3)$ be the normalized
Frobenius block, with eigenvalues $z_1,z_2,z_3$, and put
$$
D_\pp
:=
\sum_{i<j}|z_i-z_j|^2
=
9-|\operatorname{tr}y_\pp|^2.
$$
Thus
$$
\pi_{\mathcal{J}}(x)
=
\#\{\pp\notin S:N\pp\le x,\ D_\pp=8\}.
$$
The conjugacy invariant set $\{g\in G,|\mathrm{tr}(g)|^2=1\}$ has Sato-Tate measure $0$.
Since
$$
|z_1+z_2+z_3|^2-1
=
\frac{(z_1+z_2)(z_1+z_3)(z_2+z_3)}
{z_1z_2z_3},
$$
the inverse image of the Sato-Tate set in the maximal torus is
$$
D_T
=
\bigcup_{1\le i<j\le3}
\{\theta_i-\theta_j\equiv\pi\pmod{2\pi}\}.
$$
This is a finite union of smooth hypersurfaces and
therefore has finite $2$-dimensional upper Minkowski
content.

Moreover, $D_T$ is invariant under scalar multiplication
$U\mapsto e^{it}U$, hence under translation by the normal
subgroup $\mathrm U(1)\subset\mathrm U(3)$.
By \cref{imp}, the improved spectral exponent is
$$
\delta_N
=
\operatorname{rank}(\mathrm{PU}(3))
+
|\Phi^+(\mathrm{PU}(3))|
=
2+3=5.
$$
\cref{cst3} gives, for
sufficiently large $x$,
$$
\pi_{\mathcal J}(x)
\ll_{\mathcal J,K}
x^{9/10}(\log x)^{-3/5},
$$
which implies the stated bound.
\end{proof}
\begin{remark}
Although the improvement over Serre's original exponent is modest, namely $1/60\approx 0.0167$, our method highlights the role of geometric symmetries of the underlying Sato-Tate set in determining the exponent of the error term. In the present work, reductions in the rank parameter arise from invariance under compact subgroup actions. It would be interesting to investigate if more general symmetries, encoded by more general and possibly singular Riemannian foliations on $\st_K(M)$, could yield further improvements.
\end{remark}

\clearpage
\appendix
\crefalias{section}{appendix}
\chapter{Riemannian geometric tools}
In this Appendix we introduce selected tools from Riemannian geometry that are relevant to our Beurling-Selberg machinery.
Our main reference for Riemannian geometry is Lee's book
\cite{Lee2018}. Since $(\M,g)$ is assumed compact and connected,
the Hopf-Rinow theorem
\cite[Theorem 6.19 and Corollary 6.22]{Lee2018}
implies that it is geodesically complete. Consequently,
the exponential map
$\exp\colon T\M\to\M$, defined in \cite[p. 128]{Lee2018},
is globally defined and smooth
\cite[Proposition 5.19]{Lee2018}.
Let 
$P\subset \M$ be a compact embedded smooth
submanifold. We denote by
$NP$
its normal bundle. Define the tubular neighborhood by $$
\mathcal N_\varepsilon P
:=
\{(p,v)\in NP:|v|_g<\varepsilon\}.
$$
\begin{theorem}\cite[Theorems 5.25 and 6.40]{Lee2018}\label{exp}
There exists $\varepsilon_0=\varepsilon(P)>0$
such that, for every $0<\varepsilon<\varepsilon_0$,
the normal exponential map induces a diffeomorphism
$$
\exp_P:
\mathcal N_\varepsilon P
\rightarrow
P^{(\varepsilon)},
\qquad
(p,v)\longmapsto\exp_p(v).
$$
Moreover, we have
$
d_g(\exp_p(v),P)=|v|_g.
$
\end{theorem}
For smooth tensor fields $S,T$ and a vector field $X$,
$$
\nabla_X(S\otimes T)
=
(\nabla_X S)\otimes T
+
S\otimes(\nabla_X T).
$$
Iterating this identity gives the covariant Leibniz formula:

\begin{lemma}[Covariant Leibniz]\label{lbt}
    Let $U\subset \M$ be an open set and $x\in U$. Let $k\geq 1$ and $F,G\in C^\infty(U)$ and $X_1,\ldots,X_k\in T_x\M$,
$$
\bigl(\nabla^k(FG)\bigr)_x(X_1,\ldots,X_k)
=
\sum_{I\subseteq\{1,\ldots,k\}}
\Bigl(
(\nabla^{|I|}F)_x
\otimes
(\nabla^{k-|I|}G)_x
\Bigr)(X_I,X_{I^c}),
$$
where $X_I=(X_{i_1},\ldots,X_{i_j})$ for
$I=\{i_1<\cdots<i_j\}$, and the empty-index terms
are interpreted as $F(x)$ and $G(x)$.
\end{lemma}
A similar result holds for the covariant derivative of the composite of two functions:
\begin{lemma}[Covariant Faà-Di-Bruno]\label{fdb}
Let $U\subset \M$ be open and $x\in U$. Let $F\in C^\infty(U)$, and
$\psi\in C^\infty(\mathbb R)$ be real valued smooth functions. For $k\geq1$ let $\Pi_k$
be the set of partitions of $\{1,\ldots,k\}$. Given $x\in\M$ and $X_1,\ldots,X_k\in T_x\M$, we have
\begin{align*}
&(\nabla^k(\psi\circ F))_x(X_1,\ldots,X_k)
\\=&
\sum_{j=1}^{k}
\sum_{\substack{\pi=\{B_1,\ldots,B_j\}\in\Pi_k\\|\pi|=j}}
\psi^{(j)}(F(x))
\left(\bigotimes_{a=1}^{j}(\nabla^{|B_a|}F)_x\right)
(X_{B_1},\ldots,X_{B_j}), 
\end{align*}
where
$X_B=(X_{i_1},\ldots,X_{i_m})$ for $B=\{i_1<\cdots<i_m\}$.
\end{lemma}
\begin{proof}
    The result follows by the same induction as the classical multivariate
Faà di Bruno formula \cite[Theorem 2.1]{ConstantineSavits1996}, using
the covariant Leibniz rule of \cref{lbt}.
\end{proof}
For a covariant $k$-tensor $A$, define its divergence by
$$
(\operatorname{div}_g A)(X_2,\ldots,X_k)
=
\sum_{a=1}^{d}
(\nabla_{e_a}A)(e_a,X_2,\ldots,X_k),
$$
where $(e_a)_a$ is a local orthonormal frame.
\begin{lemma}\label{div}
Define the operator $T_{n,g}:C^\infty(\M)
\longrightarrow\Gamma^\infty((T^*\M)^{\otimes n})$ by the formulas:
\begin{align*}
T_{2m,g}u
&:=(-1)^m(\Delta_g^m u)\,g^{\otimes m},\\
T_{2m+1,g}u
&:=(-1)^{m+1}d(\Delta_g^m u)\otimes g^{\otimes m},
\end{align*}
with $g^{\otimes0}=1$. Then, for every $n\geq1$ and every
$u\in C^\infty(\M)$,
$$
\operatorname{div}_g^{\,n}(T_{n,g}u)=\Delta_g^n u.
$$
\end{lemma}

\begin{proof}
Since $\nabla g=0$ we have for every smooth
function $f$ and $m\geq 1$,
$$
\operatorname{div}_g(fg^{\otimes m})
=
df\otimes g^{\otimes(m-1)}
,
$$
and for $m\geq 0$ we have
$$
\operatorname{div}_g(df\otimes g^{\otimes m})
=
\operatorname{tr}_g(\nabla^2f)\,g^{\otimes m}
=
-(\Delta_gf)\,g^{\otimes m}.
$$
Consequently, for $m\geq1$,
$$
\begin{aligned}
\operatorname{div}_g(T_{2m,g}u)
&=(-1)^m d(\Delta_g^m u)
       \otimes g^{\otimes(m-1)}\\
&=T_{2m-1,g}(\Delta_g u),
\end{aligned}
$$
whereas, for $m\geq0$,
$$
\begin{aligned}
\operatorname{div}_g(T_{2m+1,g}u)
&=(-1)^m(\Delta_g^{m+1}u)\,g^{\otimes m}\\
&=T_{2m,g}(\Delta_g u).
\end{aligned}
$$
Thus
$$
\operatorname{div}_g(T_{n,g}u)
=
T_{n-1,g}(\Delta_g u)
\qquad(n\geq1).
$$
Iterating this identity and using $T_{0,g}f=f$ yields
$$
\operatorname{div}_g^{\,n}(T_{n,g}u)
=
T_{0,g}(\Delta_g^n u)
=
\Delta_g^n u.
$$
\end{proof}
Put $K_{\mathrm{Ric}}:=(d-1)\kappa$, with
$K_{\mathrm{Ric}}=0$ when $d=1$.
\begin{lemma}\label{decay}
For every integer $n\ge1$, $f\in C^\infty(\M)$, and $j\ge1$,
$$
|\hat{f}(j)|
 \le C_{n,g}\lambda_j^{-n/2}
       \|\nabla^n f\|_{L^1(M)}\|\varphi_j\|_{L^\infty(\M)},
$$
where
$$
 C_{n,g}={d^{\lfloor n/2\rfloor/2}}
 \begin{cases}
  1, & n\ \text{even},\\
  \displaystyle\sqrt{\frac{2e}{\pi}
                     \left(1+\frac{K_{\mathrm{Ric}}}{\lambda_1}\right)},
       & n\ \text{odd}.
 \end{cases}
$$
\end{lemma}

\begin{proof}
Let $T_{n,g}:C^\infty(\M)
\longrightarrow\Gamma^\infty((T^*\M)^{\otimes n})$ be the operator of \cref{div} so that $\Delta_g^n=\mathrm{div}^n\circ T_{n,g}$.
Fix $j\geq 1$.

By \cite[Corollary 1.3]{ArnaudonThalmaierWang2020},
$$
 \|\nabla\varphi_j\|_\inft
 \le \sqrt{\frac{2e}{\pi}(\lambda_j+K_{\mathrm{Ric}})}\,
       \|\varphi_j\|_\inft.
$$
Using $|g^{\otimes m}|_g=d^{m/2}$ and $\lambda_j\ge\lambda_1$ gives
$\|T_{n,g}\varphi_j\|_\inft
 \le C_{n,g}\lambda_j^{n/2}\|\varphi_j\|_\inft$.
For smooth $f\in C^\infty(\M)$, integration by parts yields
$$
 \lambda_j^n\int_{\M} f\overline{\varphi_j}\,d\s
 =(-1)^n\int_{\M}
   \langle\nabla^n f,\overline{T_{n,g}\varphi_j}\rangle_g\,d\s.
$$
Taking the absolute values gives 
\begin{align*}
  \left|\int_{\M} f\overline{\varphi_j}\,d\s\right|&\leq \frac{1}{\lambda_j^n}\|\nabla^nf\|_{L^1(\M)}\|T_{n,g}\varphi_j\|_{L^\infty(\M)}\\&\leq C_{n,g}\frac{1}{\lambda_j^{n/2}}\|\nabla^nf\|_{L^1(\M)}\|\varphi_j\|_{L^\infty(\M)}.  
\end{align*}
\end{proof}
\chapter{Review of Fraenkel's regularized distance}\label[appendix]{reg}
The notion of regularized distance goes back to Stein \cite{Stein1970}. Given a non-empty closed set $F\subset \R^d,$ Stein constructs a function $\rho_F\in C^\infty(\R^d\setminus F)$ such that $\rho_F(x)\asymp d_{\R^d}(x,F)$ and $|D^\alpha\rho_F(x)|\ll d_{\R^d}(x,F)^{1-|\alpha|}$ for every multi-index $\alpha$. Fix a closed subset $F\subset \bb{R}^d$. Fraenkel gave a different construction that allows
$\rho_F$ to be chosen arbitrarily close to
$d_{\mathbb R^d}(\,\cdot\,,F)$. Applying Fraenkel's smoothing construction \cite[Theorem 1]{Fraenkel1979} to the Lipschitz function $x\rightarrow d_{\bb{R}^d}(x,F)$ yields the following:
\begin{theorem}\label{rgd}
Let $F\subset\mathbb{R}^{d}$ be nonempty and closed. Define the mollifier
$$
J(y)=
\begin{cases}
c_{0}\exp\left(-\dfrac{1}{1-|y|^{2}}\right),& |y|<1,\\
0,& |y|\geq1,
\end{cases}
$$
where $c_0$ is chosen so that
$\int_{\mathbb R^d}J(y)\,dy=1$.
Writing $J(y)=J_{*}(|y|)$, set
$$
J_{1}(y):=(d-1)J_{*}(|y|)+|y|J_{*}'(|y|).
$$
For every $\xi\in(0,1)$ define the function
$\rho_{F,\xi}\in C^{\infty}(\mathbb{R}^{d}\setminus F)$ by the formula
$$
\rho_{F,\xi}(x)
:=
\xi^{-d}
\int_{\mathbb{R}^d\setminus F}
d_{\R^d}(y,F)^{1-d}
J\left(\frac{y-x}{\xi d_{\R^d}(y,F
)}\right)dy,
\qquad x\in\mathbb{R}^d\setminus F.
$$
Then we have
$$
(1+\xi)^{-2}d_{\e}(x,F)
\leq \rho_{F,\xi}(x)
\leq(1-\xi)^{-2}d_{\e}(x,F),
$$
and for every multi-index $|\alpha|\geq1$ we have
\begin{equation}\label{der}
|D^\alpha\rho_{F,\xi}(x)|
\leq P_\alpha\bigl(\xi d_{\bb{R}^d}(x,F)\bigr)^{1-|\alpha|},    
\end{equation}
where if $\alpha=\beta+\gamma$ with $|\gamma|=1$ one may take
$$
P_\alpha
=
\int_{B(0,1)}
|D^\beta J_1(y)|
\frac{(1+|y|)^{|\beta|}}{1-|y|}\,dy .
$$
\end{theorem}
We state and prove a Hausdorff stability lemma for Fraenkel's distance which will be a key ingredient in our Beurling-Selberg theory.
\begin{lemma}\label{hsd}
    Let $A,B\subset\mathbb R^d$ be nonempty closed sets such that
$d^{\e}_H(A,B)<\infty$. Fix $\xi\in(0,1)$, and let
$\rho_{A,\xi}$ and $\rho_{B,\xi}$ be the Fraenkel regularized
distances constructed using the same mollifier $J$. Given a multi-index $\alpha$ define $k=|\alpha|$.
    Define the constants
    \begin{align*}
        M_{\alpha,J}:&=\sup_{|u|\leq 1}|(d+k-1)D^\alpha J(u)+u\cdot\nabla D^\alpha J(u)|,\\
        C_{\alpha,d,J,\xi}:&=\max\left \{4\xi^{1-k}P_\alpha,2^k3^d\omega_d\frac{(1+\xi)^{d+k}}{\xi^k(1-\xi)^d}M_{\alpha,J},\frac{3}{(1-\xi)^2}\right\},
    \end{align*}
with the convention that $P_0=0$.
    Then for $x\in \e\setminus (A\cup B)$ we have
    $$|D^\alpha\rho_{A,\xi}(x)-D^{\alpha}\rho_{B,\xi}(x)|\leq C_{\alpha,d,J,\xi}\frac{d_H(A,B)}{\min(d_{\e}(x,A),d_{\e}(x,B))^k}.$$
    \end{lemma}
\begin{proof}
Let $x\notin A\cup B$.
    Define $r=\min(d_{\R^d}(x,A),d_{\R^d}(x,B))$ and $h=d_H(A,B)$. Recall that for $F\in \{A,B\}$ and $\xi\in (0,1)$ we have $$
\rho_{F,\xi}(x)
=
\xi^{-d}
\int_{\mathbb{R}^d\setminus F}
d_{\R^d}(y,F)^{1-d}
J\left(\frac{y-x}{\xi d_{\R^d}(y,F
)}\right)dy.
$$
Differentiating under the integral sign with respect to $x$ gives
$$
D^\alpha\rho_{F,\xi}(x)
=
(-1)^k\xi^{-d-k}
\int_{\mathbb{R}^d\setminus F}
d_{\R^d}(y,F)^{1-d-k}
D^\alpha J\left(\frac{y-x}{\xi d_{\R^d}(y,F
)}\right)dy.
$$
Define $H_\alpha:\mathbb R^d\times[0,\infty)\to\mathbb R$ by
$$
H_\alpha(w,t):=
\begin{cases}
(-1)^k\xi^{-d-k}t^{1-d-k}
D^\alpha J\left(\dfrac{w}{\xi t}\right),&t>0,\\
0,&t=0.
\end{cases}
$$
For $F\in\{A,B\}$, we can therefore write
$$
D^\alpha\rho_{F,\xi}(x)
=\int_{\mathbb R^d}
  H_\alpha(y-x,d_{\mathbb R^d}(y,F))\,dy.
$$
We have $|\partial_t H_\alpha(w,t)|\leq \xi^{-d-k}t^{-d-k}M_{\alpha,J}\ind_{|w|\leq \xi t}.$

Assume that $h<r/2$. Let $z\in \R^d$. Let $t$
lie between $d_{\mathbb{R}^d}(z,A)$ and
$d_{\mathbb{R}^d}(z,B)$. Choose $F\in\{A,B\}$ such that
$$
d_{\mathbb{R}^d}(z,F)
=
\max\left\{
d_{\mathbb{R}^d}(z,A),
d_{\mathbb{R}^d}(z,B)
\right\}.
$$
We have
$$
0
\le
d_{\mathbb{R}^d}(z,F)-t
\le
\left|
d_{\mathbb{R}^d}(z,A)-d_{\mathbb{R}^d}(z,B)
\right|
\le
h.
$$
Since $d_{\mathbb R^d}(x,F)\leq r+h$, whenever
$|z-x|\leq\xi t$ we have
$$
|z-x|
\leq\xi t
\leq\xi d_{\mathbb R^d}(z,F)
\leq\xi\bigl(r+h+|z-x|\bigr).
$$
Consequently,
$$
|z-x|\leq\frac{\xi(r+h)}{1-\xi}
<R_\xi:=\frac{3\xi r}{2(1-\xi)}.
$$
Thus, for every $t$ between $d_{\mathbb R^d}(z,A)$
and $d_{\mathbb R^d}(z,B)$,
$\partial_tH_\alpha(z-x,t)$ vanishes when
$z\notin B_{\mathbb R^d}(x,R_\xi)$.
Moreover, whenever $|z-x|\leq\xi t$, the preceding bound
$d_{\mathbb R^d}(z,F)-t\leq h$ gives
$$
r
\leq d_{\mathbb R^d}(x,F)
\leq d_{\mathbb R^d}(z,F)+|z-x|
\leq t+h+\xi t.
$$
Hence
$$
t\geq\frac{r-h}{1+\xi}
>\frac{r}{2(1+\xi)}.
$$
Using the mean-value theorem, we conclude that
$$
\begin{aligned}
&\left|
D^\alpha\rho_{A,\xi}(x)
-
D^\alpha\rho_{B,\xi}(x)
\right|\\
&\quad \leq h\sup_{z\in B_{\e}(x,R_\xi),s\geq \frac{r}{2(1+\xi)}}|\partial_tH_{\alpha}(z-x,s)|\mathrm{vol}_{\e}(B_{\e}(x,R_\xi))
\\
&\quad\leq
h\,
\omega_d
\left(\frac{3r\xi}{2(1-\xi)}\right)^d
\xi^{-d-k}
\left(\frac{2(1+\xi)}{r}\right)^{d+k}
M_{\alpha,J}                                                \\
&\quad=
2^k3^d\omega_d
\frac{(1+\xi)^{d+k}}
     {\xi^k(1-\xi)^d}
M_{\alpha,J}
\frac{h}{r^k}.
\end{aligned}
$$
Suppose now that $h\geq r/2$ and $k\geq 1$. Fraenkel's derivative
estimate gives, for $F\in\{A,B\}$,
$$
|D^\alpha\rho_{F,\xi}(x)|
\leq
P_\alpha\bigl(\xi d_{\e}(x,F)\bigr)^{1-k}.\, 
$$
Since $d_{\R^d}(x,F)\geq r$, it follows that
$$
\begin{aligned}
\left|
D^\alpha\rho_{A,\xi}(x)
-
D^\alpha\rho_{B,\xi}(x)
\right|
&\leq
2P_\alpha\xi^{1-k}r^{1-k} \\
&\leq
4P_\alpha\xi^{1-k}
\frac{h}{r^k}.
\end{aligned}
$$
Finally, suppose that $h\geq r/2$ and $k=0$. We have
\begin{align*}
    |\rho_{A,\xi}(x)-\rho_{B,\xi}(x)|&\leq \frac{1}{(1-\xi)^2} \max\{d_{\mathbb R^d}(x,A),d_{\mathbb R^d}(x,B)\}\\
    &\leq \frac{1}{(1-\xi)^2}(h+r)\\
    &\leq \frac{3}{(1-\xi)^2}h
\end{align*}
Combining the two cases proves the statement.
\end{proof}
\begin{lemma}\label{Stein}
     Suppose
$
0<r<
\min\left\{
D_{\mathbb M},\frac{r_{\mathrm{inj}}}{2},
\frac{1}{\sqrt{2K}}
\right\},
$
where $1/\sqrt{2K}:=\infty$ if $K=0$. Let $p\in M$.
Fix a linear isometry $B:\mathbb R^d\to T_p\M$, and set
$$
\psi:=\exp_p\circ B:
B_{\mathbb R^d}(0,2r)\longrightarrow B_g(p,2r).
$$
For every non-empty closed set $F\subsetneq \M$, define the closed set
$$
E(F):=\psi^{-1}\bigl(F\cap B_g(p,2r)\bigr)
       \cup\bigl(\mathbb R^d\setminus B_{\mathbb R^d}(0,2r)\bigr).
$$
Set
$
\kappa^-(r):=\frac{1}{1+2Kr^2},$ and $
\kappa^+(r):=\frac{1}{1-2Kr^2}.
$
Then:
\begin{enumerate}
\item For every $x\in B_g(p,r)$,
$$
\frac{r}{D_{\mathbb M}}\kappa^-(r)d_g(x,F)
\leq d_{\e}\bigl(\psi^{-1}(x),E(F)\bigr)
\leq \kappa^+(r)d_g(x,F).
$$
If moreover $x\in F^{(r)}$, then
$$
d_{\e}\bigl(\psi^{-1}(x),E(F)\bigr)
\geq \kappa^-(r)d_g(x,F).
$$

\item For every non-empty closed set $F'\subsetneq M$,
$$
d^{\e}_H\bigl(E(F),E(F')\bigr)
\leq \kappa^+(r)d_H^g(F,F').
$$
\end{enumerate}
\end{lemma}
\begin{proof}
Suppose that $
0<r<
\min\left\{
D_{\mathbb M},\frac{r_{\mathrm{inj}}}{2},
\frac{1}{\sqrt{2K}}
\right\}
$.
By \cite[Lemma 5]{DyerVegterWintraecken2015}, for
$|v|<2r$ and $w\in \e$ we have
$$
(1-2Kr^2)|w|
\leq |(d\psi)_v(w)|_g
\leq (1+2Kr^2)|w|.
$$
Since $2r<r_{\mathrm{inj}}$, the map $\psi$ is a
diffeomorphism, and for $x\in B_g(p,2r)$ we have
$
|\psi^{-1}(x)|=d_g(p,x)
$.
Fix $x\in B_g(p,r)$. Choose $q\in F$ such that
$d_g(x,q)=d_g(x,F)$. Let
$\gamma$ be a minimizing geodesic from $x$ to $q$. 
If $\gamma\subset B_g(p,2r)$ then
$$
d_{\e}(\psi^{-1}(x),E(F))
\leq |\psi^{-1}(x)-\psi^{-1}(q)|
\leq \kappa^+(r)d_g(x,q).
$$
Otherwise, let $z$ be the first point where $\gamma$ meets
$\partial B_g(p,2r)$. Then
\begin{align*}
    d_{\e}(\psi^{-1}(x),E(F))&\leq d_{\R^d}(\psi^{-1}(x),\R^d\setminus B_{\R^d}(0,2r))\\
&\leq 2r-|\psi^{-1}(x)|\\&
=2r-d_g(p,x)\\&
\leq d_g(x,z)\\&
\leq d_g(x,q).
\end{align*}

This proves the upper bound.

For the lower bound, given
$w\in\psi^{-1}(F\cap B_g(p,2r))$ we have
$$
d_g(x,\psi(w))
\leq (1+2Kr^2)|\psi^{-1}(x)-w|.
$$
Consequently,
$$
\begin{aligned}
d_{\e}(\psi^{-1}(x),E(F))
&\geq
\min\left\{\kappa^-(r)d_g(x,F),\,2r-|\psi^{-1}(x)|\right\}\\
&\geq
\min\left\{\kappa^-(r)d_g(x,F),\,r\right\}\\
&\geq
\frac{r}{D_{\mathbb M}}\kappa^-(r)d_g(x,F).
\end{aligned}
$$
The last step uses $d_g(x,F)\leq D_{\mathbb M}$,
$r<D_{\mathbb M}$, and $\kappa^-(r)\leq1$.
If $x\in F^{(r)}$, then $d_g(x,F)<r$, so the same
minimum is at least $\kappa^-(r)d_g(x,F)$.

For the second assertion, put $h:=d_H^g(F,F')$ and fix
$u\in E(F)$. If $u\notin B_{\mathbb R^d}(0,2r)$, then
$u\in E(F')$, so $d_{\e}(u,E(F'))=0$.

Otherwise, $x:=\psi(u)\in F\cap B_g(p,2r)$.
Since $F'$ is compact, choose $y\in F'$ such that
$$
d_g(x,y)=d_g(x,F')\leq h.
$$
Let $\gamma$ be a minimizing geodesic from $x$ to $y$.
If $\gamma\subset B_g(p,2r)$, then
$$
d_{\e}(u,E(F'))
\leq |u-\psi^{-1}(y)|
\leq \kappa^+(r)d_g(x,y)
\leq \kappa^+(r)h.
$$
Otherwise, if $z$ is the first point where $\gamma$
meets $\partial B_g(p,2r)$, then
$$
d_{\e}(u,E(F'))
\leq 2r-|u|
\leq d_g(x,z)
\leq d_g(x,y)
\leq h
\leq \kappa^+(r)h.
$$
Taking the supremum over $u\in E(F)$, and then
interchanging $F$ and $F'$, proves the assertion.
\end{proof}

\chapter{Computation of Minkowski contents}\label[appendix]{content}
For a compact connected smooth Riemannian manifold
$(\M,g)$ of dimension $d$, one has 
$$
\mathcal{H}^q_g(E)
\ll_{\M}
\underline{\mathcal{M}}^q_g(E);
$$
see \cite[Section 5.5, p. 80]{Mattila1995} for the Euclidean bound. The Riemannian variant is proved by covering $E$ by finitely many bi-Lipschitz charts. We record several classes of Minkowski-measurable sets
whose Minkowski content agrees with the associated
Hausdorff measure.

\begin{enumerate}
\item Consider a compact smooth submanifold $P\subset M$ of dimension $q<d$. By Gray's tube formula \cite[Theorem 9.23]{GrayTubes} one gets
\begin{equation}\label{gray}
    \s\bigl(P^{(\r)}\bigr)
=
\omega_{d-q}\mathcal H_g^q(P)\r^{d-q}
+O_{P,g}\left(\r^{d-q+2}\right).
\end{equation}
Thus, $P$ is $q$-dimensional Minkowski measurable and
$$
\mathcal M_g^q(P)=\mathcal H_g^q(P).
$$

\item Suppose that $A\subset\M$ is compact and $q$-rectifiable in
Federer's sense: there exist a bounded set $B\subset\mathbb R^q$
and a Lipschitz map $f:B\to\M$ such that $A=f(B)$.
Federer's theorem, localized in Riemannian charts
\cite[3.2.39, 3.2.46]{Federer1969}, implies, with the standard
normalizations, that
$$
    \mathcal M_g^q(A)=\mathcal H_g^q(A)<\infty.
$$
Consequently, $A$ is $q$-dimensional Minkowski measurable.

\item Assume that $(\M,g)$ is a real-analytic Riemannian manifold and that $A\subset \M$ is compact and subanalytic. Let $m=\dim_H A$ be the Hausdorff dimension. By the analytic uniformization theorem \cite[Theorem 0.1]{bierstone1988semianalytic} there exist a compact real-analytic manifold $Y$ of dimension $m$ and a proper real-analytic map
$$
    \pi:Y\rightarrow M
$$
such that $\pi(Y)=A$. Since $Y$ is compact, $\pi$ is Lipschitz, hence $A$ is $m$-rectifiable. So
$$
    \mathcal M_g^m(A)=\mathcal H_g^m(A).
$$

\item Let $A\subset \M$ be a compact set of positive reach, and assume that $\s(A)=0$. According to the tube formula of Solanes-Trillo \cite[Corollary 4.7]{SolanesTrillo2025}, the set $A$ is Minkowski measurable and
$$
    \mathcal M_g^{d-1}(A)=\mathcal H_g^{d-1}(A).
$$

\end{enumerate}
As $\s(\M)=1$ as per our conventions, if $\overline{\mathcal M}_g^q(A)<\infty$ then $$\sup_{0<r<r_0}\frac{\s(A^{(r)})}{r^{d-q}}<\infty,$$ for all $r_0>0$.

Assume that $(\M,g)$ is a compact real-analytic Riemannian
manifold of dimension $d$, and let $(N,h)$ be another
real-analytic Riemannian manifold, possibly with boundary.
Let $A\subset\M\times N$ be a subanalytic set, and assume
that, for every $y\in N$, the fiber
$$
A_y:=\{x\in\M:(x,y)\in A\}
$$
is a possibly empty proper compact subset of $\M$.
Define the family of sets $\mathcal A:=(A_y)_{y\in N}$. 

\begin{proposition}\label{analytictube}
Suppose moreover that $A$ is compact and nonempty, and put
$$
q:=\sup_{y\in N}\dim_{\mathrm H}A_y.
$$
If $q<d$, then $$
0<\overline{\mathcal M}_g^{q}(\mathcal A)<\infty.
$$
\end{proposition}

\begin{proof}
Let
$$
\iota_\M:\M\hookrightarrow\mathbb R^{n_1},
\qquad
\iota_N:N\hookrightarrow\mathbb R^{n_2}
$$
be analytic embeddings, and set
$$
\iota:=(\iota_\M,\iota_N):
\M\times N\hookrightarrow\mathbb R^{n_1+n_2}.
$$
Then
$$
\widetilde A:=\iota(A)
$$
is a compact subanalytic subset of
$\mathbb R^{n_1+n_2}$, and hence globally subanalytic;
see \cite[Definition 1.1.6]{Valette2025Subanalytic}.
Likewise, $\iota_\M(\M)$ is globally subanalytic.
Let
$$
T:=\left\{
t\in\mathbb R^{n_2}:
\widetilde A_t\ne\varnothing
\right\}.
$$
Define
$$
\widehat A
:=
\widetilde A
\cup
\left(
\{0\}\times(\mathbb R^{n_2}\setminus T)
\right).
$$
Then $\widehat A$ is globally subanalytic and
$$
\widehat A_t
=
\begin{cases}
\widetilde A_t,&t\in T,\\
\{0\},&t\notin T.
\end{cases}
$$
Furthermore, the fibers of $\widehat A$ have uniformly
bounded diameters, and
$$
\sup_{t\in\mathbb R^{n_2}}
\dim_{\mathrm H}\widehat A_t=q.
$$
By \cite[Proposition 4.3.4, (4.3.3)]{Valette2025Subanalytic},
applied with the constant second family
$\iota_\M(\M)$, there exist $C>0$ and $r_0>0$ such that
$$
\sup_{t\in\mathbb R^{n_2}}
\mathcal H_{\mathbb R^{n_1}}^d
\left(
(\widehat A_t)^{(r)}
\cap\iota_\M(\M)
\right)
\le Cr^{d-q},
\qquad 0<r<r_0.
$$

Since $\M$ is compact and $\iota_\M$ is a
$C^1$-embedding, the Riemannian metric $g$ and the metric
induced on $\iota_\M(\M)$ by the Euclidean metric are
uniformly equivalent. Consequently, there exist
$C_1,C_2>0$ such that, for every $y\in N$ and every
sufficiently small $r>0$,
$$
\iota_\M\bigl((A_y)^{(r)}\bigr)
\subseteq
\bigl(\widetilde A_{\iota_N(y)}\bigr)^{(C_1r)}
\cap\iota_\M(\M)
$$
and
$$
\s\bigl((A_y)^{(r)}\bigr)
\le
C_2\,
\mathcal H_{\mathbb R^{n_1}}^d
\left(
\bigl(\widetilde A_{\iota_N(y)}\bigr)^{(C_1r)}
\cap\iota_\M(\M)
\right).
$$
Since $\widetilde A_t\subseteq\widehat A_t$, it follows
that
$$
\begin{aligned}
\sup_{y\in N}\s\bigl((A_y)^{(r)}\bigr)
&\le
C_2
\sup_{t\in\mathbb R^{n_2}}
\mathcal H_{\mathbb R^{n_1}}^d
\left(
(\widehat A_t)^{(C_1r)}
\cap\iota_\M(\M)
\right)
\\
&\le C_2CC_1^{d-q}r^{d-q}.
\end{aligned}
$$
In other words,
$$
\sup_{y\in N}\s\bigl((A_y)^{(r)}\bigr)
\ll_{A,\M,N}r^{d-q}.
$$

Finally, since the fiber dimensions are integers, there
exists $y_0\in N$ with $\dim_{\mathrm H}A_{y_0}=q$.
By Whitney stratification, $A_{y_0}$ contains a
$q$-dimensional analytic stratum. Choosing a relatively
compact smooth stratum with positive
$\mathcal H_g^q$-measure gives
$$
\limsup_{r\to0^+}
\frac{\s\bigl((A_{y_0})^{(r)}\bigr)}
{r^{d-q}}>0.
$$
Therefore,
$$
0<\overline{\mathcal M}_g^{q}(\mathcal A)<\infty.
$$
This completes the proof.
\end{proof}

Let $m\geq1$, let $(K,h)$ be a compact Riemannian manifold,
possibly with boundary, and let
$$
F=(F_1,\ldots,F_m):\M\times K\longrightarrow\mathbb R^m
$$
be smooth. For $t\in K$, define
$$
A_t:=\{x\in\M:F_i(x,t)=0,\ 1\leq i\leq m\},
$$
and
$$
A=\{(x,t)\in \M\times K:x\in A_t\}
=\bigcup_{t\in K}A_t\times\{t\}.
$$
We assume that $A_t\neq\varnothing$ for all $t\in K$, and that for all $(x,t)\in A$, the vectors
$$
\nabla_xF_1(x,t),\ldots,\nabla_xF_m(x,t)
$$
are linearly independent in $T_x\M$.
Ehresmann's fibration Theorem implies that 
$$
0<\inf_{t\in K}\mathcal H_g^{d-m}(A_t)
\leq\sup_{t\in K}\mathcal H_g^{d-m}(A_t)<\infty.
$$
We shall prove that the family $(A_t)_{t\in K}$ has finite uniform Minkowski content. 
\begin{proposition}\label{tube3}
There exists $r_0>0$, independent of
$t\in K$, such that, uniformly in $t$,
$$
    \s\bigl(A_t^{(r)}\bigr)
    =
    \omega_mr^m\mathcal H_g^{d-m}(A_t)
    +O(r^{m+2}),
    \qquad 0<r<r_0.
$$
In particular the family $\mathcal{A}=(A_t)_{t\in K}$ satisfies the uniform Minkowski estimate
$$
    \underline{\mathcal M}_g^{\,d-m}(\mathcal A)
    =
    \inf_{y\in K}\mathcal H_g^{d-m}(A_y),
    \qquad
    \overline{\mathcal M}_g^{\,d-m}(\mathcal A)
    =
    \sup_{y\in K}\mathcal H_g^{d-m}(A_y).
$$
and
$$
    0<
    \underline{\mathcal M}_g^{\,d-m}(\mathcal A)
    \leq
    \overline{\mathcal M}_g^{\,d-m}(\mathcal A)
    <\infty.
$$
\end{proposition}

\begin{proof}
Since the vectors
$\bigl(\nabla_xF_i(x,t)\bigr)_{1\leq i\leq m}$ are linearly independent
on $A$, the differential
$
    D_xF(x,t):T_x\M\longrightarrow\mathbb R^m
$
is surjective at every $(x,t)\in A$. Hence $A$ is a smooth submanifold
of $\M\times K$ of codimension $m$, and each fiber $A_t$ is a
smooth embedded submanifold of $\M$ of codimension $m$. Moreover,
$$
    N_pA_t
    =
    \operatorname{span}
    \bigl\{\nabla_xF_1(p,t),\ldots,\nabla_xF_m(p,t)\bigr\}.
$$

Consider the fiberwise normal bundle
$$
    \mathcal N A
    :=
    \left\{
        (p,t,\xi):
        (p,t)\in A,\ \xi\in N_pA_t
    \right\}
$$
and the normal exponential map
$$
\begin{aligned}
    \mathcal E:\mathcal N A&\longrightarrow M\times K,\\
    (p,t,\xi)&\longmapsto\bigl(\exp_p(\xi),t\bigr).
\end{aligned}
$$
Let $(p,t,0)\in \mathcal{N}A$. The differential $d\mathcal{E}_{(p,t,0)}$ is invertible and thus by the inverse function theorem, $\mathcal E$ is a local smooth diffeomorphism in a neighborhood of $(p,t,0)$. Since
$A$ is compact, this neighborhood may be chosen uniformly: there exists
$r>0$ such that $\mathcal E$ is a local diffeomorphism on
$$
    \mathcal N_r A
    :=
    \{(p,t,\xi)\in\mathcal N A:|\xi|_g<r\}.
$$

We claim that, after possibly decreasing $r$, the map
$\mathcal E$ is injective on $\mathcal N_r A$. Otherwise, there
would exist sequences
$
    (p_j,t_j,\xi_j)\neq(q_j,t_j,\zeta_j)$ as $
    |\xi_j|_g+|\zeta_j|_g\rightarrow 0,
$
such that
$
    \exp_{p_j}(\xi_j)=\exp_{q_j}(\zeta_j).
$
Since
$
    d_g(p_j,q_j)
    \leq |\xi_j|_g+|\zeta_j|_g\rightarrow0,
$
compactness of $A$ implies the existence of two subsequences converging to the same point $(p,t,0)$ of the zero
section. This contradicts the local injectivity of $\mathcal E$ near
$(p,t,0)$. Thus there exists $r_0>0$ such that
$$
    \mathcal E:\mathcal N_{r_0}A
    \longrightarrow
    \bigl\{(x,t)\in M\times K:d_g(x,A_t)<r_0\bigr\}
$$
is a smooth diffeomorphism.

Let $J_t(p,\xi)$ denote the Jacobian of the exponential map with
respect to the product of the induced measure $d\sigma_t(p)$ on $A_t$
and the Euclidean measure $d\xi$ on $N_pA_t$. Along the zero section,
$
    J_t(p,0)=1.
$
The usual Fermi-coordinate expansion, together with compactness of $A$,
gives
$$
    J_t(p,\xi)=1+\l H_t(p),\xi\rr_g+ O(|\xi|^2_g)
$$
uniformly in $(p,t)\in A$. Therefore, for $0<r<r_0$,
\begin{align*}
    \s\bigl(A_t^{(r)}\bigr)
    &=
    \int_{A_t}
    \int_{\{\,\xi\in N_pA_t:\,|\xi|_g<r\,\}}
        J_t(p,\xi)\,d\xi\,d\sigma_t(p)\\
    &=
    \omega_mr^m\mathcal H_g^{d-m}(A_t)
    +O\left(
        r^{m+2}\mathcal H_g^{d-m}(A_t)
    \right).
\end{align*}

and hence
$$
    \s\bigl(A_t^{(r)}\bigr)
    =
    \omega_mr^m\mathcal H_g^{d-m}(A_t)
    +O(r^{m+2}).
$$
This proves the asserted expansion uniformly in $t$.
\end{proof}
\begin{remark}
Our method also applies to families of sets $(A_i)_{i\in I}$
with uniformly positive reach and null measure, that is, satisfying
$\inf_{i\in I}\mathrm{reach}(A_i)>0$ in the sense of
\cite[Definition 4.2]{SolanesTrillo2025} following Solanes-Trillo tube formula \cite[Corollary 4.7]{SolanesTrillo2025} and \cite[Proposition 4.6]{SolanesTrillo2025} which gives an analogue of the tubular exponential morphism. However, we do not pursue
this extension here.
\end{remark}
\chapter{Lie groups and symmetric spaces}\label[appendix]{symmetric}
Our main references are Helgason's works
\cite{Helgason2000,Helgason2001} for the basic theory of
symmetric spaces and \cite[Sections 2-3]{RemlingRosler2015}
for the formalism of Heckman-Opdam Jacobi polynomials.
 
Assume that $G$ is a semisimple simply connected compact Lie group. Let $\sigma:G\to G$ be an involutive automorphism and $G^\sigma$ the subgroup of fixed points. Let $(G^\sigma)^0\subset K\subset G^\sigma$ be a compact subgroup. Denote by $X:=G/K$ the associated compact symmetric space and by $\mu_X$ the pushforward of the normalized Haar measure along the projection $G\rightarrow X$.
We recall some basic facts on analysis on compact symmetric spaces.
The differential $d\sigma$ induces an eigenspace decomposition of the Lie algebra:
$\mathfrak{g}=\mathfrak{l}\oplus\mathfrak{p}$, with $d\sigma|_{\mathfrak{l}}=\mathrm{Id}$ and $d\sigma|_{\mathfrak{p}}=-\mathrm{Id}$. Let $\mathfrak{a}\subset\mathfrak{p}$ be a maximal abelian subspace and $F:=\exp(\mathfrak{a})/(\exp(\mathfrak{a})\cap K)$ the associated flat torus in $X$. Let $\Sigma$ be the resulting restricted root system, $W_K:=N_K(\mathfrak a)/Z_K(\mathfrak a)$ the restricted Weyl group, and $\Sigma^+$ a choice of positive roots.
Given $\alpha\in\Sigma$, let $m_\alpha$ denote the corresponding root multiplicity. Similarly to the case of conjugacy classes, we have a homeomorphism $F/W_K\simeq K\setminus X$
\cite[Chapter VII, Theorems 8.3-8.6]{Helgason2001}. 
Given a spherical representation $\pi_\lambda$ of restricted highest weight $\lambda\in X^*(F)^+$, we define the spherical function
$\f_\lambda(H)=\l\pi_\lambda(\exp H)v_K,v_K\rr$
for a unit vector $v_K\in V_\lambda^K$.
For $\lambda\in X^*(F)^+$, define the orbit sum $M_\lambda:F\to\bb C$ by
$$
M_\lambda((\exp H)K)
:=
\sum_{\mu\in W_K\lambda}e^{i\l\mu,H\rr}.
$$
Define the Jacobi polynomials
$$
P_\lambda
:=
\sum_{\substack{\mu\in X^*(F)^+\\\mu\leq\lambda}}
c_{\lambda,\mu}M_\mu,
$$
where $\leq$ denotes the dominance order associated with $2\Sigma$, and the coefficients $(c_{\lambda,\mu})$ are chosen so that $c_{\lambda,\lambda}=1$ and
$$
\l P_\lambda,M_\mu\rr_{L^2(F,J_K\,d\theta)}=0
\qquad(\mu\prec\lambda).
$$
The Jacobi polynomials $(P_\lambda)_\lambda$ form an orthogonal basis of $L^2(F,J_K\,d\theta)^{W_K}$.
We collect in \cref{even} explicit
realizations of spherical functions in terms of Heckman-Opdam
Jacobi polynomials for compact symmetric spaces with even
restricted root multiplicities. Our main references are
\cite{Heckman1991,VanDiejenLapointeMorse2004,
BransonOlafssonPasquale2005Compact,
BransonOlafssonPasquale2005Jacobi}.
\begingroup
\renewcommand{\arraystretch}{1.3}
\setlength{\tabcolsep}{10pt}

\begin{xltabular}{\linewidth}{@{}p{0.29\linewidth} X@{}}
\toprule
Restricted root data & Normalized spherical function\\
\midrule
\endfirsthead

\toprule
Restricted root data & Normalized spherical function\\
\midrule
\endhead

\addlinespace[0.5em]
\multicolumn{2}{@{}p{\linewidth}@{}}{%
  \textbf{Group type}\qquad
  $(H\times H)/\operatorname{diag}H$
}\\*[0.8em]

$\Sigma=\Sigma_H$, \quad $m_\alpha=2$
&
$\displaystyle
  \varphi_\kappa(h)
  =\frac{\chi_\kappa(h)}{\chi_\kappa(1)}
$
\par\smallskip
Weyl characters.
\\
\addlinespace[1em]
\midrule

\addlinespace[0.5em]
\multicolumn{2}{@{}p{\linewidth}@{}}{%
  \textbf{Type A\,II}\qquad
  $\mathrm{SU}(2n)/\mathrm{Sp}(n)$
}\\*[0.8em]

$\Sigma=A_{n-1}$, \quad $m_\alpha=4$
&
$\displaystyle
  \varphi_\kappa(z)
  =\frac{P_\kappa^{(1/2)}(z)}
         {P_\kappa^{(1/2)}(1^n)}
$
\par\smallskip
Jack polynomials with parameter $\alpha=\tfrac12$.
\\
\addlinespace[1em]
\midrule

\addlinespace[0.5em]
\multicolumn{2}{@{}p{\linewidth}@{}}{%
  \textbf{Type E\,IV}\qquad
  $E_6/F_4$
}\\*[0.8em]

$\Sigma=A_2$, \quad $m_\alpha=8$
&
$\displaystyle
  \varphi_\kappa(z)
  =\frac{P_\kappa^{(1/4)}(z)}
         {P_\kappa^{(1/4)}(1^3)}
$
\par\smallskip
Jack polynomials with parameter $\alpha=\tfrac14$.
\\
\addlinespace[1em]
\midrule

\addlinespace[0.5em]
\multicolumn{2}{@{}p{\linewidth}@{}}{%
  \textbf{Odd-dimensional spheres}\qquad
  $S^{2r+1}
   =\mathrm{Spin}(2r+2)/\mathrm{Spin}(2r+1)$
}\\*[0.8em]

$\Sigma=A_1$, \quad $m_\alpha=2r$
&
$\displaystyle
  \varphi_\ell(\theta)
  =\frac{C_\ell^{\,r}(\cos\theta)}
         {C_\ell^{\,r}(1)}
$
\par\medskip
$\displaystyle
  ={}_2F_1\left(
    \begin{matrix}
      -\ell,\ \ell+2r\\
      r+\frac12
    \end{matrix}
    ;\sin^2\frac{\theta}{2}
  \right).
$
\par\smallskip
Here $C_\ell^r$ denotes the Gegenbauer polynomial of degree
$\ell$ and parameter $r$, and ${}_2F_1$ denotes the Gauss
hypergeometric function.
\\
\addlinespace[0.8em]
\bottomrule

\caption{Spherical functions for compact symmetric spaces with
even restricted root multiplicities.}
\label{even}\\
\end{xltabular}

\smallskip

\endgroup
\chapter{Manifold submetries and singular foliations}\label[appendix]{what}
In this appendix, we collect the geometric and analytic
properties of Riemannian submetries and possibly singular
Riemannian foliations used throughout the article.
We recall the role of basic mean curvature, examine the
behavior of Minkowski content under Riemannian submetries,
and establish the regularity and spectral properties of
the averaging operator onto basic functions.
These results provide the necessary framework for extending
our spectral approximation theory to singular leaf spaces.
\begin{theorem}\cite[Theorem B]{LinMendesRadeschi2025}\label{f}
    If $(\M,g)$ is a compact normal homogeneous space (in particular a compact Lie group), any Riemannian submetry with compact fibers has basic mean curvature.
   If $\ff$ is a (possibly singular) Riemannian foliation with compact leaves on such a manifold, then the projection map into the leaf space $\sigma:\M \rightarrow \M/\ff$ is a manifold submetry with basic mean curvature.
\end{theorem}
\begin{remark}\label{saturated}
\begin{enumerate}
    \item 
   Let $\sigma:\M\to X$ be a Riemannian submetry, and equip $X$ with
$\mu_X:=\sigma_*\operatorname{vol}_g$. For $B\subset X$, put $A:=\sigma^{-1}(B)$. Then
$$
\sigma^{-1}(B^{(r)})=A^{(r)},
\qquad
\sigma^{-1}(\partial B)=\partial A.
$$
Consequently, using the same reference dimension $d=\dim\M$,
the upper and lower Minkowski contents of $\partial B$ with
respect to $\mu_X$ equal those of $\partial A$ with respect to
$\operatorname{vol}_g$.
    \item 
    The map $A\rightarrow \sigma(A)$ induces a bijection between $\sigma$-saturated sets of $\M$ and subsets of $X$.
    \end{enumerate}
\end{remark}
\begin{lemma}\cite[Lemma 26]{LinMendesRadeschi2025}\label{no}
Let $\sigma\colon \M\to X$ be a Riemannian submetry with compact fibers
and basic mean curvature. Then the following hold:
\begin{enumerate}
    \item[(a)] The averaging operator is well-defined and is the
    $L^{2}$-orthogonal projection onto the subspace of $L^{2}(\M)$
    consisting of the elements with a basic representative.

    \item[(b)] For every $f\in C^{\infty}(\M)$, the element
    $P_\sigma f\in L^{2}(\M)$ has a smooth representative.

    \item[(c)] The operators
    $\Delta_g,P_\sigma\colon C^{\infty}(\M)\to C^{\infty}(\M)$
    commute.
\end{enumerate}
\end{lemma}
We notice that \cite[Lemma 26]{LinMendesRadeschi2025} provides only a smooth representative for the basic projection. As remarked by \cite{LytchakRadeschi2018} in the case of spherical foliations it is possible to define $P_\sigma f\in C^\infty(\M)$ pointwise. As we seek mainly smooth approximations, we shall give a rigorous proof of that fact in our Riemannian setting.

\begin{lemma}\label{basic}
Let $\sigma:\M\rightarrow X$ be a Riemannian submetry
with basic mean curvature and compact smooth fibers.
Then
$$
P_\sigma C^\infty(\M)\subset C^\infty(\M,\sigma)
$$
pointwise.
\end{lemma}

\begin{proof}
Fix $x\in\M$, and write $L_p:=\sigma^{-1}(p)$ for the
fiber containing $x$, of codimension $q$.
By \cref{exp}, for sufficiently small
$\varepsilon>0$, there exists a tubular diffeomorphism
$$
\psi_p:N_\varepsilon L_p\rightarrow L_p^{(\varepsilon)},
\qquad
\psi_p(z,0)=z.
$$
Denote its Jacobian by $J$, so that $J(z,0)=1$.
For $v\in C^\infty(\M)$, the change-of-variables formula
gives
$$
\int_{L_p^{(\varepsilon)}}v\,d\s
=
\int_{N_\varepsilon L_p}
(v\circ\psi_p)(z,t)J(z,t)
\,d\operatorname{vol}_{g|L_p}(z)\,dt.
$$
Moreover, the leading term in Gray's tube formula as in \cref{gray} gives
$$
\s\bigl(L_p^{(\varepsilon)}\bigr)
=
\omega_q\varepsilon^q
\operatorname{vol}(L_p)
+
O(\varepsilon^{q+2}),
$$
where $\omega_q$ is the volume of the Euclidean unit
ball in $\mathbb R^q$.

Let $f\in C^\infty(\M)$, and let $u\in C^\infty(\M)$
be the representative of $P_\sigma f$ obtained above.
On the dense regular part $\M_{reg}$, $u$ is basic, by continuity and
the submetry property, it is basic on every fiber.

Since $L_p^{(\varepsilon)}$ is saturated, its indicator
is basic. Orthogonality of $P_\sigma$ therefore gives
\begin{align*}
u(x)
&=
\frac1{\operatorname{vol}(L_p)}
\int_{L_p}u\,d\operatorname{vol}_{g|L_p}
\\
&=
\lim_{\varepsilon\to0^+}
\frac1{\s(L_p^{(\varepsilon)})}
\int_{L_p^{(\varepsilon)}}u\,d\s
\\
&=
\lim_{\varepsilon\to0^+}
\frac1{\s(L_p^{(\varepsilon)})}
\int_{L_p^{(\varepsilon)}}P_\sigma f\,d\s
\\
&=
\lim_{\varepsilon\to0^+}
\frac1{\s(L_p^{(\varepsilon)})}
\int_{L_p^{(\varepsilon)}}f\,d\s
\\
&=
\lim_{\varepsilon\to0^+}
\frac{
\displaystyle
\int_{N_\varepsilon L_p}
(f\circ\psi_p)J
\,d\operatorname{vol}_{g|L_p}\,dt
}{
\omega_q\varepsilon^q\operatorname{vol}(L_p)
}
\\
&=\lim_{\varepsilon\to0^+}
\frac1{\s(L_p^{(\varepsilon)})}
\int_{L_p^{(\varepsilon)}}f\,d\s
\\
&=
\frac1{\operatorname{vol}(L_p)}
\int_{L_p}f\,d\operatorname{vol}_{g|L_p}
\\
&=(P_\sigma f)(x).
\end{align*}
Thus the fiberwise average agrees everywhere with its
smooth basic representative.
\end{proof}
We give examples of Riemannian submetries arising from
possibly singular Riemannian foliations with closed leaves. Let $G$ be a connected compact Lie group equipped with the $\mathrm{Ad}$-invariant metric. For an irreducible unitary representation $\pi_\lambda$
and an orthonormal basis $(e_i)$, put
$$
m_{ij}^\lambda(g)
:=\sqrt{\dim V_\lambda}\,
\langle\pi_\lambda(g)e_i,e_j\rangle.
$$


\par\medskip
\noindent
\begingroup
\small
\setlength{\tabcolsep}{3pt}
\renewcommand{\arraystretch}{1.3}

\captionof{table}{Basic projections and their projected coefficients}
\label{tableofbasic}

\noindent
\begin{tabularx}{\linewidth}{@{}XX@{}}
\toprule
Basic projection & Projected coefficients\\
\midrule

\multicolumn{2}{@{}p{\linewidth}@{}}{%
  \textbf{Conjugacy}\quad
  $L_g=\{hgh^{-1}:h\in G\}$%
}\\[2mm]
$\displaystyle
 \begin{aligned}
 (P_{\mathrm{Conj}}f)(g)
 &=\int_G f(hgh^{-1})\\
 &\qquad d\mu_G(h).
 \end{aligned}
$
&
$\displaystyle
 P_{\mathrm{Conj}}m_{ij}^\lambda
 =\frac{\delta_{ij}}{\sqrt{\dim V_\lambda}}\chi_\lambda.
$
\\[3mm]
\midrule

\multicolumn{2}{@{}p{\linewidth}@{}}{%
  \textbf{Conjugacy and normal subgroup translation}%
}\\
\multicolumn{2}{@{}p{\linewidth}@{}}{%
  $L_g=\{nhgh^{-1}:n\in N,\ h\in G\},
   \qquad N\trianglelefteq G$%
}\\[2mm]
$\displaystyle
 \begin{aligned}
 (P_{\mathrm{Conj},N}f)(g)
 &=\int_N\int_G f(nhgh^{-1})\\
 &\qquad d\mu_G(h)\,d\mu_N(n).
 \end{aligned}
$
&
If $N\subseteq\ker\pi_\lambda$,
\par\smallskip
$\displaystyle
 P_{\mathrm{Conj},N}m_{ij}^\lambda
 =\frac{\delta_{ij}}{\sqrt{\dim V_\lambda}}\chi_\lambda;
$
\par\smallskip
otherwise the projection is zero.
\\[3mm]
\midrule

\multicolumn{2}{@{}p{\linewidth}@{}}{%
  \textbf{Twisted conjugacy}\quad
  $L_g=\{hg\theta(h)^{-1}:h\in G\}$%
}\\[2mm]
$\displaystyle
 \begin{aligned}
 (P_\theta f)(g)
 &=\int_G f(hg\theta(h)^{-1})\\
 &\qquad d\mu_G(h).
 \end{aligned}
$
&
If $A_\lambda:\pi_\lambda\simeq\pi_\lambda\circ\theta$ is an intertwiner,
\par\smallskip
$\displaystyle
 P_\theta m_{ij}^\lambda
 =\frac{\overline{(A_\lambda)_{ji}}}
        {\sqrt{\dim V_\lambda}}\chi_\lambda^\theta;
$
\par\smallskip
otherwise the projection is zero.
\\[3mm]
\bottomrule
\end{tabularx}

\endgroup
\par\medskip

We record some explicit examples of Riemannian submetries
arising from foliations on the sphere.
\begin{example}\label{foliationexamples}
Equip the source spheres with their unit round metrics.
\begin{enumerate}
\item \textbf{Latitude foliations}:
For $d\geq2$, the action of the special orthogonal group
$\mathrm{SO}(d)$ on $\bb{S}^d$ induces a submetry
$$
\sigma:\S^d\to[-1,1],\qquad \sigma(x)=x_{d+1}.
$$
The leaf space has the quotient distance
$$
d_X(s,t)=|\arccos s-\arccos t|.
$$
The regular leaves are $(d-1)$-spheres, while
the poles are singular leaves. The basic eigenfunctions
are Gegenbauer polynomials, with eigenvalues
$\lambda_k=k(k+d-1)$, $k\geq0$;
see \cite{AlexandrinoBettiol2015}.

\item \textbf{Products of spheres}:
The action of $\mathrm{SO}(p)\times\mathrm{SO}(q)$
on $\bb{S}^{p+q-1}\subset\R^p\times\R^q$, for $p,q\geq2$,
induces a submetry
$$
\sigma(u,v)=|u|^2\in[0,1].
$$
The leaf space carries the quotient distance
$$
d_X(s,t)=|\arccos\sqrt{s}-\arccos\sqrt{t}|.
$$
The regular leaves are diffeomorphic to
$\S^{p-1}\times\S^{q-1}$.
The basic eigenfunctions are Jacobi polynomials,
with eigenvalues
$$
\lambda_k=2k(2k+p+q-2),\qquad k\geq0;
$$
see \cite{AlexandrinoBettiol2015,Dadok1985}.

\item \textbf{Hopf fibrations}:
The classical Hopf fibrations, equipped with their
quotient metrics, are$$
\bb{S}^{2m+1}\to\mathbb{CP}^m,\qquad
\S^{4m+3}\to\mathbb{HP}^m,\qquad
\S^{15}\to \bb{O}\mathbb{P}^1.
$$
Their basic eigenvalues are, respectively,$$
4k(k+m),\qquad
4k(k+2m+1),\qquad
4k(k+7),
\qquad k\geq0,
$$
where $m\geq1$.
See \cite{BaudoinCho2021,OrneaPartonPiccinniVuletescu2013}.
\end{enumerate}
\end{example}

We conclude by recalling the basic Weyl law used in
the application to Shimura varieties.
\begin{theorem}[{\cite[Theorem A]{LinMendesRadeschi2025Weyl}}]
\label{basicweyl}
Let $\mathcal F$ be a possibly singular
Riemannian foliation with closed leaves and basic mean
curvature. Put $X:=\M/\mathcal F$ equipped with the
pushforward $\mu_X$ of the volume measure and set
$m:=\dim X$.
Then the basic spectral counting function satisfies
$$
N_{\mathcal F}(T)
\sim
\frac{\omega_m\mathrm{vol}(X_{\mathrm{reg}})}{(2\pi)^m}
T^{m/2}.
$$
In particular,
$$
N_{\mathcal F}(T)
\ll_{\M,\mathcal F}(1+T)^{m/2}.
$$
\end{theorem}

\chapter{Motivic $L$-functions and character sums}\label[appendix]{formalism}
Let $M\in\mathcal{M}_{\sim}(K,E)$ be a motive over the number
field $K$ with coefficients in the number field $E$, and let
$(\rho_{M,\lambda})_\lambda$
be its $\lambda$-adic realization. We recall Deligne's formalism of the invariants of the motivic $L$-function
\cite[Section 5.2(a), p. 328]{DeligneValeursL}.

Assume that the realizations $(\rho_{M,\lambda})_\lambda$
form a linearly strongly compatible system.
For each finite place $v$ of $K$, choose $\lambda$
of different residue characteristic and write
$$
\operatorname{WD}_v(M)=(r_{M,v},N_{M,v})
$$
for the associated Frobenius-semisimple
Weil--Deligne representation.
Let $V_{M,v}$ be its underlying vector space, and put
$n_M:=\dim_E H_B(M)$ and $q_v:=Nv$.

Let $I_v$ be the inertia subgroup, and
let $(G_{K_v}^{u})_{u\geq0}$ be the upper-numbering
ramification filtration. The Swan conductor is
$$
\operatorname{Sw}(r_{M,v})
:=
\int_0^\infty
\operatorname{codim} V_{M,v}^{G_{K_v}^{u}}\,du.
$$
The local conductor exponent is
$$
a_v(M)
:=
\operatorname{Sw}(r_{M,v})
+n_M
-\dim\bigl((\ker N_{M,v})^{I_v}\bigr).
$$
We set
$$
\mathfrak f(M)
:=
\prod_{v\nmid\infty}\mathfrak{p}_v^{a_v(M)}
$$
and
$$
q_M
:=
|d_K|^{n_M}\,
\mathrm N_{K/\mathbb Q}\mathfrak{f}(M),
$$
where $d_K$ is the discriminant of $K$.
Let $v$ be a finite place of $K$, and choose $\lambda$ whose residue
characteristic differs from that of $v$. Denote by $I_v$ the inertia
subgroup of a decomposition group at $v$, and by $F_v$ a geometric
Frobenius acting on $H_\lambda(M)^{I_v}$. Strong compatibility implies that
$$
P_v(M,T)
:=
\det\left(1-F_vT\mid (\ker N_{M,v})^{r_{M,v}(I_v)}\right)
\in E[T]
$$
is independent of $\lambda$. For an embedding $\iota\colon E\hookrightarrow\C$ define
\begin{align*}
L_v(\iota,M,s)&:=\iota\left(P_v(M,q_v^{-s})\right)^{-1},\\
L(\iota,M,s)&:=\prod_{v\nmid\infty}L_v(\iota,M,s).
\end{align*}

Let $v\mid\infty$, choose an embedding
$\tau:K\hookrightarrow\mathbb C$ inducing $v$, and set
$V:=H_{B,\tau}(M)\otimes_{E,\iota}\mathbb C$.
Write $$V=\bigoplus_{p,q}V^{p,q},$$
 define $h_v^{p,q}(\iota,M):=\dim_\bb{C}V^{p,q}$. If $v$ is complex, then
\begin{equation*}
L_v(\iota,M,s):=\prod_{p,q}\Gamma_{\bb C}\left(s-\min\{p,q\}\right)^{h_v^{p,q}(\iota,M)}. 
\end{equation*}
If $v$ is real, then complex conjugation gives an involution $F_{\infty,v}$ on $V$ such that
\begin{equation*}
F_{\infty,v}V^{p,q}=V^{q,p}.
\end{equation*}
For ${\varepsilon}\in\{0,1\}$, we define
\begin{equation*}
h_v^{p,p,{\varepsilon}}(\iota,M):=\dim_\bb{C}\left\{x\in V^{p,p}:F_{\infty,v}x=(-1)^{p+{\varepsilon}}x\right\}.
\end{equation*}
Thus,
\begin{align*}
L_v(\iota,M,s)&:=\prod_{p<q}\Gamma_{\bb C}(s-p)^{h_v^{p,q}(\iota,M)}
\prod_p\prod_{{\varepsilon}=0}^{1}\Gamma_{\R}(s+{\varepsilon}-p)^{h_v^{p,p,{\varepsilon}}(\iota,M)},
\end{align*}
where $\Gamma_{\R}(s):=\pi^{-s/2}\Gamma(s/2)$ and $\Gamma_{\bb C}(s):=2(2\pi)^{-s}\Gamma(s)$. Hence, define
\begin{equation*}
L_\infty(\iota,M,s):=\prod_{v\mid\infty}L_v(\iota,M,s).
\end{equation*}
The completed motivic $L$-function is 
\begin{equation*}
\Lambda(\iota,M,s):=L_\infty(\iota,M,s)L(\iota,M,s).
\end{equation*}
Write:
$$
L_\infty\left(\iota,M,s+\frac{w}{2}\right)
=
\prod_j\Gamma_{\mathbb R}(s+\mu_{j,\iota,M})
\prod_k\Gamma_{\mathbb C}(s+\nu_{k,\iota,M}),
$$
with nonnegative gamma shifts appearing with respective multiplicities.
The real shifts are $0$ or $1$, while the complex shifts are
given by the corresponding quantities $|q-p|/2$.
Define
$$
H_\infty(M)
:=
\max_{\iota:E\hookrightarrow\mathbb C}
\max\left(
\{0\}
\cup\{\mu_{j,\iota,M}\}_j
\cup\{\nu_{k,\iota,M}\}_k
\right).
$$
Fix $T\geq0$. For a complex place $v$, put
$$
\mathcal Q_{\iota,M,v}(T)
:=
\prod_{p,q}
\left(3+T+\frac{|q-p|}{2}\right)^{2h_v^{p,q}(\iota,M)}.
$$
For a real place $v$, put
$$
\begin{aligned}
\mathcal Q_{\iota,M,v}(T)
:={}&
\prod_{p<q}
\left(3+T+\frac{q-p}{2}\right)^{2h_v^{p,q}(\iota,M)}
\\[3pt]
&\times
\prod_{\varepsilon=0}^{1}
(3+T+\varepsilon)^{h_v^{w/2,w/2,\varepsilon}(\iota,M)}.
\end{aligned}
$$
If $w$ is odd, the diagonal Hodge numbers are understood to be zero.
Set
$$
\mathcal Q_{\iota,M}(T)
:=q_M\prod_{v\mid\infty}\mathcal Q_{\iota,M,v}(T),
\qquad
Q(\iota,M):=\mathcal Q_{\iota,M}(0).
$$
Let $r_1$ and $r_2$ denote the numbers of real and complex places
of $K$, respectively, so that $n_K=r_1+2r_2$.
Recall that, at every complex place,
$$
\sum_{p,q}h_v^{p,q}(\iota,M)=n_M,
$$
while, at every real place,
$$
2\sum_{p<q}h_v^{p,q}(\iota,M)
+
\sum_{\varepsilon=0}^{1}
h_v^{w/2,w/2,\varepsilon}(\iota,M)
=n_M.
$$

Using the elementary inequality
$$
\frac{3+T+x}{3+x}\leq1+\frac{T}{3},
$$
for $T,x\geq 0$, we obtain
$$
\mathcal Q_{\iota,M}(T)
\leq
Q(\iota,M)\left(1+\frac{T}{3}\right)^{n_Kn_M}.
$$

Moreover,

\begin{equation}
\begin{aligned}
\log Q(\iota,M)
&\leq
\log q_M
+n_M(r_1+2r_2)\log\bigl(3+H_\infty(M)\bigr)
\\
&\leq
\log q_M
+2n_Mn_K\log\bigl(3+H_\infty(M)\bigr).
\end{aligned}
\end{equation}

\begin{lemma}\label{tensor}
Let $M,N\in\mathcal M_\sim(K,E)$ be pure motives of ranks
$n_M$ and $n_N$, respectively, whose $\lambda$-adic realizations
form strongly compatible systems. Then
$$
q_{M\otimes N}
\mid
q_M^{n_N}q_N^{n_M}.
$$
Moreover, for every embedding $\iota:E\hookrightarrow\mathbb C$
and every $T\geq1$,
$$
\log\mathcal Q_{\iota,M\otimes N}(T)
\leq
2n_N\log\mathcal Q_{\iota,M}(T)
+
2n_M\log\mathcal Q_{\iota,N}(T).
$$
\end{lemma}

\begin{proof}
 
Recall that given two finite-dimensional $E$-vector spaces $V,W$ with the action of a finite group $H$
one has $V^H\otimes W^H\subseteq(V\otimes W)^H$
$$
\operatorname{codim}(V\otimes W)^H
\leq
\dim W\,\operatorname{codim}V^H
+
\dim V\,\operatorname{codim}W^H.
$$
Applying this estimate to the higher ramification groups yields
\def \codim {\mathrm{codim}}
\begin{align*}
    \operatorname{Sw}(r_{M,v}\otimes r_{N,v})=&\int_{0}^\infty \codim((V_{M,v}\otimes V_{N,v})^{G^u_{K_v}})du\\&\leq \int_{0}^\infty n_M \codim(V_{N,v}^{G^u_{K_v}})+n_N\codim(V_{M,v}^{G^u_{K_v}})du\\
    &= n_N\Sw(r_{M,v})+n_M\Sw(r_{N,v}).
\end{align*}

Similarly we have the inclusion
$$
(\ker N_{M,v})^{I_v}
\otimes
(\ker N_{N,v})^{I_v}
\subseteq
(\ker N_{M,v}\otimes1+1\otimes N_{N,v})^{I_v},
$$
which implies

\begin{align*}
&n_Mn_N-
\dim\bigl(\ker(N_{M,v}\otimes1+1\otimes N_{N,v})\bigr)^{I_v}\\
&\quad\leq
n_N\bigl(n_M-\dim(\ker N_{M,v})^{I_v}\bigr)
+n_M\bigl(n_N-\dim(\ker N_{N,v})^{I_v}\bigr).
\end{align*}

Using the conductor formula and combining these inequalities,
one finds
$$
a_v(M\otimes N)\leq n_N\,a_v(M)+n_M\,a_v(N).
$$
Thus
$$
\mathfrak{f}(M\otimes N)
\mid
\mathfrak{f}(M)^{n_N}\mathfrak{f}(N)^{n_M}.
$$
In particular,
$$
\mathfrak{a}
:=
\mathfrak{f}(M)^{n_N}\mathfrak{f}(N)^{n_M}
\mathfrak{f}(M\otimes N)^{-1}
$$
is an integral ideal, and
$$
\frac{q_M^{n_N} q_N^{n_M}}{q_{M\otimes N}}
=
|d_K|^{n_Mn_N}\,
\mathrm N_{K/\mathbb Q}\mathfrak{a}
\in\mathbb N.
$$
This shows the claimed divisibility.

Consider now the archimedean factors. For fixed $\iota$ and $T\geq1$,
set $A:=3+T$. The Hodge decomposition of a tensor product gives
$$
h_v^{p,q}(\iota,M\otimes N)
=
\sum_{\substack{a+r=p\\b+s=q}}
h_v^{a,b}(\iota,M)\,
h_v^{r,s}(\iota,N).
$$
At a complex place, we have
$$
3+T+\frac{|(a+r)-(b+s)|}{2}
\leq
\left(3+T+\frac{|a-b|}{2}\right)
\left(3+T+\frac{|r-s|}{2}\right).
$$
Combining this with the dimensional identities gives
$$
\begin{aligned}
\mathcal Q_{\iota,M\otimes N,v}(T)
&\leq
\prod_{a,b,r,s}
\left[
\left(3+T+\frac{|a-b|}{2}\right)
\left(3+T+\frac{|r-s|}{2}\right)
\right]^{
2h_v^{a,b}(\iota,M)h_v^{r,s}(\iota,N)}
\\
&=
\mathcal Q_{\iota,M,v}(T)^{n_N}
\mathcal Q_{\iota,N,v}(T)^{n_M}.
\end{aligned}
$$

At a real place, put $A:=3+T$.
Given a Hodge block $B$ in the realization of a motive $M$,
let $\mathscr Q_M(B)$ denote its contribution to
$\mathcal Q_{\iota,M,v}(T)$.
For a diagonal block $D_{\varepsilon}$ of parity
${\varepsilon}\in\{0,1\}$ and an off-diagonal block $O_r$
of normalized gap $r>0$, one has
$$
\begin{array}{c|c|c}
B&\dim B&\mathscr Q_M(B)\\ \hline
D_{\varepsilon}&1&A+{\varepsilon}\\
O_r&2&(A+r)^2
\end{array}
$$
We claim that, for any two elementary summands $B,C$,
$$
\mathscr Q_{M\otimes N}(B\otimes C)
\leq
\mathscr Q_M(B)^{2\dim C}
\mathscr Q_N(C)^{2\dim B}.
$$

For two diagonal summands, let
$\overline{{\varepsilon}+\eta}\in\{0,1\}$ denote reduction
modulo $2$. Then
$$
\begin{aligned}
\mathscr Q_{M\otimes N}(D_{\varepsilon}\otimes D_\eta)
&=3+T+\overline{{\varepsilon}+\eta}\\
&\leq(3+T+{\varepsilon})^2(3+T+\eta)^2,
\end{aligned}
$$
and the claim holds.

For a diagonal summand and an off-diagonal summand,
$$
\mathscr Q(D_{\varepsilon}\otimes O_r)
=(3+T+r)^2
\leq
(3+T+{\varepsilon})^4(3+T+r)^4,
$$
so the claim holds in this case as well.

For two off-diagonal summands with different gaps
$r\ne\eta$, the tensor product consists of two
off-diagonal terms with gaps $r+\eta$ and
$|r-\eta|$. We have
$$
3+T+r+\eta,\quad 3+T+|r-\eta|
\leq(3+T+r)(3+T+\eta),
$$
and thus
$$
\begin{aligned}
\mathscr Q(O_r\otimes O_\eta)
&=(A+r+\eta)^2(A+|r-\eta|)^2\\
&\leq(A+r)^4(A+\eta)^4\\
&\leq(A+r)^8(A+\eta)^8.
\end{aligned}
$$
Thus the desired claim holds in this case as well.

Finally,
$$
\mathscr Q(O_r\otimes O_r)
=(3+T+2r)^2(3+T)(3+T+1).
$$
Since $r\geq\tfrac12$,
$$
(3+T)(3+T+1)\leq(3+T+r)^2,
\qquad
3+T+2r\leq(3+T+r)^2.
$$
Consequently,
$$
\mathscr Q(O_r\otimes O_r)
\leq(3+T+r)^6
\leq(3+T+r)^{16},
$$
which proves the claim in the remaining case.

As a result, for real places $v$,
\begin{align*}
\mathcal Q_{\iota,M\otimes N,v}(T)
&\leq
\prod_{i,j}
\mathscr Q_M(B_i)^{2\dim C_j}
\mathscr Q_N(C_j)^{2\dim B_i}\\
&=
\mathcal Q_{\iota,M,v}(T)^{2n_N}
\mathcal Q_{\iota,N,v}(T)^{2n_M}.
\end{align*}

Combining this with the finite conductor bound and the
complex-place estimate, and using $q_M,q_N\geq1$, gives
$$
\begin{aligned}
\mathcal Q_{\iota,M\otimes N}(T)
&\leq
q_M^{n_N}q_N^{n_M}
\prod_{v\mid\infty}
\mathcal Q_{\iota,M,v}(T)^{2n_N}
\mathcal Q_{\iota,N,v}(T)^{2n_M}\\
&\leq
\mathcal Q_{\iota,M}(T)^{2n_N}
\mathcal Q_{\iota,N}(T)^{2n_M}.
\end{aligned}
$$
Taking logarithms proves the assertion.
\end{proof}
The argument of the next proposition follows the classical proof of the character sums for Artin
$L$-functions. Since these estimates are a key input
to our motivic equidistribution theorems, we
sketch the proof in order to keep track of the dependence on the
relevant motivic invariants.
\begin{proposition}\label{sumc}
Let $M\in\mathcal M_\sim(K,E)$ be a pure motive of weight $w$
and degree $n_M$ whose $\lambda$-adic realizations form linearly strongly
compatible. Let $S_M$ be the set of finite places where
the corresponding Weil-Deligne datum is ramified, and define
$$
S(M):=n_M|S_M|.
$$
Fix an embedding $\iota:E\hookrightarrow\mathbb C$.
Assume the $L$-function formalism of Conjecture \ref{grh} for $M$ and the motivic Ramanujan conjecture for $M^\vee$.
Put
$$
d:=
\begin{cases}
\displaystyle
\dim_E\operatorname{Hom}_{\mathcal M_\sim(K,E)}
(\mathbf1,M(w/2)), & w\in2\mathbb Z,\\
0, & w\notin2\mathbb Z.
\end{cases}
$$
For $\pp\notin S_M$, define
$$
\lambda_{\iota,M}(\pp)
:=
(N\pp)^{-w/2}
\iota\left(
\operatorname{Tr}
\bigl(\operatorname{Frob}_{\pp}\mid H_\lambda(M)\bigr)
\right),
$$
where $\lambda$ has residue characteristic different from that
of $\pp$, and Frobenius is geometric.
Write $Q(\iota,M,T):=\mathcal Q_{\iota,M}(T)$ and
$Q(\iota,M):=Q(\iota,M,0)$.

Then, for every $x\geq2$,
$$
\sum_{\substack{N\pp\leq x\\\pp\notin S_M}}
\lambda_{\iota,M}(\pp)
=
d\,\Li(x)
+
O\left(
x^{1/2}\bigl(1+d+\log Q(\iota,M,x)\bigr)
\right)
+
O(S(M)).
$$
Moreover, for every $f\in C_c^\infty((0,\infty))$ and $x\geq2$,
\begin{align*}
\sum_{\pp\notin S_M}
\lambda_{\iota,M}(\pp)
f\left(\frac{N\pp}{x}\right)
={}&
d\int_2^\infty\frac{f(t/x)}{\log t}\,dt
\\
&+
O_f\left(
\frac{x^{1/2}}{\log x}
\bigl(1+d+\log Q(\iota,M)\bigr)
\right)
+
O_f\left({S(M)}\right).
\end{align*}
The implied constants in the first estimate are absolute,
and those in the second depend only on $\|f^{(j)}\|_{L^\infty((0,\infty))}$ for $j\in \{0,1,2\}$ and the support of $f$.
\end{proposition}

\begin{proof}
We follow the classical $L$-function argument of
\cite[Proposition 4.1]{Murty1985}.
Put
$$
\mathcal L(s):=L(\iota,M,s+w/2),
\qquad
D:=n_Kn_M,
\qquad
\kappa(T):=1+d+\log \mathcal{Q}_{\iota,M}(T).
$$
Recall that
$$
D\log3\leq\log Q(\iota,M),
\qquad
\kappa(T)\leq\kappa(0)+D\log(1+T/3).
$$
By Hypothesis \ref{grh}, $\mathcal L(s)$ has a pole
of order $d$ at $s=1$.
Write the normalized finite Euler factors as
$$
L_{\pp}(\iota,M,s+w/2)
=
\prod_{j=1}^{n_{\pp}}
\left(1-\alpha_{\pp,j}(N\pp)^{-s}\right)^{-1},
$$
with $n_{\pp}\leq n_M$, and put
$$
b_{\pp}(m):=\sum_{j=1}^{n_{\pp}}\alpha_{\pp,j}^{\,m}.
$$
The motivic Ramanujan conjecture implies that
$|b_{\pp}(m)|\leq n_M$.
Let $\Lambda(n)$ denote the classical von Mangoldt
function. For $\Re s>1$, we have
$$
-\frac{\mathcal L'}{\mathcal L}(s)
=
\sum_{\pp}\sum_{m\geq1}
b_{\pp}(m)\frac{\log N\pp}{(N\pp)^{ms}}
=
\sum_{n\geq1}\frac{b(n)}{n^s},
$$
where $|b(n)|\leq D\Lambda(n)$.
Define
$$
\psi_M(x):=\sum_{n\leq x}b(n),
\qquad
\psi_{M,1}(x):=\sum_{n\leq x}b(n)(x-n).
$$
Taking $h:=x^{1/2}$, the bound on $b(n)$ gives
\begin{equation}\label{psi}
\frac{\psi_{M,1}(x+h)-\psi_{M,1}(x)}h
=
\psi_M(x)+O(Dx^{1/2}\log x).
\end{equation}

Recall the classical identity
$$
\frac{1}{2\pi i}
\int_{2-i\infty}^{2+i\infty}
\frac{x^{s+1}n^{-s}}{s(s+1)}\,ds
=
\begin{cases}
x-n,&n<x,\\
0,&n\geq x.
\end{cases}
$$
Absolute convergence on $\Re s=2$ therefore gives
$$
\psi_{M,1}(x)
=
\frac1{2\pi i}
\int_{2-i\infty}^{2+i\infty}
-\frac{\mathcal L'}{\mathcal L}(s)
\frac{x^{s+1}}{s(s+1)}\,ds.
$$
Shifting the contour along suitable heights avoiding
the zero ordinates yields
\begin{align*}
\psi_{M,1}(x)
={}&
\frac1{2\pi i}
\int_{-1/4-i\infty}^{-1/4+i\infty}
-\frac{\mathcal L'}{\mathcal L}(s)
\frac{x^{s+1}}{s(s+1)}\,ds\\
&+\frac{dx^2}{2}
-\sum_\rho\frac{x^{\rho+1}}{\rho(\rho+1)}
+R_0(x),
\end{align*}
where $\rho=\frac12+i\gamma$ runs through the
nontrivial zeros, counted with multiplicity, and
$R_0(x)$ is the corresponding residue at $0$.

Put
$$
\gamma_M(s):=L_\infty(\iota,M,s+w/2),
\qquad
\mathcal L^\vee(s):=L(\iota,M^\vee,s-w/2).
$$
The normalized archimedean factors and conductors
of $M$ and $M^\vee$ agree.
By the functional equation,
\begin{align*}
-\frac{\mathcal L'}{\mathcal L}(-1/4+it)
={}&
\log q_M
+\frac{\gamma_M'}{\gamma_M}(-1/4+it)
+\frac{\gamma_M'}{\gamma_M}(5/4-it)\\
&+\frac{(\mathcal L^\vee)'}{\mathcal L^\vee}(5/4-it).
\end{align*}
By the motivic Ramanujan conjecture for $M^\vee$,
$$
\left|
\frac{(\mathcal L^\vee)'}{\mathcal L^\vee}(5/4-it)
\right|
\ll D.
$$
By Stirling’s formula, we have, for $z\in\{-1/4,5/4\}$,
$$
\left|
\frac{\gamma_M'}{\gamma_M}(z+it)
\right|
\ll\kappa(0)+D\log(2+|t|).
$$
Thus
$$
\left|
\frac{\mathcal L'}{\mathcal L}(-1/4+it)
\right|
\ll\kappa(0)+D\log(2+|t|).
$$
Consequently,
\begin{align*}
&\left|
\int_{-1/4-i\infty}^{-1/4+i\infty}
-\frac{\mathcal L'}{\mathcal L}(s)
\frac{x^{s+1}}{s(s+1)}\,ds
\right|\\
&\qquad\ll
x^{3/4}
\int_{\R}
\frac{\kappa(0)+D\log(2+|t|)}{1+t^2}\,dt
\ll x^{3/4}\kappa(0).
\end{align*}

Write the Laurent expansion
$$
\frac{\mathcal L'}{\mathcal L}(s)
=
\frac{m_0}{s}+c_0+O(s),
\qquad
m_0:=\operatorname{ord}_{s=0}\mathcal L(s).
$$
Let $a_\infty$ be the pole order of $\gamma_M$ at $0$.
Each gamma factor has at most a simple pole there, and
their number, counted with multiplicity, is at most $D$.
Thus $0\leq a_\infty\leq D$.
Since $F_M$ has order $-d$ at $0$, we obtain
$$
m_0-a_\infty=-d,
\qquad
|m_0|\leq D+d.
$$

To bound $c_0$, put
$$
\Xi_M(s):=s^d(1-s)^dF_M(s).
$$
Hadamard factorization gives
$$
\frac{\Xi_M'}{\Xi_M}(0)
=
\frac{\Xi_M'}{\Xi_M}(2)
-\sum_\rho\frac{2}{\rho(2-\rho)}.
$$
The Euler-product and gamma-factor bounds give
$\Xi_M'/\Xi_M(2)\ll\kappa(0)$.
Under the Riemann hypothesis,
$|\rho(2-\rho)|\asymp1+\gamma^2$.
Thus \eqref{sum-zero-square}, established below, implies
$$
\frac{\Xi_M'}{\Xi_M}(0)\ll\kappa(0).
$$
Write
$$
\frac{\gamma_M'}{\gamma_M}(s)
=
-\frac{a_\infty}{s}+b_\infty+O(s).
$$
The gamma-factor expansions give
$b_\infty\ll\kappa(0)$.
Comparing finite parts at $0$ therefore yields
$$
c_0=
\frac{\Xi_M'}{\Xi_M}(0)
+d-\frac12\log q_M-b_\infty
\ll\kappa(0).
$$

Since
$$
\frac{x^{s+1}}{s(s+1)}
=
\frac{x}{s}+x(\log x-1)+O_x(s),
$$
we have
$$
R_0(x)=-m_0x\log x+(m_0-c_0)x.
$$
Hence
$$
|R_0(x)|
\ll
x\bigl((D+d)\log x+\kappa(0)\bigr).
$$
It follows that
$$
\psi_{M,1}(x)
=
\frac{dx^2}{2}
-
\sum_\rho\frac{x^{\rho+1}}{\rho(\rho+1)}
+
O\left(
x\bigl((D+d)\log x+\kappa(0)\bigr)
\right).
$$

As in the case of Artin $L$-functions, Hadamard
factorization, the functional equation and the
gamma-factor bounds give
$$
\#\{\rho:|\gamma-t|\leq1\}
\ll\kappa(2+|t|).
$$
Stieltjes integration therefore gives, for $T\geq2$,
\begin{align}
\sum_\rho\frac1{1+\gamma^2}
&\ll\kappa(0),\label{sum-zero-square}\\
\sum_{|\gamma|\leq T}\frac1{1+|\gamma|}
&\ll\log(2T)\kappa(T),\label{sum-zero-first}\\
\sum_{|\gamma|>T}\frac1{1+\gamma^2}
&\ll\frac{\kappa(T)}{T}.
\label{sum-zero-tail}
\end{align}
For each nontrivial zero,
$$
\left|
\frac{(x+h)^{\rho+1}-x^{\rho+1}}
{h\rho(\rho+1)}
\right|
\ll
\min\left\{
\frac{x^{1/2}}{1+|\gamma|},
\frac{x}{1+\gamma^2}
\right\}.
$$
Taking $T:=\max\{2,x^{1/2}\}$, splitting the zero sum
at $T$, and using \cref{psi}, we obtain
$$
\psi_M(x)
=
dx+O\bigl(x^{1/2}\log x\,\kappa(x)\bigr).
$$
Abel summation now gives
$$
\sum_{\substack{\pp,\ m\geq1\\(N\pp)^m\leq x}}
\frac{b_{\pp}(m)}m
=
d\,\Li(x)+O\bigl(x^{1/2}\kappa(x)\bigr).
$$
The terms with $m\geq2$ contribute $O(Dx^{1/2})$.
Here we use the elementary bounds
$$
\sum_{N\pp\leq y}\log N\pp\ll n_Ky,
\qquad
\sum_{m\geq2}\ \sum_{(N\pp)^m\leq y}\log N\pp
\ll n_Ky^{1/2}.
$$
Finally, the terms with $m=1$ and $\pp\in S_M$ contribute
at most $n_M|S_M|=S(M)$.
Since $D\ll\kappa(x)$, this proves the first assertion.

For the smoothed estimate, we use
\cite[Theorem 5.11]{IwaniecKowalski}, following
\cite[Lemma 7.1]{ZywinaLangTrotter}, keeping track of
the dependence on the test function.
Choose $0<a<b$ with $\operatorname{supp}f\subset(a,b)$,
and first suppose that
$$
x\geq x_f:=\max\{2,4/a,a^{-2}\}.
$$
Define
$$
g_x(u):=\frac{f(u)}{\log(xu)}
$$
on $(a,b)$ and extend it by zero outside this interval.
For every integer $j\geq0$,
$$
\|g_x^{(j)}\|_\infty
\ll_{f,j}\frac1{\log x}.
$$
Consequently, its Mellin transform satisfies, for every
integer $A\geq0$,
$$
\widehat g_x(\sigma+it)
:=
\int_0^\infty g_x(u)u^{\sigma+it}\,\frac{du}{u}
\ll_{f,A}
\frac1{\log x\,(1+|t|)^A},
\qquad
-\frac14\leq\sigma\leq2.
$$

Similarly, Mellin inversion and the
explicit formula give
\begin{align*}
&\sum_{\pp}\sum_{m\geq1}
b_{\pp}(m)\log N\pp\,
g_x\left(\frac{(N\pp)^m}{x}\right)\\
&\qquad=
dx\widehat g_x(1)
-
\sum_\rho x^\rho\widehat g_x(\rho)
+
O_f\left(\frac{\kappa(0)}{\log x}\right).
\end{align*}
Indeed, the residue at $0$ is
$-m_0\widehat g_x(0)=O_f((D+d)/\log x)$,
and the integral on $\Re s=-1/4$ is
$O_f(x^{-1/4}\kappa(0)/\log x)$.

By the motivic Riemann hypothesis, the Mellin bound
with $A=2$, and \eqref{sum-zero-square},
$$
\sum_\rho
\left|x^\rho\widehat g_x(\rho)\right|
\ll_f
\frac{x^{1/2}\kappa(0)}{\log x}.
$$
Moreover,
$$
x\widehat g_x(1)
=
\int_2^\infty\frac{f(t/x)}{\log t}\,dt.
$$
Since
$$
\log N\pp\,
g_x\left(\frac{(N\pp)^m}{x}\right)
=
\frac1m f\left(\frac{(N\pp)^m}{x}\right),
$$
we obtain
\begin{align*}
\sum_{\pp}\sum_{m\geq1}
\frac{b_{\pp}(m)}m
f\left(\frac{(N\pp)^m}{x}\right)
={}&
d\int_2^\infty\frac{f(t/x)}{\log t}\,dt\\
&+
O_f\left(
\frac{x^{1/2}\kappa(0)}{\log x}
\right).
\end{align*}

On the support of the summands,
$\log((N\pp)^m)\geq\frac12\log x$.
Thus the terms with $m\geq2$ contribute
$$
\ll_f
\frac{n_M}{\log x}
\sum_{m\geq2}\ \sum_{(N\pp)^m\leq bx}\log N\pp
\ll_f
\frac{Dx^{1/2}}{\log x}.
$$
Finally, the exceptional primes contribute
$$
\left|
\sum_{\pp\in S_M}
b_{\pp}(1)f\left(\frac{N\pp}{x}\right)
\right|
\leq
\|f\|_\infty n_M|S_M|
=
O_f(S(M)).
$$
Given that
$$
|\widehat g_x(\sigma+it)|
\ll_{a,b}
\frac{\|f\|_\infty+\|f'\|_\infty+\|f''\|_\infty}
{\log x\,(1+|t|)^2},
$$
The implied constant depends only on the $\|f^{(j)}\|_{L^\infty((0,\infty))}$ for $j\in \{0,1,2\}$ and the support of $f$.
This proves the second assertion for $x\geq x_f$.

For $2\leq x<x_f$, the prime sum and the main term
are bounded by $O_f(D+d)\ll_f\kappa(0)$.
Since $x^{1/2}/\log x$ is bounded below for $x\geq2$,
the estimate follows.
\end{proof}
We conclude this section by computing the Hodge-circle
constant $C_H$ in several cases, using the Hodge data
available in the literature. The entries in Table3 use the Betti trace metric
$$
\|X\|^2=-\frac12\operatorname{Re}\operatorname{Tr}(X^2).
$$
For this metric,
$$
C_H=2\pi\max_\tau
\left(\frac12\sum_{p,q}
h_\tau^{p,q}(M)(p-q)^2\right)^{1/2},
$$
where $\tau$ ranges over the embeddings $K\hookrightarrow\C$.
\clearpage
\begingroup
\small
\renewcommand{\arraystretch}{1.3}
\setlength{\tabcolsep}{10pt}
\setlength{\LTcapwidth}{\linewidth}

\begin{xltabular}{\linewidth}{@{}
  >{\raggedright\arraybackslash}p{0.32\linewidth}
  >{\centering\arraybackslash}X
  >{\centering\arraybackslash}p{0.20\linewidth}
@{}}
\toprule
\textbf{Motive}
& \textbf{Hodge data}
& $\boldsymbol{\frac{C_H}{2\pi}}$
\\
\midrule
\endfirsthead

\toprule
\textbf{Motive}
& \textbf{Hodge data}
& $\boldsymbol{\frac{C_H}{2\pi}}$
\\
\midrule
\endhead

Modular motive $M_f$,\newline weight $k-1$
&
$h^{k-1,0}=h^{0,k-1}=1$
&
$k-1$
\\ \addlinespace

$h^1(A)$,\newline $\dim A=g$
&
$h^{1,0}=h^{0,1}=g$
&
$\sqrt{g}$
\\ \addlinespace

$h^m(A)$,\newline $1\le m\le2g-1$
&
$h^{p,m-p}=
  \binom{g}{p}\binom{g}{m-p}$
&
$\sqrt{g\binom{2g-2}{m-1}}$
\\ \addlinespace

$h^0(A)$ or $h^{2g}(A)$
&
Diagonal types
&
$0$
\\ \addlinespace

K3-type motive,\newline rank $d$
&
$\begin{gathered}
h^{2,0}=h^{0,2}=1\\
h^{1,1}=d-2
\end{gathered}$
&
$2$
\\ \addlinespace

Rank-two motive,\newline $a\ne b$
&
$h^{a,b}=h^{b,a}=1$
&
$|a-b|$
\\ \addlinespace

Effective pure motive,\newline weight $w$, rank $d$
&
$\begin{gathered}
p,q\ge0\\
p+q=w
\end{gathered}$
&
$\le w\sqrt{d/2}$
\\ \addlinespace

Pure motive,\newline rank $d$
&
$\begin{gathered}
h^{p,q}\ne0\\
\Longrightarrow |p-q|\le L
\end{gathered}$
&
$\le L\sqrt{d/2}$
\\ \addlinespace

Tate twist of \newline Artin motives
&
Diagonal types
&
$0$
\\
\bottomrule

\caption{Hodge-circle constants for selected examples of motives}\\
\end{xltabular}
\endgroup
\chapter{Some Lie-theoretic estimates}\label[appendix]{ST}
In this section, we collect explicit formulas for the
Sato-Tate measures of the classical compact groups
and their small-radius asymptotics
on geodesic balls. Our main references for the density
formulas are \cite[Section 3.1]{Meckes2019}
and \cite{KatzSarnak1999}. We use the Euclidean angular metric on $(\mathbb{R}/2\pi\Z)^{\mathrm{rank}(G)}$.
\begingroup
\small
\renewcommand{\arraystretch}{1.3}
\setlength{\tabcolsep}{10pt}

\begin{xltabular}{\linewidth}{|
  >{\centering\arraybackslash}p{0.16\linewidth}|
  >{\centering\arraybackslash}X|
  >{\centering\arraybackslash}X|}
\hline
$G$ & $d\mu_{\st}$
& $\mu_{\st}(B_g(0,\varepsilon))$
\\[6pt] \hline
\endfirsthead

\hline
$G$ & $d\mu_{\st}$
& $\mu_{\st}(B_g(0,\varepsilon))$
\\[6pt] \hline
\endhead

$\mathrm{SU}(2)$
&
$\displaystyle
\frac{1}{\pi}\sin^2\theta\,d\theta
$
&
$\displaystyle
\begin{gathered}
\frac{2}{3\pi}\varepsilon^3\\
{}+O(\varepsilon^5)
\end{gathered}
$
\\[8pt] \hline

$\mathrm{USp}(2g)$
&
$\displaystyle
\begin{gathered}
\frac{2^{g(g-1)}}{g!\pi^g}
\prod_{j=1}^{g}\sin^2\theta_j\\
{}\times
\prod_{i<j}(\cos\theta_i-\cos\theta_j)^2\,d\theta
\end{gathered}
$
&
$\displaystyle
\begin{gathered}
C_g^1\varepsilon^{g(2g+1)}\\
{}+O_g(\varepsilon^{g(2g+1)+2})
\end{gathered}
$
\\[8pt] \hline

$\mathrm{U}(n)$
&
$\displaystyle
\begin{gathered}
\frac{1}{n!(2\pi)^n}\\
{}\times
\prod_{i<j}|e^{i\theta_i}-e^{i\theta_j}|^2\,d\theta
\end{gathered}
$
&
$\displaystyle
\begin{gathered}
C_n^4\varepsilon^{n^2}\\
{}+O_n(\varepsilon^{n^2+2})
\end{gathered}
$
\\[8pt] \hline

$\mathrm{SO}(2n+1)$
&
$\displaystyle
\begin{gathered}
\frac{2^{n(n-1)}}{n!\pi^n}
\prod_{j=1}^{n}\sin^2(\theta_j/2)\\
{}\times
\prod_{i<j}(\cos\theta_i-\cos\theta_j)^2\,d\theta
\end{gathered}
$
&
$\displaystyle
\begin{gathered}
C_n^2\varepsilon^{n(2n+1)}\\
{}+O_n(\varepsilon^{n(2n+1)+2})
\end{gathered}
$
\\[8pt] \hline

$\mathrm{SO}(2n)$,
$n\geq2$
&
$\displaystyle
\begin{gathered}
\frac{2^{n^2-3n+1}}{n!\pi^n}\\
{}\times
\prod_{i<j}(\cos\theta_i-\cos\theta_j)^2\,d\theta
\end{gathered}
$
&
$\displaystyle
\begin{gathered}
C_n^3\varepsilon^{n(2n-1)}\\
{}+O_n(\varepsilon^{n(2n-1)+2})
\end{gathered}
$
\\[8pt] \hline
\caption{Sato-Tate densities on the full angular torus
and small-radius asymptotics at the identity.}
\end{xltabular}
\endgroup
The constants above are given by
\begin{align*}
C_g^1
&=\frac{\prod_{j=1}^{g}(2j)!}
 {g!(2\pi)^{g/2}2^{g(2g+1)/2}
  \Gamma(1+g(2g+1)/2)},\\
C_n^2
&=\frac{\prod_{j=1}^{n}(2j)!}
 {n!(2\pi)^{n/2}2^{2n+n(2n+1)/2}
  \Gamma(1+n(2n+1)/2)},\\
C_n^3
&=\frac{\prod_{j=1}^{n-1}(2j)!}
 {(2\pi)^{n/2}2^{n-1+n(2n-1)/2}
  \Gamma(1+n(2n-1)/2)},\\
C_n^4
&=\frac{\prod_{j=1}^{n}j!}
 {n!(2\pi)^{n/2}2^{n^2/2}\Gamma(1+n^2/2)}.
\end{align*}
We conclude this section with examples of the computation of the Lie-theoretic parameters appearing in our Chebotarev-Sato-Tate estimate and extensively used in \cref{equid,sato2}.
\begingroup
\small
\renewcommand{\arraystretch}{1.5}
\setlength{\tabcolsep}{10pt}

\begin{xltabular}{\linewidth}{|
  c|
  >{\centering\arraybackslash}X|
  >{\centering\arraybackslash}X|
  >{\centering\arraybackslash}X|}
\hline
$G$ & $\operatorname{rank}G$ & $|\Phi^+|$ & $\delta_G$
\\ \hline
\endfirsthead

\hline
$G$ & $\dim T=\operatorname{rank}G$ & $|\Phi^+|$ & $\delta_G$
\\ \hline
\endhead

$\mathrm{SU}(2)$
& $1$ & $1$ & $2$
\\ \hline

$\mathrm{USp}(2g)$
& $g$ & $g^2$ & $g(g+1)$
\\ \hline

$\mathrm{U}(n)$
& $n$ & $\frac{n(n-1)}{2}$ & $\frac{n(n+1)}{2}$
\\ \hline

$\mathrm{SO}(2n+1)$
& $n$ & $n^2$ & $n(n+1)$
\\ \hline

$\mathrm{SO}(2n)$, $n\geq2$
& $n$ & $n(n-1)$ & $n^2$
\\ \hline

\caption{Dimensions of maximal tori, numbers of positive roots,
and $\delta_G=\operatorname{rank}G+|\Phi^+|$.}\\
\end{xltabular}
\endgroup

\backmatter
\bibliographystyle{amsplain}
\bibliography{references}
\end{document}